\documentclass[12pt]{amsart}

\usepackage[paperheight=11in,paperwidth=8.5in,left=1in,top=1.25in,right=1in,bottom=1.25in,]{geometry}
\usepackage[linktocpage=false,colorlinks=true,linkcolor=black,citecolor=black,
filecolor=black,urlcolor=black,pagebackref=false,pdfstartpage={1},pdfstartview={FitH},pdftitle={},
pdfauthor={},pdfsubject={},pdfcreator={},pdfproducer={},pdfkeywords={}]{hyperref}

\usepackage[lofdepth,lotdepth]{subfig}

\usepackage{afterpage}
\usepackage{amsfonts}
\usepackage{amsmath}
\allowdisplaybreaks[4]
\usepackage{amssymb}
\usepackage{amsthm}
\usepackage{bbm}
\usepackage[justification=centering]{caption}
\usepackage[capitalize,nameinlink,noabbrev,nosort]{cleveref}
\usepackage{enumerate}
\usepackage{enumitem}
\usepackage{float}
\restylefloat{figure}
\usepackage{graphicx}
\usepackage{kbordermatrix}
\usepackage{mathrsfs}
\usepackage{mathtools}
\usepackage[framemethod=tikz,ntheorem,xcolor]{mdframed}
\usepackage{parcolumns}
\usepackage{tabularx}
\usepackage{tikz}\pdfpageattr{/Group <</S /Transparency /I true /CS /DeviceRGB>>}
\usetikzlibrary{arrows,calc,intersections,matrix,positioning,through}
\usepackage{tkz-euclide}
\usetkzobj{all}
\usepackage{tikz-cd}
\tikzset{commutative diagrams/.cd,every label/.append style = {font = \normalsize}}
\usepackage{tkz-graph}
\usetikzlibrary{arrows,%
                shapes,positioning}

\tikzset{
  dep u/.style={insert path={-- ++(0,15) node{}}},
  dep r/.style={insert path={-- ++(15,0) node{}}},
  dep d/.style={insert path={-- ++(0,-15) node{}}},
  dep l/.style={insert path={-- ++(-15,0) node{}}},
  recurse lattice path/.code args={#1#2}{
    \ifx#1.\else\tikzset{dep #1,recurse lattice path=#2}\fi
  },
  lattice path/.style={recurse lattice path=#1.}
}

\usetikzlibrary{external}
\usepackage[all]{xy}
\usepackage[aligntableaux=center]{ytableau}
\usepackage{epstopdf}

\usepackage{array}
\usepackage{color}

\numberwithin{equation}{section}

\newtheorem{thm}[equation]{Theorem}
\AfterEndEnvironment{thm}{\noindent\ignorespaces}
\newtheorem{cor}[equation]{Corollary}
\AfterEndEnvironment{cor}{\noindent\ignorespaces}
\newtheorem{lem}[equation]{Lemma}
\AfterEndEnvironment{lem}{\noindent\ignorespaces}
\newtheorem{prop}[equation]{Proposition}
\AfterEndEnvironment{prop}{\noindent\ignorespaces}
\newtheorem{conj}[equation]{Conjecture}
\AfterEndEnvironment{conj}{\noindent\ignorespaces}
\newtheorem{prob}[equation]{Problem}
\AfterEndEnvironment{prob}{\noindent\ignorespaces}
\newtheorem{question}[equation]{Question}
\AfterEndEnvironment{question}{\noindent\ignorespaces}

\AfterEndEnvironment{exercise}{\noindent\ignorespaces}

\theoremstyle{definition}
\newtheorem{defn}[equation]{Definition}

\AfterEndEnvironment{defn}{\noindent\ignorespaces}
\AfterEndEnvironment{definition}{\noindent\ignorespaces}
\newtheorem*{pf_no_qed}{Proof}
\newenvironment{pf}[1][]{\begin{pf_no_qed}[#1]\pushQED{\qed}}{\popQED\end{pf_no_qed}}
\AfterEndEnvironment{pf}{\noindent\ignorespaces}

\AfterEndEnvironment{proof}{\noindent\ignorespaces}
\newtheorem{eg_no_qed}[equation]{Example}
\newenvironment{eg}[1][]{\begin{eg_no_qed}[#1]\pushQED{\qed}}{\popQED\end{eg_no_qed}}
\AfterEndEnvironment{eg}{\noindent\ignorespaces}

\AfterEndEnvironment{example}{\noindent\ignorespaces}

\AfterEndEnvironment{ex}{\noindent\ignorespaces}
\newtheorem{rmk}[equation]{Remark}
\AfterEndEnvironment{rmk}{\noindent\ignorespaces}
\newtheorem{remark}[equation]{Remark}
\AfterEndEnvironment{remark}{\noindent\ignorespaces}

\theoremstyle{remark}
\newtheorem*{claim}{Claim}
\AfterEndEnvironment{claim}{\noindent\ignorespaces}
\newtheorem*{claimpf_no_qed}{Proof of Claim}
\newenvironment{claimpf}[1][]{\begin{claimpf_no_qed}[#1]\pushQED{\qed}}{\popQED\end{claimpf_no_qed}}
\AfterEndEnvironment{claimpf}{\noindent\ignorespaces}

\usepackage{comment}

\font\pipefont=lcircle10
\def\elbow{\smash{\raise3pt\hbox{\pipefont\rlap{\rlap{\char'014}\char'016}}}}
\def\textelbow{\;\;\,\elbow\;\;\,}
\def\halfelbow{\smash{\raise2pt\hbox{\pipefont\rlap{\rlap{\rlap{\char'015}\phantom{\char'017}}}}}}
\def\cross{\smash{\lower5pt\hbox{\rlap{\vrule height16pt}}\raise3pt\hbox{\rlap{\hskip-8pt \vrule height0.4pt depth0pt width16pt}}}}
\def\textcross{\;\;\,\cross\;\;\,}
\ytableausetup{boxsize=16pt}

\newcommand{\C}{\mathbb{C}}

\newcommand{\on}{\operatorname}

\newcommand{\Mat}{\on{Mat}}
\newcommand{\Gr}{\on{Gr}}

\newcommand{\w}{\on{weight}}

\newcommand{\dom}{\on{dom}}
\newcommand{\Dom}{\on{Dom}}
\newcommand{\var}{\on{var}}
\newcommand{\varbar}{\overline{\on{var}}}
\newcommand{\sign}{\on{sign}}

\newcommand{\Le}{\textup{\protect\scalebox{-1}[1]{L}}}

\newcommand{\R}{\mathbb{R}}

\DeclareMathOperator{\down}{\mathsf{down}}

\newcommand{\rf}[1]{\hyperref[#1]{(\ref*{#1})}}

\DeclareMathOperator{\shadow}{shadow}

\DeclareMathOperator{\supp}{supp}

\DeclareMathOperator{\touch}{touch}

\DeclareMathOperator{\up}{\mathsf{up}}

\title{Decompositions of amplituhedra}

\author{Steven N.\ Karp}
\address{Department of Mathematics, University of Michigan}
\email{\href{mailto:snkarp@umich.edu}{snkarp@umich.edu}}

\author{Lauren K.\ Williams}
\address{Department of Mathematics, University of California, Berkeley}
\email{\href{mailto:williams@math.berkeley.edu}{williams@math.berkeley.edu}}

\author{Yan X Zhang}
\address{Department of Mathematics and Statistics, San Jos\'{e} State University}
\email{\href{mailto:yan.x.zhang@sjsu.edu}{yan.x.zhang@sjsu.edu}}

\begin{document}

\makeatletter
\tikzset{
    fading speed/.code={
        \pgfmathtruncatemacro\tikz@startshading{50-(100-#1)*0.25}
        \pgfmathtruncatemacro\tikz@endshading{50+(100-#1)*0.25}
        \pgfdeclareverticalshading[%
            tikz@axis@top,tikz@axis@middle,tikz@axis@bottom%
        ]{axis#1}{100bp}{%
            color(0bp)=(tikz@axis@bottom);
            color(\tikz@startshading)=(tikz@axis@bottom);
            color(50bp)=(tikz@axis@middle);
            color(\tikz@endshading)=(tikz@axis@top);
            color(100bp)=(tikz@axis@top)
        }
        \tikzset{shading=axis#1}
    }
}
\makeatother

\tikzexternaldisable

\begin{abstract}
The \emph{(tree) amplituhedron} $\mathcal{A}_{n,k,m}$ is the image
in the Grassmannian $\Gr_{k,k+m}$ of the totally nonnegative Grassmannian
$\Gr_{k,n}^{\geq 0}$, under a (map induced by a) linear map which is totally
positive. It was introduced by
Arkani-Hamed and Trnka in 2013 in order to give a geometric basis for the
computation of scattering amplitudes in $\mathcal{N}=4$
supersymmetric Yang-Mills theory.  
In the case relevant to physics ($m=4$), there is a collection of 
recursively-defined 
$4k$-dimensional \emph{BCFW cells} in 
$\Gr_{k,n}^{\geq 0}$,
whose images conjecturally 
``triangulate" the amplituhedron---that is, their images are disjoint
and cover a dense subset of 
$\mathcal{A}_{n,k,4}$.
In this paper, we approach this problem by first giving an explicit
(as opposed to recursive) description of the BCFW cells.  We then 
develop sign-variational tools which we use to 
prove that when $k=2$, the images of these cells are disjoint in $\mathcal{A}_{n,k,4}$.  We also conjecture that for arbitrary even $m$,
there is a decomposition of the amplituhedron $\mathcal{A}_{n,k,m}$
involving precisely $M(k, n-k-m, \frac{m}{2})$ top-dimensional cells
(of dimension $km$),
where $M(a,b,c)$ is the number of plane partitions contained in an
$a \times b \times c$ box.  This agrees with the fact that when $m=4$, the number of BCFW cells is the Narayana
number $N_{n-3,k+1} = \frac{1}{n-3}\binom{n-3}{k+1}\binom{n-3}{k}$.
\end{abstract}

\maketitle
\setcounter{tocdepth}{1}
\tableofcontents

\section{Introduction}\label{sec:intro}

\noindent The totally nonnegative Grassmannian 
$\Gr_{k,n}^{\ge 0}$ 
is the subset of the real Grassmannian $\Gr_{k,n}$ 
where all Pl\"ucker coordinates are nonnegative.  Following seminal
work of Lusztig \cite{lusztig}, as well as by Fomin and Zelevinsky \cite{FZ}, 
Postnikov initiated the combinatorial study of $\Gr_{k,n}^{\geq 0}$ and its cell decomposition \cite{postnikov}.  In particular, Postnikov 
showed how the cells in the cell decomposition are naturally indexed by 
combinatorial objects including \emph{decorated permutations}, 
 \emph{\Le -diagrams}, 
equivalence classes of \emph{plabic graphs}.
Since then the totally nonnegative Grassmannian has found applications in 
diverse contexts such as mirror symmetry \cite{MarshRietsch}, 
soliton solutions to the KP equation \cite{KodamaWilliams}, 
and scattering amplitudes for $\mathcal{N}=4$
supersymmetric Yang-Mills theory 
\cite{abcgpt}.

Building on \cite{abcgpt}, 
Arkani-Hamed and Trnka \cite{arkani-hamed_trnka} recently introduced a beautiful new 
mathematical object called the \emph{(tree) amplituhedron}, which 
is the image of the totally nonnegative Grassmannian under a particular map.

\begin{defn}\label{def:amp}
For $a \le b$, define $\Mat_{a,b}^{>0}$ as the set of real $a \times b$ matrices whose $a\times a$ minors are all positive. 
Let $Z \in \Mat_{k+m,n}^{>0}$, where
$m\ge 0$ is fixed with $k+m \leq n$. 
Then $Z$ induces a map
$$\tilde{Z}:\Gr_{k,n}^{\ge 0} \to \Gr_{k,k+m}$$ defined by 
$$\tilde{Z}(\langle v_1,\dots, v_k \rangle) := \langle Z(v_1),\dots, Z(v_k) \rangle,$$
where $\langle v_1,\dots,v_k\rangle$ is an element of 
$\Gr_{k,n}^{\ge 0}$ written as the span of $k$ basis vectors.\footnote{The
fact that $Z$ has positive maximal minors ensures that $\tilde{Z}$
is well defined
\cite{arkani-hamed_trnka}.  
See \cite[Theorem 4.2]{karp} for a characterization of when 
a matrix $Z$ gives rise to a well-defined map $\tilde{Z}$.}
The \emph{(tree) amplituhedron} $\mathcal{A}_{n,k,m}(Z)$ is defined to be the image
$\tilde{Z}(\Gr_{k,n}^{\ge 0})$ inside $\Gr_{k,k+m}$.  
\end{defn}

In special cases the amplituhedron recovers familiar objects. If $Z$ is a square matrix, i.e.\ 
$k+m=n$, then $\mathcal{A}_{n,k,m}(Z)$ is isomorphic to 
the totally nonnegative Grassmannian. If $k=1$, then it follows from 
\cite{Sturmfels} that $\mathcal{A}_{n,1,m}(Z)$ is a {\itshape cyclic polytope} in projective space $\mathbb{P}^m$.

While the amplituhedron $\mathcal{A}_{n,k,m}(Z)$ is an interesting mathematical object for any $m$, the case relevant
 to physics is $m=4$. In this case, it provides a geometric basis for the computation of {\itshape scattering amplitudes} in $\mathcal{N}=4$ supersymmetric Yang-Mills theory. These amplitudes are complex numbers related to the probability of observing a certain scattering process of $n$ particles.
It is expected that the leading-order term of such amplitudes can be expressed as the integral of a canonical form on the amplituhedron $\mathcal{A}_{n,k,4}(Z)$.  This 
statement is closely related to the following
 conjecture of Arkani-Hamed and Trnka \cite{arkani-hamed_trnka}.

\begin{conj}\label{AHTConj}
Let $Z \in \Mat_{k+4,n}^{>0}$, and 
let $\mathcal{C}_{n,k,4}$ be the collection of \emph{BCFW cells} in 
$\Gr_{k,n}^{\geq 0}$. Then the images under $\tilde{Z}$ of the 
cells $\mathcal{C}_{n,k,4}$ ``triangulate" the $m=4$ amplituhedron, 
i.e.\ they are 
pairwise disjoint, and together they cover a dense subset of 
the amplituhedron $\mathcal{A}_{n,k,4}(Z)$.
\end{conj}

More specifically, 
the BCFW recurrence
\cite{BCF, BCFW}, of Britto, Cachazo, Feng, and Witten,
provides one way to compute scattering amplitudes.
Translated into the Grassmannian formulation of 
\cite{abcgpt}, the terms in the BCFW recurrence can 
be identified with a collection of $4k$-dimensional cells in 
$\Gr_{k,n}^{\geq 0}$ which we refer to as the \emph{BCFW cells} $\mathcal{C}_{n,k,4}$. 
If the images of these BCFW cells in 
$\mathcal{A}_{n,k,4}(Z)$ fit together in a nice way, then we can combine the canonical form coming from each term to obtain the canonical form on $\mathcal{A}_{n,k,4}(Z)$.

In this paper, we make a first step towards understanding 
\cref{AHTConj}.  The 
BCFW cells are defined 
recursively  in 
terms of plabic graphs (see \cref{sec:recursive}), and moreover there is a `shift by 2' applied
at the end of the recursion (see 
\cref{def:BCFW-permutations}).  Hence proving anything about how the 
images of the BCFW cells fit together is not at all straightforward from 
the definitions.
To approach \cref{AHTConj}, we start by 
giving an explicit, non-recursive description of the BCFW cells.
Namely, we index the BCFW cells in $\Gr_{k,n}^{\geq 0}$ 
by pairs of noncrossing lattice paths inside a $k\times (n-k-4)$ rectangle, and 
associate a \Le -diagram to each pair of lattice paths, from which 
we can read off the corresponding cell (see  
\cref{thm:BCFW}, proved in \cref{sec:bijection}).  We then use these \Le -diagrams to understand the case $k=2$:
we derive an elegant description of basis vectors for elements of each BCFW cell in terms of `dominoes' (\cref{thm:4-dominoes}), and show that the images of distinct BCFW cells are disjoint in the amplituhedron $\mathcal{A}_{n,2,4}(Z)$ (\cref{thm:disjointness=k=2}). The proof uses classical results about sign variation, along with some new tools particularly suited to our problem. 

We expect that our techniques may be helpful in understanding the case of arbitrary $k$, and we make a step in this direction in an appendix to 
this paper, which is joint with Hugh Thomas.
We use a bijection between BCFW cells and Dyck paths to associate 
a conjectural 
`domino basis' to each element of a 
BCFW cell (\cref{conj:domino}).
We leave the proof of the conjecture and the analysis of 
how general BCFW cells fit together to future work.

As a warmup to our study of BCFW cells in the case $m=4$, in \cref{sec:m2} we develop an analogous story in the case $m=2$. Namely, we give a BCFW-style recursion on plabic graphs, describe the resulting cells of $\Gr_{k,n}^{\ge 0}$ using lattice paths and domino bases, and prove disjointness of the images of these cells inside the $m=2$ amplituhedron $\mathcal{A}_{n,k,2}(Z)$.

It was observed (e.g.\ in \cite{abcgpt}) that 
the number of BCFW cells $|\mathcal{C}_{n,k,4}|$ is the Narayana number
$\frac{1}{n-3}\binom{n-3}{k+1}\binom{n-3}{k}$.
Motivated by this fact, as well as known results about 
decompositions of amplituhedra
when $m=2$ or $k=1$, we conjecture that when $m$ is even, there is a 
decomposition of $\mathcal{A}_{n,k,m}(Z)$ which 
involves precisely $M(k, n-k-m, \frac{m}{2})$ top-dimensional 
cells (see \cref{mainconj}). Here $M(a,b,c)$ denotes the number of collections of $c$ noncrossing lattice paths inside an $a\times b$ rectangle, or equivalently, 
the number of plane partitions which fit inside an $a\times b\times c$ box. See \cref{sec:Nabc} for other combinatorial interpretations
of $M(a,b,c)$, as well as \cref{oddm} for a possible extension of 
\cref{mainconj} to odd $m$.

\textsc{Acknowledgements:}
We are grateful to Nima Arkani-Hamed, Jacob Bourjaily, Greg Kuperberg, Hugh Thomas, and Jaroslav Trnka for helpful conversations.
This material is based upon work supported by the National Science Foundation
under agreement No.\ DMS-1128155 and No.\ DMS-1600447.  Any opinions,
findings and conclusions or recommendations expressed in this material
are those of the authors and do not necessarily reflect  the
views of the National Science Foundation.

\section{Background on the totally nonnegative Grassmannian}\label{sec:background}

\noindent The {\itshape (real) Grassmannian} $\Gr_{k,n}$ (for $0\le k \le n$) is the space of all $k$-dimensional subspaces of $\R^n$.  An element of
$\Gr_{k,n}$ can be viewed as a $k\times n$ matrix of rank $k$ modulo invertible row operations, whose rows give a basis for the $k$-dimensional subspace.

Let $[n]$ denote $\{1,\dots,n\}$, and $\binom{[n]}{k}$ denote the set of all $k$-element subsets of $[n]$. Given $V\in\Gr_{k,n}$ represented by a $k\times n$ matrix $A$, for $I\in \binom{[n]}{k}$ we let $\Delta_I(V)$ be the $k\times k$ minor of $A$ using the columns $I$. The $\Delta_I(V)$ do not depend on our choice of matrix $A$ (up to simultaneous rescaling by a nonzero constant), and are called the {\itshape Pl\"{u}cker coordinates} of $V$.

\begin{defn}[{\cite[Section~3]{postnikov}}]\label{def:positroid}
We say that $V\in\Gr_{k,n}$ is {\itshape totally nonnegative} if $\Delta_I(V)\ge 0$ for all $I\in\binom{[n]}{k}$, and {\itshape totally positive} if $\Delta_I(V) > 0$ for all $I\in\binom{[n]}{k}$. The set of all totally nonnegative $V\in\Gr_{k,n}$ is the {\it totally nonnegative Grassmannian} $\Gr_{k,n}^{\ge 0}$, and the set of all totally positive $V$ is the {\itshape totally positive Grassmannian} $\Gr_{k,n}^{>0}$. For $M\subseteq \binom{[n]}{k}$, let $S_{M}$ be
the set of $V\in\Gr_{k,n}^{\geq 0}$ with the prescribed collection of Pl\"{u}cker coordinates strictly positive (i.e.\ $\Delta_I(V)>0$ for all $I\in M$), and the remaining Pl\"{u}cker coordinates
equal to zero (i.e.\ $\Delta_J(V)=0$ for all $J\in\binom{[n]}{k}\setminus M$). If $S_M\neq\emptyset$, we call $M$ a \emph{positroid} and $S_M$ its \emph{positroid cell}.
\end{defn}

Each positroid cell $S_{M}$ is indeed a topological cell \cite[Theorem 6.5]{postnikov}, and moreover, the positroid cells of $\Gr_{k,n}^{\ge 0}$ glue together to form a CW complex \cite{PSW}.

\subsection{Combinatorial objects parameterizing cells}

In \cite{postnikov}, Postnikov defined several families of combinatorial objects which are in bijection with cells of the totally nonnegative Grassmannian, including \emph{decorated permutations}, \emph{\Le -diagrams}, and equivalence classes of \emph{reduced plabic graphs}. In this section, we introduce these objects, and give bijections between them.  This will give us a canonical way to label each positroid by a decorated permutation, a \Le-diagram, and an equivalence class of plabic graphs.
\begin{defn}\label{defn:decperm}
A \emph{decorated permutation} of $[n]$ is a bijection $\pi : [n] \to [n]$ whose fixed points are each colored either black or white. We denote a black fixed point $i$ by $\pi(i) = \underline{i}$, and a white fixed point $i$ by $\pi(i) = \overline{i}$.
An \emph{anti-excedance} of the decorated permutation $\pi$ is an element $i \in [n]$ such that either $\pi^{-1}(i) > i$ or $\pi(i)=\overline{i}$. 
\end{defn}
Postnikov showed that the positroids for $\Gr_{k,n}^{\ge 0}$ are indexed by decorated permutations of $[n]$ with exactly $k$ anti-excedances \cite[Section 16]{postnikov}.

\begin{figure}[ht]
\begin{center}
$$
\begin{tikzpicture}[baseline=(current bounding box.center)]
\pgfmathsetmacro{\scalar}{1.6};
\pgfmathsetmacro{\unit}{\scalar*0.922/1.6};
\draw[thick](0,0)rectangle(6*\unit,-4*\unit);
\foreach \x in {1,...,5}{
\draw[thick](\x*\unit-\unit,0)rectangle(\x*\unit,-\unit);}
\foreach \x in {1,...,5}{
\draw[thick](\x*\unit-\unit,-\unit)rectangle(\x*\unit,-2*\unit);}
\foreach \x in {1,...,3}{
\draw[thick](\x*\unit-\unit,-2*\unit)rectangle(\x*\unit,-3*\unit);}
\foreach \x in {1,...,2}{
\draw[thick](\x*\unit-\unit,-3*\unit)rectangle(\x*\unit,-4*\unit);}
\node[inner sep=0]at(0.5*\unit,-0.5*\unit){\scalebox{\scalar}{$+$}};
\node[inner sep=0]at(1.5*\unit,-0.5*\unit){\scalebox{\scalar}{$0$}};
\node[inner sep=0]at(2.5*\unit,-0.5*\unit){\scalebox{\scalar}{$+$}};
\node[inner sep=0]at(3.5*\unit,-0.5*\unit){\scalebox{\scalar}{$0$}};
\node[inner sep=0]at(4.5*\unit,-0.5*\unit){\scalebox{\scalar}{$+$}};
\node[inner sep=0]at(0.5*\unit,-1.5*\unit){\scalebox{\scalar}{$+$}};
\node[inner sep=0]at(1.5*\unit,-1.5*\unit){\scalebox{\scalar}{$+$}};
\node[inner sep=0]at(2.5*\unit,-1.5*\unit){\scalebox{\scalar}{$+$}};
\node[inner sep=0]at(3.5*\unit,-1.5*\unit){\scalebox{\scalar}{$+$}};
\node[inner sep=0]at(4.5*\unit,-1.5*\unit){\scalebox{\scalar}{$+$}};
\node[inner sep=0]at(0.5*\unit,-2.5*\unit){\scalebox{\scalar}{$0$}};
\node[inner sep=0]at(1.5*\unit,-2.5*\unit){\scalebox{\scalar}{$0$}};
\node[inner sep=0]at(2.5*\unit,-2.5*\unit){\scalebox{\scalar}{$0$}};
\node[inner sep=0]at(0.5*\unit,-3.5*\unit){\scalebox{\scalar}{$+$}};
\node[inner sep=0]at(1.5*\unit,-3.5*\unit){\scalebox{\scalar}{$+$}};
\end{tikzpicture}\qquad\qquad\qquad
\begin{tikzpicture}[baseline=(current bounding box.center)]
\pgfmathsetmacro{\unit}{0.922};
\useasboundingbox(0,0)rectangle(6*\unit,-4*\unit);
\coordinate (vstep)at(0,-0.24*\unit);
\coordinate (hstep)at(0.17*\unit,0);
\draw[thick](0,0)--(6*\unit,0) (0,0)--(0,-4*\unit);
\node[inner sep=0]at(0,0){\scalebox{1.6}{\begin{ytableau}
\none \\
\none \\
\none \\
\none \\
\none & \none & \none & \none & \none & \cross & \elbow & \cross & \elbow & \cross \\
\none & \none & \none & \none & \none & \elbow & \elbow & \elbow & \elbow & \elbow \\
\none & \none & \none & \none & \none & \cross & \cross & \cross \\
\none & \none & \none & \none & \none & \elbow & \elbow
\end{ytableau}}};
\node[inner sep=0]at($(5.5*\unit,0)+(vstep)$){$1$};
\node[inner sep=0]at($(5*\unit,-0.5*\unit)+(hstep)$){$2$};
\node[inner sep=0]at($(5*\unit,-1.5*\unit)+(hstep)$){$3$};
\node[inner sep=0]at($(4.5*\unit,-2*\unit)+(vstep)$){$4$};
\node[inner sep=0]at($(3.5*\unit,-2*\unit)+(vstep)$){$5$};
\node[inner sep=0]at($(3*\unit,-2.5*\unit)+(hstep)$){$6$};
\node[inner sep=0]at($(2.5*\unit,-3*\unit)+(vstep)$){$7$};
\node[inner sep=0]at($(2*\unit,-3.5*\unit)+(hstep)$){$8$};
\node[inner sep=0]at($(1.5*\unit,-4*\unit)+(vstep)$){$9$};
\node[inner sep=0]at($(0.5*\unit,-4*\unit)+(vstep)$){$10$};
\node[inner sep=0]at($(5.5*\unit,0)-(vstep)$){$1$};
\node[inner sep=0]at($(4.5*\unit,0)-(vstep)$){$4$};
\node[inner sep=0]at($(3.5*\unit,0)-(vstep)$){$5$};
\node[inner sep=0]at($(2.5*\unit,0)-(vstep)$){$7$};
\node[inner sep=0]at($(1.5*\unit,0)-(vstep)$){$9$};
\node[inner sep=0]at($(0.5*\unit,0)-(vstep)$){$10$};
\node[inner sep=0]at($(0,-0.5*\unit)-(hstep)$){$2$};
\node[inner sep=0]at($(0,-1.5*\unit)-(hstep)$){$3$};
\node[inner sep=0]at($(0,-2.5*\unit)-(hstep)$){$6$};
\node[inner sep=0]at($(0,-3.5*\unit)-(hstep)$){$8$};
\end{tikzpicture}
$$
\caption{A \Le -diagram $D$ of type $(4,10)$ with shape $\lambda=(5,5,3,2)$, and its corresponding pipe dream $P$ with decorated permutation $\pi_D = (\underline{1},5,4,9,7,\overline{6},2,10,3,8)$.}
\label{fig:Le}
\end{center}
\end{figure}
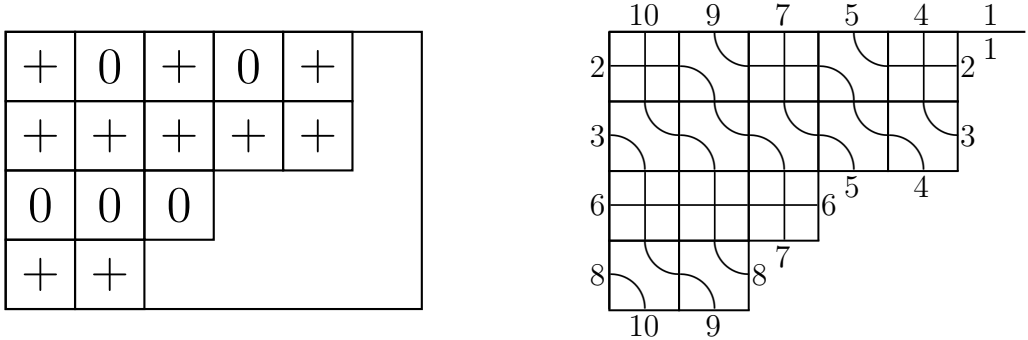

Now we introduce certain fillings of Young diagrams with the symbols $0$ and $+$, called {\it $\oplus$-diagrams}, and associate a decorated permutation to each such diagram. Postnikov \cite[Section 20]{postnikov} showed that a special subset of these diagrams, called {\it \Le -diagrams}, are in bijection with decorated permutations. We introduce the more general $\oplus$-diagrams here, since in \cref{sec:BCFW} we will use a distinguished subset of them to index the BCFW cells of $\Gr_{k,n}^{\ge 0}$.
\begin{defn}\label{def:oplus}
Fix $0 \le k \le n$. Given a partition $\lambda$, we let $Y_{\lambda}$ denote the Young diagram of $\lambda$. A {\itshape $\oplus$-diagram of type (k,n)} is a filling $D$ of a Young diagram $Y_\lambda$ fitting inside a $k\times (n-k)$ rectangle with the symbols $0$ and $+$ (such that each box of $Y$ is filled with exactly one symbol). We call $\lambda$ the {\itshape shape} of $D$. (See \cref{fig:Le}.)

We associate a decorated permutation $\pi$ of $[n]$ to $D$ as follows (see \cite[Section 19]{postnikov}).
\begin{enumerate}[label=\arabic*.]
\item Replace each $+$ in $D$ with an elbow $\textelbow$, and each $0$ in $D$ with a cross $\textcross$.
\item View the southeast border of $Y_\lambda$ as a lattice path of $n$ steps from the northeast corner to the southwest corner of the $k \times (n-k)$ rectangle, and label its edges by $1, \dots, n$.
\item Label each edge of the northwest border of $Y_{\lambda}$ with the label of its opposite edge on the southeast border. This gives a {\itshape pipe dream} $P$ associated to $D$.
\item Read off the decorated permutation $\pi$ from $P$ by following the `pipes' from the southeast border to the northwest border $Y_\lambda$. If the pipe originating at $i$ ends at $j$, we set $\pi(i) := j$. If $\pi(i)=i$, then either $i$ labels two horizontal edges or two vertical edges of $P$. In the former case, we set $\pi(i):=\underline{i}$, and in the latter case, we set $\pi(i):=\overline{i}$.
\end{enumerate}
\cref{fig:Le} illustrates this procedure. We denote the pipe dream $P$ by $P(D)$, and the decorated permutation $\pi$ by $\pi_D$. Note that the anti-excedances of $\pi$ correspond to the vertical steps of the southeast border of $Y_\lambda$, so $\pi$ has exactly $k$ anti-excedances. We denote the corresponding positroid cell of $\Gr_{k,n}^{\ge 0}$ (see \cite[Section 16]{postnikov}, and also \cref{network_param}) by $S_\pi$ or $S_D$.
\end{defn}

\begin{defn}\label{def:le_reduced}
Let $D$ be a $\oplus$-diagram, and $P$ its associated pipe dream from \cref{def:oplus}.

(i) We say that $D$ is {\itshape reduced} if no two pipes of $P$ cross twice.

(ii) We say that $D$ is a {\itshape \Le-diagram} (or Le-diagram) if it avoids the \Le -configuration, i.e.\ it has no $0$ with both a $+$ above it in the same column and a $+$ to its left in the same row. Equivalently (see \cite[Lemma 19.3]{postnikov}), any pair of pipes which cross do not later touch or cross, when read from southeast to northwest. (Two pipes {\itshape cross} if they form a cross $\textcross$, and {\itshape touch} if they form an elbow $\textelbow$.) So, \Le -diagrams are reduced.
\end{defn}

Postnikov showed that \Le -diagrams correspond to decorated permutations. Lam and Williams later showed how to transform any reduced $\oplus$-diagram into a \Le -diagram by using certain moves. We state these results.
\begin{lem}[{\cite[Corollary 20.1 and Theorem 6.5]{postnikov}}]\label{D_to_pi}
The map $D\mapsto\pi_D$ from \cref{def:oplus} is a bijection from the set of \Le -diagrams of type $(k,n)$ to the set of decorated permutations of $[n]$ with $k$ anti-excedances. Therefore, \Le -diagrams of type $(k,n)$ index the cells of $\Gr_{k,n}^{\ge 0}$. The dimension of the positroid cell $S_D$ indexed by $D$ is the number of $+$'s in $D$.
\end{lem}

\begin{lem}[\cite{LamWilliams}]\label{le_moves}
Let $D$ be a $\oplus$-diagram. Then $D$ is reduced if and only if $D$ can be transformed into a \Le -diagram $D'$ by a sequence of {\it \Le -moves}:
$$
\begin{array}{ccc}
\begin{ytableau}+\end{ytableau} & \cdots & \begin{ytableau}+\end{ytableau} \\
\vdots & \ddots & \vdots \\
\begin{ytableau}+\end{ytableau} & \cdots & \begin{ytableau}0\end{ytableau} \\
\end{array}\qquad\mapsto\qquad
\begin{array}{ccc}
\begin{ytableau}0\end{ytableau} & \cdots & \begin{ytableau}+\end{ytableau} \\
\vdots & \ddots & \vdots \\
\begin{ytableau}+\end{ytableau} & \cdots & \begin{ytableau}+\end{ytableau} \\
\end{array}.
$$
In this picture, the four boxes on each side denote the corners of a rectangle whose height and width are both at least $2$, and whose other boxes (aside from the corners) all contain $0$'s. \Le -moves preserve the decorated permutation of a $\oplus$-diagram. Hence $D$ indexes the cell $S_{D'}$ with decorated permutation $\pi_D = \pi_{D'}$, and the dimension of $S_{D'}$ is the number of $+$'s in $D$.
\end{lem}
For example, here is a sequence of \Le -moves which transforms a reduced $\oplus$-diagram into a \Le -diagram:
$$
\begin{ytableau}
+ & 0 & + & + & 0 & + \\
+ & 0 & 0 & 0 & +
\end{ytableau}\quad\mapsto\quad
\begin{ytableau}
0 & 0 & + & + & 0 & + \\
+ & 0 & + & 0 & +
\end{ytableau}\quad\mapsto\quad
\begin{ytableau}
0 & 0 & 0 & + & 0 & + \\
+ & 0 & + & + & +
\end{ytableau}\;.
$$
\begin{pf}
We explain how to deduce this result from the work of Lam and Williams. They showed \cite[Section 5]{LamWilliams} that {\itshape any} $\oplus$-diagram can be transformed into a \Le -diagram using the \Le -moves above, as well as the {\itshape uncrossing moves}
$$
\begin{array}{ccc}
\begin{ytableau}0\end{ytableau} & \cdots & \begin{ytableau}+\end{ytableau} \\
\vdots & \ddots & \vdots \\
\begin{ytableau}+\end{ytableau} & \cdots & \begin{ytableau}0\end{ytableau} \\
\end{array}\qquad\mapsto\qquad
\begin{array}{ccc}
\begin{ytableau}+\end{ytableau} & \cdots & \begin{ytableau}+\end{ytableau} \\
\vdots & \ddots & \vdots \\
\begin{ytableau}+\end{ytableau} & \cdots & \begin{ytableau}+\end{ytableau} \\
\end{array}.
$$
Note that a \Le -move in the statement of the lemma takes a reduced $\oplus$-diagram to a reduced $\oplus$-diagram, and an uncrossing move above can only be performed on a non-reduced $\oplus$-diagram. These observations imply the result.
\end{pf}

\begin{defn}
A {\it plabic graph}\footnote{``Plabic'' stands for {\itshape planar bi-colored}.}  is an undirected planar graph $G$ drawn inside a disk
(considered modulo homotopy)
with $n$ {\it boundary vertices} on the boundary of the disk,
labeled $1,\dots,n$ in clockwise order, as well as some {\it internal vertices}. Each boundary vertex is incident to a single edge, and each internal vertex is colored either black or white. If a boundary vertex is incident to a leaf (a vertex of degree $1$), we refer to that leaf as a \emph{lollipop}.
\end{defn}

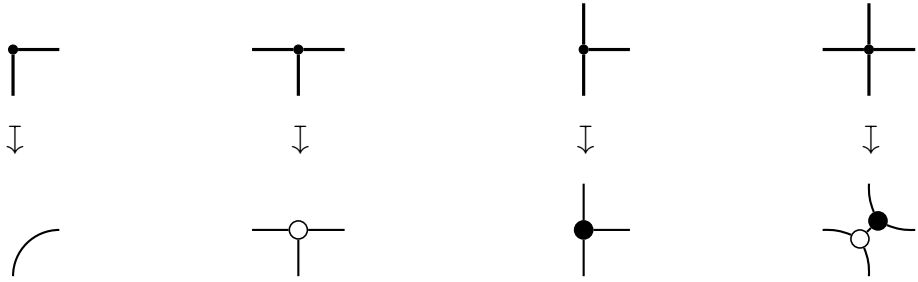
\begin{figure}[ht]
\begin{center}
\begin{tikzpicture}[baseline=(current bounding box.center),scale=0.8]
\tikzstyle{out1}=[inner sep=0,minimum size=2.4mm,circle,draw=black,fill=black,semithick]
\tikzstyle{in1}=[inner sep=0,minimum size=2.4mm,circle,draw=black,fill=white,semithick]
\tikzstyle{hookvertex}=[inner sep=0,minimum size=1.2mm,circle,draw=black,fill=black]
\pgfmathsetmacro{\unit}{1.5};
\pgfmathsetmacro{\side}{0.8};
\node[hookvertex](uc)at(0,\unit){};
\node[inner sep=0](un)at(0,\unit+\side){};
\node[inner sep=0](ue)at(\side,\unit){};
\node[inner sep=0](us)at(0,\unit-\side){};
\node[inner sep=0](uw)at(-\side,\unit){};
\node[inner sep=0](dc)at(0,-\unit){};
\node[inner sep=0](dn)at(0,-\unit+\side){};
\node[inner sep=0](de)at(\side,-\unit){};
\node[inner sep=0](ds)at(0,-\unit-\side){};
\node[inner sep=0](dw)at(-\side,-\unit){};
\node[inner sep=0]at(0,0){\rotatebox[origin=c]{-90}{$\mapsto$}};
\path[very thick](uc)edge(ue) (uc)edge(us);
\path[thick](de)edge[bend right=45](ds);
\end{tikzpicture}\qquad\qquad\qquad
\begin{tikzpicture}[baseline=(current bounding box.center),scale=0.8]
\tikzstyle{out1}=[inner sep=0,minimum size=2.4mm,circle,draw=black,fill=black,semithick]
\tikzstyle{in1}=[inner sep=0,minimum size=2.4mm,circle,draw=black,fill=white,semithick]
\tikzstyle{hookvertex}=[inner sep=0,minimum size=1.2mm,circle,draw=black,fill=black]
\pgfmathsetmacro{\unit}{1.5};
\pgfmathsetmacro{\side}{0.8};
\node[hookvertex](uc)at(0,\unit){};
\node[inner sep=0](un)at(0,\unit+\side){};
\node[inner sep=0](ue)at(\side,\unit){};
\node[inner sep=0](us)at(0,\unit-\side){};
\node[inner sep=0](uw)at(-\side,\unit){};
\node[in1](dc)at(0,-\unit){};
\node[inner sep=0](dn)at(0,-\unit+\side){};
\node[inner sep=0](de)at(\side,-\unit){};
\node[inner sep=0](ds)at(0,-\unit-\side){};
\node[inner sep=0](dw)at(-\side,-\unit){};
\node[inner sep=0]at(0,0){\rotatebox[origin=c]{-90}{$\mapsto$}};
\path[very thick](uc)edge(ue) (uc)edge(us) (uc)edge(uw);
\path[thick](dc)edge(de) (dc)edge(ds) (dc)edge(dw);
\end{tikzpicture}\qquad\qquad\qquad
\begin{tikzpicture}[baseline=(current bounding box.center),scale=0.8]
\tikzstyle{out1}=[inner sep=0,minimum size=2.4mm,circle,draw=black,fill=black,semithick]
\tikzstyle{in1}=[inner sep=0,minimum size=2.4mm,circle,draw=black,fill=white,semithick]
\tikzstyle{hookvertex}=[inner sep=0,minimum size=1.2mm,circle,draw=black,fill=black]
\pgfmathsetmacro{\unit}{1.5};
\pgfmathsetmacro{\side}{0.8};
\node[hookvertex](uc)at(0,\unit){};
\node[inner sep=0](un)at(0,\unit+\side){};
\node[inner sep=0](ue)at(\side,\unit){};
\node[inner sep=0](us)at(0,\unit-\side){};
\node[inner sep=0](uw)at(-\side,\unit){};
\node[out1](dc)at(0,-\unit){};
\node[inner sep=0](dn)at(0,-\unit+\side){};
\node[inner sep=0](de)at(\side,-\unit){};
\node[inner sep=0](ds)at(0,-\unit-\side){};
\node[inner sep=0](dw)at(-\side,-\unit){};
\node[inner sep=0]at(0,0){\rotatebox[origin=c]{-90}{$\mapsto$}};
\path[very thick](uc)edge(ue) (uc)edge(us) (uc)edge(un);
\path[thick](dc)edge(dn) (dc)edge(de) (dc)edge(ds);
\end{tikzpicture}\qquad\qquad\qquad
\begin{tikzpicture}[baseline=(current bounding box.center),scale=0.8]
\tikzstyle{out1}=[inner sep=0,minimum size=2.4mm,circle,draw=black,fill=black,semithick]
\tikzstyle{in1}=[inner sep=0,minimum size=2.4mm,circle,draw=black,fill=white,semithick]
\tikzstyle{hookvertex}=[inner sep=0,minimum size=1.2mm,circle,draw=black,fill=black]
\pgfmathsetmacro{\unit}{1.5};
\pgfmathsetmacro{\side}{0.8};
\pgfmathsetmacro{\shift}{0.15};
\node[hookvertex](uc)at(0,\unit){};
\node[inner sep=0](un)at(0,\unit+\side){};
\node[inner sep=0](ue)at(\side,\unit){};
\node[inner sep=0](us)at(0,\unit-\side){};
\node[inner sep=0](uw)at(-\side,\unit){};
\node[out1](dout1)at(\shift,-\unit+\shift){};
\node[in1](din1)at(-\shift,-\unit-\shift){};
\node[inner sep=0](dn)at(0,-\unit+\side){};
\node[inner sep=0](de)at(\side,-\unit){};
\node[inner sep=0](ds)at(0,-\unit-\side){};
\node[inner sep=0](dw)at(-\side,-\unit){};
\node[inner sep=0]at(0,0){\rotatebox[origin=c]{-90}{$\mapsto$}};
\path[very thick](uc)edge(ue) (uc)edge(us) (uc)edge(uw) (uc)edge(un);
\path[thick](dout1)edge[bend left=12](dn) (dout1)edge[bend right=12](de) (dout1)edge(din1) (din1)edge[bend left=12](ds) (din1)edge[bend right=12](dw);
\end{tikzpicture}
\caption{Local substitutions for going from the hook diagram $H(D)$ to the plabic graph $G(D)$.}
\label{fig:local}
\end{center}
\end{figure}

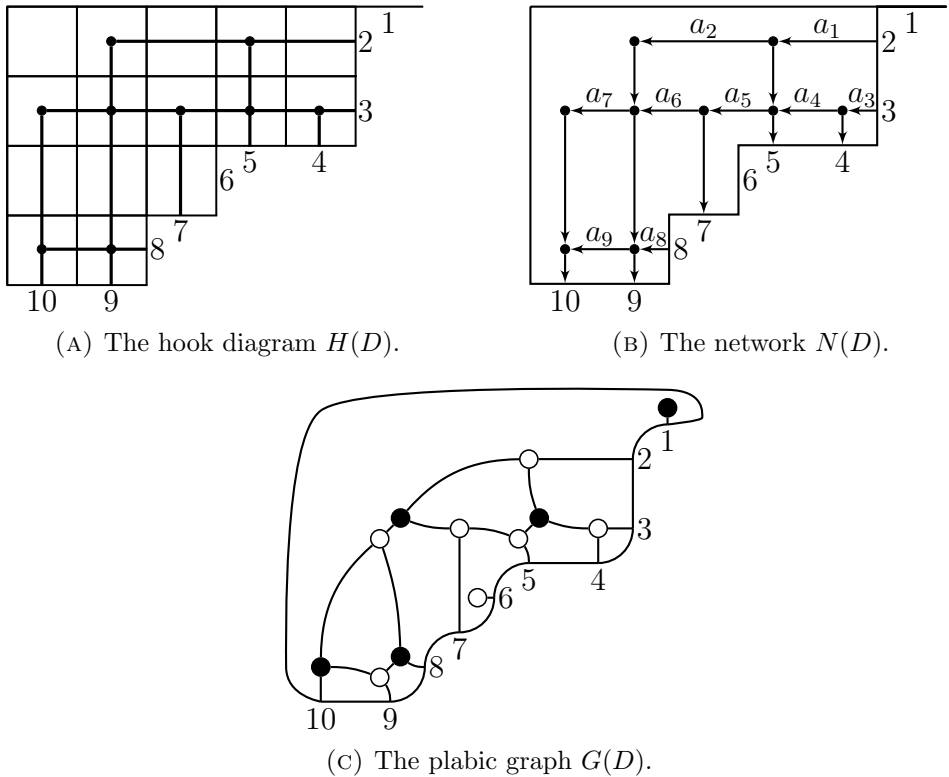
\begin{figure}[ht]
\centering
\subfloat[][The hook diagram $H(D)$.]{
\begin{tikzpicture}[baseline=(current bounding box.center)]
\tikzstyle{out1}=[inner sep=0,minimum size=1.2mm,circle,draw=black,fill=black]
\tikzstyle{in1}=[inner sep=0,minimum size=1.2mm,circle,draw=black,fill=white]
\pgfmathsetmacro{\unit}{0.922};
\useasboundingbox(0,0)rectangle(6.5*\unit,-4.5*\unit);
\coordinate (vstep)at(0,-0.24*\unit);
\coordinate (hstep)at(0.17*\unit,0);
\coordinate (vepsilon)at(0,-0.02*\unit);
\coordinate (hepsilon)at(0.02*\unit,0);
\draw[thick](0,0)--(6*\unit,0) (0,0)--(0,-4*\unit);
\node[inner sep=0]at(0,0){\scalebox{1.6}{\begin{ytableau}
\none \\
\none \\
\none \\
\none \\
\none & \none & \none & \none & \none & & & & & \\
\none & \none & \none & \none & \none & & & & & \\
\none & \none & \none & \none & \none & & & \\
\none & \none & \none & \none & \none & &
\end{ytableau}}};
\node[inner sep=0]at($(5.5*\unit,0)+(vstep)$){$1$};
\node[inner sep=0]at($(5*\unit,-0.5*\unit)+(hstep)$){$2$};
\node[inner sep=0]at($(5*\unit,-1.5*\unit)+(hstep)$){$3$};
\node[inner sep=0]at($(4.5*\unit,-2*\unit)+(vstep)$){$4$};
\node[inner sep=0]at($(3.5*\unit,-2*\unit)+(vstep)$){$5$};
\node[inner sep=0]at($(3*\unit,-2.5*\unit)+(hstep)$){$6$};
\node[inner sep=0]at($(2.5*\unit,-3*\unit)+(vstep)$){$7$};
\node[inner sep=0]at($(2*\unit,-3.5*\unit)+(hstep)$){$8$};
\node[inner sep=0]at($(1.5*\unit,-4*\unit)+(vstep)$){$9$};
\node[inner sep=0]at($(0.5*\unit,-4*\unit)+(vstep)$){$10$};
\node[inner sep=0](b1)at($(5.5*\unit,0)+(vepsilon)$){};
\node[inner sep=0](b2)at($(5*\unit,-0.5*\unit)+(hepsilon)$){};
\node[inner sep=0](b3)at($(5*\unit,-1.5*\unit)+(hepsilon)$){};
\node[inner sep=0](b4)at($(4.5*\unit,-2*\unit)+(vepsilon)$){};
\node[inner sep=0](b5)at($(3.5*\unit,-2*\unit)+(vepsilon)$){};
\node[inner sep=0](b6)at($(3*\unit,-2.5*\unit)+(hepsilon)$){};
\node[inner sep=0](b7)at($(2.5*\unit,-3*\unit)+(vepsilon)$){};
\node[inner sep=0](b8)at($(2*\unit,-3.5*\unit)+(hepsilon)$){};
\node[inner sep=0](b9)at($(1.5*\unit,-4*\unit)+(vepsilon)$){};
\node[inner sep=0](b10)at($(0.5*\unit,-4*\unit)+(vepsilon)$){};
\node[out1](i12)at($(2*\unit,-1*\unit)+(-0.5*\unit,0.5*\unit)$){};
\node[out1](i14)at($(4*\unit,-1*\unit)+(-0.5*\unit,0.5*\unit)$){};
\node[out1](i21)at($(1*\unit,-2*\unit)+(-0.5*\unit,0.5*\unit)$){};
\node[out1](i22)at($(2*\unit,-2*\unit)+(-0.5*\unit,0.5*\unit)$){};
\node[out1](i23)at($(3*\unit,-2*\unit)+(-0.5*\unit,0.5*\unit)$){};
\node[out1](i24)at($(4*\unit,-2*\unit)+(-0.5*\unit,0.5*\unit)$){};
\node[out1](i25)at($(5*\unit,-2*\unit)+(-0.5*\unit,0.5*\unit)$){};
\node[out1](i41)at($(1*\unit,-4*\unit)+(-0.5*\unit,0.5*\unit)$){};
\node[out1](i42)at($(2*\unit,-4*\unit)+(-0.5*\unit,0.5*\unit)$){};
\path[very thick](b2.center)edge(i12) (b3.center)edge(i21) (b8.center)edge(i41) (i25)edge(b4.center) (i14)edge(b5.center) (i23)edge(b7.center) (i12)edge(b9.center) (i21)edge(b10.center);
\end{tikzpicture}}\qquad
\subfloat[][The network $N(D)$.]{
\begin{tikzpicture}[baseline=(current bounding box.center)]
\tikzstyle{out1}=[inner sep=0,minimum size=1.2mm,circle,draw=black,fill=black]
\tikzstyle{in1}=[inner sep=0,minimum size=1.2mm,circle,draw=black,fill=white]
\pgfmathsetmacro{\unit}{0.922};
\useasboundingbox(0,0)rectangle(6.5*\unit,-4.5*\unit);
\coordinate (vstep)at(0,-0.24*\unit);
\coordinate (hstep)at(0.17*\unit,0);
\coordinate (vepsilon)at(0,-0.02*\unit);
\coordinate (hepsilon)at(0.02*\unit,0);
\draw[thick](0,0)--(6*\unit,0)--(5*\unit,0)--(5*\unit,-2*\unit)--(3*\unit,-2*\unit)--(3*\unit,-3*\unit)--(2*\unit,-3*\unit)--(2*\unit,-4*\unit)--(0,-4*\unit)--(0,0);
\node[inner sep=0]at($(5.5*\unit,0)+(vstep)$){$1$};
\node[inner sep=0]at($(5*\unit,-0.5*\unit)+(hstep)$){$2$};
\node[inner sep=0]at($(5*\unit,-1.5*\unit)+(hstep)$){$3$};
\node[inner sep=0]at($(4.5*\unit,-2*\unit)+(vstep)$){$4$};
\node[inner sep=0]at($(3.5*\unit,-2*\unit)+(vstep)$){$5$};
\node[inner sep=0]at($(3*\unit,-2.5*\unit)+(hstep)$){$6$};
\node[inner sep=0]at($(2.5*\unit,-3*\unit)+(vstep)$){$7$};
\node[inner sep=0]at($(2*\unit,-3.5*\unit)+(hstep)$){$8$};
\node[inner sep=0]at($(1.5*\unit,-4*\unit)+(vstep)$){$9$};
\node[inner sep=0]at($(0.5*\unit,-4*\unit)+(vstep)$){$10$};
\node[inner sep=0](b1)at($(5.5*\unit,0)+(vepsilon)$){};
\node[inner sep=0](b2)at($(5*\unit,-0.5*\unit)+(hepsilon)$){};
\node[inner sep=0](b3)at($(5*\unit,-1.5*\unit)+(hepsilon)$){};
\node[inner sep=0](b4)at($(4.5*\unit,-2*\unit)+(vepsilon)$){};
\node[inner sep=0](b5)at($(3.5*\unit,-2*\unit)+(vepsilon)$){};
\node[inner sep=0](b6)at($(3*\unit,-2.5*\unit)+(hepsilon)$){};
\node[inner sep=0](b7)at($(2.5*\unit,-3*\unit)+(vepsilon)$){};
\node[inner sep=0](b8)at($(2*\unit,-3.5*\unit)+(hepsilon)$){};
\node[inner sep=0](b9)at($(1.5*\unit,-4*\unit)+(vepsilon)$){};
\node[inner sep=0](b10)at($(0.5*\unit,-4*\unit)+(vepsilon)$){};
\node[out1](i12)at($(2*\unit,-1*\unit)+(-0.5*\unit,0.5*\unit)$){};
\node[out1](i14)at($(4*\unit,-1*\unit)+(-0.5*\unit,0.5*\unit)$){};
\node[out1](i21)at($(1*\unit,-2*\unit)+(-0.5*\unit,0.5*\unit)$){};
\node[out1](i22)at($(2*\unit,-2*\unit)+(-0.5*\unit,0.5*\unit)$){};
\node[out1](i23)at($(3*\unit,-2*\unit)+(-0.5*\unit,0.5*\unit)$){};
\node[out1](i24)at($(4*\unit,-2*\unit)+(-0.5*\unit,0.5*\unit)$){};
\node[out1](i25)at($(5*\unit,-2*\unit)+(-0.5*\unit,0.5*\unit)$){};
\node[out1](i41)at($(1*\unit,-4*\unit)+(-0.5*\unit,0.5*\unit)$){};
\node[out1](i42)at($(2*\unit,-4*\unit)+(-0.5*\unit,0.5*\unit)$){};
\path[-latex',thick](b2)edge node[above=-3pt]{$a_1$}(i14) (i14)edge node[above=-3pt]{$a_2$}(i12) (b3)edge node[above=-3pt]{$a_3$}(i25) (i25)edge node[above=-3pt]{$a_4$}(i24) (i24)edge node[above=-3pt]{$a_5$}(i23) (i23)edge node[above=-3pt]{$a_6$}(i22) (i22)edge node[above=-3pt]{$a_7$}(i21) (b8)edge node[above=-3pt]{$a_8$}(i42) (i42)edge node[above=-3pt]{$a_9$}(i41) (i25)edge(b4) (i14)edge(i24) (i24)edge(b5) (i23)edge(b7) (i12)edge(i22) (i22)edge(i42) (i42)edge(b9) (i21)edge(i41) (i41)edge(b10);
\end{tikzpicture}}\\
\subfloat[][The plabic graph $G(D)$.]{
\begin{tikzpicture}[baseline=(current bounding box.center)]
\tikzstyle{out1}=[inner sep=0,minimum size=2.4mm,circle,draw=black,fill=black,semithick]
\tikzstyle{in1}=[inner sep=0,minimum size=2.4mm,circle,draw=black,fill=white,semithick]
\pgfmathsetmacro{\unit}{0.922};
\useasboundingbox(-0.5*\unit,0.5*\unit)rectangle(6.5*\unit,-4.5*\unit);
\coordinate (vstep)at(0,-0.24*\unit);
\coordinate (hstep)at(0.17*\unit,0);
\coordinate (vloll)at(0,0.24*\unit);
\coordinate (hloll)at(-0.24*\unit,0);
\coordinate (dstep)at(0.15*\unit,0.15*\unit);
\node[inner sep=0]at($(5.5*\unit,0)+(vstep)$){$1$};
\node[inner sep=0]at($(5*\unit,-0.5*\unit)+(hstep)$){$2$};
\node[inner sep=0]at($(5*\unit,-1.5*\unit)+(hstep)$){$3$};
\node[inner sep=0]at($(4.5*\unit,-2*\unit)+(vstep)$){$4$};
\node[inner sep=0]at($(3.5*\unit,-2*\unit)+(vstep)$){$5$};
\node[inner sep=0]at($(3*\unit,-2.5*\unit)+(hstep)$){$6$};
\node[inner sep=0]at($(2.5*\unit,-3*\unit)+(vstep)$){$7$};
\node[inner sep=0]at($(2*\unit,-3.5*\unit)+(hstep)$){$8$};
\node[inner sep=0]at($(1.5*\unit,-4*\unit)+(vstep)$){$9$};
\node[inner sep=0]at($(0.5*\unit,-4*\unit)+(vstep)$){$10$};
\coordinate(v05h)at(5.5*\unit,0*\unit);
\coordinate(v0h5)at(5*\unit,-0.5*\unit);
\coordinate(v1h5)at(5*\unit,-1.5*\unit);
\coordinate(v24h)at(4.5*\unit,-2*\unit);
\coordinate(v23h)at(3.5*\unit,-2*\unit);
\coordinate(v2h3)at(3*\unit,-2.5*\unit);
\coordinate(v32h)at(2.5*\unit,-3*\unit);
\coordinate(v3h2)at(2*\unit,-3.5*\unit);
\coordinate(v41h)at(1.5*\unit,-4*\unit);
\coordinate(v40h)at(0.5*\unit,-4*\unit);
\node[out1](v0e5h)at($(5.5*\unit,0*\unit)+(vloll)$){};
\node[in1](v0h3h)at(3.5*\unit,-0.5*\unit){};
\node[inner sep=0](v0h2)at(2*\unit,-0.5*\unit){};
\node[inner sep=0](v11h)at(1.5*\unit,-1*\unit){};
\node[in1](v1h4h)at(4.5*\unit,-1.5*\unit){};
\node[out1](v1h3hout1)at($(3.5*\unit,-1.5*\unit)+(dstep)$){};
\node[in1](v1h3hin1)at($(3.5*\unit,-1.5*\unit)-(dstep)$){};
\node[in1](v1h2h)at(2.5*\unit,-1.5*\unit){};
\node[out1](v1h1hout1)at($(1.5*\unit,-1.5*\unit)+(dstep)$){};
\node[in1](v1h1hin1)at($(1.5*\unit,-1.5*\unit)-(dstep)$){};
\node[inner sep=0](v1h1)at(1*\unit,-1.5*\unit){};
\node[inner sep=0](v20h)at(0.5*\unit,-2*\unit){};
\node[in1](v2h3e)at($(3*\unit,-2.5*\unit)+(hloll)$){};
\node[out1](v3h1hout1)at($(1.5*\unit,-3.5*\unit)+(dstep)$){};
\node[in1](v3h1hin1)at($(1.5*\unit,-3.5*\unit)-(dstep)$){};
\node[out1](v3h0h)at(0.5*\unit,-3.5*\unit){};
\node[inner sep=0](border0)at(6*\unit,0.1*\unit){};
\node[inner sep=0](border1)at(5.5*\unit,0.5*\unit){};
\node[inner sep=0](border2)at(0.5*\unit,0.2*\unit){};
\node[inner sep=0](border3)at(0,-3.5*\unit){};
\path[thick](v05h)edge[bend right=45](v0h5) (v0h5)edge(v1h5) (v1h5)edge[bend left=45](v24h) (v24h)edge(v23h) (v23h)edge[bend right=45](v2h3) (v2h3)edge[bend left=45](v32h) (v32h)edge[bend right=45](v3h2) (v3h2)edge[bend left=45](v41h) (v41h)edge(v40h) (v05h)edge(v0e5h) (v0h5)edge(v0h3h) (v0h3h)edge[bend right=24](v1h1hout1) (v0h3h)edge[bend right=10](v1h3hout1) (v1h5)edge(v1h4h) (v1h4h)edge(v24h) (v1h4h)edge[bend left=10](v1h3hout1) (v1h3hout1)edge(v1h3hin1) (v1h3hin1)edge[bend left=16](v23h) (v1h3hin1)edge[bend right=10](v1h2h) (v1h2h)edge[bend left=10](v1h1hout1) (v1h2h)edge(v32h) (v1h1hout1)edge(v1h1hin1) (v1h1hin1)edge[bend right=24](v3h0h) (v1h1hin1)edge[bend left=8](v3h1hout1) (v2h3)edge(v2h3e) (v3h2)edge[bend left=16](v3h1hout1) (v3h1hout1)edge(v3h1hin1) (v3h1hin1)edge[bend right=10](v3h0h) (v3h1hin1)edge[bend left=16](v41h) (v3h0h)edge(v40h);
\draw[thick]plot[smooth,tension=0.36]coordinates{(v05h) (border0) (border1) (border2)  (border3) (v40h)};
\end{tikzpicture}}
\caption{The hook diagram, network,
 and plabic graph associated to the \Le -diagram $D$ from \cref{fig:Le}.}
\label{fig:plabic}
\end{figure}

The following construction of Postnikov associates a hook diagram, network, and plabic graph to any \Le -diagram.
\begin{defn}[{\cite[Sections 6 and 20]{postnikov}}]\label{def:Le-plabic}
Let $D$ be a \Le-diagram of type $(k,n)$. We define the {\itshape hook diagram} $H(D)$ of $D$, the {\itshape network} $N(D)$ of $D$, and the {\itshape plabic graph} $G(D)$ of $D$. (See \cref{fig:plabic} for examples.)

To construct $H(D)$, we delete the $0$'s of $D$, and replace each $+$ with a vertex. From each vertex we construct a hook which goes east and south, to the border of the Young diagram of $D$. We label the edges of the southeast border of $D$ by $1, \dots, n$ from northeast to southwest.

To construct $N(D)$, we direct the edges of $H(D)$ west and south. Let $E$ be the set of horizontal edges of $N(D)$. To each such edge $e\in E$, we associate a variable $a_e$.

To construct $G(D)$ from $H(D)$, we place boundary vertices $1, \dots, n$ along the southeast border. Then we replace the local region around each internal vertex as in \cref{fig:local}, and add a black (respectively, white) lollipop for each black (respectively, white) fixed point of the decorated permutation $\pi_D$.
\end{defn}

More generally, each \Le-diagram $D$ 
is associated with a family of \emph{reduced plabic graphs}
consisting of $G(D)$ together with other plabic graphs which can be obtained
from $G(D)$ by certain \emph{moves}; see \cite[Section 12]{postnikov}.

From the plabic graph constructed in \cref{def:Le-plabic} (and more generally from a reduced plabic graph $G$), we can read off the corresponding decorated permutation $\pi_G$ as follows.
\begin{defn}\label{def:rules}
Let $G$ be a reduced plabic graph as above with boundary vertices $1,\dots, n$. For each boundary vertex $i\in [n]$, we follow a path along the edges of $G$ starting at $i$, turning (maximally) right at every internal black vertex, and (maximally) left at every internal white vertex. This path ends at some boundary vertex $\pi(i)$. By \cite[Section 13]{postnikov}, the fact that $G$ is reduced implies that each fixed point of $\pi$ is attached to a lollipop; we color each fixed point by the color of its lollipop. In this way we obtain the \emph{decorated permutation} $\pi_G = \pi$ of $G$. We say that $G$ is of {\itshape type} $(k,n)$, where $k$ is the number of anti-excedances of $\pi_G$. By \cite[Definition 11.5]{postnikov}\footnote{There is a typo in Postnikov's formula: $k+(n-k)$ should read $k-(n-k)$. Our formula looks different than his, but is equivalent.}, we have
\begin{align}\label{k-statistic}
k = \#\text{edges} - \#\text{black vertices} - \sum_{\text{white vertices }v}(\deg(v)-1).
\end{align}
\end{defn}
We invite the reader to verify that the decorated permutation of the plabic graph $G$ of \cref{fig:plabic} is $\pi_G = (\underline{1},5,4,9,7,\overline{6},2,10,3,8)$. By \eqref{k-statistic}, we calculate $k = 21 - 5 - 12 = 4$.

We now explain how to parameterize elements of the cell $S_D$ associated to a \Le -diagram $D$ from its network $N(D)$.
\begin{thm}[{\cite[Section 6]{postnikov}}]\label{network_param}
Let $D$ be a \Le -diagram of type $(k,n)$, with associated cell $S_D\subseteq\Gr_{k,n}^{\ge 0}$ and network $N(D)$ from \cref{def:Le-plabic}. Let $E$ denote the set of horizontal edges of $N(D)$. We obtain a parameterization $\mathbb{R}_{>0}^E\to S_D$ of $S_D$, as follows. Let $s_1 < \dots < s_k$ be the labels of the vertical edges of the southeast border of $N(D)$. Also, to any directed path $p$ of $N(D)$, we define its {\itshape weight} $w_p$ to be the product of the edge variables $a_e$ over all horizontal edges $e$ in $p$. We let $A$ be the $k\times n$ matrix with rows indexed by $\{s_1,\dots, s_k\}$ whose $(s_i,j)$-entry equals
$$
(-1)^{|\{i'\in [k]: s_i < s_{i'} < j\}|}\sum_{p:s_i\to j}w_p,
$$
where the sum is over all directed paths $p$ from $s_i$ to  $j$. Then the map sending $(a_e)_{e\in E}\in\mathbb{R}_{>0}^E$ to the element of $\Gr_{k,n}^{\ge 0}$ represented by $A$ is a homeomorphism from 
$\R_{>0}^E$ to  $S_D$.
\end{thm}
For example, the network in \cref{fig:plabic} gives the parameterization
$$
\kbordermatrix{
& 1 & 2 & 3 & 4 & 5 & 6 & 7 & 8 & 9 & 10 \cr
2 & 0 & 1 & 0 & 0 & -a_1 & 0 & a_1a_5 & 0 & -a_1(a_2+a_5a_6) & -a_1(a_2+a_5a_6)(a_7+a_9) \cr
3 & 0 & 0 & 1 & a_3 & a_3a_4 & 0 & -a_3a_4a_5 & 0 & a_3a_4a_5a_6 & a_3a_4a_5a_6(a_7+a_9) \cr
6 & 0 & 0 & 0 & 0 & 0 & 1 & 0 & 0 & 0 & 0 \cr
8 & 0 & 0 & 0 & 0 & 0 & 0 & 0 & 1 & a_8 & a_8a_9}
$$
of $S_D\subseteq\Gr_{4,10}^{\ge 0}$, where $D$ is the \Le -diagram in \cref{fig:Le}.

\section{Background on sign variation}

\noindent In this section we provide some background on sign variation,
and state \cref{criterion}, which will be useful later for proving that two cells have disjoint image in the amplituhedron.

\begin{defn}\label{def:sign_vector}
A \emph{sign vector} is a vector with entries in $\{0, +, -\}$. For a vector $v\in\mathbb{R}^n$, we let $\sign(v)\in\{0,+,-\}^n$ be the sign vector naturally associated to $v$. For example, $\sign(1,0,-4,2) = (+,0,-,+)$.
\end{defn}

\begin{defn}\label{def:var}
Given $v\in\mathbb{R}^n$, we let $\var(v)$ be the number of times $v$ changes sign, when viewed as a sequence of $n$ numbers and ignoring any zeros. Also let
$$
\varbar(v) := \max\{\var(w) : \text{$w\in\mathbb{R}^n$ such that $w_i = v_i$ for all $i\in [n]$ with $v_i\neq 0$}\},
$$
i.e.\ $\varbar(v)$ is the maximum number of times $v$ changes sign after we choose a sign for each zero coordinate. Note that we can apply $\var$ and $\overline{\var}$ to sign vectors. For example, $\var(+,0,0,+,-) = 1$ and $\varbar(+,0,0,+,-) = 3$.
\end{defn}

The following result of Gantmakher and Krein characterizes total positivity in $\Gr_{k,n}$ using sign variation.
\begin{thm}[{\cite[Theorems V.3, V.7, V.1, V.6]{gantmakher_krein_50}}]\label{thm:G-K}
Let $V\in\Gr_{k,n}$. \\
(i) $V\in\Gr_{k,n}^{\ge 0}\iff\var(v)\le k-1\text{ for all }v\in V\setminus\{0\}\iff\overline{\var}(w)\ge k\text{ for all }w\in V^\perp\setminus\{0\}$. \\
(ii) $V\in\Gr_{k,n}^{>0}\iff\overline{\var}(v)\le k-1\text{ for all }v\in V\setminus\{0\}\iff\var(w)\ge k\text{ for all }w\in V^\perp\setminus\{0\}$.
\end{thm}

The following lemma gives a sign-variational characterization for when two elements of $\Gr_{k,n}^{\ge 0}$ correspond to the same point of the amplituhedron $\mathcal{A}_{n,k,m}(Z)$.  We will apply it repeatedly in \cref{sec:k=2-disjointness} to show that the images of the BCFW cells are disjoint in the $k=2, m=4$ amplituhedron $\mathcal{A}_{n,2,4}(Z)$.
\begin{lem}\label{criterion}
Let $Z\in \Mat_{k+m,n}^{>0}$, where $k,m,n\ge 0$ satisfy $k+m \le n$, and $V,V'\in\Gr_{k,n}^{\ge 0}$. Then $\tilde{Z}(V) = \tilde{Z}(V')$ if and only if for all $v\in V$, there exists a unique $v'\in V'$ such that ${Z}(v) = {Z}(v')$. We call $v'$ the \emph{matching vector} for $v$. Note that by \cref{thm:G-K}(ii), we either have $v = v'$ or $\var(v-v')\ge k+m$.
\end{lem}

\begin{pf}
($\Rightarrow$): Suppose that $\tilde{Z}(V) = \tilde{Z}(V')$, i.e.\ $\{Z(v) : v\in V\} = \{Z(v') : v'\in V'\}$. Then for any $v\in V$, there exists $v'\in V'$ with $Z(v) = Z(v')$. Since $\dim(\tilde{Z}(V')) = k$ (see \cref{def:amp}), we have $\ker(Z)\cap V' = \{0\}$, which implies that $v'$ is unique.

($\Leftarrow$): 
Suppose that for all $v\in V$, there exists a unique $v'\in V'$ such that 
$Z(v) = Z(v')$.  Then 
$\tilde Z(V) \subseteq \tilde Z(V')$, and since 
$\dim(\tilde Z(V)) = \dim(\tilde Z(V'))=k$, we get
$\tilde Z(V) = \tilde Z(V')$.
\end{pf}

\section{Warmup: the \texorpdfstring{$m=2$}{m=2} amplituhedron}\label{sec:m2}

\noindent In this section we focus on the $m=2$ amplituhedron
 $\mathcal{A}_{n,k,2}(Z)$, as a warmup to our study of the $m=4$ amplituhedron in the following sections.
In analogy to the $m=4$ case, the outline of this section is as follows.
\begin{enumerate}[label=\arabic*.]
\item In \cref{sec:m=2_BCFW}, we give a recurrence on plabic graphs, which (after a `shift by $1$') produces a collection $\mathcal{C}_{n,k,2}$ of $2k$-dimensional cells of $\Gr_{k,n}^{\ge 0}$, whose images should triangulate $\mathcal{A}_{n,k,2}(Z)$. (This is analogous to the $m=4$ BCFW recursion; see \cref{sec:BCFW}.)
\item In \cref{sec:m=2_diagrams}, we describe the \Le -diagrams of these $2k$-dimensional cells, indexed by lattice paths inside a $k\times (n-k-2)$ rectangle. (This is analogous to our description for $m=4$ in terms of pairs of noncrossing lattice paths inside a $k\times (n-k-4)$ rectangle; see \cref{sec:lattice-paths} and \cref{sec:bijection}.)
\item In \cref{sec:m=2_bases}, we find a nice basis of any subspace $V$ coming from our cells $\mathcal{C}_{n,k,2}$. (This is analogous to our `domino bases' for $m=4$; see \cref{sec:k=2-bases} for the case $k=2$, and \cref{conj:domino} for a conjectural generalization to all $k$.)
\item In \cref{sec:m=2_disjointness}, we prove that the images of the cells $\mathcal{C}_{n,k,2}$ in $\mathcal{A}_{n,k,2}(Z)$ are disjoint, using sign variation arguments. (In \cref{sec:k=2-disjointness}, we will prove disjointness for $m=4$ BCFW cells when $k=2$, using more intricate arguments along the same lines.) 
\end{enumerate}

While we will not need the results of this section to handle the $m=4$ case, the ideas and techniques used here will hopefully give the reader a flavor of our arguments for $m=4$. We will also cite the $m=2$ case as evidence for our conjecture in \cref{sec:numerology}.

\begin{rmk}\label{m=2_ATT}
Arkani-Hamed, Thomas, and Trnka \cite[Section 7]{ATT} consider the same collection of cells $\mathcal{C}_{n,k,2}$, up to a cyclic shift. They define their cells in terms of bases as in item 3 above, and do not consider a recurrence on plabic graphs or the \Le -diagrams. They show both that the images of these cells in $\mathcal{A}_{n,k,2}(Z)$ are disjoint, and that their union covers a dense subset of $\mathcal{A}_{n,k,2}(Z)$ (we will not prove the latter fact here). We remark that their arguments also employ sign variation.
\end{rmk}

\subsection{A BCFW-style recursion for \texorpdfstring{$m=2$}{m=2}}\label{sec:m=2_BCFW}

We start by giving a recursion on plabic graphs which 
produces, for each pair $(k,n)$ with $1 \leq k \leq n-1$,
a collection $\mathcal{C}_{n,k,2}$ of $2k$-dimensional positroid cells 
of $\Gr_{k,n}^{\geq 0}$.  This recursion is an analogue of the well-known \emph{BCFW recursion}, which produces the
BCFW cells $\mathcal{C}_{n,k,4}$
for the $m=4$ amplituhedron (see \cref{sec:BCFW}).

\begin{defn}
We define a local operation on plabic graphs which we call
\emph{splitting}. Let $G$ be a plabic graph of type 
$(k,n)$ (see \cref{def:rules}), and $i\in [n]$ a boundary vertex incident to a unique edge $e$. We {\itshape split $e$ at $i$} by locally replacing $e$ by a vertex $v$ with three incident edges as follows, while leaving the rest of the graph unchanged:
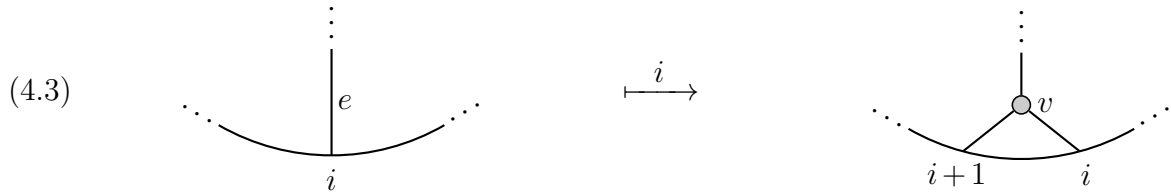
\begin{align}\label{split}
\quad\begin{tikzpicture}[baseline=(current bounding box.center)]
\tikzstyle{out1}=[inner sep=0,minimum size=2.4mm,circle,draw=black,fill=black,semithick]
\tikzstyle{in1}=[inner sep=0,minimum size=2.4mm,circle,draw=black,fill=white,semithick]
\tikzstyle{ambiguous}=[inner sep=0,minimum size=2.4mm,circle,draw=black,fill=black!20,semithick]
\pgfmathsetmacro{\radius}{3.0};
\pgfmathsetmacro{\l}{0.72};
\draw[thick](-120:\radius)arc(-120:-60:\radius);
\node[inner sep=0,rotate=150]at($(-120:\radius)+(150:0.36)$){$\cdots$};
\node[inner sep=0,rotate=30]at($(-60:\radius)+(30:0.36)$){$\cdots$};
\node[inner sep=0](i)at(-90:\radius){};
\node at(-90:\radius+0.32){$i$};
\node[inner sep=0](top)at(-90:\radius-2*\l){};
\node[inner sep=0,rotate=90]at($(top)+(90:0.36)+(0.013,0)$){$\cdots$};
\path[thick](top)edge node[right=-2pt]{$e$}(i.center);
\end{tikzpicture}\qquad\qquad\xmapsto{\mbox{\;\, $i$ \;\,}}\qquad\qquad
\quad\begin{tikzpicture}[baseline=(current bounding box.center)]
\tikzstyle{out1}=[inner sep=0,minimum size=2.4mm,circle,draw=black,fill=black,semithick]
\tikzstyle{in1}=[inner sep=0,minimum size=2.4mm,circle,draw=black,fill=white,semithick]
\tikzstyle{ambiguous}=[inner sep=0,minimum size=2.4mm,circle,draw=black,fill=black!20,semithick]
\pgfmathsetmacro{\radius}{3.0};
\pgfmathsetmacro{\l}{0.72};
\draw[thick](-120:\radius)arc(-120:-60:\radius);
\node[inner sep=0,rotate=150]at($(-120:\radius)+(150:0.36)$){$\cdots$};
\node[inner sep=0,rotate=30]at($(-60:\radius)+(30:0.36)$){$\cdots$};
\node[inner sep=0](i)at(-75:\radius){};
\node at(-75:\radius+0.32){$i$};
\node[inner sep=0](i+1)at(-105:\radius){};
\node at(-105:\radius+0.32){$i\hspace*{-1pt}+\hspace*{-1pt}1$};
\node[inner sep=0](top)at(-90:\radius-2*\l){};
\node[inner sep=0,rotate=90]at($(top)+(90:0.36)+(0.013,0)$){$\cdots$};
\node[ambiguous](v)at(-90:\radius-\l){};
\node[inner sep=0]at($(v)+(0:0.32)$){$v$};
\path[thick](v)edge(i.center) edge(i+1.center) edge(top);
\end{tikzpicture}
\end{align}
(The color of $v$ may be either black or white.)
This produces a plabic graph $G'$ of type 
$(k',n+1)$, whose boundary vertices are labeled by $[n+1]$, so that 
$i$ and $i+1$ are labeled as above
(i.e.\ labels $1,\dots,i-1$ from $G$ remain unchanged,
and labels $i+1,\dots,n$ from $G$ get increased by $1$).
Note that by \eqref{k-statistic}, $k' = k$ if 
$v$ is white, while $k' = k+1$ if $v$ is black.
\end{defn}

\begin{defn}[Recursion for cells when $m=2$]\label{def:BCFWm2}
For positive integers $k$ and $n$ such that $1 \leq k \leq n-1$, 
we recursively define a set $\mathcal{\tilde{G}}_{n,k,2}$ of plabic graphs of 
type $(k,n)$ as follows.\footnote{In our notation $\mathcal{\tilde{G}}_{n,k,2}$, the subscript $2$ is to remind us that these graphs correspond to the $m=2$ amplituhedron, and the tilde is to remind us that these graphs do {\bfseries not} directly label the cells $\mathcal{C}_{n,k,2}$. Rather, we must first `shift by $1$'; see \cref{def:m-is-2-permutations}.}
\begin{enumerate}
\item If $n=2$ and $k=1$, then $\mathcal{\tilde{G}}_{n,k,2}$ contains a unique element, the plabic graph
$\begin{tikzpicture}[baseline=(current bounding box.center)]
\tikzstyle{out1}=[inner sep=0,minimum size=2.4mm,circle,draw=black,fill=black,semithick]
\tikzstyle{in1}=[inner sep=0,minimum size=2.4mm,circle,draw=black,fill=white,semithick]
\pgfmathsetmacro{\radius}{0.48};
\draw[thick](0,0)circle[radius=\radius];
\foreach \x in {1,...,2}{
\node[inner sep=0](b\x)at(180-\x*180:\radius){};
\node at(180-\x*180:\radius+.24){$\x$};}
\path[thick](b1.center)edge(b2.center);
\end{tikzpicture}$.
\item For $n\ge 3$, 
$\mathcal{\tilde{G}}_{n,k,2}$ 
is the set of all plabic graphs obtained either by splitting a plabic graph in $\mathcal{\tilde{G}}_{n-1,k,2}$ 
at $n-1$ with $v$ white,
or by splitting a plabic graph in 
$\mathcal{\tilde{G}}_{n-1,k-1,2}$ at $n-1$ with $v$ black. (See \cref{fig:BCFW-graphs-m=2}.)
\end{enumerate}
\end{defn}
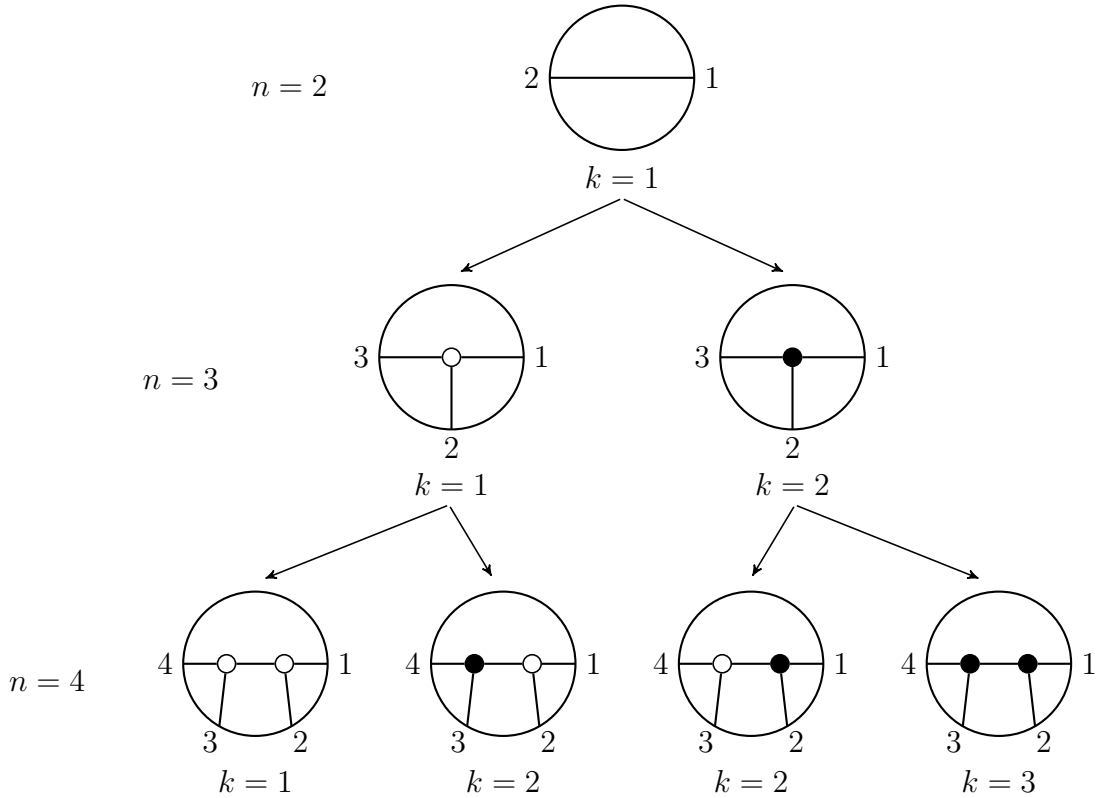
\begin{figure}[ht]
\begin{center}
$n=2$\qquad\qquad\qquad\begin{tikzpicture}[baseline=(current bounding box.center)]
\tikzstyle{out1}=[inner sep=0,minimum size=2.4mm,circle,draw=black,fill=black,semithick]
\tikzstyle{in1}=[inner sep=0,minimum size=2.4mm,circle,draw=black,fill=white,semithick]\
\pgfmathsetmacro{\radius}{0.96};
\draw[thick](0,0)circle[radius=\radius];
\foreach \x in {1,...,2}{
\node[inner sep=0](b\x)at(180-\x*180:\radius){};
\node at(180-\x*180:\radius+.24){$\x$};}\
\path[thick](b1.center)edge(b2.center);
\node[inner sep=0]at(0,-\radius-.36){$k=1$};
\end{tikzpicture}\qquad\qquad\qquad\phantom{$n=2$}
\end{center}
\begin{center}
\begin{tikzpicture}[baseline=(current bounding box.center),
  to/.style={->,>=stealth',shorten >=1pt,semithick}]
  \node[inner sep=0](in1)at(0,1){};
  \node[inner sep=0](out1)at(-2.2,0){};
  \node[inner sep=0](out2)at(2.2,0){};
  \draw[to](in1)edge(out1);
  \draw[to](in1)edge(out2);
\end{tikzpicture}
\end{center}
\begin{center}
$n=3$\qquad\qquad\begin{tikzpicture}[baseline=(current bounding box.center)]
\tikzstyle{out1}=[inner sep=0,minimum size=2.4mm,circle,draw=black,fill=black,semithick]
\tikzstyle{in1}=[inner sep=0,minimum size=2.4mm,circle,draw=black,fill=white,semithick]
\pgfmathsetmacro{\radius}{0.96};
\draw[thick](0,0)circle[radius=\radius];
\foreach \x in {1,...,3}{
\node[inner sep=0](b\x)at(90-\x*90:\radius){};
\node at(90-\x*90:\radius+.24){$\x$};}
\node[in1](c1)at(0,0){};
\path[thick](b1.center)edge(c1) (b2.center)edge(c1) (b3.center)edge(c1);
\node[inner sep=0]at(0,-\radius-.72){$k=1$};
\end{tikzpicture}\qquad\qquad\begin{tikzpicture}[baseline=(current bounding box.center)]
\tikzstyle{out1}=[inner sep=0,minimum size=2.4mm,circle,draw=black,fill=black,semithick]
\tikzstyle{in1}=[inner sep=0,minimum size=2.4mm,circle,draw=black,fill=white,semithick]
\pgfmathsetmacro{\radius}{0.96};
\draw[thick](0,0)circle[radius=\radius];
\foreach \x in {1,...,3}{
\node[inner sep=0](b\x)at(90-\x*90:\radius){};
\node at(90-\x*90:\radius+.24){$\x$};}
\node[out1](c1)at(0,0){};
\path[thick](b1.center)edge(c1) (b2.center)edge(c1) (b3.center)edge(c1);
\node[inner sep=0]at(0,-\radius-.72){$k=2$};
\end{tikzpicture}\qquad\qquad\phantom{$n=3$}
\end{center}
\begin{tikzpicture}[baseline=(current bounding box.center),to/.style={->,>=stealth',shorten >=1pt,semithick}]
  \node[inner sep=0](in1)at(-2.3,1){};
  \node[inner sep=0](in2)at(2.3,1){};
  \node[inner sep=0](out1)at(-4.8,0){};
  \node[inner sep=0](out2)at(-1.7,0){};
  \node[inner sep=0](out4)at(1.7,0){};
  \node[inner sep=0](out5)at(4.8,0){};
  \path[to](in1)edge(out1);
  \path[to](in1)edge(out2);
  \path[to](in2)edge(out5);
  \path[to](in2)edge(out4);
\end{tikzpicture}
\begin{center}
\makebox[\textwidth]{
$n=4$\qquad\begin{tikzpicture}[baseline=(current bounding box.center)]
\tikzstyle{out1}=[inner sep=0,minimum size=2.4mm,circle,draw=black,fill=black,semithick]
\tikzstyle{in1}=[inner sep=0,minimum size=2.4mm,circle,draw=black,fill=white,semithick]
\pgfmathsetmacro{\radius}{0.96};
\draw[thick](0,0)circle[radius=\radius];
\foreach \x in {1,...,4}{
\node[inner sep=0](b\x)at(60-\x*60:\radius){};
\node at(60-\x*60:\radius+.24){$\x$};}
\node[in1](c1)at(\radius*0.4,0){};
\node[in1](c2)at(-\radius*0.4,0){};
\path[thick](b1.center)edge(c1) (b2.center)edge(c1) (b3.center)edge(c2) (b4.center)edge(c2) (c1)edge(c2);
\node[inner sep=0]at(0,-\radius-.60){$k=1$};
\end{tikzpicture}\quad
\begin{tikzpicture}[baseline=(current bounding box.center)]
\tikzstyle{out1}=[inner sep=0,minimum size=2.4mm,circle,draw=black,fill=black,semithick]
\tikzstyle{in1}=[inner sep=0,minimum size=2.4mm,circle,draw=black,fill=white,semithick]
\pgfmathsetmacro{\radius}{0.96};
\draw[thick](0,0)circle[radius=\radius];
\foreach \x in {1,...,4}{
\node[inner sep=0](b\x)at(60-\x*60:\radius){};
\node at(60-\x*60:\radius+.24){$\x$};}
\node[in1](c1)at(\radius*0.4,0){};
\node[out1](c2)at(-\radius*0.4,0){};
\path[thick](b1.center)edge(c1) (b2.center)edge(c1) (b3.center)edge(c2) (b4.center)edge(c2) (c1)edge(c2);
\node[inner sep=0]at(0,-\radius-.60){$k=2$};
\end{tikzpicture}\quad
\begin{tikzpicture}[baseline=(current bounding box.center)]
\tikzstyle{out1}=[inner sep=0,minimum size=2.4mm,circle,draw=black,fill=black,semithick]
\tikzstyle{in1}=[inner sep=0,minimum size=2.4mm,circle,draw=black,fill=white,semithick]
\pgfmathsetmacro{\radius}{0.96};
\draw[thick](0,0)circle[radius=\radius];
\foreach \x in {1,...,4}{
\node[inner sep=0](b\x)at(60-\x*60:\radius){};
\node at(60-\x*60:\radius+.24){$\x$};}
\node[out1](c1)at(\radius*0.4,0){};
\node[in1](c2)at(-\radius*0.4,0){};
\path[thick](b1.center)edge(c1) (b2.center)edge(c1) (b3.center)edge(c2) (b4.center)edge(c2) (c1)edge(c2);
\node[inner sep=0]at(0,-\radius-.60){$k=2$};
\end{tikzpicture}\quad
\begin{tikzpicture}[baseline=(current bounding box.center)]
\tikzstyle{out1}=[inner sep=0,minimum size=2.4mm,circle,draw=black,fill=black,semithick]
\tikzstyle{in1}=[inner sep=0,minimum size=2.4mm,circle,draw=black,fill=white,semithick]
\pgfmathsetmacro{\radius}{0.96};
\draw[thick](0,0)circle[radius=\radius];
\foreach \x in {1,...,4}{
\node[inner sep=0](b\x)at(60-\x*60:\radius){};
\node at(60-\x*60:\radius+.24){$\x$};}
\node[out1](c1)at(\radius*0.4,0){};
\node[out1](c2)at(-\radius*0.4,0){};
\path[thick](b1.center)edge(c1) (b2.center)edge(c1) (b3.center)edge(c2) (b4.center)edge(c2) (c1)edge(c2);
\node[inner sep=0]at(0,-\radius-.60){$k=3$};
\end{tikzpicture}\qquad\phantom{$n=4$}}
\end{center}
\caption{The first two steps of the $m=2$ recursion from \cref{def:BCFWm2}, giving the graphs $\mathcal{\tilde{G}}_{n,k,2}$ for $n\le 4$.}
\label{fig:BCFW-graphs-m=2}
\end{figure}

\begin{figure}[ht]
\begin{tikzpicture}[baseline=(current bounding box.center)]
\tikzstyle{out1}=[inner sep=0,minimum size=2.4mm,circle,draw=black,fill=black,semithick]
\tikzstyle{in1}=[inner sep=0,minimum size=2.4mm,circle,draw=black,fill=white,semithick]
\tikzstyle{ambiguous}=[inner sep=0,minimum size=2.4mm,circle,draw=black,fill=black!20,semithick]
\pgfmathsetmacro{\radius}{2.40};
\pgfmathsetmacro{\s}{0.60};
\pgfmathsetmacro{\angle}{16};
\draw[thick](-45:\radius)arc(-45:225:\radius);
\foreach \x in {1,...,3}{
\node[inner sep=0](b\x)at(\angle-\x*\angle:\radius){};
\node at(\angle-\x*\angle:\radius+.24){$\x$};}
\foreach \x in {4,...,6}{
\node[inner sep=0](b\x)at(180+6*\angle-\x*\angle:\radius){};}
\node[inner sep=0]at(180+2*\angle:\radius+.48){$n\hspace*{-1pt}-\hspace*{-1pt}2$};
\node[inner sep=0]at(180+\angle:\radius+.48){$n\hspace*{-1pt}-\hspace*{-1pt}1$};
\node[inner sep=0]at(180:\radius+.24){$n$};
\node[inner sep=0]at($(-90:\radius*cos{45})+(0.03,0)$){$\cdots$};
\node[ambiguous](v1)at($(b1)+(-\s,0)$){};
\node[ambiguous](v2)at($(v1)+(-5*\s/3,0)$){};
\node[inner sep=0](v2l)at($(v2)+(-2*\s/3,0)$){};
\node[ambiguous](v4)at($(b6)+(\s,0)$){};
\node[ambiguous](v3)at($(v4)+(5*\s/3,0)$){};
\node[inner sep=0](v3r)at($(v3)+(2*\s/3,0)$){};
\node[inner sep=0]at($(v1)+(0,0.32)$){$v_1$};
\node[inner sep=0]at($(v2)+(0,0.32)$){$v_2$};
\node[inner sep=0]at($(v3)+(0,0.32)$){$v_{n-3}$};
\node[inner sep=0]at($(v4)+(0,0.32)$){$v_{n-2}$};
\node[inner sep=0]at(0.03,0){$\cdots$};
\path[thick](b1.center)edge(v1) (b2.center)edge(v1) (b3.center)edge(v2) (b4.center)edge(v3) (b5.center)edge(v4) (b6.center)edge(v4) (v1)edge(v2) (v2)edge(v2l) (v3r)edge(v3) (v3)edge(v4);
\end{tikzpicture}
\caption{An arbitrary element of $\mathcal{\tilde{G}}_{n,k,2}$, where precisely $k-1$ of the vertices $v_1, \dots, v_{n-2}$ are black.}
\label{fig:tree}
\end{figure}
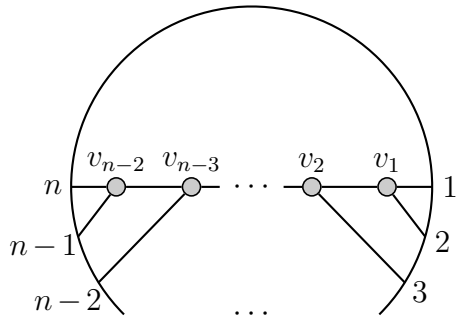

Alternatively, the plabic graphs in $\mathcal{\tilde{G}}_{n,k,2}$ are precisely those shown in \cref{fig:tree}, where $k-1$ of the vertices $v_1, \dots, v_{n-2}$ are black, and the rest are white. Therefore $|\mathcal{\tilde{G}}_{n,k,2}| = \binom{n-2}{k-1}$.

\begin{defn}\label{def:m-is-2-permutations}
Let $k,n\ge 0$ satisfy $k\le n-2$, and 
$c_n := (n \;\; n{-}1 \; \cdots \; 2 \;\; 1)$ be the long cycle in the symmetric group on $[n]$. We define a collection $\Pi_{n,k,2}$ of decorated permutations of type $(k,n)$ by
$$
\Pi_{n,k,2}:= \{c_n\pi_G : G\in \mathcal{\tilde{G}}_{n,k+1,2}\},
$$
where above we color any fixed points of $c_n\pi_G$ black. 
(Note that if $G$ is a plabic graph coming from the recursion in \cref{def:BCFWm2}, then $\pi_G$ has no fixed points, so indeed multiplying $\pi_G$ by $c_n$ on the left decreases the number of anti-excedances by $1$.) 
We let $\mathcal{C}_{n,k,2}$ denote the collection of 
cells $S_\pi\subseteq\Gr_{k,n}^{\ge 0}$ corresponding to 
permutations $\pi\in\Pi_{n,k,2}$.
\end{defn}

\subsection{\texorpdfstring{\Le}{Le}-diagrams for \texorpdfstring{$m=2$}{m=2}}\label{sec:m=2_diagrams}

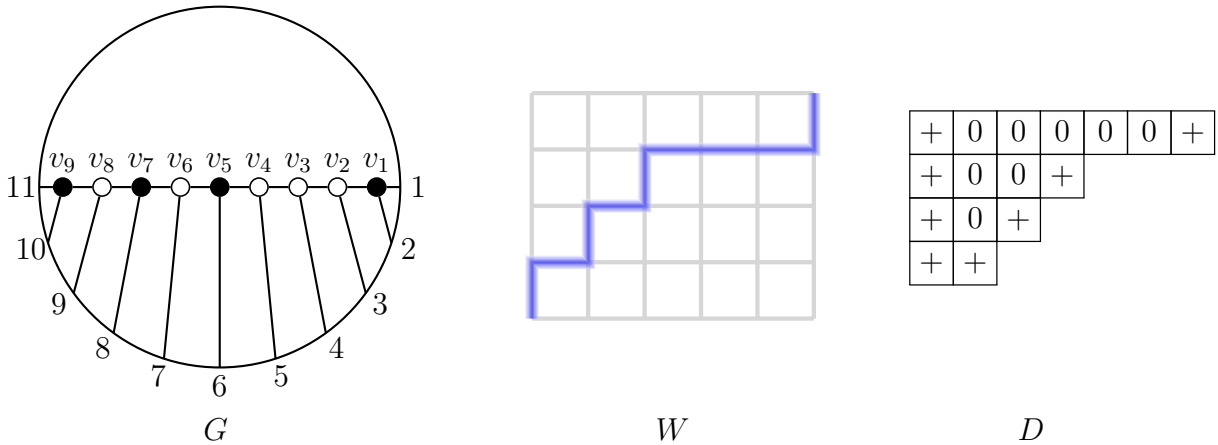
\begin{figure}[ht]
\begin{tabular}{ccc}
\begin{tikzpicture}[baseline=(current bounding box.center)]
\tikzstyle{out1}=[inner sep=0,minimum size=2.4mm,circle,draw=black,fill=black,semithick]
\tikzstyle{in1}=[inner sep=0,minimum size=2.4mm,circle,draw=black,fill=white,semithick]
\tikzstyle{ambiguous}=[inner sep=0,minimum size=2.4mm,circle,draw=black,fill=black!20,semithick]
\pgfmathsetmacro{\n}{11};
\pgfmathsetmacro{\radius}{2.40};
\pgfmathsetmacro{\s}{2*5*\radius/(5*\n-9)};
\pgfmathsetmacro{\angle}{180/(\n-1)};
\draw[thick](0,0)circle[radius=\radius];
\foreach \x in {1,...,11}{
\node[inner sep=0](b\x)at(\angle-\x*\angle:\radius){};
\node at(\angle-\x*\angle:\radius+.24){$\x$};}
\node[out1](v1)at($(b1)+(-3*\s/5,0)$){};
\node[in1](v2)at($(v1)+(-\s,0)$){};
\node[in1](v3)at($(v2)+(-\s,0)$){};
\node[in1](v4)at($(v3)+(-\s,0)$){};
\node[out1](v5)at($(v4)+(-\s,0)$){};
\node[in1](v6)at($(v5)+(-\s,0)$){};
\node[out1](v7)at($(v6)+(-\s,0)$){};
\node[in1](v8)at($(v7)+(-\s,0)$){};
\node[out1](v9)at($(v8)+(-\s,0)$){};
\foreach \x in {1,...,9}{
\node[inner sep=0]at($(v\x)+(0,0.32)$){$v_{\x}$};}
\path[thick](b1.center)edge(v1) (v1)edge(b2.center) edge(v2) (v2)edge(b3.center) edge(v3) (v3)edge(b4.center) edge(v4) (v4)edge(b5.center) edge(v5) (v5)edge(b6.center) edge(v6) (v6)edge(b7.center) edge(v7) (v7)edge(b8.center) edge(v8) (v8)edge(b9.center) edge(v9) (v9)edge(b10.center) edge(b11.center);
\end{tikzpicture}& \qquad
\tikzexternalenable\begin{tikzpicture}[baseline=(current bounding box.center)]
\tikzstyle{hu}=[top color=red!10,bottom color=red!10,middle color=red,opacity=0.70]
\tikzstyle{vu}=[left color=red!10,right color=red!10,middle color=red,opacity=0.70]
\tikzstyle{hl}=[top color=blue!10,bottom color=blue!10,middle color=blue,opacity=0.55]
\tikzstyle{vl}=[left color=blue!10,right color=blue!10,middle color=blue,opacity=0.55]
\pgfmathsetmacro{\u}{0.75};
\pgfmathsetmacro{\w}{0.12*\u};
\coordinate(h)at(-\u,0);
\coordinate(v)at(0,-\u);
\draw[step=\u,color=black!16,ultra thick](0,0)grid(5*\u,4*\u);
\node[inner sep=0](l1)at(5*\u,4*\u){};
\node[inner sep=0](l2)at($(l1)+(v)$){};
\node[inner sep=0](l3)at($(l2)+(h)$){};
\node[inner sep=0](l4)at($(l3)+(h)$){};
\node[inner sep=0](l5)at($(l4)+(h)$){};
\node[inner sep=0](l6)at($(l5)+(v)$){};
\node[inner sep=0](l7)at($(l6)+(h)$){};
\node[inner sep=0](l8)at($(l7)+(v)$){};
\node[inner sep=0](l9)at($(l8)+(h)$){};
\node[inner sep=0](l10)at($(l9)+(v)$){};
\begin{scope}
\clip($(l1.center)+(-\w,0)$)--++($1*(v)+(0,\w)$)--++(2*\w,-2*\w)--++($-1*(v)+(0,\w)$);
\path[vl]($(l1.center)+(\w,\w)$)rectangle($(l2.center)+(-\w,-\w)$);
\end{scope}
\begin{scope}
\clip($(l2.center)+(-\w,\w)$)--++($3*(h)$)--++(2*\w,-2*\w)--++($-3*(h)$);
\path[hl]($(l2.center)+(\w,\w)$)rectangle($(l5.center)+(-\w,-\w)$);
\end{scope}
\begin{scope}
\clip($(l5.center)+(-\w,\w)$)--++($1*(v)$)--++(2*\w,-2*\w)--++($-1*(v)$);
\path[vl]($(l5.center)+(\w,\w)$)rectangle($(l6.center)+(-\w,-\w)$);
\end{scope}
\begin{scope}
\clip($(l6.center)+(-\w,\w)$)--++($1*(h)$)--++(2*\w,-2*\w)--++($-1*(h)$);
\path[hl]($(l6.center)+(\w,\w)$)rectangle($(l7.center)+(-\w,-\w)$);
\end{scope}
\begin{scope}
\clip($(l7.center)+(-\w,\w)$)--++($1*(v)$)--++(2*\w,-2*\w)--++($-1*(v)$);
\path[vl]($(l7.center)+(\w,\w)$)rectangle($(l8.center)+(-\w,-\w)$);
\end{scope}
\begin{scope}
\clip($(l8.center)+(-\w,\w)$)--++($1*(h)$)--++(2*\w,-2*\w)--++($-1*(h)$);
\path[hl]($(l8.center)+(\w,\w)$)rectangle($(l9.center)+(-\w,-\w)$);
\end{scope}
\begin{scope}
\clip($(l9.center)+(-\w,\w)$)--++($1*(v)+(0,-\w)$)--++(2*\w,0)--++($-1*(v)+(0,-\w)$);
\path[vl]($(l9.center)+(\w,\w)$)rectangle($(l10.center)+(-\w,-\w)$);
\end{scope}
\end{tikzpicture}\tikzexternaldisable\qquad &
\qquad\begin{ytableau}
+ & 0 & 0 & 0 & 0 & 0 & + \\
+ & 0 & 0 & + \\
+ & 0 & + \\
+ & +
\end{ytableau}\\
$G$\rule[13pt]{0pt}{0pt} & \qquad$W$\qquad & $D$
\end{tabular}\vspace*{-3pt}
\caption{A plabic graph $G\in\mathcal{\tilde{G}}_{11,5,2}$, the corresponding lattice path $W$ inside a $4\times 5$ rectangle, and the corresponding \Le -diagram $D\in\mathcal{D}_{11,4,2}$. We have $\pi_D = (2,11,\underline{3},\underline{4},6,1,8,5,10,7,9) = c_{11}\pi_G$.}
\label{fig:m=2_bijection}
\end{figure}

\begin{defn}\label{Dnk2}
Let $\mathcal{D}_{n,k,2}$ denote the set of all \Le -diagrams $D$ of type $(k,n)$ such that:
\begin{itemize}
\item each of the $k$ rows of $D$ contains at least $2$ boxes;
\item the leftmost and rightmost entry in each row of $D$ is $+$, and all other entries are $0$.
\end{itemize}
(See \cref{fig:m=2_bijection}.) In particular, $D$ has precisely $2k$ $+$'s, and hence indexes a cell $S_D$ of $\Gr_{k,n}^{\ge 0}$ of dimension $2k$. Also note that the elements of $\mathcal{D}_{n,k,2}$ are naturally indexed by lattice paths inside a $k\times (n-k-2)$ rectangle\footnote{i.e.\ lattice paths moving from northeast to southwest by unit steps west and south; see \cref{def:noncrossing}.}; the lattice path corresponding to $D$ is the southeast border of the Young diagram obtained by deleting the two leftmost columns of $D$.
\end{defn}

We can verify that in fact $\mathcal{D}_{n,k,2}$ indexes the cells $\mathcal{C}_{n,k,2}$, via the following bijection.
\begin{prop}
Given $G\in\mathcal{\tilde{G}}_{n,k+1,2}$ as in \cref{fig:tree}, let $W$ by the lattice path inside a $k\times (n-k-2)$ rectangle given by reading the vertices $v_1, \dots, v_{n-2}$ of $G$, moving west if $v_i$ is white and south if $v_i$ is black. Let $D\in\mathcal{D}_{n,k,2}$ be the \Le -diagram corresponding to $W$, as in \cref{Dnk2}. (See \cref{fig:m=2_bijection}.) Then
$$
\pi_D = c_n\pi_G,
$$
where $c_n := (n \;\; n{-}1 \; \cdots \; 2 \;\; 1)$, $\pi_D$ is defined in \cref{def:oplus}, and all fixed points of $c_n\pi_G$ are colored black. In particular, $\mathcal{D}_{n,k,2}$ indexes the cells $\mathcal{C}_{n,k,2}$.
\end{prop}
We leave the proof as an exercise to the reader, as a warmup to \cref{thm:shift_map_holds}.
\begin{rmk}\label{m=2_to_m=1}
Let $\mathcal{D}_{n,k,1}$ be a set of \Le -diagrams of type $(k,n)$ in bijection with $\mathcal{D}_{n+1,k,2}$, where the element of $\mathcal{D}_{n,k,1}$ corresponding to $D\in\mathcal{D}_{n+1,k,2}$ is formed by deleting the leftmost column of $D$. For example,
$$
\begin{ytableau}
+ & 0 & 0 & 0 & 0 & 0 & + \\
+ & 0 & 0 & + \\
+ & 0 & + \\
+ & +
\end{ytableau}\;\in\mathcal{D}_{11,4,2}\qquad\longleftrightarrow\qquad
\begin{ytableau}
0 & 0 & 0 & 0 & 0 & + \\
0 & 0 & + \\
0 & + \\
+
\end{ytableau}\;\in\mathcal{D}_{10,4,1}\;.
$$
Let $\mathcal{C}_{n,k,1}$ be the cells of $\Gr_{k,n}^{\ge 0}$ indexed by $\mathcal{D}_{n,k,1}$. In \cite{karpwilliams}, two of us showed that the images of $\mathcal{C}_{n,k,1}$ inside the $m=1$ amplituhedron $\mathcal{A}_{n,k,1}(Z)$ induce a triangulation (in fact, a cell decomposition) of $\mathcal{A}_{n,k,1}(Z)$. In fact, one of our motivations for studying this cell decomposition came from discovering the $m=2$ BCFW-like recursion and the corresponding diagrams $\mathcal{D}_{n,k,2}$. In light of the bijection $\mathcal{D}_{n+1,k,2}\to\mathcal{D}_{n,k,1}$, it is natural to ask whether the diagrams corresponding to the BCFW cells for $m=4$ can be similarly truncated to give diagrams inducing a triangulation of the $m=3$ amplituhedron. In \cref{sec:m=3} we show that this fails, at least for the most obvious such truncation.
\end{rmk}

\subsection{Domino bases for \texorpdfstring{$m=2$}{m=2}}\label{sec:m=2_bases}

We now describe a nice basis of any subspace coming from a cell of $\mathcal{C}_{n,k,2}$. We remark that this is the same basis appearing in \cite[(7.5)]{ATT}, up to a cyclic shift (which preserves the positivity of $\Gr_{k,n}$ and $Z$). Therefore, our collection of cells $\mathcal{C}_{n,k,2}$ is the same as Arkani-Hamed, Trnka, and Thomas' up to a cyclic shift.
\begin{defn}\label{def:domino_short}
We say that $d\in\mathbb{R}^n\setminus\{0\}$ is a {\itshape domino} if either $d$ has exactly two nonzero components, which are adjacent and have the same sign, or $d$ has exactly one nonzero component, which is component $n$.
\end{defn}
For example, the vectors $(0,0,0,6,1,0)$ and $(0,0,0,0,0,-2)$ in $\mathbb{R}^6$ are dominoes.
\begin{lem}\label{domino_m=2}
Given $D\in\mathcal{D}_{n,k,2}$, label the edges of the southeast border of $D$ by $1, \dots, n$ from northeast to southwest, and let $s_1 < \dots < s_k$ be the labels of the vertical steps. Also let $V\in\Gr_{k,n}$. Then $V\in S_D$ if and only if $V$ has a basis $\{v^{(1)}, \dots, v^{(k)}\}$ such that for all $i\in [k]$,
\begin{itemize}
\item $v^{(i)}_j = 0$ for all $j\neq s_i, s_i+1, n$; and
\item $v^{(i)}_{s_i}$, $v^{(i)}_{s_i+1}$, and $(-1)^{k-i}v^{(i)}_n$ are all positive.
\end{itemize}
\end{lem}
Note that such a vector $v^{(i)}$ is a sum of two dominoes.
\begin{pf}
Let $A$ be the $k\times n$ matrix from \cref{network_param}, which 
represents a general element of $S_D$. For $i = 1, \dots, k-1$, we perform the following row operations on $A$ to obtain a new $k\times n$ matrix $A'$, 
which represents the same element of $\Gr_{k,n}^{\geq 0}$
that $A$ does:
\begin{itemize}
\item if $s_{i+1} > s_i+1$, we leave row $s_i$ unchanged;
\item if $s_{i+1} = s_i+1$, we add the appropriate scalar multiple of row $s_{i+1}$ to row $s_i$ in order to make the $(s_i,j)$-entry zero, where $j > s_i$ is minimum such that $j\neq s_1, \dots, s_k$.
\end{itemize}
For example, for the \Le -diagram
$$
D = \begin{tikzpicture}[baseline=(current bounding box.center)]
\pgfmathsetmacro{\unit}{0.922};
\useasboundingbox(0,0)rectangle(3.5*\unit,-3.3*\unit);
\node[inner sep=0]at(0,0){\scalebox{1.6}{\begin{ytableau}
\none \\
\none \\
\none \\
\none & \none & \none & + & 0 & + \\
\none & \none & \none & + & 0 & + \\
\none & \none & \none & + & +
\end{ytableau}}};
\end{tikzpicture}\text{ with associated network }
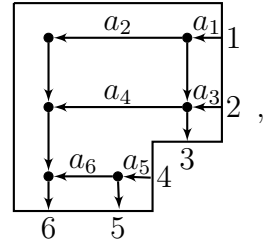
\begin{tikzpicture}[baseline=(current bounding box.center)]
\tikzstyle{out1}=[inner sep=0,minimum size=1.2mm,circle,draw=black,fill=black]
\tikzstyle{in1}=[inner sep=0,minimum size=1.2mm,circle,draw=black,fill=white]
\pgfmathsetmacro{\unit}{0.922};
\useasboundingbox(-0.5*\unit,0)rectangle(3.5*\unit,-3.3*\unit);
\coordinate (vstep)at(0,-0.24*\unit);
\coordinate (hstep)at(0.17*\unit,0);
\coordinate (vepsilon)at(0,-0.02*\unit);
\coordinate (hepsilon)at(0.02*\unit,0);
\draw[thick](0,0)--(3*\unit,0)--(3*\unit,-2*\unit)--(2*\unit,-2*\unit)--(2*\unit,-3*\unit)--(0,-3*\unit)--(0,0);
\node[inner sep=0]at($(3*\unit,-0.5*\unit)+(hstep)$){$1$};
\node[inner sep=0]at($(3*\unit,-1.5*\unit)+(hstep)$){$2$};
\node[inner sep=0]at($(2.5*\unit,-2*\unit)+(vstep)$){$3$};
\node[inner sep=0]at($(2*\unit,-2.5*\unit)+(hstep)$){$4$};
\node[inner sep=0]at($(1.5*\unit,-3*\unit)+(vstep)$){$5$};
\node[inner sep=0]at($(0.5*\unit,-3*\unit)+(vstep)$){$6$};
\node[inner sep=0](b1)at($(3*\unit,-0.5*\unit)+(hepsilon)$){};
\node[inner sep=0](b2)at($(3*\unit,-1.5*\unit)+(hepsilon)$){};
\node[inner sep=0](b3)at($(2.5*\unit,-2*\unit)+(vepsilon)$){};
\node[inner sep=0](b4)at($(2*\unit,-2.5*\unit)+(vepsilon)$){};
\node[inner sep=0](b5)at($(1.5*\unit,-3*\unit)+(hepsilon)$){};
\node[inner sep=0](b6)at($(0.5*\unit,-3*\unit)+(vepsilon)$){};
\node[out1](i11)at($(1*\unit,-1*\unit)+(-0.5*\unit,0.5*\unit)$){};
\node[out1](i13)at($(3*\unit,-1*\unit)+(-0.5*\unit,0.5*\unit)$){};
\node[out1](i21)at($(1*\unit,-2*\unit)+(-0.5*\unit,0.5*\unit)$){};
\node[out1](i23)at($(3*\unit,-2*\unit)+(-0.5*\unit,0.5*\unit)$){};
\node[out1](i31)at($(1*\unit,-3*\unit)+(-0.5*\unit,0.5*\unit)$){};
\node[out1](i32)at($(2*\unit,-3*\unit)+(-0.5*\unit,0.5*\unit)$){};
\path[-latex',thick](b1)edge node[above=-3pt]{$a_1$}(i13) (i13)edge node[above=-3pt]{$a_2$}(i11) (b2)edge node[above=-3pt]{$a_3$}(i23) (i23)edge node[above=-3pt]{$a_4$}(i21) (b4)edge node[above=-3pt]{$a_5$}(i32) (i32)edge node[above=-3pt]{$a_6$}(i31)
(i13)edge(i23) (i11)edge(i21) (i23)edge(b3) (i21)edge(i31) (i32)edge(b5) (i31)edge(b6);
\end{tikzpicture},
$$
the parameterization matrices $A$ and $A'$ are
$$
A = \kbordermatrix{
& 1 & 2 & 3 & 4 & 5 & 6 \cr
1 & 1 & 0 & -a_1 & 0 & 0 & a_1(a_2+a_4) \cr
2 & 0 & 1 & a_3 & 0 & 0 & -a_3a_4 \cr
4 & 0 & 0 & 0 & 1 & a_5 & a_5a_6}\;\;\rightsquigarrow\;\;
A' = \kbordermatrix{
& 1 & 2 & 3 & 4 & 5 & 6 \cr
1 & 1 & a_1/a_3 & 0 & 0 & 0 & a_1a_2 \cr
2 & 0 & 1 & a_3 & 0 & 0 & -a_3a_4 \cr
4 & 0 & 0 & 0 & 1 & a_5 & a_5a_6}.
$$

Given $V\in S_D$, we can specialize the edge variables to positive real numbers so that $A'$ represents $V$. Then we may take $v^{(1)}, \dots, v^{(k)}$ to be the rows of $A'$. Conversely, given any $v^{(1)}, \dots, v^{(k)}$ as in the statement of the lemma, we can find positive values of the edge variables so that the rows of $A'$ are $v^{(1)}, \dots, v^{(k)}$, after we scale row $s_i$ by $v^{(i)}_{s_i}$.
\end{pf}

\subsection{Disjointness for \texorpdfstring{$m=2$}{m=2}}\label{sec:m=2_disjointness}

We show that the images of the cells $\mathcal{C}_{n,k,2}$ inside the $m=2$ amplituhedron $\mathcal{A}_{n,k,2}(Z)$ are disjoint. This is a warmup to \cref{sec:k=2-disjointness}, where we prove disjointness for $m=4$ when $k=2$.
\begin{lem}\label{domino_sum}
If $v\in\mathbb{R}^n$ is a sum of $k\ge 1$ dominoes, then $\var(v)\le k-1$.
\end{lem}
We leave the proof as an exercise to the reader. We will provide a more general version of \cref{domino_sum} in \cref{lem:obo}.
\begin{prop}\label{m=2_disjointness_result}
Let $Z \in \Mat_{k+2,n}^{>0}$.
Then $\tilde{Z}$ maps the cells $\mathcal{C}_{n,k,2}$ of $\Gr_{2,n}^{\ge 0}$ injectively into the amplituhedron $\mathcal{A}_{n,k,2}(Z)$, and their images are pairwise disjoint.
\end{prop}

\begin{pf}
We must show that given $D,D'\in\mathcal{D}_{n,k,2}$ and $V\in S_D, V'\in S_{D'}$ such that $\tilde{Z}(V) = \tilde{Z}(V')$, we have $V = V'$. Let $\{v^{(1)}, \dots, v^{(k)}\}$ be the distinguished basis of $V$ from \cref{domino_m=2}.
\begin{claim}
For any $i\in [k]$ and $v'\in V'$, we have $\var(v^{(i)}-v')\le k+1$.
\end{claim}
\begin{claimpf}
Let $\{w^{(1)}, \dots, w^{(k)}\}$ be the distinguished basis of $V'$ from \cref{domino_m=2}, and write $v' = \sum_{i=1}^kc_iw^{(i)}$ for some $c_1, \dots, c_k\in\mathbb{R}$. Then $v^{(i)} - v' = v^{(i)} - \sum_{i=1}^kc_iw^{(i)}$. The right-hand side can be written as  a sum of $k+2$ or fewer dominoes, 
so $\var(v^{(i)}-v')\le k+1$ by \cref{domino_sum}.
\end{claimpf}
For $i\in [k]$, let $v'^{(i)}\in V'$ be the matching vector for $v^{(i)}$, as in \cref{criterion}. Then by the claim and \cref{criterion}, we have $v^{(i)} = v'^{(i)}$, i.e.\ $v^{(i)}\in V'$. Hence $V\subseteq V'$, and so $V = V'$.
\end{pf}

\section{Binary trees and BCFW plabic graphs}\label{sec:BCFW}

\noindent In the case $m=4$ (which is the case of interest in physics), there is a distinguished collection of  $4k$-dimensional cells of $\Gr_{k,n}^{\geq 0}$ called \emph{BCFW cells}. These cells are named after Britto, Cachazo, Feng, and Witten, 
who gave recursion relations for computing scattering amplitudes \cite{BCF, BCFW}. These relations were translated in \cite{abcgpt} into a recursion on plabic graphs, which we define in \cref{sec:recursive}. The resulting graphs are easily seen to be in bijection with complete binary trees (well-known Catalan objects), as we explain in \cref{sec:trees-to-graphs}. However, these plabic graphs do not literally label the BCFW cells in the sense of \cref{def:rules}; rather, we must perform a `shift by $2$' to the decorated permutation corresponding to the graph, which makes it difficult to explicitly describe the BCFW cells from the BCFW recursion.\footnote{Bai and He \cite[(3.3)]{bai_he_15} found a recursion on the plabic graphs which do literally index the BCFW cells. However, their recursion is more complicated and does not obviously correspond to a family of Catalan objects, so we prefer to work with the original version of the recursion.} In \cref{sec:lattice-paths}, we will use a different family of Catalan objects, namely pairs of noncrossing lattice paths inside a rectangle, to explicitly describe the BCFW cells in terms of $\oplus$-diagrams.

\subsection{The recursive description of BCFW cells from plabic graphs}\label{sec:recursive}

The {\itshape BCFW recursion} \cite{BCF, BCFW} is a recursive operation that produces, 
for each pair $(k,n)$ with $0 \le k \le n-4$, a set $\mathcal{C}_{n,k,4}$ of $4k$-dimensional 
positroid cells of $\Gr_{k,n}^{\ge 0}$, 
called the {\itshape $(k,n)$-BCFW cells}.\footnote{In fact, there are many different ways to 
carry out the BCFW recursion, which each lead to a possibly different set of positroid 
cells of $\Gr_{k,n}^{\ge 0}$. For example one can cyclically permute all boundary labels of the 
plabic graphs one obtains in the recursion.  In this paper 
we fix a canonical way to perform the recursion.} 
This operation is described in \cite[Section 17.2/16.2\footnote{These two section numbers refer, respectively, to the published book and the arXiv preprint.}]{abcgpt}.

\begin{defn}
We define a local operation on plabic graphs which we call 
{\itshape blowing up}. Let $G$ be a plabic graph of type $(k,n)$ 
(see \cref{def:rules}),
and $i\in [n]$ a boundary vertex incident to a unique internal vertex $v$, which has degree $3$. We {\itshape blow up $G$ at $i$} by locally replacing $v$ by a square, as follows, while leaving the rest of the graph unchanged:
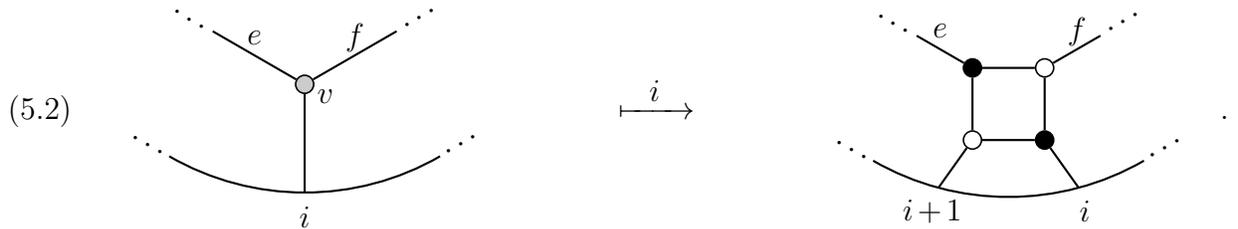
\begin{align}\label{blowup}
\quad\begin{tikzpicture}[baseline=(current bounding box.center)]
\tikzstyle{out1}=[inner sep=0,minimum size=2.4mm,circle,draw=black,fill=black,semithick]
\tikzstyle{in1}=[inner sep=0,minimum size=2.4mm,circle,draw=black,fill=white,semithick]
\tikzstyle{ambiguous}=[inner sep=0,minimum size=2.4mm,circle,draw=black,fill=black!20,semithick]
\pgfmathsetmacro{\radius}{3.6};
\pgfmathsetmacro{\l}{1.44};
\draw[thick](-120:\radius)arc(-120:-60:\radius);
\node[inner sep=0,rotate=150]at($(-120:\radius)+(150:0.36)$){$\cdots$};
\node[inner sep=0,rotate=30]at($(-60:\radius)+(30:0.36)$){$\cdots$};
\node[inner sep=0](i)at(-90:\radius){};
\node at(-90:\radius+0.32){$i$};
\node[ambiguous](v)at(-90:\radius-\l){};
\node[inner sep=0]at($(v)+(-30:0.32)$){$v$};
\node[inner sep=0](left)at($(v)+(150:\l)$){};
\node[inner sep=0,rotate=150]at($(left)+(150:0.36)$){$\cdots$};
\node[inner sep=0](right)at($(v)+(30:\l)$){};
\node[inner sep=0,rotate=30]at($(right)+(30:0.36)$){$\cdots$};
\path[thick](v)edge(i.center) (v)edge node[above]{$e$}(left) (v)edge node[above=-3pt]{$f$}(right);
\end{tikzpicture}\qquad\qquad\xmapsto{\mbox{\;\, $i$ \;\,}}\qquad\qquad
\begin{tikzpicture}[baseline=(current bounding box.center)]
\tikzstyle{out1}=[inner sep=0,minimum size=2.4mm,circle,draw=black,fill=black,semithick]
\tikzstyle{in1}=[inner sep=0,minimum size=2.4mm,circle,draw=black,fill=white,semithick]
\tikzstyle{ambiguous}=[inner sep=0,minimum size=2.4mm,circle,draw=black,fill=black!20,semithick]
\pgfmathsetmacro{\radius}{3.6};
\pgfmathsetmacro{\l}{1.44};
\pgfmathsetmacro{\s}{0.96};
\draw[thick](-120:\radius)arc(-120:-60:\radius);
\node[inner sep=0,rotate=150]at($(-120:\radius)+(150:0.36)$){$\cdots$};
\node[inner sep=0,rotate=30]at($(-60:\radius)+(30:0.36)$){$\cdots$};
\node[inner sep=0](i)at(-75:\radius){};
\node at(-75:\radius+0.32){$i$};
\node[inner sep=0](i+1)at(-105:\radius){};
\node at(-105:\radius+0.32){$i\hspace*{-1pt}+\hspace*{-1pt}1$};
\node[inner sep=0](v)at(-90:\radius-\l){};
\node[out1](a)at($(v)+(150:{\s/1.7320508075688772935274463415058723669428052538104})$){};
\node[in1](b)at($(a)+(\s,0)$){};
\node[out1](c)at($(a)+(\s,-\s)$){};
\node[in1](d)at($(a)+(0,-\s)$){};
\node[inner sep=0](left)at($(v)+(150:\l)$){};
\node[inner sep=0,rotate=150]at($(left)+(150:0.36)$){$\cdots$};
\node[inner sep=0](right)at($(v)+(30:\l)$){};
\node[inner sep=0,rotate=30]at($(right)+(30:0.36)$){$\cdots$};
\path[thick](a)edge(b) (b)edge(c) (c)edge(d) (d)edge(a) (c)edge(i.center) (d)edge(i+1.center) (a)edge node[above]{$e$}(left) (b)edge node[above=-3pt]{$f$}(right);
\end{tikzpicture}\quad.
\end{align}
(The color of $v$ may be black or white.) This produces a plabic graph $G'$ of type $(k',n+1)$, whose boundary vertices are labeled by $[n+1]$ so that $i$ and $i+1$ are labeled as above (i.e.\ labels $1,\dots, i-1$ from $G$ remain unchanged, and labels $i+1, \dots, n$ from $G$ get increased by $1$). In the blowup, the orientation of the square is important: its vertices alternate in color, with $i$ incident to a black vertex and $i+1$ to a white vertex. 
Note that by \eqref{k-statistic}, $k' = k$ if $v$ is black, while $k' = k+1$ if $v$ is white.
\end{defn}

\begin{defn}[BCFW recursion]\label{def:BCFW}
For positive integers $k$ and $n$ such that $2 \leq k \leq n-2$, we recursively define a set $\mathcal{\tilde{G}}_{n,k,4}$ of plabic graphs of type $(k,n)$ as follows.\footnote{In our notation $\mathcal{\tilde{G}}_{n,k,4}$, the tilde is to remind us that these graphs do {\bfseries not} directly label the $m=4$ BCFW cells. Rather, we must first `shift by $2$'; see \cref{def:BCFW-permutations}.}
\begin{enumerate}[label=\arabic*.]
\item If $n=4$ and $k=2$, then $\mathcal{\tilde{G}}_{n,k,4}$ contains a unique element, the plabic graph
$$
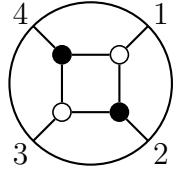
\quad\begin{tikzpicture}[baseline=(current bounding box.center)]
\tikzstyle{out1}=[inner sep=0,minimum size=2.4mm,circle,draw=black,fill=black,semithick]
\tikzstyle{in1}=[inner sep=0,minimum size=2.4mm,circle,draw=black,fill=white,semithick]
\pgfmathsetmacro{\radius}{1.08};
\draw[thick](0,0)circle[radius=\radius];
\foreach \x in {1,...,4}{
\node[inner sep=0](b\x)at(135-\x*90:\radius){};
\node at(135-\x*90:\radius+.24){$\x$};}
\node[in1](i1)at(45:\radius/2){};
\node[out1](i2)at(-45:\radius/2){};
\node[in1](i3)at(-135:\radius/2){};
\node[out1](i4)at(135:\radius/2){};
\path[thick](b1.center)edge(i1) (b2.center)edge(i2) (i3)edge(b3.center) (i4)edge(b4.center) (i1)edge(i2)edge(i4) (i2)edge(i3) (i3)edge(i4);
\end{tikzpicture}\quad.
$$
\item For $n\ge 5$, $\mathcal{\tilde{G}}_{n,k,4}$ is the set of all plabic graphs obtained either by blowing up a plabic graph in $\mathcal{\tilde{G}}_{n-1,k,4}$ at some $i\neq 1,n-1$ incident to a black vertex, or by blowing up a plabic graph in $\mathcal{\tilde{G}}_{n-1,k-1,4}$ at some $i\neq 1,n-1$ incident to a white vertex. (We may obtain such a plabic graph multiple times in this way, but we only count it once in $\mathcal{\tilde{G}}_{n,k,4}$.) We emphasize that we never blow up at the first or last boundary positions.
\end{enumerate}
See \cref{fig:BCFW-graphs} for a depiction of the graphs
$\mathcal{\tilde{G}}_{n,k,4}$ for $n\le 6$. Note that we can read the $k$ statistic from each graph as the number of black vertices incident to the boundary.
\end{defn}
\begin{figure}[ht]
\begin{center}
$n=4$\qquad\qquad\qquad\begin{tikzpicture}[baseline=(current bounding box.center)]
      \tikzstyle{out1}=[inner sep=0,minimum size=2.4mm,circle,draw=black,fill=black,semithick]
      \tikzstyle{in1}=[inner sep=0,minimum size=2.4mm,circle,draw=black,fill=white,semithick]
      \pgfmathsetmacro{\radius}{1.2};
      \draw[thick](0,0)circle[radius=\radius];
      \foreach \x in {1,...,4}{
        \node[inner sep=0](b\x)at(135-\x*90:\radius){};
        \node at(135-\x*90:\radius+.24){$\x$};}
      \node[in1](i1)at(45:\radius/2){};
      \node[out1](i2)at(-45:\radius/2){};
      \node[in1](i3)at(-135:\radius/2){};
      \node[out1](i4)at(135:\radius/2){};
      \path[thick](b1.center)edge(i1) (b2.center)edge(i2) (i3)edge(b3.center) (i4)edge(b4.center) (i1)edge(i2)edge(i4) (i2)edge(i3) (i3)edge(i4);
\node[inner sep=0]at(0,-\radius-.36){$k=2$};
\end{tikzpicture}\qquad\qquad\qquad\phantom{$n=4$}
\end{center}
\begin{tikzpicture}[baseline=(current bounding box.center),
  to/.style={->,>=stealth',shorten >=1pt,semithick}]
  \node[inner sep=0](in1)at(0,1.5){};
  \node[inner sep=0](out1)at(-2.3,0){};
  \node[inner sep=0](out2)at(2.3,0){};
  \draw[to](in1)edge node[left=4pt]{$2$}(out1);
  \draw[to](in1)edge node[right=4pt]{$3$}(out2);
\end{tikzpicture}
\begin{center}
$n=5$\qquad\qquad\begin{tikzpicture}[baseline=(current bounding box.center)]
    \tikzstyle{out1}=[inner sep=0,minimum size=1.8mm,circle,draw=black,fill=black,semithick]
    \tikzstyle{in1}=[inner sep=0,minimum size=1.8mm,circle,draw=black,fill=white,semithick]
    \pgfmathsetmacro{\radius}{1.20};
    \pgfmathsetmacro{\radiuss}{.55};
    \pgfmathsetmacro{\radiust}{.36};
    \pgfmathsetmacro{\radiusu}{.24};
    \pgfmathsetmacro{\shift}{.36};
    \draw[thick](0,0)circle[radius=\radius];
    \foreach \x in {1,...,5}{
      \node[inner sep=0](b\x)at(126-\x*72:\radius){};
      \node at(126-\x*72:\radius+.24){$\x$};}
    \node[inner sep=0](s0)at(90:\shift){};
    \node[out1](s1)at($(s0)+(135:\radiuss)$){};
    \node[in1](s2)at($(s0)+(45:\radiuss)$){};
    \node[inner sep=0](s3)at($(s0)+(-45:\radiuss)$){};
    \node[in1](s4)at($(s0)+(-135:\radiuss)$){};
    \node[out1](t1)at($(s3)+(180:\radiust)$){};
    \node[in1](t2)at($(s3)+(90:\radiust)$){};
    \node[out1](t3)at($(s3)+(0:\radiust)$){};
    \node[in1](t4)at($(s3)+(-90:\radiust)$){};
    \path[thick](b1.center)edge(s2) (b2.center)edge(t3) (b3.center)edge(t4) (b4.center)edge(s4) (b5.center)edge(s1)  (s1)edge(s2) (s1)edge(s4) (s2)edge(t2) (s4)edge(t1) (t1)edge(t2) (t1)edge(t4) (t2)edge(t3) (t3)edge(t4);
\node[inner sep=0]at(0,-\radius-.72){$k=2$};
\end{tikzpicture}\qquad\qquad\begin{tikzpicture}[baseline=(current bounding box.center)]
    \tikzstyle{out1}=[inner sep=0,minimum size=1.8mm,circle,draw=black,fill=black,semithick]
    \tikzstyle{in1}=[inner sep=0,minimum size=1.8mm,circle,draw=black,fill=white,semithick]
    \pgfmathsetmacro{\radius}{1.20};
    \pgfmathsetmacro{\radiuss}{.55};
    \pgfmathsetmacro{\radiust}{.36};
    \pgfmathsetmacro{\radiusu}{.24};
    \pgfmathsetmacro{\shift}{.36};
    \draw[thick](0,0)circle[radius=\radius];
    \foreach \x in {1,...,5}{
      \node[inner sep=0](b\x)at(126-\x*72:\radius){};
      \node at(126-\x*72:\radius+.24){$\x$};}
    \node[inner sep=0](s0)at(90:\shift){};
    \node[out1](s1)at($(s0)+(135:\radiuss)$){};
    \node[in1](s2)at($(s0)+(45:\radiuss)$){};
    \node[out1](s3)at($(s0)+(-45:\radiuss)$){};
    \node[inner sep=0](s4)at($(s0)+(-135:\radiuss)$){};
    \node[in1](t1)at($(s4)+(180:\radiust)$){};
    \node[out1](t2)at($(s4)+(90:\radiust)$){};
    \node[in1](t3)at($(s4)+(0:\radiust)$){};
    \node[out1](t4)at($(s4)+(-90:\radiust)$){};
    \path[thick](b1.center)edge(s2) (b2.center)edge(s3) (b3.center)edge(t4) (b4.center)edge(t1) (b5.center)edge(s1)  (s1)edge(s2) (s2)edge(s3) (s3)edge(t3) (s1)edge(t2) (t1)edge(t2) (t1)edge(t4) (t2)edge(t3) (t3)edge(t4);
\node[inner sep=0]at(0,-\radius-.72){$k=3$};
\end{tikzpicture}\qquad\qquad\phantom{$n=5$}
\end{center}
\begin{tikzpicture}[baseline=(current bounding box.center),to/.style={->,>=stealth',shorten >=1pt,semithick}]
  \node[inner sep=0](in1)at(-2.4,2){};
  \node[inner sep=0](in2)at(2.4,2){};
  \node[inner sep=0](out1)at(-6.4,0){};
  \node[inner sep=0](out2)at(-3.2,0){};
  \node[inner sep=0](out3l)at(-0.1,0){};
  \node[inner sep=0](out3r)at(0.1,0){};
  \node[inner sep=0](out4)at(3.2,0){};
  \node[inner sep=0](out5)at(6.4,0){};
  \path[to](in1)edge node[left=6pt]{$2$}(out1);
  \path[to](in1)edge node[left=-1pt]{$3$}(out2);
  \path[to](in1)edge node[right=2pt]{$4$}(out3l);
  \path[to](in2)edge node[right=6pt]{$4$}(out5);
  \path[to](in2)edge node[right=-1pt]{$3$}(out4);
  \path[to](in2)edge node[left=2pt]{$2$}(out3r);
\end{tikzpicture}
\begin{center}
\makebox[\textwidth]{
$n=6$\qquad\begin{tikzpicture}[baseline=(current bounding box.center)]
    \tikzstyle{out1}=[inner sep=0,minimum size=1.8mm,circle,draw=black,fill=black,semithick]
    \tikzstyle{in1}=[inner sep=0,minimum size=1.8mm,circle,draw=black,fill=white,semithick]
    \pgfmathsetmacro{\radius}{1.1};
    \pgfmathsetmacro{\radiuss}{.55};
    \pgfmathsetmacro{\radiust}{.36};
    \pgfmathsetmacro{\radiusu}{.2};
    \pgfmathsetmacro{\shift}{.36};
    \draw[thick](0,0)circle[radius=\radius];
    \foreach \x in {1,...,6}{
      \node[inner sep=0](b\x)at(180-\x*60:\radius){};
      \node at(120-\x*60:\radius+.24){$\x$};}
    \node[inner sep=0](s0)at(90:\shift){};
    \node[out1](s1)at($(s0)+(135:\radiuss)$){};
    \node[in1](s2)at($(s0)+(45:\radiuss)$){};
    \node[inner sep=0](s3)at($(s0)+(-45:\radiuss)$){};
    \node[in1](s4)at($(s0)+(-135:\radiuss)$){};
    \node[out1](t1)at($(s3)+(180:\radiust)$){};
    \node[in1](t2)at($(s3)+(90:\radiust)$){};
    \node[inner sep=0](t3)at($(s3)+(0:\radiust)$){};
    \node[in1](t4)at($(s3)+(-90:\radiust)$){};
    \node[in1](u1)at($(t3)+(135:\radiusu)$){};
    \node[out1](u2)at($(t3)+(45:\radiusu)$){};
    \node[in1](u3)at($(t3)+(-45:\radiusu)$){};
    \node[out1](u4)at($(t3)+(-135:\radiusu)$){};
    \path[thick](b1.center)edge(s1) (b2.center)edge(s2) (b3.center)edge(u2) (b4.center)edge(u3) (b5.center)edge(t4) (b6.center)edge(s4)  (s1)edge(s2) (s1)edge(s4) (s2)edge(t2) (s4)edge(t1) (t1)edge(t2) (t1)edge(t4) (t4)edge(u4) (t2)edge(u1) (u1)edge(u2) (u1)edge(u4) (u2)edge(u3) (u3)edge(u4);
\node[inner sep=0]at(0,-\radius-.60){$k=2$};
\end{tikzpicture}\hspace*{5pt}\begin{tikzpicture}[baseline=(current bounding box.center)]
    \tikzstyle{out1}=[inner sep=0,minimum size=1.8mm,circle,draw=black,fill=black,semithick]
    \tikzstyle{in1}=[inner sep=0,minimum size=1.8mm,circle,draw=black,fill=white,semithick]
    \pgfmathsetmacro{\radius}{1.1};
    \pgfmathsetmacro{\radiuss}{.55};
    \pgfmathsetmacro{\radiust}{.36};
    \pgfmathsetmacro{\radiusu}{.24};
    \pgfmathsetmacro{\shift}{.36};
    \draw[thick](0,0)circle[radius=\radius];
    \foreach \x in {1,...,6}{
      \node[inner sep=0](b\x)at(180-\x*60:\radius){};
      \node at(120-\x*60:\radius+.24){$\x$};}
    \node[inner sep=0](s0)at(90:\shift){};
    \node[out1](s1)at($(s0)+(135:\radiuss)$){};
    \node[in1](s2)at($(s0)+(45:\radiuss)$){};
    \node[inner sep=0](s3)at($(s0)+(-45:\radiuss)$){};
    \node[in1](s4)at($(s0)+(-135:\radiuss)$){};
    \node[out1](t1)at($(s3)+(180:\radiust)$){};
    \node[in1](t2)at($(s3)+(90:\radiust)$){};
    \node[out1](t3)at($(s3)+(0:\radiust)$){};
    \node[inner sep=0](t4)at($(s3)+(-90:\radiust)$){};
    \node[out1](u1)at($(t4)+(135:\radiusu)$){};
    \node[in1](u2)at($(t4)+(45:\radiusu)$){};
    \node[out1](u3)at($(t4)+(-45:\radiusu)$){};
    \node[in1](u4)at($(t4)+(-135:\radiusu)$){};
    \path[thick](b1.center)edge(s1) (b2.center)edge(s2) (b3.center)edge(t3) (b4.center)edge(u3) (b5.center)edge(u4) (b6.center)edge(s4)  (s1)edge(s2) (s1)edge(s4) (s2)edge(t2) (s4)edge(t1) (t1)edge(t2) (t1)edge(u1) (t2)edge(t3) (t3)edge(u2) (u1)edge(u2) (u1)edge(u4) (u2)edge(u3) (u3)edge(u4);
\node[inner sep=0]at(0,-\radius-.60){$k=3$};
\end{tikzpicture}\hspace*{5pt}\begin{tikzpicture}[baseline=(current bounding box.center)]
\tikzstyle{out1}=[inner sep=0,minimum size=1.8mm,circle,draw=black,fill=black,semithick]
\tikzstyle{in1}=[inner sep=0,minimum size=1.8mm,circle,draw=black,fill=white,semithick]
\pgfmathsetmacro{\radius}{1.1};
\pgfmathsetmacro{\radiuss}{.55};
\pgfmathsetmacro{\radiust}{.24};
\pgfmathsetmacro{\radiusu}{.24};
\pgfmathsetmacro{\shift}{.36};
\draw[thick](0,0)circle[radius=\radius];
\foreach \x in {1,...,6}{
\node[inner sep=0](b\x)at(180-\x*60:\radius){};
\node at(120-\x*60:\radius+.24){$\x$};}
\node[inner sep=0](s0)at(90:\shift){};
\node[out1](s1)at($(s0)+(135:\radiuss)$){};
\node[in1](s2)at($(s0)+(45:\radiuss)$){};
\node[inner sep=0](s3)at($(s0)+(-45:\radiuss)$){};
\node[inner sep=0](s4)at($(s0)+(-135:\radiuss)$){};
\node[out1](t1)at($(s3)+(180:\radiust)$){};
\node[in1](t2)at($(s3)+(90:\radiust)$){};
\node[out1](t3)at($(s3)+(0:\radiust)$){};
\node[in1](t4)at($(s3)+(-90:\radiust)$){};
\node[out1](u1)at($(s4)+(90:\radiusu)$){};
\node[in1](u2)at($(s4)+(0:\radiusu)$){};
\node[out1](u3)at($(s4)+(-90:\radiusu)$){};
\node[in1](u4)at($(s4)+(-180:\radiusu)$){};
\path[thick](b1.center)edge(s1) (b2.center)edge(s2) (b3.center)edge(t3) (b4.center)edge(t4) (b5.center)edge(u3) (b6.center)edge(u4) (s1)edge(s2) (s1)edge(u1) (s2)edge(t2) (t1)edge(t2) (t1)edge(t4) (t1)edge(u2) (t2)edge(t3) (t3)edge(t4) (u1)edge(u2) (u1)edge(u4) (u2)edge(u3) (u3)edge(u4);
\node[inner sep=0]at(0,-\radius-.60){$k=3$};
\end{tikzpicture}\hspace*{5pt}\begin{tikzpicture}[baseline=(current bounding box.center)]
\tikzstyle{out1}=[inner sep=0,minimum size=1.8mm,circle,draw=black,fill=black,semithick]
\tikzstyle{in1}=[inner sep=0,minimum size=1.8mm,circle,draw=black,fill=white,semithick]
\pgfmathsetmacro{\radius}{1.1};
\pgfmathsetmacro{\radiuss}{.55};
\pgfmathsetmacro{\radiust}{.36};
\pgfmathsetmacro{\radiusu}{.24};
\pgfmathsetmacro{\shift}{.36};
\draw[thick](0,0)circle[radius=\radius];
\foreach \x in {1,...,6}{
\node[inner sep=0](b\x)at(180-\x*60:\radius){};
\node at(120-\x*60:\radius+.24){$\x$};}
\node[inner sep=0](s0)at(90:\shift){};
\node[out1](s1)at($(s0)+(135:\radiuss)$){};
\node[in1](s2)at($(s0)+(45:\radiuss)$){};
\node[out1](s3)at($(s0)+(-45:\radiuss)$){};
\node[inner sep=0](s4)at($(s0)+(-135:\radiuss)$){};
\node[out1](t1)at($(s4)+(90:\radiust)$){};
\node[in1](t2)at($(s4)+(0:\radiust)$){};
\node[inner sep=0](t3)at($(s4)+(-90:\radiust)$){};
\node[in1](t4)at($(s4)+(-180:\radiust)$){};
\node[out1](u1)at($(t3)+(135:\radiusu)$){};
\node[in1](u2)at($(t3)+(45:\radiusu)$){};
\node[out1](u3)at($(t3)+(-45:\radiusu)$){};
\node[in1](u4)at($(t3)+(-135:\radiusu)$){};
\path[thick](b1.center)edge(s1) (b2.center)edge(s2) (b3.center)edge(s3) (b4.center)edge(u3) (b5.center)edge(u4) (b6.center)edge(t4) (s1)edge(s2) (s1)edge(t1) (s2)edge(s3) (s3)edge(t2) (t1)edge(t2) (t1)edge(t4) (t2)edge(u2) (t4)edge(u1) (u1)edge(u2) (u1)edge(u4) (u2)edge(u3) (u3)edge(u4);
\node[inner sep=0]at(0,-\radius-.60){$k=3$};
\end{tikzpicture}\hspace*{5pt}\begin{tikzpicture}[baseline=(current bounding box.center)]
\tikzstyle{out1}=[inner sep=0,minimum size=1.8mm,circle,draw=black,fill=black,semithick]
\tikzstyle{in1}=[inner sep=0,minimum size=1.8mm,circle,draw=black,fill=white,semithick]
\pgfmathsetmacro{\radius}{1.1};
\pgfmathsetmacro{\radiuss}{.55};
\pgfmathsetmacro{\radiust}{.36};
\pgfmathsetmacro{\radiusu}{.2};
\pgfmathsetmacro{\shift}{.36};
\draw[thick](0,0)circle[radius=\radius];
\foreach \x in {1,...,6}{
\node[inner sep=0](b\x)at(180-\x*60:\radius){};
\node at(120-\x*60:\radius+.24){$\x$};}
\node[inner sep=0](s0)at(90:\shift){};
\node[out1](s1)at($(s0)+(135:\radiuss)$){};
\node[in1](s2)at($(s0)+(45:\radiuss)$){};
\node[out1](s3)at($(s0)+(-45:\radiuss)$){};
\node[inner sep=0](s4)at($(s0)+(-135:\radiuss)$){};
\node[out1](t1)at($(s4)+(90:\radiust)$){};
\node[in1](t2)at($(s4)+(0:\radiust)$){};
\node[out1](t3)at($(s4)+(-90:\radiust)$){};
\node[inner sep=0](t4)at($(s4)+(-180:\radiust)$){};
\node[in1](u1)at($(t4)+(135:\radiusu)$){};
\node[out1](u2)at($(t4)+(45:\radiusu)$){};
\node[in1](u3)at($(t4)+(-45:\radiusu)$){};
\node[out1](u4)at($(t4)+(-135:\radiusu)$){};
\path[thick](b1.center)edge(s1) (b2.center)edge(s2) (b3.center)edge(s3) (b4.center)edge(t3) (u3)edge(t3) (b5.center)edge(u4) (b6.center)edge(u1) (s1)edge(s2) (s1)edge(t1) (s2)edge(s3) (s3)edge(t2) (t1)edge(t2) (t1)edge(u2) (t3)edge(t2) (u1)edge(u2) (u1)edge(u4) (u2)edge(u3) (u3)edge(u4);
\node[inner sep=0]at(0,-\radius-.60){$k=4$};
\end{tikzpicture}\qquad\phantom{$n=6$}}
\end{center}
\caption{The first two steps of the BCFW recursion. The Catalan sequence $1, 2, 5, 14, \dots$ gives the number of graphs in each row. Note that the $k$ statistic has not yet been `shifted by $2$'. Observe that the middle graph in the bottom row can be obtained via two different sequences of blowups.}
\label{fig:BCFW-graphs}
\end{figure}
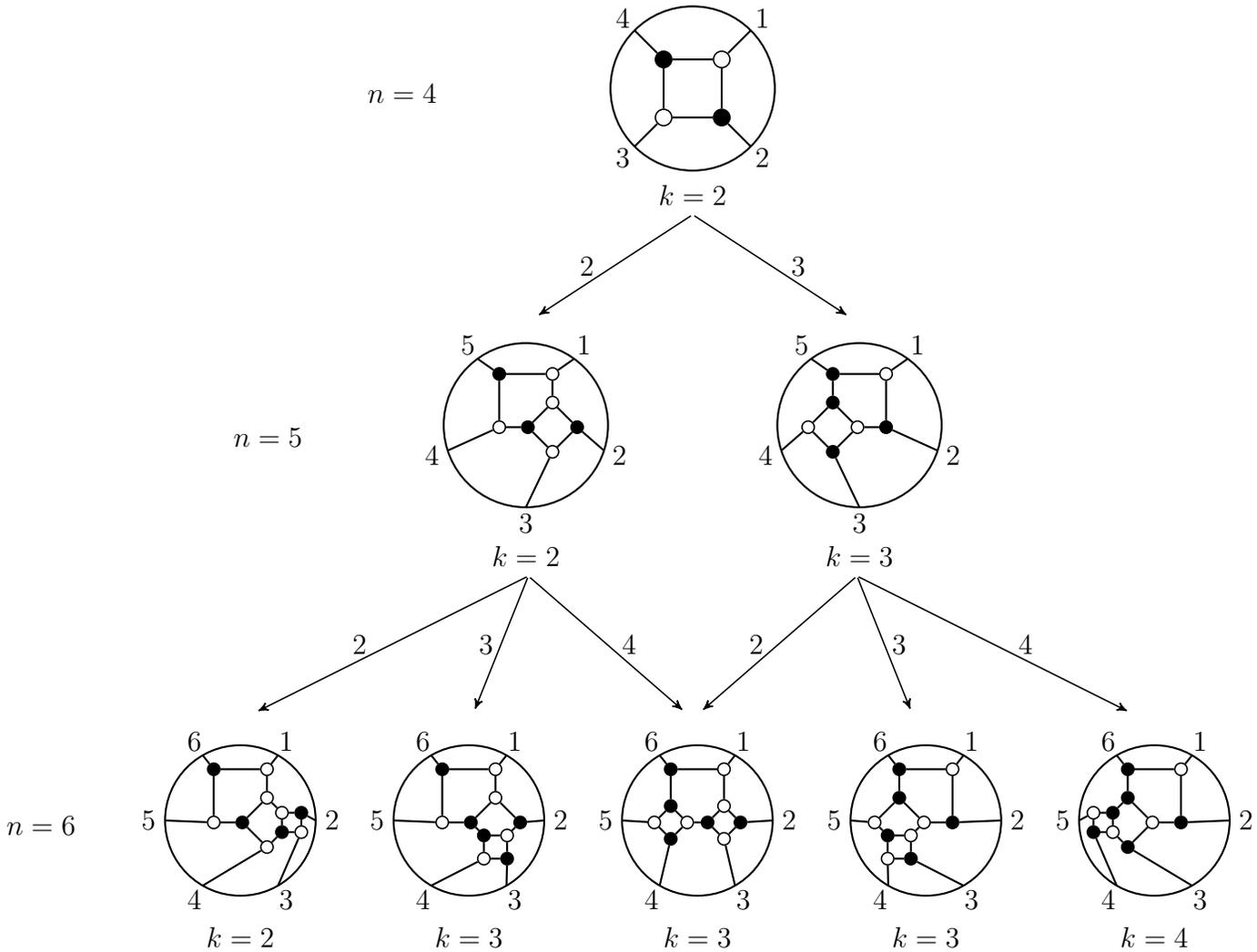

\begin{defn}\label{def:BCFW-permutations}
Let $k,n\ge 0$ satisfy $k\le n-4$, and $c_n := (n \;\; n{-}1 \; \cdots \; 2 \;\; 1)$ be the long cycle in the symmetric group on $[n]$. We define the {\itshape BCFW permutations} of type $(k,n)$ (for $m=4$) as
$$
\Pi_{n,k,4}:= \{c_n^2\pi_G : G\in \mathcal{\tilde{G}}_{n,k+2,4}\},
$$
where above we color any fixed points of $c_n^2\pi_G$ black. (Note that for any plabic graph $G$ coming from the BCFW recursion, we have $\pi_G(i)\neq i, i+1\pmod{n}$ for all $i\in [n]$, so indeed multiplying $\pi_G$ by $c_n^2$ on the left decreases the number of anti-excedances by $2$.) 
The collection
$\mathcal{C}_{n,k,4}$ of 
cells $S_\pi\subseteq\Gr_{k,n}^{\ge 0}$ corresponding to BCFW permutations $\pi\in\Pi_{n,k,4}$ are called the {\itshape $(k,n)$-BCFW cells} (for $m=4$).\footnote{Some authors use the term `BCFW cells' to refer to the images $\tilde{Z}(S_\pi)$ in the amplituhedron $\mathcal{A}_{n,k,4}(Z)$ (though in general it is not known that these images are indeed topological cells).}
\end{defn}

\begin{conj}[{\cite[Section 5]{arkani-hamed_trnka}}]\label{amp-conj}
Let $Z \in \Mat_{k+4,n}^{>0}$, where $k,n\ge 0$ satisfy $k\le n-4$. Then the images under $\tilde{Z}$ of the BCFW cells $\mathcal{C}_{n,k,4}$ ``triangulate" the $m=4$ amplituhedron, i.e.\ they are pairwise disjoint, and together they cover a dense subset of $\mathcal{A}_{n,k,4}(Z)$.
\end{conj}

Interestingly, the number of BCFW cells in $\Gr_{k,n}^{\geq 0}$ is a {\itshape Narayana number}. The Narayana numbers $N_{a,b} := \frac{1}{a}\binom{a}{b}\binom{a}{b-1}$ refine the Catalan numbers, i.e.\ $\sum_{b=1}^a N_{a,b}$ is the Catalan number $\frac{1}{a+1}\binom{2a}{a}$ \cite[A46]{stanley_catalan}.

\begin{lem}[{\cite[(17.7)/(16.8)]{abcgpt}}]\label{lem:bcfw-m=4}
For $m=4$, the number of $(k,n)$-BCFW cells is the Narayana number $N_{n-3,k+1} = \frac{1}{n-3}\binom{n-3}{k+1}\binom{n-3}{k}$.
\end{lem}

\subsection{Complete binary trees and a bijection to BCFW graphs}\label{sec:trees-to-graphs}

In this section we explain how to index the plabic graphs $\mathcal{\tilde{G}}_{n,k,4}$ coming from the BCFW recursion by {\itshape complete binary trees}. This construction also appears in a similar form in lecture notes written by Morales based on a course taught by Postnikov \cite[Figure 110]{morales}.
\begin{defn}\label{def:binary-tree}
A {\itshape complete} (or {\itshape plane}) {\itshape binary tree} $T$ is a rooted tree such that every vertex either has $2$ ordered child vertices, one joined by a horizontal (left) edge and the other by a vertical (up) edge, or $0$ child vertices. We require that the root vertex has $2$ children. (See \cref{fig:binary-tree} for an example.) We call childless vertices {\itshape leaves}, and other vertices {\itshape internal vertices}. An edge of $T$ is called {\itshape external} if it is incident to a leaf; otherwise it is called {\itshape internal}. If $T$ has $n-2$ leaves ($n\ge 4$), we label them by $2, 3, \dots, n-1$ clockwise, such that the leaf $2$ is joined to the root vertex by a path of horizontal edges (and similarly, $n-1$ is joined to the root by a path of vertical edges). We let $\mathcal{T}_{n,k,4}$ denote the set of complete binary trees $T$ with $n-2$ leaves, exactly $k+1$ of which are incident to a horizontal edge.
\end{defn}
\begin{figure}[ht]
\begin{center}
$$
\quad\begin{tikzpicture}[baseline=(current bounding box.center)]
\tikzstyle{out1}=[inner sep=0,minimum size=2.4mm,circle,draw=black,fill=black,semithick]
\tikzstyle{in1}=[inner sep=0,minimum size=2.4mm,circle,draw=black,fill=white,semithick]
\tikzstyle{vertex}=[inner sep=0,minimum size=1.2mm,circle,draw=black,fill=black,semithick]
\pgfmathsetmacro{\l}{0.48};
\pgfmathsetmacro{\d}{0.36};
\pgfmathsetmacro{\s}{1.44};
\node[vertex](r)at(0,0){};
\node[vertex](rh)at($(r)+(-\s,0)$){};
\node[vertex](rv)at($(r)+(0,\s)$){};
\node[inner sep=0](rhh)at($(rh)+(-\l,0)$){};
\node[inner sep=0](rhv)at($(rh)+(0,\l)$){};
\node[inner sep=0](rvh)at($(rv)+(-\l,0)$){};
\node[inner sep=0](rvv)at($(rv)+(0,\l)$){};
\node[inner sep=0]at($(rhh)+(-0.18,0)$){$2$};
\node[inner sep=0]at($(rhv)+(0,0.24)$){$3$};
\node[inner sep=0]at($(rvh)+(-0.18,0)$){$4$};
\node[inner sep=0]at($(rvv)+(0,0.24)$){$5$};
\path[thick](r)edge(rh) edge(rv) (rh)edge(rhh) edge(rhv) (rv)edge(rvh) edge(rvv);
\end{tikzpicture}\qquad\qquad\qquad\quad
\begin{tikzpicture}[baseline=(current bounding box.center)]
\tikzstyle{out1}=[inner sep=0,minimum size=2.4mm,circle,draw=black,fill=black,semithick]
\tikzstyle{in1}=[inner sep=0,minimum size=2.4mm,circle,draw=black,fill=white,semithick]
\tikzstyle{vertex}=[inner sep=0,minimum size=1.2mm,circle,draw=black,fill=black,semithick]
\pgfmathsetmacro{\s}{0.72};
\pgfmathsetmacro{\r}{0.36};
\pgfmathsetmacro{\l}{0.48};
\useasboundingbox(-2*\r-\s-\l,-\l-0.18)rectangle(\r+\s+\l,3*\r+2*\s+\l-0.36);
\node[in1](i1)at(0,0){};
\node[out1](i2)at($(i1)+(0,\s)$){};
\node[in1](i3)at($(i2)+(\r,\r)$){};
\node[out1](i4)at($(i3)+(\s,0)$){};
\node[out1](i5)at($(i4)+(0,-\s)$){};
\node[in1](i6)at($(i5)+(-\r,-\r)$){};
\node[inner sep=0](bn)at($(i5)+(\l,0)$){};
\node[inner sep=0](b1)at($(i6)+(0,-\l)$){};
\node[inner sep=0]at($(bn)+(0.18,0)$){$6$};
\node[inner sep=0]at($(b1)+(0,-0.24)$){$1$};
\node[out1](hur)at($(i2)+(-\r,\r)$){};
\node[in1](hul)at($(hur)+(-\s,0)$){};
\node[out1](hml)at($(hul)+(-\r,-\r)$){};
\node[in1](hbl)at($(hml)+(0,-\s)$){};
\node[inner sep=0](b2)at($(hml)+(-\l,0)$){};
\node[inner sep=0](b3)at($(hul)+(0,\l)$){};
\node[inner sep=0]at($(b2)+(-0.18,0)$){$2$};
\node[inner sep=0]at($(b3)+(0,0.24)$){$3$};
\node[out1](vbr)at($(i4)+(0,\s+2*\r)$){};
\node[in1](vmr)at($(vbr)+(-\s,0)$){};
\node[out1](vur)at($(vmr)+(-\r,-\r)$){};
\node[in1](vul)at($(vur)+(0,-\s)$){};
\node[inner sep=0](b4)at($(vur)+(-\l,0)$){};
\node[inner sep=0](b5)at($(vmr)+(0,\l)$){};
\node[inner sep=0]at($(b4)+(-0.18,0)$){$4$};
\node[inner sep=0]at($(b5)+(0,0.24)$){$5$};
\path[thick](i1)edge(i2) (i2)edge(i3) (i3)edge(i4) (i4)edge(i5) (i5)edge(i6) (i6)edge(i1) (i5)edge(bn) (i6)edge(b1) (i2)edge(hur) (hur)edge(hul) (hul)edge(hml) (hml)edge(hbl) (hbl)edge(i1) (hml)edge(b2) (hul)edge(b3) (vbr)edge(vmr) (vmr)edge(vur) (vur)edge(vul) (vul)edge(i3) (i4)edge(vbr) (vur)edge(b4) (vmr)edge(b5);
\end{tikzpicture}\qquad\qquad\qquad\quad
\begin{tikzpicture}[baseline=(current bounding box.center)]
\tikzstyle{out1}=[inner sep=0,minimum size=2.4mm,circle,draw=black,fill=black,semithick]
\tikzstyle{in1}=[inner sep=0,minimum size=2.4mm,circle,draw=black,fill=white,semithick]
\tikzstyle{vertex}=[inner sep=0,minimum size=1.2mm,circle,draw=black,fill=black,semithick]
\pgfmathsetmacro{\s}{0.72};
\pgfmathsetmacro{\r}{0.36};
\pgfmathsetmacro{\l}{0.48};
\useasboundingbox(-2*\r-\s-\l,-\l-0.18)rectangle(\r+\s+\l,3*\r+2*\s+\l-0.36);
\node[in1](i1)at(0,0){};
\node[out1](i2)at($(i1)+(0,\s)$){};
\node[in1](i3)at($(i2)+(\r,\r)$){};
\node[out1](i4)at($(i3)+(\s,0)$){};
\node[out1](i5)at($(i4)+(0,-\s)$){};
\node[in1](i6)at($(i5)+(-\r,-\r)$){};
\node[inner sep=0](bn)at($(i5)+(\l,0)$){};
\node[inner sep=0](b1)at($(i6)+(0,-\l)$){};
\node[inner sep=0]at($(bn)+(0.18,0)$){$6$};
\node[inner sep=0]at($(b1)+(0,-0.24)$){$1$};
\node[inner sep=0](hur)at($(i2)+(-\r,\r)$){};
\node[in1](hul)at($(hur)+(-\s,0)$){};
\node[out1](hml)at($(hul)+(-\r,-\r)$){};
\node[inner sep=0](hbl)at($(hml)+(0,-\s)$){};
\node[inner sep=0](b2)at($(hml)+(-\l,0)$){};
\node[inner sep=0](b3)at($(hul)+(0,\l)$){};
\node[inner sep=0]at($(b2)+(-0.18,0)$){$2$};
\node[inner sep=0]at($(b3)+(0,0.24)$){$3$};
\node[inner sep=0](vbr)at($(i4)+(0,\s+2*\r)$){};
\node[in1](vmr)at($(vbr)+(-\s,0)$){};
\node[out1](vur)at($(vmr)+(-\r,-\r)$){};
\node[inner sep=0](vul)at($(vur)+(0,-\s)$){};
\node[inner sep=0](b4)at($(vur)+(-\l,0)$){};
\node[inner sep=0](b5)at($(vmr)+(0,\l)$){};
\node[inner sep=0]at($(b4)+(-0.24,0.04)$){$4$};
\node[inner sep=0]at($(b5)+(0,0.24)$){$5$};
\node[inner sep=0](b12)at($(b1)+(-2*\s-2*\r,\r)$){};
\node[inner sep=0](b5n)at($(bn)+(-\r,2*\s+2*\r)$){};
\path[thick](i1)edge(i2) (i2)edge(i3) (i3)edge(i4) (i4)edge(i5) (i5)edge(i6) (i6)edge(i1) (i5)edge(bn.center) (i6)edge(b1.center) (i2)edge(hul) (hul)edge(hml) (hml)edge(i1) (hml)edge(b2.center) (hul)edge(b3.center) (vmr)edge(vur) (vur)edge(i3) (i4)edge(vmr) (vur)edge(b4.center) (vmr)edge(b5.center);
\draw[thick]plot[smooth cycle,tension=0.6]coordinates{(b1) (b12) (b2) (b3) (b4) (b5) (b5n) (bn)};
\end{tikzpicture}
$$
\caption{A complete binary tree $T\in\mathcal{T}_{6,1,4}$ and its plabic graph $G(T)\in\mathcal{\tilde{G}}_{6,3,4}$ (before and after contracting bivalent vertices and enclosing the graph in a disk). Note that $G(T)$ is the middle graph in the bottom row of \cref{fig:BCFW-graphs}.}
\label{fig:binary-tree}
\end{center}
\end{figure}
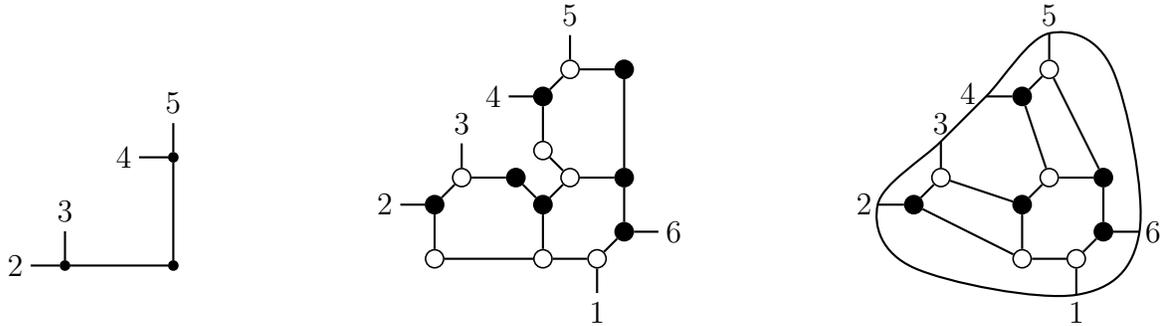
\begin{defn}\label{def:tree-to-graph}
Given a complete binary tree $T\in\mathcal{T}_{n,k,4}$, we define a plabic graph $G(T)$ with $n$ boundary vertices, as follows.
\begin{itemize}
\item We replace the root vertex of $T$ with a face, as shown:
$$
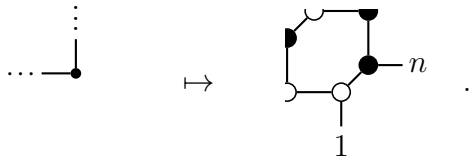
\begin{tikzpicture}[baseline=(current bounding box.center)]
\tikzstyle{out1}=[inner sep=0,minimum size=2.4mm,circle,draw=black,fill=black,semithick]
\tikzstyle{in1}=[inner sep=0,minimum size=2.4mm,circle,draw=black,fill=white,semithick]
\tikzstyle{vertex}=[inner sep=0,minimum size=1.2mm,circle,draw=black,fill=black,semithick]
\pgfmathsetmacro{\l}{0.48};
\pgfmathsetmacro{\d}{0.36};
\useasboundingbox(-\l-\d,-\l-\d-0.48)rectangle(0.48,\l+\d);
\node[vertex](v)at(0,0){};
\node[inner sep=0](N)at(0,\l){};
\node[inner sep=0](E)at(\l,0){};
\node[inner sep=0](S)at(0,-\l){};
\node[inner sep=0](W)at(-\l,0){};
\path[thick](v)edge(N) edge(W);
\node[inner sep=0]at($(N)+(0,\d)$){$\vdots$};
\node[inner sep=0,rotate=90]at($(W)+(-\d,0)$){$\vdots$};
\end{tikzpicture}\qquad\mapsto\qquad
\begin{tikzpicture}[baseline=(current bounding box.center)]
\tikzstyle{out1}=[inner sep=0,minimum size=2.4mm,circle,draw=black,fill=black,semithick]
\tikzstyle{in1}=[inner sep=0,minimum size=2.4mm,circle,draw=black,fill=white,semithick]
\tikzstyle{vertex}=[inner sep=0,minimum size=1.2mm,circle,draw=black,fill=black,semithick]
\pgfmathsetmacro{\s}{0.72};
\pgfmathsetmacro{\r}{0.36};
\pgfmathsetmacro{\l}{0.48};
\pgfmathsetmacro{\epsilon}{0.02};
\clip(-\epsilon,-\l-0.6)rectangle(\s+\r+\l+0.54,\s+\r+\epsilon);
\node[in1](i1)at(0,0){};
\node[out1](i2)at($(i1)+(0,\s)$){};
\node[in1](i3)at($(i2)+(\r,\r)$){};
\node[out1](i4)at($(i3)+(\s,0)$){};
\node[out1](i5)at($(i4)+(0,-\s)$){};
\node[in1](i6)at($(i5)+(-\r,-\r)$){};
\node[inner sep=0](bn)at($(i5)+(\l,0)$){};
\node[inner sep=0](b1)at($(i6)+(0,-\l)$){};
\node[inner sep=0]at($(bn)+(0.18,0)$){$n$};
\node[inner sep=0]at($(b1)+(0,-0.24)$){$1$};
\path[thick](i1)edge(i2) (i2)edge(i3) (i3)edge(i4) (i4)edge(i5) (i5)edge(i6) (i6)edge(i1) (i5)edge(bn) (i6)edge(b1);
\end{tikzpicture}\;\;.\vspace*{-4pt}
$$
(As we will explain shortly, the `half-vertices' on the right are intentional.)
\item We replace each internal vertex of $T$ with a face, as shown:
$$
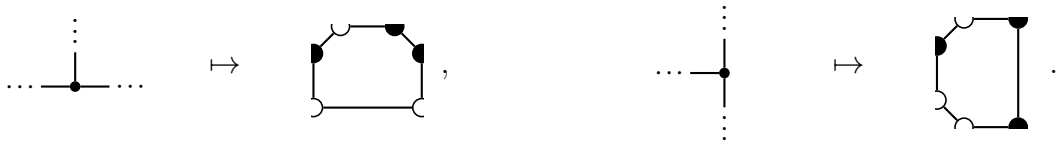
\begin{tikzpicture}[baseline=(current bounding box.center)]
\tikzstyle{out1}=[inner sep=0,minimum size=2.4mm,circle,draw=black,fill=black,semithick]
\tikzstyle{in1}=[inner sep=0,minimum size=2.4mm,circle,draw=black,fill=white,semithick]
\tikzstyle{vertex}=[inner sep=0,minimum size=1.2mm,circle,draw=black,fill=black,semithick]
\pgfmathsetmacro{\l}{0.48};
\pgfmathsetmacro{\d}{0.36};
\useasboundingbox(-\l-\d,-0.48)rectangle(\l+\d,\l+\d);
\node[vertex](v)at(0,0){};
\node[inner sep=0](N)at(0,\l){};
\node[inner sep=0](E)at(\l,0){};
\node[inner sep=0](S)at(0,-\l){};
\node[inner sep=0](W)at(-\l,0){};
\path[thick](v)edge(N) edge(W) edge(E);
\node[inner sep=0]at($(N)+(0,\d)$){$\vdots$};
\node[inner sep=0,rotate=-90]at($(E)+(\d,0)$){$\vdots$};
\node[inner sep=0,rotate=90]at($(W)+(-\d,0)$){$\vdots$};
\end{tikzpicture}\qquad\mapsto\qquad
\begin{tikzpicture}[baseline=(current bounding box.center)]
\tikzstyle{out1}=[inner sep=0,minimum size=2.4mm,circle,draw=black,fill=black,semithick]
\tikzstyle{in1}=[inner sep=0,minimum size=2.4mm,circle,draw=black,fill=white,semithick]
\tikzstyle{vertex}=[inner sep=0,minimum size=1.2mm,circle,draw=black,fill=black,semithick]
\pgfmathsetmacro{\s}{0.72};
\pgfmathsetmacro{\r}{0.36};
\pgfmathsetmacro{\epsilon}{0.02};
\clip(-\s-2*\r-\epsilon,\s+\r+\epsilon)rectangle(\epsilon,-0.18);
\node[in1](br)at(0,0){};
\node[out1](mr)at($(br)+(0,\s)$){};
\node[out1](ur)at($(mr)+(-\r,\r)$){};
\node[in1](ul)at($(ur)+(-\s,0)$){};
\node[out1](ml)at($(ul)+(-\r,-\r)$){};
\node[in1](bl)at($(ml)+(0,-\s)$){};
\path[thick](br)edge(mr) (mr)edge(ur) (ur)edge(ul) (ul)edge(ml) (ml)edge(bl) (bl)edge(br);
\end{tikzpicture}\;\;,\qquad\qquad\qquad
\begin{tikzpicture}[baseline=(current bounding box.center)]
\tikzstyle{out1}=[inner sep=0,minimum size=2.4mm,circle,draw=black,fill=black,semithick]
\tikzstyle{in1}=[inner sep=0,minimum size=2.4mm,circle,draw=black,fill=white,semithick]
\tikzstyle{vertex}=[inner sep=0,minimum size=1.2mm,circle,draw=black,fill=black,semithick]
\pgfmathsetmacro{\l}{0.48};
\pgfmathsetmacro{\d}{0.36};
\node[vertex](v)at(0,0){};
\node[inner sep=0](N)at(0,\l){};
\node[inner sep=0](E)at(\l,0){};
\node[inner sep=0](S)at(0,-\l){};
\node[inner sep=0](W)at(-\l,0){};
\path[thick](v)edge(N) edge(W) edge(S);
\node[inner sep=0]at($(N)+(0,\d)$){$\vdots$};
\node[inner sep=0,rotate=180]at($(S)+(0,-\d)$){$\vdots$};
\node[inner sep=0,rotate=90]at($(W)+(-\d,0)$){$\vdots$};
\end{tikzpicture}\qquad\mapsto\qquad
\begin{tikzpicture}[baseline=(current bounding box.center)]
\tikzstyle{out1}=[inner sep=0,minimum size=2.4mm,circle,draw=black,fill=black,semithick]
\tikzstyle{in1}=[inner sep=0,minimum size=2.4mm,circle,draw=black,fill=white,semithick]
\tikzstyle{vertex}=[inner sep=0,minimum size=1.2mm,circle,draw=black,fill=black,semithick]
\pgfmathsetmacro{\s}{0.72};
\pgfmathsetmacro{\r}{0.36};
\pgfmathsetmacro{\epsilon}{0.02};
\clip(-\s-\r-\epsilon,-\s-2*\r-\epsilon)rectangle(0.18,\epsilon);
\node[out1](br)at(0,0){};
\node[in1](mr)at($(br)+(-\s,0)$){};
\node[out1](ur)at($(mr)+(-\r,-\r)$){};
\node[in1](ul)at($(ur)+(0,-\s)$){};
\node[in1](ml)at($(ul)+(\r,-\r)$){};
\node[out1](bl)at($(ml)+(\s,0)$){};
\path[thick](br)edge(mr) (mr)edge(ur) (ur)edge(ul) (ul)edge(ml) (ml)edge(bl) (bl)edge(br);
\end{tikzpicture}\;\;.
$$
\item We replace each leaf of $T$, as shown:
$$
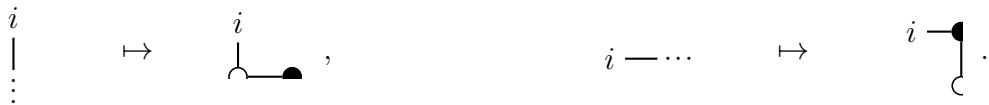
\quad\begin{tikzpicture}[baseline=(current bounding box.center)]
\tikzstyle{out1}=[inner sep=0,minimum size=2.4mm,circle,draw=black,fill=black,semithick]
\tikzstyle{in1}=[inner sep=0,minimum size=2.4mm,circle,draw=black,fill=white,semithick]
\tikzstyle{vertex}=[inner sep=0,minimum size=1.2mm,circle,draw=black,fill=black,semithick]
\pgfmathsetmacro{\l}{0.48};
\pgfmathsetmacro{\d}{0.36};
\node[inner sep=0](v)at(0,0){};
\node[inner sep=0](N)at(0,\l){};
\node[inner sep=0](E)at(\l,0){};
\node[inner sep=0](S)at(0,-\l){};
\node[inner sep=0](W)at(-\l,0){};
\node[inner sep=0]at($(v)+(0,0.24)$){$i$};
\path[thick](v)edge(S);
\node[inner sep=0,rotate=180]at($(S)+(0,-\d)$){$\vdots$};
\end{tikzpicture}\qquad\mapsto\qquad
\begin{tikzpicture}[baseline=(current bounding box.center)]
\tikzstyle{out1}=[inner sep=0,minimum size=2.4mm,circle,draw=black,fill=black,semithick]
\tikzstyle{in1}=[inner sep=0,minimum size=2.4mm,circle,draw=black,fill=white,semithick]
\tikzstyle{vertex}=[inner sep=0,minimum size=1.2mm,circle,draw=black,fill=black,semithick]
\pgfmathsetmacro{\s}{0.72};
\pgfmathsetmacro{\r}{0.36};
\pgfmathsetmacro{\l}{0.48};
\pgfmathsetmacro{\epsilon}{0.02};
\useasboundingbox(-0.18,\l+0.6)rectangle(\s+0.18,-0.6);
\clip(-0.18,\l+0.6)rectangle(\s+0.18,-\epsilon);
\node[in1](l)at(0,0){};
\node[out1](r)at($(l)+(\s,0)$){};
\node[inner sep=0](b)at($(l)+(0,\l)$){};
\node[inner sep=0]at($(b)+(0,0.24)$){$i$};
\path[thick](l)edge(r) edge(b);
\end{tikzpicture}\;\;,\qquad\qquad\qquad\qquad
\begin{tikzpicture}[baseline=(current bounding box.center)]
\tikzstyle{out1}=[inner sep=0,minimum size=2.4mm,circle,draw=black,fill=black,semithick]
\tikzstyle{in1}=[inner sep=0,minimum size=2.4mm,circle,draw=black,fill=white,semithick]
\tikzstyle{vertex}=[inner sep=0,minimum size=1.2mm,circle,draw=black,fill=black,semithick]
\pgfmathsetmacro{\l}{0.48};
\pgfmathsetmacro{\d}{0.36};
\node[inner sep=0](v)at(0,0){};
\node[inner sep=0](N)at(0,\l){};
\node[inner sep=0](E)at(\l,0){};
\node[inner sep=0](S)at(0,-\l){};
\node[inner sep=0](W)at(-\l,0){};
\node[inner sep=0]at($(v)+(-0.18,0)$){$i$};
\path[thick](v)edge(E);
\node[inner sep=0,rotate=-90]at($(E)+(\d,0)$){$\vdots$};
\end{tikzpicture}\qquad\mapsto\qquad
\begin{tikzpicture}[baseline=(current bounding box.center)]
\tikzstyle{out1}=[inner sep=0,minimum size=2.4mm,circle,draw=black,fill=black,semithick]
\tikzstyle{in1}=[inner sep=0,minimum size=2.4mm,circle,draw=black,fill=white,semithick]
\tikzstyle{vertex}=[inner sep=0,minimum size=1.2mm,circle,draw=black,fill=black,semithick]
\pgfmathsetmacro{\s}{0.72};
\pgfmathsetmacro{\r}{0.36};
\pgfmathsetmacro{\l}{0.48};
\pgfmathsetmacro{\epsilon}{0.02};
\clip(-\l-0.6,0.18)rectangle(\epsilon,-\s-0.18);
\node[out1](r)at(0,0){};
\node[in1](l)at($(r)+(0,-\s)$){};
\node[inner sep=0](b)at($(r)+(-\l,0)$){};
\node[inner sep=0]at($(b)+(-0.18,0)$){$i$};
\path[thick](r)edge(l) edge(b);
\end{tikzpicture}\;\;.
$$
In both of the resulting local pictures of a plabic graph above, the vertex not incident to $i$ is incident to exactly one other vertex (which appears just outside this local picture), of the same color; we contract these two vertices. (This is to avoid having any bivalent vertices in the resulting plabic graph.)
\item We draw a curve through $1, \dots, n$, enclosing the resulting graph in a disk.
\end{itemize}
This gives a plabic graph $G(T)$. See \cref{fig:binary-tree} for an example. (Above, we have depicted only `half vertices' of $G(T)$, since each internal vertex of $G(T)$, aside from those incident to $1$ and $n$, comes from two vertices of $T$.)
\end{defn}

\begin{lem}\label{trees_to_graphs}
The map $T\mapsto G(T)$ from \cref{def:tree-to-graph} gives a bijection $\mathcal{T}_{n,k,4}\to\mathcal{\tilde{G}}_{n,k+2,4}$.
\end{lem}

\begin{pf}
It follows from \cref{def:tree-to-graph} that adding child vertices to a leaf $i$ of $T$ corresponds to blowing up $G(T)$ at $i$. That is, the BCFW recursion acts on complete binary trees by adding children.
\end{pf}
Note that it is straightforward to recover $T$ from $G(T)$: the graph formed by the internal edges of $T$ is dual to the graph formed by the internal faces of $G(T)$.

\section{Pairs of noncrossing lattice paths and BCFW \texorpdfstring{\Le}{Le}-diagrams}\label{sec:lattice-paths}

\noindent In this section, we index the $m=4$ BCFW cells by pairs of noncrossing lattice paths inside a rectangle. We explain how to obtain a $\oplus$-diagram of a BCFW cell from such a pair.

\begin{defn}\label{def:noncrossing}
Fix $a,b\in\mathbb{N}$. A {\itshape lattice path $W$} inside an $a\times b$ rectangle is a path that moves from the northeast corner to the southwest corner, taking unit steps west and south. We represent $W$ by a word of length $a+b$ on the alphabet $\{H,V\}$, with exactly $a$ letters $V$ (corresponding to the vertical steps) and $b$ letters $H$ (corresponding to the horizontal steps). For example, the upper lattice path $W_U$ in \cref{fig:PathsToLeDiagram} is given by $W_U = HHVHHVVH$.

Let $\mathcal{L}_{n,k,4}$ denote the set of all pairs $(W_U, W_L)$ of {\itshape noncrossing} lattice paths inside a $k\times (n-k-4)$ rectangle, where $W_U$ denotes the upper path and $W_L$ denotes the lower path. That is, $W_L$ is weakly below $W_U$ (but the two paths are allowed to overlap); see \cref{fig:PathsToLeDiagram}. In terms of words, this means that for any $i\in [n-4]$, there are at least as many $V$'s among the first $i$ letters of $W_L$ as among the first $i$ letters of $W_U$.
\end{defn}

We give an injective map from elements of $\mathcal{L}_{n,k,4}$ to
$4k$-dimensional positroid cells of $\Gr_{k,n}^{\geq 0}$.  We will prove that the collection of cells in its image are precisely the $(k,n)$-BCFW cells.
\afterpage{\begin{figure}[p!]
\begin{align*}
& (W_U, W_L) =
\tikzexternalenable\begin{tikzpicture}[baseline=(current bounding box.center)]
\tikzstyle{hu}=[top color=red!10,bottom color=red!10,middle color=red,opacity=0.70]
\tikzstyle{vu}=[left color=red!10,right color=red!10,middle color=red,opacity=0.70]
\tikzstyle{hl}=[top color=blue!10,bottom color=blue!10,middle color=blue,opacity=0.55]
\tikzstyle{vl}=[left color=blue!10,right color=blue!10,middle color=blue,opacity=0.55]
\pgfmathsetmacro{\u}{0.75};
\pgfmathsetmacro{\w}{0.12*\u};
\coordinate(h)at(-\u,0);
\coordinate(v)at(0,-\u);
\draw[step=\u,color=black!16,ultra thick](0,0)grid(5*\u,3*\u);
\node[inner sep=0](l1)at(5*\u,3*\u){};
\node[inner sep=0](l2)at($(l1)+(v)$){};
\node[inner sep=0](l3)at($(l2)+(v)$){};
\node[inner sep=0](l4)at($(l3)+(h)$){};
\node[inner sep=0](l5)at($(l4)+(v)$){};
\node[inner sep=0](l6)at($(l5)+(h)$){};
\node[inner sep=0](l7)at($(l6)+(h)$){};
\node[inner sep=0](l8)at($(l7)+(h)$){};
\node[inner sep=0](l9)at($(l8)+(h)$){};
\node[inner sep=0](u1)at(5*\u,3*\u){};
\node[inner sep=0](u2)at($(u1)+(h)$){};
\node[inner sep=0](u3)at($(u2)+(h)$){};
\node[inner sep=0](u4)at($(u3)+(v)$){};
\node[inner sep=0](u5)at($(u4)+(h)$){};
\node[inner sep=0](u6)at($(u5)+(h)$){};
\node[inner sep=0](u7)at($(u6)+(v)$){};
\node[inner sep=0](u8)at($(u7)+(v)$){};
\node[inner sep=0](u9)at($(u8)+(h)$){};
\begin{scope}
\clip($(u1.center)+(0,\w)$)--++($2*(h)+(-\w,0)$)--++(2*\w,-2*\w)--++($-2*(h)+(-\w,0)$);
\path[hu]($(u1.center)+(\w,\w)$)rectangle($(u3.center)+(-\w,-\w)$);
\end{scope}
\begin{scope}
\clip($(u3.center)+(-\w,\w)$)--++($1*(v)$)--++(2*\w,-2*\w)--++($-1*(v)$);
\path[vu]($(u3.center)+(\w,\w)$)rectangle($(u4.center)+(-\w,-\w)$);
\end{scope}
\begin{scope}
\clip($(u4.center)+(-\w,\w)$)--++($2*(h)$)--++(2*\w,-2*\w)--++($-2*(h)$);
\path[hu]($(u4.center)+(\w,\w)$)rectangle($(u6.center)+(-\w,-\w)$);
\end{scope}
\begin{scope}
\clip($(u6.center)+(-\w,\w)$)--++($2*(v)$)--++(2*\w,-2*\w)--++($-2*(v)$);
\path[vu]($(u6.center)+(\w,\w)$)rectangle($(u8.center)+(-\w,-\w)$);
\end{scope}
\begin{scope}
\clip($(u8.center)+(-\w,\w)$)--++($1*(h)+(\w,0)$)--++(0,-2*\w)--++($-1*(h)+(\w,0)$);
\path[hu]($(u8.center)+(\w,\w)$)rectangle($(u9.center)+(-\w,-\w)$);
\end{scope}
\begin{scope}
\clip($(l1.center)+(-\w,0)$)--++($2*(v)+(0,\w)$)--++(2*\w,-2*\w)--++($-2*(v)+(0,\w)$);
\path[vl]($(l1.center)+(\w,\w)$)rectangle($(l3.center)+(-\w,-\w)$);
\end{scope}
\begin{scope}
\clip($(l3.center)+(-\w,\w)$)--++($1*(h)$)--++(2*\w,-2*\w)--++($-1*(h)$);
\path[hl]($(l3.center)+(\w,\w)$)rectangle($(l4.center)+(-\w,-\w)$);
\end{scope}
\begin{scope}
\clip($(l4.center)+(-\w,\w)$)--++($1*(v)$)--++(2*\w,-2*\w)--++($-1*(v)$);
\path[vl]($(l4.center)+(\w,\w)$)rectangle($(l5.center)+(-\w,-\w)$);
\end{scope}
\begin{scope}
\clip($(l5.center)+(-\w,\w)$)--++($4*(h)+(\w,0)$)--++(0,-2*\w)--++($-4*(h)+(\w,0)$);
\path[hl]($(l5.center)+(\w,\w)$)rectangle($(l9.center)+(-\w,-\w)$);
\end{scope}
\end{tikzpicture}\tikzexternaldisable\qquad\qquad
Y_U = \ydiagram{3,1,1}\qquad\qquad
Y_L = \ydiagram{5,5,4} \\
& \quad\textnormal{\bfseries Step 1.}\qquad \begin{ytableau}
\none \\
\none \\
+ &   &   &   &   &   &   &   & + \\
+ &   &   &   &   &   &   &   & + \\
+ &   &   &   &   &   &   & + 
\end{ytableau} \\[-12pt]
& \quad\textnormal{\bfseries Step 2.}\qquad \begin{ytableau}
\none \\
\none \\
+ &   &   &   &   &   & 0 &   & + \\
+ &   &   &   &   &   & 0 &   & + \\
+ &   &   &   &   &   & 0 & + 
\end{ytableau} \\[-12pt]
& \quad\textnormal{\bfseries Step 2.}\qquad 
\begin{ytableau}
\none \\
\none \\
+ &   &   &   &   &   & 0 & 0 & + \\
+ &   &   &   &   &   & 0 &   & + \\
+ &   &   &   &   &   & 0 & + 
\end{ytableau} \\[-12pt]
& \quad\textnormal{\bfseries Step 2.}\qquad \begin{ytableau}
\none \\
\none \\
+ &   &   &   &   & 0 & 0 & 0 & + \\
+ &   &   &   &   &   & 0 &   & + \\
+ &   &   &   &   &   & 0 & + 
\end{ytableau} \\[-12pt]
& \quad\textnormal{\bfseries Step 3.}\qquad \begin{ytableau}
\none \\
\none \\
+ &   &   & + & + & 0 & 0 & 0 & + \\
+ &   &   &   &   & + & 0 & + & + \\
+ &   &   &   & + & + & 0 & + 
\end{ytableau} \\[-12pt]
& \quad\textnormal{\bfseries Step 4.}\qquad \begin{ytableau}
\none \\
\none \\
+ & 0 & 0 & + & + & 0 & 0 & 0 & + \\
+ & 0 & 0 & 0 & 0 & + & 0 & + & + \\
+ & 0 & 0 & 0 & + & + & 0 & + 
\end{ytableau}
\end{align*}
\caption{The map $\Omega_{\mathcal{L} \mathcal{D}}(W_U, W_L)$ from \cref{def:latticetocell} takes a pair of lattice paths $(W_U, W_L)\in\mathcal{L}_{n,k,4}$ to a reduced $\oplus$-diagram in $\mathcal{D}_{n,k,4}$, corresponding to a $4k$-dimensional cell of $\Gr_{k,n}^{\ge 0}$. Here $k=3$ and $n=12$.}
\label{fig:PathsToLeDiagram}
\end{figure}\clearpage}
\begin{defn}\label{def:latticetocell}
Given $(W_U, W_L)\in\mathcal{L}_{n,k,4}$, let $Y_U$ be the Young diagram inside a $k\times (n-k-4)$ rectangle whose southeast border is $Y_U$, and similarly define $Y_L$. We associate a $\oplus$-diagram $D$ of type $(k,n)$ (recall \cref{def:oplus}) to $(W_U, W_L)$ as follows. (See \cref{fig:PathsToLeDiagram} for an example.)
\begin{enumerate}[leftmargin=*,widest=Step 1]
\item[Step 1.] The Young diagram of $D$ is obtained from $Y_L$ by adding $m=4$ extra columns at the left of height $k$. Place a $+$ at the far left end and far right end of each row of $D$, leaving the other boxes empty.
\item[Step 2.] Consider each column of $Y_U$ in turn, reading the columns from left to right. For each column of $Y_U$ of height $i$, place a top-justified column of $i$ $0$'s in $D$ as far right as possible.
\item[Step 3.] In each row of $D$, place two $+$'s as far to the right as possible.
\item[Step 4.] Fill any remaining empty boxes of $D$ with a $0$.
\end{enumerate}
We denote $D$ by 
$\Omega_{\mathcal{L} \mathcal{D}}(W_U, W_L)$.  
This defines an injection 
$\Omega_{\mathcal{L} \mathcal{D}}(W_U, W_L)$  
from $\mathcal{L}_{n,k,4}$ to the set of $\oplus$-diagrams of type $(k,n)$. Let $\mathcal{D}_{n,k,4}$ be the image of 
$\Omega_{\mathcal{L} \mathcal{D}}(W_U, W_L)$. 
 We will show in \cref{lem:reduced} that each $D\in\mathcal{D}_{n,k,4}$ is a reduced $\oplus$-diagram with exactly $4k$ $+$'s. Hence by \cref{le_moves}, $D$ corresponds to a $4k$-dimensional cell of $\Gr_{k,n}^{\ge 0}$, and we can use \Le -moves to find the \Le -diagram of $D$. We will also show in \cref{lem:reduced} that these $4k$-dimensional cells are all distinct.
\end{defn}

\begin{thm}\label{thm:BCFW}
The $\oplus$-diagrams $\mathcal{D}_{n,k,4}$ index the $(k,n)$-BCFW cells $\mathcal{C}_{n,k,4}$.
\end{thm}
\cref{thm:BCFW} follows from \cref{trees_to_graphs} and \cref{thm:shift_map_holds}, the latter of which we will prove in \cref{sec:bijection}. First, we show that the $\oplus$-diagrams in $\mathcal{D}_{n,k,4}$ are reduced and represent distinct positroid cells.
\begin{lem}\label{lem:reduced}
The $\oplus$-diagrams in $\mathcal{D}_{n,k,4}$ are reduced, and each correspond to a  distinct $4k$-dimensional cell of $\Gr_{k,n}^{\ge 0}$.
\end{lem}

\begin{pf}
We show that each $D\in\mathcal{D}_{n,k,4}$ is reduced. Let $D'$ be the $\oplus$-diagram of type $(k,n-1)$ obtained from $D$ by deleting the leftmost column. Then $D'$ is a \Le -diagram: each $0$ added in Step 2 of \cref{def:latticetocell} has no $+$ above it in the same column, and each $0$ added in Step 4 has no $+$ to its left in the same row (except for the $+$ at the left end of the row, which has been deleted to form $D'$ from $D$). Hence $D'$ is reduced, so its pipe dream $P(D')$ has no double crossings. We form $P(D)$ from $P(D')$ by adding a column of elbows at the left, which introduces no new crossings. Hence $D$ is reduced. Since $D$ has $4k$ $+$'s, it corresponds to a $4k$-dimensional cell of $\Gr_{k,n}^{\ge 0}$.

Now we show that the diagrams in $\mathcal{D}_{n,k,4}$ index distinct cells of $\Gr_{k,n}^{\ge 0}$. Suppose that $D_1, D_2\in\mathcal{D}_{n,k,4}$ index the same cell of $\Gr_{k,n}^{\ge 0}$, and as above let $D_1'$ and $D_2'$ be the \Le -diagrams of type $(k,n-1)$ formed by deleting the leftmost column of $D_1$ and $D_2$, respectively. In the construction of the pipe dreams $P(D_1)$ and $P(D_2)$ from \cref{def:oplus}, the edges of the west border of the Young diagram in each pipe dream are labeled with the anti-excedances of the corresponding decorated permutation. Since $\pi_{D_1} = \pi_{D_2}$, the pipe dreams $P(D_1)$ and $P(D_2)$ have the same shape and their edges are labeled in the same way. Hence $\pi_{D_1'} = \pi_{D_2'}$, whence $D_1' = D_2'$ by \cref{D_to_pi}. This implies $D_1 = D_2$.
\end{pf}

\section{From binary trees to pairs of lattice paths}\label{sec:bijection}

\noindent In this section we prove \cref{thm:BCFW}. 
Our strategy is to construct a map 
$\Omega_{\mathcal{T} \mathcal{L}}: \mathcal{T}_{n,k,4} \to 
\mathcal{L}_{n,k,4}$ 
which takes a complete binary tree $T$ to a pair of noncrossing lattice paths inside a $k\times (n-k-4)$ rectangle, such that the decorated permutation of the 
$\oplus$-diagram $\Omega_{\mathcal{L} \mathcal{D}}(\Omega_{\mathcal{T} \mathcal{L}}(T))$ equals 
$c_n^2 \pi_{G(T)}$ (recall \cref{def:BCFW-permutations}).
\begin{defn}\label{def:tree-to-paths}
Given a complete binary tree $T$, we let $r_H$ and $r_V$ denote the horizontal and vertical child vertices of the root $r$ of $T$. We let $T_H$ denote the subtree of $T$ rooted at $r$ obtained from $T$ by deleting all children of $r_V$. Also, if $r_H$ is an internal vertex of $T_H$, we let $T'_H$ be the subtree of $T$ rooted at $r_H$ formed from $T_H$ by deleting $r$ and its two incident edges. We similarly define the subtrees $T_V$ and $T'_V$ of $T$, by switching the roles of $r_H$ and $r_V$.

We associate two lattice paths $W_U(T), W_L(T)$ to $T$ by the following recursive definition. (Recall from \cref{def:noncrossing} that we identify a lattice path of length $l$ with its corresponding word in $\{H,V\}^l$.) Below, $\cdot$ denotes concatenation of words.
\begin{itemize}
\item If $r_H$ and $r_V$ are both leaves, then $W_U(T)$ and $W_L(T)$ are the empty words.
\item If $r_H$ is not a leaf and $r_V$ is a leaf, then
$$
\qquad W_U(T) := H\cdot W_U(T'_H), \qquad W_L(T) := W_L(T'_H)\cdot H.
$$
\item If $r_H$ is a leaf and $r_V$ is not a leaf, then
$$
\qquad W_U(T) := V\cdot W_U(T'_V), \qquad W_L(T) := V\cdot W_L(T'_V).
$$
\item Otherwise,
$$
\qquad W_U(T) = W_U(T_H)\cdot W_U(T_V), \qquad W_L(T) = W_L(T_H)\cdot W_L(T_V).
$$
\end{itemize}
Note that in fact $W_U(T) = W_U(T_H)\cdot W_U(T_V)$ and $W_L(T) = W_L(T_H)\cdot W_L(T_V)$ for all $T$.

We let $\Omega_{\mathcal{T} \mathcal{L}}: \mathcal{T}_{n,k,4} \to \mathcal{L}_{n,k,4}$ be the map which sends $T$ to $(W_U(T), W_L(T))$. It follows from the definitions of $W_U(T)$ and $W_L(T)$ that they are both words in $\{H,V\}^{n-4}$ with precisely $k$ $V$'s, and moreover that for any $i\in [n-4]$, there are at least as many $V$'s among the first $i$ letters of $W_L(T)$ as among the first $i$ letters of $W_U(T)$. Hence the pair of lattice paths $(W_U(T), W_L(T))$ is indeed noncrossing, and represents an element of $\mathcal{L}_{n,k,4}$.
\end{defn}

For example, for the complete binary tree $T$ in \cref{fig:finalV-example}, we have $W_U(T) = HVHVH$ and $W_L(T) = HVHHV$.
\begin{rmk}\label{readoff}
Given $T\in\mathcal{T}_{n,k,4}$, we can alternatively find $W_U(T), W_L(T)\in\{H,V\}^{n-4}$ as follows. We obtain $W_U(T)$ by reading the internal edges of $T$ in a depth-first search starting at the root, preferentially reading horizontal edges over vertical edges; we record an $H$ for each horizontal edge and a $V$ for each vertical edge. We obtain $W_L(T)$ by reading the leaves $3, 4, \dots, n-2$ of $T$ in order, recording an $H$ for each leaf incident to a vertical edge, and a $V$ for each leaf incident to a horizontal edge.
\end{rmk}

\begin{prop}\label{prop:lattice_bijection}
The map $\Omega_{\mathcal{T}\mathcal{L}}:\mathcal{T}_{n,k,4}\to\mathcal{L}_{n,k,4}$ is a bijection.
\end{prop}

\begin{pf}
It follows from \cref{lem:bcfw-m=4} and \cref{prop:objects}(1) that $\mathcal{T}_{n,k,4}$ and $\mathcal{L}_{n,k,4}$ have the same cardinality, namely the Narayana number $N_{n-3,k+1}$.\footnote{For further references on the enumeration of $\mathcal{L}_{n,k,4}$, see \cite[pp.\ 66-67]{stanley_catalan}.} Hence it suffices to show that $\Omega_{\mathcal{T}\mathcal{L}}$ is injective. We can prove this by induction on $n$, using the recursive definitions of $W_U(T)$ and $W_L(T)$. The key observation is that when we regard the pair of lattice paths $(W_U(T), W_L(T))$ as the concatenation of the pairs $(W_U(T_H), W_L(T_H))$ and $(W_U(T_V), W_L(T_V))$, we recover where $(W_U(T_V), W_L(T_V))$ begins as the first occurrence of overlapping vertical steps in $(W_U(T), W_L(T))$. (If there are no overlapping vertical steps, then $r_V$ is a leaf.) To see this, note that $W_U(T_V)$ and $W_L(T_V)$ both begin with a vertical step (if they have any steps at all); conversely, the paths $W_U(T_H)$ and $W_L(T_H)$ do not have any overlapping vertical steps, since otherwise the paths $W_U(T'_H)$ and $W_L(T'_H)$ would cross each other.
\end{pf}

Now we state and prove the main result of this section.
\begin{thm}\label{thm:shift_map_holds}
For $T\in\mathcal{T}_{n,k,4}$, we have
$$
\pi_D = c_n^2\pi_{G(T)},
$$
where $D = \Omega_{\mathcal{L} \mathcal{D}}(\Omega_{\mathcal{T} \mathcal{L}}(T))\in\mathcal{D}_{n,k,4}$ and $c_n := (n \;\; n{-}1 \; \cdots \; 2 \;\; 1)$ is the long cycle in the symmetric group on $[n]$. (As in \cref{def:BCFW-permutations}, we color the fixed points of $c_n^2\pi_{G(T)}$ black.)
\end{thm}
As we have already noted, this implies \cref{thm:BCFW}.
\begin{pf}
We proceed by induction on $n$. Given $T\in\mathcal{T}_{n,k,4}$, let $D := \Omega_{\mathcal{L} \mathcal{D}}(\Omega_{\mathcal{T} \mathcal{L}}(T))\in\mathcal{D}_{n,k,4}$ be its associated $\oplus$-diagram. Note that every row of $D$ contains at least one $+$, so $\pi_D$ has no white fixed points. Hence it suffices to show that the equality $\pi_D = c_n^2\pi_{G(T)}$ holds for (undecorated) permutations. If $n=4$, then $T$ is necessarily the tree $\begin{tikzpicture}[baseline=(current bounding box.center)]
\tikzstyle{out1}=[inner sep=0,minimum size=2.4mm,circle,draw=black,fill=black,semithick]
\tikzstyle{in1}=[inner sep=0,minimum size=2.4mm,circle,draw=black,fill=white,semithick]
\tikzstyle{vertex}=[inner sep=0,minimum size=1.2mm,circle,draw=black,fill=black,semithick]
\pgfmathsetmacro{\l}{0.48};
\pgfmathsetmacro{\d}{0.36};
\pgfmathsetmacro{\s}{1.44};
\node[vertex](r)at(0,0){};
\node[inner sep=0](rh)at($(r)+(-\l,0)$){};
\node[inner sep=0](rv)at($(r)+(0,\l)$){};
\path[thick](r)edge(rh) edge(rv);
\end{tikzpicture}\;$, and $D$ is the empty $\oplus$-diagram inside a $0\times 4$ rectangle. Hence $\pi_{G(T)} = 3412$ and $\pi_D = 1234 = c_4^2\pi_{G(T)}$. This proves the base case.

Now suppose that $n\ge 5$ and \cref{thm:shift_map_holds} holds for smaller values of $n$. We consider two different cases, depending on whether the last letter of $W_U(T)$ is $V$ or $H$.

{\bfseries Case 1: $W_U(T)$ ends in $V$.}
Recall from \cref{readoff} that we obtain $W_U(T)$ by reading the internal edges of $T$ in a depth-first search. Let $e$ be the internal edge of $T$ corresponding to the last letter of $W_U(T)$, and let $w$ be the child vertex incident to $e$ (i.e.\ $w$ is a child of the other vertex incident to $e$). Then the children of $w$ are both leaves, and are labeled by $p$ and $p+1$ for some $2 \le p \le n-2$. Let $T'\in\mathcal{T}_{n-1,k-1,4}$ be obtained from $T$ by replacing $e$ and its children by a vertical boundary edge incident to a leaf (which will be labeled by $p$). Since $W_U(T)$ ends in $V$, we obtain $(W_U(T'), W_L(T'))$ from $(W_U(T), W_L(T))$ by deleting the last vertical step of each path. Let $D'\in\mathcal{D}_{n-1,k-1,4}$ be the associated $\oplus$-diagram. (See \cref{fig:finalV-example}.) By the induction hypothesis, we have $\pi_{D'} = c_{n-1}^2\pi_{G(T')}$. (We will regard permutations of $[n-1]$ as permutations of $[n]$ which fix $n$.)
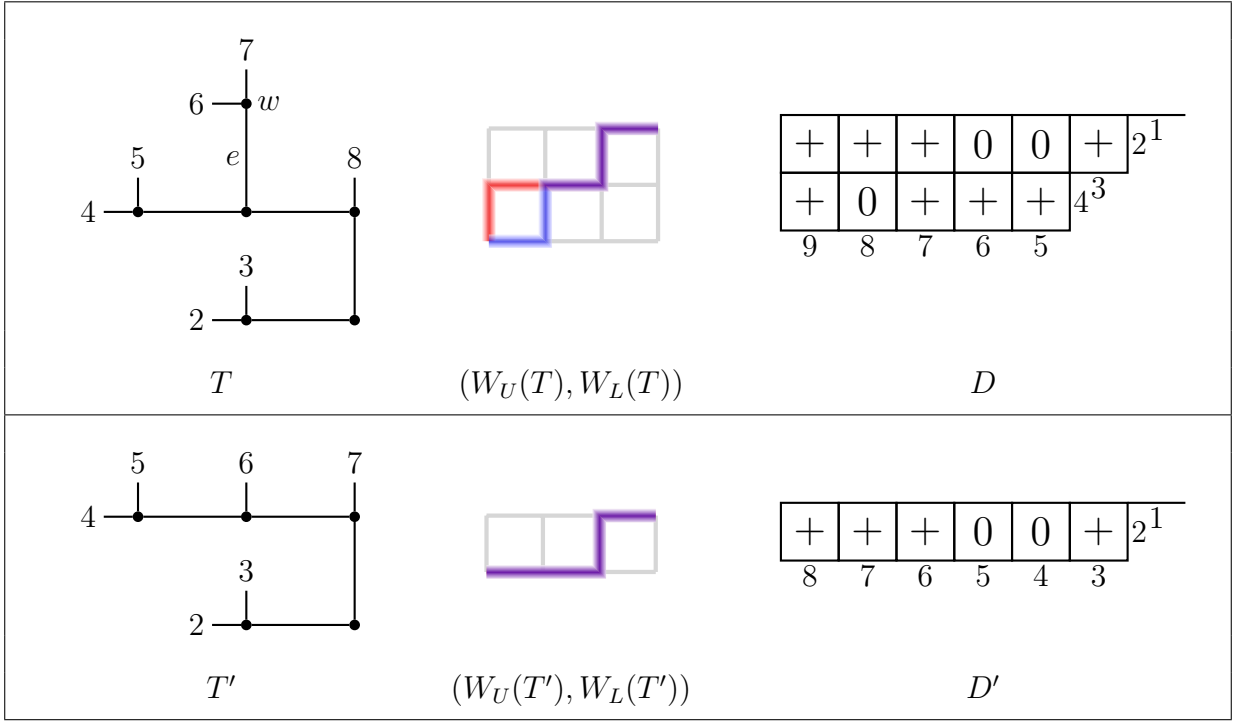
\begin{figure}[htb]
\begin{center}
\begin{tabular}{|ccc|}
\hline & & \\
\qquad\begin{tikzpicture}[baseline=(current bounding box.center)]
\tikzstyle{out1}=[inner sep=0,minimum size=2.4mm,circle,draw=black,fill=black,semithick]
\tikzstyle{in1}=[inner sep=0,minimum size=2.4mm,circle,draw=black,fill=white,semithick]
\tikzstyle{vertex}=[inner sep=0,minimum size=1.2mm,circle,draw=black,fill=black,semithick]
\pgfmathsetmacro{\l}{0.48};
\pgfmathsetmacro{\d}{0.36};
\pgfmathsetmacro{\s}{1.44};
\node[vertex](r)at(0,0){};
\node[vertex](rh)at($(r)+(-\s,0)$){};
\node[vertex](rv)at($(r)+(0,\s)$){};
\node[vertex](rvh)at($(rv)+(-\s,0)$){};
\node[vertex](rvhh)at($(rvh)+(-\s,0)$){};
\node[vertex](rvhv)at($(rvh)+(0,\s)$){};
\node[inner sep=0](rhh)at($(rh)+(-\l,0)$){};
\node[inner sep=0](rhv)at($(rh)+(0,\l)$){};
\node[inner sep=0](rvv)at($(rv)+(0,\l)$){};
\node[inner sep=0](rvhhh)at($(rvhh)+(-\l,0)$){};
\node[inner sep=0](rvhhv)at($(rvhh)+(0,\l)$){};
\node[inner sep=0](rvhvh)at($(rvhv)+(-\l,0)$){};
\node[inner sep=0](rvhvv)at($(rvhv)+(0,\l)$){};
\node[inner sep=0]at($(rhh)+(-0.18,0)$){$2$};
\node[inner sep=0]at($(rhv)+(0,0.24)$){$3$};
\node[inner sep=0]at($(rvhhh)+(-0.18,0)$){$4$};
\node[inner sep=0]at($(rvhhv)+(0,0.24)$){$5$};
\node[inner sep=0]at($(rvhvh)+(-0.18,0)$){$6$};
\node[inner sep=0]at($(rvhvv)+(0,0.24)$){$7$};
\node[inner sep=0]at($(rvv)+(0,0.24)$){$8$};
\node[inner sep=0]at($(rvhv)+(0.30,0)$){$w$};
\path[thick](r)edge(rh) edge(rv) (rh)edge(rhh) edge(rhv) (rv)edge(rvh) edge(rvv) (rvh)edge(rvhh) edge node[left=-2pt]{$e$}(rvhv) (rvhh)edge(rvhhh) edge(rvhhv) (rvhv)edge(rvhvh) edge(rvhvv);
\end{tikzpicture}\qquad
&
\qquad\tikzexternalenable\begin{tikzpicture}[baseline=(current bounding box.center)]
\tikzstyle{hu}=[top color=red!10,bottom color=red!10,middle color=red,opacity=0.70]
\tikzstyle{vu}=[left color=red!10,right color=red!10,middle color=red,opacity=0.70]
\tikzstyle{hl}=[top color=blue!10,bottom color=blue!10,middle color=blue,opacity=0.55]
\tikzstyle{vl}=[left color=blue!10,right color=blue!10,middle color=blue,opacity=0.55]
\pgfmathsetmacro{\u}{0.75};
\pgfmathsetmacro{\w}{0.12*\u};
\coordinate(h)at(-\u,0);
\coordinate(v)at(0,-\u);
\draw[step=\u,color=black!16,ultra thick](0,0)grid(3*\u,2*\u);
\node[inner sep=0](l1)at(3*\u,2*\u){};
\node[inner sep=0](l2)at($(l1)+(h)$){};
\node[inner sep=0](l3)at($(l2)+(v)$){};
\node[inner sep=0](l4)at($(l3)+(h)$){};
\node[inner sep=0](l5)at($(l4)+(v)$){};
\node[inner sep=0](l6)at($(l5)+(h)$){};
\node[inner sep=0](u1)at(3*\u,2*\u){};
\node[inner sep=0](u2)at($(u1)+(h)$){};
\node[inner sep=0](u3)at($(u2)+(v)$){};
\node[inner sep=0](u4)at($(u3)+(h)$){};
\node[inner sep=0](u5)at($(u4)+(h)$){};
\node[inner sep=0](u6)at($(u5)+(v)$){};
\begin{scope}
\clip($(u1.center)+(0,\w)$)--++($1*(h)+(-\w,0)$)--++(2*\w,-2*\w)--++($-1*(h)+(-\w,0)$);
\path[hu]($(u1.center)+(\w,\w)$)rectangle($(u2.center)+(-\w,-\w)$);
\end{scope}
\begin{scope}
\clip($(u2.center)+(-\w,\w)$)--++($1*(v)$)--++(2*\w,-2*\w)--++($-1*(v)$);
\path[vu]($(u2.center)+(\w,\w)$)rectangle($(u3.center)+(-\w,-\w)$);
\end{scope}
\begin{scope}
\clip($(u3.center)+(-\w,\w)$)--++($2*(h)$)--++(2*\w,-2*\w)--++($-2*(h)$);
\path[hu]($(u3.center)+(\w,\w)$)rectangle($(u5.center)+(-\w,-\w)$);
\end{scope}
\begin{scope}
\clip($(u5.center)+(-\w,\w)$)--++($1*(v)+(0,-\w)$)--++(2*\w,0)--++($-1*(v)+(0,-\w)$);
\path[vu]($(u5.center)+(\w,\w)$)rectangle($(u6.center)+(-\w,-\w)$);
\end{scope}
\begin{scope}
\clip($(l1.center)+(0,\w)$)--++($1*(h)+(-\w,0)$)--++(2*\w,-2*\w)--++($-1*(h)+(-\w,0)$);
\path[hl]($(l1.center)+(\w,\w)$)rectangle($(l2.center)+(-\w,-\w)$);
\end{scope}
\begin{scope}
\clip($(l2.center)+(-\w,\w)$)--++($1*(v)$)--++(2*\w,-2*\w)--++($-1*(v)$);
\path[vl]($(l2.center)+(\w,\w)$)rectangle($(l3.center)+(-\w,-\w)$);
\end{scope}
\begin{scope}
\clip($(l3.center)+(-\w,\w)$)--++($1*(h)$)--++(2*\w,-2*\w)--++($-1*(h)$);
\path[hl]($(l3.center)+(\w,\w)$)rectangle($(l4.center)+(-\w,-\w)$);
\end{scope}
\begin{scope}
\clip($(l4.center)+(-\w,\w)$)--++($1*(v)$)--++(2*\w,-2*\w)--++($-1*(v)$);
\path[vl]($(l4.center)+(\w,\w)$)rectangle($(l5.center)+(-\w,-\w)$);
\end{scope}
\begin{scope}
\clip($(l5.center)+(-\w,\w)$)--++($1*(h)+(\w,0)$)--++(0,-2*\w)--++($-1*(h)+(\w,0)$);
\path[hl]($(l5.center)+(\w,\w)$)rectangle($(l6.center)+(-\w,-\w)$);
\end{scope}
\end{tikzpicture}\tikzexternaldisable\qquad
&
\qquad\begin{tikzpicture}[baseline=(current bounding box.center)]
\tikzstyle{out1}=[inner sep=0,minimum size=1.2mm,circle,draw=black,fill=black]
\tikzstyle{in1}=[inner sep=0,minimum size=1.2mm,circle,draw=black,fill=white]
\pgfmathsetmacro{\scalar}{4/3};
\pgfmathsetmacro{\unit}{\scalar*0.922/1.6};
\coordinate (vstep)at(0,-0.26*\unit);
\coordinate (hstep)at(0.20*\unit,0);
\draw[thick](6*\unit,0)--(7*\unit,0);
\foreach \x in {1,...,6}{
\draw[thick](\x*\unit-\unit,0)rectangle(\x*\unit,-\unit);}
\foreach \x in {1,...,5}{
\draw[thick](\x*\unit-\unit,-\unit)rectangle(\x*\unit,-2*\unit);}
\node[inner sep=0]at(0.5*\unit,-0.5*\unit){\scalebox{\scalar}{$+$}};
\node[inner sep=0]at(1.5*\unit,-0.5*\unit){\scalebox{\scalar}{$+$}};
\node[inner sep=0]at(2.5*\unit,-0.5*\unit){\scalebox{\scalar}{$+$}};
\node[inner sep=0]at(3.5*\unit,-0.5*\unit){\scalebox{\scalar}{$0$}};
\node[inner sep=0]at(4.5*\unit,-0.5*\unit){\scalebox{\scalar}{$0$}};
\node[inner sep=0]at(5.5*\unit,-0.5*\unit){\scalebox{\scalar}{$+$}};
\node[inner sep=0]at(0.5*\unit,-1.5*\unit){\scalebox{\scalar}{$+$}};
\node[inner sep=0]at(1.5*\unit,-1.5*\unit){\scalebox{\scalar}{$0$}};
\node[inner sep=0]at(2.5*\unit,-1.5*\unit){\scalebox{\scalar}{$+$}};
\node[inner sep=0]at(3.5*\unit,-1.5*\unit){\scalebox{\scalar}{$+$}};
\node[inner sep=0]at(4.5*\unit,-1.5*\unit){\scalebox{\scalar}{$+$}};
\node[inner sep=0]at($(6.5*\unit,0)+(vstep)$){$1$};
\node[inner sep=0]at($(6*\unit,-0.5*\unit)+(hstep)$){$2$};
\node[inner sep=0]at($(5.5*\unit,-1*\unit)+(vstep)$){$3$};
\node[inner sep=0]at($(5*\unit,-1.5*\unit)+(hstep)$){$4$};
\node[inner sep=0]at($(4.5*\unit,-2*\unit)+(vstep)$){$5$};
\node[inner sep=0]at($(3.5*\unit,-2*\unit)+(vstep)$){$6$};
\node[inner sep=0]at($(2.5*\unit,-2*\unit)+(vstep)$){$7$};
\node[inner sep=0]at($(1.5*\unit,-2*\unit)+(vstep)$){$8$};
\node[inner sep=0]at($(0.5*\unit,-2*\unit)+(vstep)$){$9$};
\end{tikzpicture}\quad\qquad
\\ & & \\ \qquad$T$\rule[-8pt]{0pt}{0pt}\qquad & \qquad$(W_U(T), W_L(T))$\qquad & \qquad$D$\quad\qquad \\ \hline & & \\
\qquad\begin{tikzpicture}[baseline=(current bounding box.center)]
\tikzstyle{out1}=[inner sep=0,minimum size=2.4mm,circle,draw=black,fill=black,semithick]
\tikzstyle{in1}=[inner sep=0,minimum size=2.4mm,circle,draw=black,fill=white,semithick]
\tikzstyle{vertex}=[inner sep=0,minimum size=1.2mm,circle,draw=black,fill=black,semithick]
\pgfmathsetmacro{\l}{0.48};
\pgfmathsetmacro{\d}{0.36};
\pgfmathsetmacro{\s}{1.44};
\node[vertex](r)at(0,0){};
\node[vertex](rh)at($(r)+(-\s,0)$){};
\node[vertex](rv)at($(r)+(0,\s)$){};
\node[vertex](rvh)at($(rv)+(-\s,0)$){};
\node[vertex](rvhh)at($(rvh)+(-\s,0)$){};
\node[inner sep=0](rhh)at($(rh)+(-\l,0)$){};
\node[inner sep=0](rhv)at($(rh)+(0,\l)$){};
\node[inner sep=0](rvv)at($(rv)+(0,\l)$){};
\node[inner sep=0](rvhhh)at($(rvhh)+(-\l,0)$){};
\node[inner sep=0](rvhhv)at($(rvhh)+(0,\l)$){};
\node[inner sep=0](rvhv)at($(rvh)+(0,\l)$){};
\node[inner sep=0]at($(rhh)+(-0.18,0)$){$2$};
\node[inner sep=0]at($(rhv)+(0,0.24)$){$3$};
\node[inner sep=0]at($(rvhhh)+(-0.18,0)$){$4$};
\node[inner sep=0]at($(rvhhv)+(0,0.24)$){$5$};
\node[inner sep=0]at($(rvhv)+(0,0.24)$){$6$};
\node[inner sep=0]at($(rvv)+(0,0.24)$){$7$};
\path[thick](r)edge(rh) edge(rv) (rh)edge(rhh) edge(rhv) (rv)edge(rvh) edge(rvv) (rvh)edge(rvhh) edge(rvhv) (rvhh)edge(rvhhh) edge(rvhhv);
\end{tikzpicture}\qquad
&
\qquad\tikzexternalenable\begin{tikzpicture}[baseline=(current bounding box.center)]
\tikzstyle{hu}=[top color=red!10,bottom color=red!10,middle color=red,opacity=0.70]
\tikzstyle{vu}=[left color=red!10,right color=red!10,middle color=red,opacity=0.70]
\tikzstyle{hl}=[top color=blue!10,bottom color=blue!10,middle color=blue,opacity=0.55]
\tikzstyle{vl}=[left color=blue!10,right color=blue!10,middle color=blue,opacity=0.55]
\pgfmathsetmacro{\u}{0.75};
\pgfmathsetmacro{\w}{0.12*\u};
\coordinate(h)at(-\u,0);
\coordinate(v)at(0,-\u);
\draw[step=\u,color=black!16,ultra thick](0,0)grid(3*\u,1*\u);
\node[inner sep=0](l1)at(3*\u,1*\u){};
\node[inner sep=0](l2)at($(l1)+(h)$){};
\node[inner sep=0](l3)at($(l2)+(v)$){};
\node[inner sep=0](l4)at($(l3)+(h)$){};
\node[inner sep=0](l5)at($(l4)+(h)$){};
\node[inner sep=0](u1)at(3*\u,1*\u){};
\node[inner sep=0](u2)at($(u1)+(h)$){};
\node[inner sep=0](u3)at($(u2)+(v)$){};
\node[inner sep=0](u4)at($(u3)+(h)$){};
\node[inner sep=0](u5)at($(u4)+(h)$){};
\begin{scope}
\clip($(u1.center)+(0,\w)$)--++($1*(h)+(-\w,0)$)--++(2*\w,-2*\w)--++($-1*(h)+(-\w,0)$);
\path[hu]($(u1.center)+(\w,\w)$)rectangle($(u2.center)+(-\w,-\w)$);
\end{scope}
\begin{scope}
\clip($(u2.center)+(-\w,\w)$)--++($1*(v)$)--++(2*\w,-2*\w)--++($-1*(v)$);
\path[vu]($(u2.center)+(\w,\w)$)rectangle($(u3.center)+(-\w,-\w)$);
\end{scope}
\begin{scope}
\clip($(u3.center)+(-\w,\w)$)--++($2*(h)+(\w,0)$)--++(0,-2*\w)--++($-2*(h)+(\w,0)$);
\path[hu]($(u3.center)+(\w,\w)$)rectangle($(u5.center)+(-\w,-\w)$);
\end{scope}
\begin{scope}
\clip($(l1.center)+(0,\w)$)--++($1*(h)+(-\w,0)$)--++(2*\w,-2*\w)--++($-1*(h)+(-\w,0)$);
\path[hl]($(l1.center)+(\w,\w)$)rectangle($(l2.center)+(-\w,-\w)$);
\end{scope}
\begin{scope}
\clip($(l2.center)+(-\w,\w)$)--++($1*(v)$)--++(2*\w,-2*\w)--++($-1*(v)$);
\path[vl]($(l2.center)+(\w,\w)$)rectangle($(l3.center)+(-\w,-\w)$);
\end{scope}
\begin{scope}
\clip($(l3.center)+(-\w,\w)$)--++($2*(h)+(\w,0)$)--++(0,-2*\w)--++($-2*(h)+(\w,0)$);
\path[hl]($(l3.center)+(\w,\w)$)rectangle($(l5.center)+(-\w,-\w)$);
\end{scope}
\end{tikzpicture}\tikzexternaldisable\qquad
& 
\qquad\begin{tikzpicture}[baseline=(current bounding box.center)]
\tikzstyle{out1}=[inner sep=0,minimum size=1.2mm,circle,draw=black,fill=black]
\tikzstyle{in1}=[inner sep=0,minimum size=1.2mm,circle,draw=black,fill=white]
\pgfmathsetmacro{\scalar}{4/3};
\pgfmathsetmacro{\unit}{\scalar*0.922/1.6};
\coordinate (vstep)at(0,-0.26*\unit);
\coordinate (hstep)at(0.20*\unit,0);
\draw[thick](6*\unit,0)--(7*\unit,0);
\foreach \x in {1,...,6}{
\draw[thick](\x*\unit-\unit,0)rectangle(\x*\unit,-\unit);}
\node[inner sep=0]at(0.5*\unit,-0.5*\unit){\scalebox{\scalar}{$+$}};
\node[inner sep=0]at(1.5*\unit,-0.5*\unit){\scalebox{\scalar}{$+$}};
\node[inner sep=0]at(2.5*\unit,-0.5*\unit){\scalebox{\scalar}{$+$}};
\node[inner sep=0]at(3.5*\unit,-0.5*\unit){\scalebox{\scalar}{$0$}};
\node[inner sep=0]at(4.5*\unit,-0.5*\unit){\scalebox{\scalar}{$0$}};
\node[inner sep=0]at(5.5*\unit,-0.5*\unit){\scalebox{\scalar}{$+$}};
\node[inner sep=0]at($(6.5*\unit,0)+(vstep)$){$1$};
\node[inner sep=0]at($(6*\unit,-0.5*\unit)+(hstep)$){$2$};
\node[inner sep=0]at($(5.5*\unit,-1*\unit)+(vstep)$){$3$};
\node[inner sep=0]at($(4.5*\unit,-1*\unit)+(vstep)$){$4$};
\node[inner sep=0]at($(3.5*\unit,-1*\unit)+(vstep)$){$5$};
\node[inner sep=0]at($(2.5*\unit,-1*\unit)+(vstep)$){$6$};
\node[inner sep=0]at($(1.5*\unit,-1*\unit)+(vstep)$){$7$};
\node[inner sep=0]at($(0.5*\unit,-1*\unit)+(vstep)$){$8$};
\end{tikzpicture}\quad\qquad
\\ & & \\ \qquad$T'$\rule[-8pt]{0pt}{0pt}\qquad & \qquad$(W_U(T'), W_L(T'))$\qquad & \qquad$D'$\quad\qquad \\ \hline
\end{tabular}
\caption{An example of Case 1 in the proof of \cref{thm:shift_map_holds}. Here $k=2$, $n=9$, $p=6$.}
\label{fig:finalV-example}
\end{center}
\end{figure}

In terms of plabic graphs, $G(T)$ is obtained from $G(T')$ by blowing up at $p$. Let us introduce the intermediate graph $G''$, obtained from $G(T')$ by inserting a lollipop in between boundary vertices $p-1$ and $p$, and increasing the labels of the boundary vertices $p, p+1, \dots, n-1$ of $G(T')$ by $1$. In terms of permutations, we have
$$
\pi_{G''} = (p \;\; p{+}1 \; \cdots \; n)\pi_{G(T')}(n \;\; n{-}1 \; \cdots \; p).
$$
Then $G(T)$ is obtained from $G(T')$ via $G''$ as follows (cf.\ \eqref{blowup}):
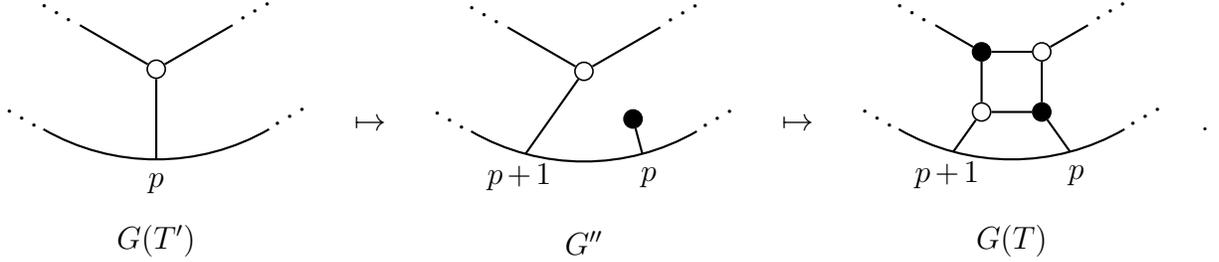
$$
\quad\begin{tikzpicture}[baseline=(current bounding box.center)]
\tikzstyle{out1}=[inner sep=0,minimum size=2.4mm,circle,draw=black,fill=black,semithick]
\tikzstyle{in1}=[inner sep=0,minimum size=2.4mm,circle,draw=black,fill=white,semithick]
\tikzstyle{ambiguous}=[inner sep=0,minimum size=2.4mm,circle,draw=black,fill=black!20,semithick]
\pgfmathsetmacro{\radius}{3.0};
\pgfmathsetmacro{\l}{1.20};
\draw[thick](-120:\radius)arc(-120:-60:\radius);
\node[inner sep=0,rotate=150]at($(-120:\radius)+(150:0.36)$){$\cdots$};
\node[inner sep=0,rotate=30]at($(-60:\radius)+(30:0.36)$){$\cdots$};
\node[inner sep=0](i)at(-90:\radius){};
\node at(-90:\radius+0.32){$p$};
\node at(-90:\radius+1.08){$G(T')$};
\node[in1](v)at(-90:\radius-\l){};
\node[inner sep=0](left)at($(v)+(150:\l)$){};
\node[inner sep=0,rotate=150]at($(left)+(150:0.36)$){$\cdots$};
\node[inner sep=0](right)at($(v)+(30:\l)$){};
\node[inner sep=0,rotate=30]at($(right)+(30:0.36)$){$\cdots$};
\path[thick](v)edge(i.center) (v)edge(left) (v)edge(right);
\end{tikzpicture}\quad\mapsto\quad
\begin{tikzpicture}[baseline=(current bounding box.center)]
\tikzstyle{out1}=[inner sep=0,minimum size=2.4mm,circle,draw=black,fill=black,semithick]
\tikzstyle{in1}=[inner sep=0,minimum size=2.4mm,circle,draw=black,fill=white,semithick]
\tikzstyle{ambiguous}=[inner sep=0,minimum size=2.4mm,circle,draw=black,fill=black!20,semithick]
\pgfmathsetmacro{\radius}{3.0};
\pgfmathsetmacro{\l}{1.20};
\draw[thick](-120:\radius)arc(-120:-60:\radius);
\node[inner sep=0,rotate=150]at($(-120:\radius)+(150:0.36)$){$\cdots$};
\node[inner sep=0,rotate=30]at($(-60:\radius)+(30:0.36)$){$\cdots$};
\node[inner sep=0](i)at(-75:\radius){};
\node at(-75:\radius+0.32){$p$};
\node[inner sep=0](i+1)at(-105:\radius){};
\node at(-105:\radius+0.32){$p\hspace*{-1pt}+\hspace*{-1pt}1$};
\node at(-90:\radius+1.08){$G''$};
\node[out1](lollipop)at(-75:\radius-0.48){};
\node[in1](v)at(-90:\radius-\l){};
\node[inner sep=0](left)at($(v)+(150:\l)$){};
\node[inner sep=0,rotate=150]at($(left)+(150:0.36)$){$\cdots$};
\node[inner sep=0](right)at($(v)+(30:\l)$){};
\node[inner sep=0,rotate=30]at($(right)+(30:0.36)$){$\cdots$};
\path[thick](i.center)edge(lollipop) (v)edge(i+1.center) (v)edge(left) (v)edge(right);
\end{tikzpicture}\quad\mapsto\quad
\begin{tikzpicture}[baseline=(current bounding box.center)]
\tikzstyle{out1}=[inner sep=0,minimum size=2.4mm,circle,draw=black,fill=black,semithick]
\tikzstyle{in1}=[inner sep=0,minimum size=2.4mm,circle,draw=black,fill=white,semithick]
\tikzstyle{ambiguous}=[inner sep=0,minimum size=2.4mm,circle,draw=black,fill=black!20,semithick]
\pgfmathsetmacro{\radius}{3.0};
\pgfmathsetmacro{\l}{1.20};
\pgfmathsetmacro{\s}{0.80};
\draw[thick](-120:\radius)arc(-120:-60:\radius);
\node[inner sep=0,rotate=150]at($(-120:\radius)+(150:0.36)$){$\cdots$};
\node[inner sep=0,rotate=30]at($(-60:\radius)+(30:0.36)$){$\cdots$};
\node[inner sep=0](i)at(-75:\radius){};
\node at(-75:\radius+0.32){$p$};
\node[inner sep=0](i+1)at(-105:\radius){};
\node at(-105:\radius+0.32){$p\hspace*{-1pt}+\hspace*{-1pt}1$};
\node at(-90:\radius+1.08){$G(T)$};
\node[inner sep=0](v)at(-90:\radius-\l){};
\node[out1](a)at($(v)+(150:{\s/1.7320508075688772935274463415058723669428052538104})$){};
\node[in1](b)at($(a)+(\s,0)$){};
\node[out1](c)at($(a)+(\s,-\s)$){};
\node[in1](d)at($(a)+(0,-\s)$){};
\node[inner sep=0](left)at($(v)+(150:\l)$){};
\node[inner sep=0,rotate=150]at($(left)+(150:0.36)$){$\cdots$};
\node[inner sep=0](right)at($(v)+(30:\l)$){};
\node[inner sep=0,rotate=30]at($(right)+(30:0.36)$){$\cdots$};
\path[thick](a)edge(b) (b)edge(c) (c)edge(d) (d)edge(a) (c)edge(i.center) (d)edge(i+1.center) (a)edge(left) (b)edge(right);
\end{tikzpicture}\quad.
$$
Since $e$ corresponds to the last letter of $W_U(T)$, and this letter is $V$, the boundary vertices $p+1, p+2, \dots, n-1$ of $G(T)$ are each incident to a white vertex. Hence $\pi_{G(T)}(n) = p$, and we see that $\pi_{G(T)} = \pi_{G''}(p \;\; p{+}1 \;\; n)$. Thus
\begin{align}\label{case_1_G_perm}
\begin{aligned}
\pi_{G(T)} &= (p \;\; p{+}1 \; \cdots \; n)\pi_{G(T')}(n \;\; n{-}1 \; \cdots \; p)(p \;\; p{+}1 \;\; n) \\
&= (p \;\; p{+}1 \; \cdots \; n)\pi_{G(T')}(n{-}1 \;\; n{-}2 \; \cdots \; p{+}1).
\end{aligned}
\end{align}

On the other hand, by \cref{def:latticetocell}, $D$ is obtained from $D'$ by appending a new row, as follows:
$$
D = \quad
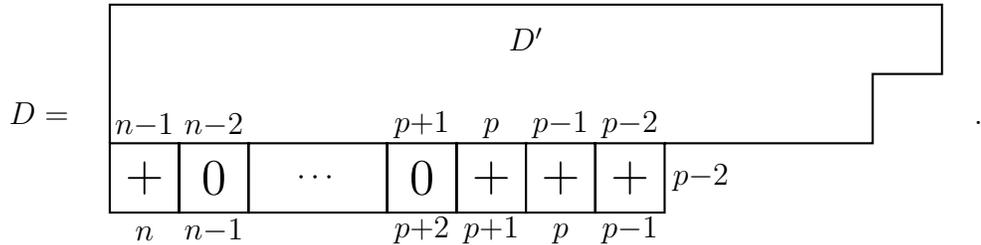
\begin{tikzpicture}[baseline=(current bounding box.center)]
\pgfmathsetmacro{\scalar}{1.6};
\pgfmathsetmacro{\unit}{\scalar*0.922/1.6};
\coordinate (vstep)at(0,-0.24);
\coordinate (hstep)at(0.48,0);
\draw[thick](0,0)--(0,2*\unit)--(12*\unit,2*\unit)--(12*\unit,\unit)--(11*\unit,\unit)--(11*\unit,0)--(8*\unit,0) (0,0)rectangle(\unit,-\unit) (\unit,0)rectangle(2*\unit,-\unit) (2*\unit,0)rectangle(4*\unit,-\unit) (4*\unit,0)rectangle(5*\unit,-\unit) (5*\unit,0)rectangle(6*\unit,-\unit) (6*\unit,0)rectangle(7*\unit,-\unit) (7*\unit,0)rectangle(8*\unit,-\unit);
\node[inner sep=0]at(6*\unit,1.5*\unit){$D'$};
\node[inner sep=0]at(0.5*\unit,-0.5*\unit){\scalebox{\scalar}{$+$}};
\node[inner sep=0]at(1.5*\unit,-0.5*\unit){\scalebox{\scalar}{$0$}};
\node[inner sep=0]at(3*\unit,-0.5*\unit){$\cdots$};
\node[inner sep=0]at(4.5*\unit,-0.5*\unit){\scalebox{\scalar}{$0$}};
\node[inner sep=0]at(5.5*\unit,-0.5*\unit){\scalebox{\scalar}{$+$}};
\node[inner sep=0]at(6.5*\unit,-0.5*\unit){\scalebox{\scalar}{$+$}};
\node[inner sep=0]at(7.5*\unit,-0.5*\unit){\scalebox{\scalar}{$+$}};
\node[inner sep=0]at($(8*\unit,-0.5*\unit)+(hstep)$){$p{-}2$};
\node[inner sep=0]at($(7.5*\unit,-1*\unit)+(vstep)$){$p{-}1$};
\node[inner sep=0]at($(6.5*\unit,-1*\unit)+(vstep)$){$\phantom{1}p\phantom{1}$};
\node[inner sep=0]at($(5.5*\unit,-1*\unit)+(vstep)$){$p{+}1$};
\node[inner sep=0]at($(4.5*\unit,-1*\unit)+(vstep)$){$p{+}2$};
\node[inner sep=0]at($(1.5*\unit,-1*\unit)+(vstep)$){$n{-}1$};
\node[inner sep=0]at($(0.5*\unit,-1*\unit)+(vstep)$){$\phantom{1}n\phantom{1}$};
\node[inner sep=0]at($(7.5*\unit,0)-(vstep)$){$p{-}2$};
\node[inner sep=0]at($(6.5*\unit,0)-(vstep)$){$p{-}1$};
\node[inner sep=0]at($(5.5*\unit,0)-(vstep)$){$\phantom{1}p\phantom{1}$};
\node[inner sep=0]at($(4.5*\unit,0)-(vstep)$){$p{+}1$};
\node[inner sep=0]at($(1.5*\unit,0)-(vstep)$){$n{-}2$};
\node[inner sep=0]at($(0.5*\unit,0)-(vstep)$){$n{-}1$};
\end{tikzpicture}\quad.
$$
Above we have given the labels of the southeast borders of $D$ and $D'$ starting at $p-2$. We can then verify from \cref{def:oplus} that
$$
\pi_D = (p{-}2 \;\; p{-}1 \; \cdots \; n)\pi_{D'}(n{-}1 \;\; n{-}2 \; \cdots \; p{+}1).
$$
Putting this together with the induction hypothesis $\pi_{D'} = c_{n-1}^2\pi_{G(T')}$ and \eqref{case_1_G_perm}, we obtain
\begin{align*}
\pi_D &= (p{-}2 \;\; p{-}1 \; \cdots \; n)c_{n-1}^2\pi_{G(T')}(n{-}1 \;\; n{-}2 \; \cdots \; p{+}1) \\
&= c_n^2(p \;\; p{+}1 \; \cdots \; n)\pi_{G(T')}(n{-}1 \;\; n{-}2 \; \cdots \; p{+}1) = c_n^2\pi_{G(T)}.
\end{align*}

{\bfseries Case 2: $W_U(T)$ ends in $H$.}
Suppose that $W_U(T)$ ends in precisely $s$ $H$'s, where $s\ge 1$. Let $e_{n-s-3}, e_{n-s-2}, \dots, e_{n-4}$ be the internal edges of $T$ corresponding to the last $s$ letters of $W_U(T)$, when we read the internal edges in a depth-first search as in \cref{readoff}.
\begin{figure}[htb]
\begin{center}
\begin{tabular}{|cc|}
\hline & \\
\;\begin{tikzpicture}[baseline=(current bounding box.center)]
\tikzstyle{out1}=[inner sep=0,minimum size=2.4mm,circle,draw=black,fill=black,semithick]
\tikzstyle{in1}=[inner sep=0,minimum size=2.4mm,circle,draw=black,fill=white,semithick]
\tikzstyle{vertex}=[inner sep=0,minimum size=1.2mm,circle,draw=black,fill=black,semithick]
\pgfmathsetmacro{\l}{0.48};
\pgfmathsetmacro{\d}{0.40};
\pgfmathsetmacro{\dbelow}{0.96};
\pgfmathsetmacro{\s}{1.32};
\node[vertex](r)at(0,0){};
\node[inner sep=0](eu)at($(r)+(0,-\dbelow)$){};
\node[inner sep=0,rotate=-90](e)at($(eu)+(0,-\d)$){$\cdots$};
\node[vertex](rh)at($(r)+(-\s,0)$){};
\node[inner sep=0](rv)at($(r)+(0,\l)$){};
\node[inner sep=0](rhv)at($(rh)+(0,\l)$){};
\node[inner sep=0](dr)at($(rh)+(-\d,0)$){};
\node[inner sep=0](d)at($(dr)+(-\d,0)$){};
\node[inner sep=0]at($(d)+(0.03,0)$){$\cdots$};
\node[inner sep=0](dl)at($(d)+(-\d,0)$){};
\node[vertex](rnew)at($(dl)+(-\d,0)$){};
\node[vertex](rnewh)at($(rnew)+(-\s,0)$){};
\node[inner sep=0](rnewv)at($(rnew)+(0,\l)$){};
\node[vertex](rnewhh)at($(rnewh)+(-\s,0)$){};
\node[inner sep=0](rnewhv)at($(rnewh)+(0,\l)$){};
\node[inner sep=0](rnewhhh)at($(rnewhh)+(-\l,0)$){};
\node[inner sep=0](rnewhhv)at($(rnewhh)+(0,\l)$){};
\node[inner sep=0]at($(rnewhhh)+(-0.18,0)$){$p$};
\node[inner sep=0]at($(rnewhhv)+(0,0.24)$){$p{+}1$};
\node[inner sep=0]at($(rnewhv)+(0,0.24)$){$p{+}2$};
\node[inner sep=0]at($(rnewv)+(0,0.24)$){$p{+}3$};
\node[inner sep=0]at($(rhv)+(0,0.24)$){$p{+}s$};
\node[inner sep=0]at($(rv)+(0,0.24)$){$p{+}s{+}1$};
\path[thick](r)edge node[below=-2pt]{$e_{n-s-3}$}(rh) edge(rv) edge(eu) (rh)edge(dr) edge(rhv) (rnew)edge(dl) edge node[below=-2pt]{$e_{n-5}$}(rnewh) edge(rnewv) (rnewh)edge node[below=-2pt]{$e_{n-4}$}(rnewhh) edge(rnewhv) (rnewhh)edge(rnewhhh) edge(rnewhhv);
\end{tikzpicture}\quad
&
\quad\begin{tikzpicture}[baseline=(current bounding box.center)]
\tikzstyle{out1}=[inner sep=0,minimum size=2.4mm,circle,draw=black,fill=black,semithick]
\tikzstyle{in1}=[inner sep=0,minimum size=2.4mm,circle,draw=black,fill=white,semithick]
\tikzstyle{vertex}=[inner sep=0,minimum size=1.2mm,circle,draw=black,fill=black,semithick]
\pgfmathsetmacro{\rstep}{0.54};
\pgfmathsetmacro{\vstep}{0.36};
\pgfmathsetmacro{\dotsstep}{0.66};
\pgfmathsetmacro{\rdots}{0.30};
\pgfmathsetmacro{\rdotsbelow}{0.96};
\pgfmathsetmacro{\vdots}{\rdots*\vstep/\rstep};
\pgfmathsetmacro{\finalrstep}{0.84};
\pgfmathsetmacro{\legstep}{0.30};
\pgfmathsetmacro{\bm}{0.72};
\pgfmathsetmacro{\bedge}{0.42};
\node[out1](m0)at(0,0){};
\node[out1](m2)at($(m0)+(2*\rstep,0)$){};
\node[out1](m4)at($(m2)+(2*\rstep,0)$){};
\node[out1](m6)at($(m4)+(2*\rstep,0)$){};
\node[out1](m8)at($(m6)+(2*\dotsstep,0)$){};
\node[out1](m10)at($(m8)+(2*\rstep,0)$){};
\node[in1](b2)at($(m2)+(0,-\bm)$){};
\node[in1](b4)at($(m4)+(0,-\bm)$){};
\node[in1](b6)at($(m6)+(0,-\bm)$){};
\node[in1](b8)at($(m8)+(0,-\bm)$){};
\node[in1](b10)at($(m10)+(0,-\bm)$){};
\node[inner sep=0](b7l)at($(b6)+(\rdots,0)$){};
\node[inner sep=0](b7r)at($(b8)+(-\rdots,0)$){};
\node[inner sep=0](b7)at($(b6)+(\dotsstep,0)$){};
\node[inner sep=0]at($(b7)+(0.03,0)$){$\cdots$};
\node[in1](t1)at($(m0)+(\rstep,\vstep)$){};
\node[in1](t3)at($(m2)+(\rstep,\vstep)$){};
\node[in1](t5)at($(m4)+(\rstep,\vstep)$){};
\node[inner sep=0](t7l)at($(m6)+(\rdots,\vdots)$){};
\node[inner sep=0](t7r)at($(m8)+(-\rdots,\vdots)$){};
\node[inner sep=0](t7)at($(m6)+(\dotsstep,\vstep/2)$){};
\node[inner sep=0]at($(t7)+(0.03,0)$){$\cdots$};
\node[in1](t9)at($(m8)+(\rstep,\vstep)$){};
\node[in1](t11)at($(m10)+(\finalrstep,\vstep)$){};
\node[in1](legl)at($(b10)+(\legstep,-\legstep)$){};
\node[out1](legr)at($(b10)+(\finalrstep,-\legstep)$){};
\node[inner sep=0](leglb)at($(legl)+(0,-\rdotsbelow)$){};
\node[inner sep=0](legrb)at($(legr)+(0,-\rdotsbelow)$){};
\node[inner sep=0]at($(leglb)+(-0.21,0.27)$){$f_1$};
\node[inner sep=0]at($(legrb)+(0.23,0.27)$){$f_2$};
\node[inner sep=0,rotate=-90]at($(leglb)+(0,-0.30)$){$\cdots$};
\node[inner sep=0,rotate=-90]at($(legrb)+(0,-0.30)$){$\cdots$};
\node[inner sep=0](bound0)at($(m0)+(-\bedge,0)$){};
\node[inner sep=0](bound1)at($(t1)+(0,\bedge)$){};
\node[inner sep=0](bound2)at($(t3)+(0,\bedge)$){};
\node[inner sep=0](bound3)at($(t5)+(0,\bedge)$){};
\node[inner sep=0](bound4)at($(bound3)+(\rstep+\dotsstep,0)$){$\cdots$};
\node[inner sep=0](bound4l)at($(bound4)+(-0.40,0)$){};
\node[inner sep=0](bound4r)at($(bound4)+(0.32,0)$){};
\node[inner sep=0](bounds)at($(t9)+(0,\bedge)$){};
\node[inner sep=0](bounds+1)at($(t11)+(0,\bedge)$){};
\node[inner sep=0](bound01)at($(bound1)+(-0.60,-0.24)$){};
\node[inner sep=0](bound-1)at($(bound0)+(-0.12,-\bm-\legstep-\rdotsbelow)$){};
\node[inner sep=0,rotate=-90]at($(bound-1)+(0,-0.30)$){$\cdots$};
\node[inner sep=0](bounds+2)at($(legrb)+(0.72,0)$){};
\node[inner sep=0,rotate=-90]at($(bounds+2)+(0,-0.30)$){$\cdots$};
\node[inner sep=0](bounds+12)at($(bounds+1)+(0.66,-0.44)$){};
\node[inner sep=0]at($(bound0)+(-0.18,0)$){$p$};
\node[inner sep=0]at($(bound1)+(0,0.24)$){$p{+}1$};
\node[inner sep=0]at($(bound2)+(0,0.24)$){$p{+}2$};
\node[inner sep=0]at($(bound3)+(0,0.24)$){$p{+}3$};
\node[inner sep=0]at($(bounds)+(0,0.24)$){$p{+}s$};
\node[inner sep=0]at($(bounds+1)+(0,0.24)$){$p{+}s{+}1$};
\path[thick](m0)edge(bound0.center) edge(t1) edge(b2) (m2)edge(t1) edge(t3) edge(b2) (m4)edge(t3) edge(t5) edge(b4) (m6)edge(t5) edge(t7l) edge(b6) (m8)edge(t7r) edge(t9) edge(b8) (m10)edge(t9) edge(t11) edge(b10) (b2)edge(b4) (b4)edge(b6) (b6)edge(b7l) (b8)edge(b7r) edge(b10) (t1)edge(bound1.center) (t3)edge(bound2.center) (t5)edge(bound3.center) (t9)edge(bounds.center) (t11)edge(bounds+1.center) edge(legr) (legl)edge(b10) edge(legr) edge(leglb) (legr)edge(legrb);
\draw[thick]plot[smooth,tension=0.5]coordinates{(bound-1) (bound0) (bound01) (bound1) (bound2) (bound3) (bound4l)};
\draw[thick]plot[smooth,tension=0.5]coordinates{(bound4r) (bounds) (bounds+1) (bounds+12) (bounds+2)};
\end{tikzpicture}\quad
\\ \qquad$T$\rule[-8pt]{0pt}{0pt}\qquad & \qquad$G(T)$\qquad \\ \hline & \\
\;\begin{tikzpicture}[baseline=(current bounding box.center)]
\tikzstyle{out1}=[inner sep=0,minimum size=2.4mm,circle,draw=black,fill=black,semithick]
\tikzstyle{in1}=[inner sep=0,minimum size=2.4mm,circle,draw=black,fill=white,semithick]
\tikzstyle{vertex}=[inner sep=0,minimum size=1.2mm,circle,draw=black,fill=black,semithick]
\pgfmathsetmacro{\l}{0.48};
\pgfmathsetmacro{\d}{0.40};
\pgfmathsetmacro{\dbelow}{0.96};
\node[vertex](r)at(0,0){};
\node[inner sep=0](eu)at($(r)+(0,-\dbelow)$){};
\node[inner sep=0,rotate=-90](e)at($(eu)+(0,-\d)$){$\cdots$};
\node[inner sep=0](rh)at($(r)+(-\l,0)$){};
\node[inner sep=0](rv)at($(r)+(0,\l)$){};
\node[inner sep=0]at($(rh)+(-0.18,0)$){$p$};
\node[inner sep=0]at($(rv)+(0,0.24)$){$p{+}1$};
\path[thick](r)edge(rh) edge(rv) edge(eu);
\end{tikzpicture}\quad
&
\begin{tikzpicture}[baseline=(current bounding box.center)]
\tikzstyle{out1}=[inner sep=0,minimum size=2.4mm,circle,draw=black,fill=black,semithick]
\tikzstyle{in1}=[inner sep=0,minimum size=2.4mm,circle,draw=black,fill=white,semithick]
\tikzstyle{vertex}=[inner sep=0,minimum size=1.2mm,circle,draw=black,fill=black,semithick]
\pgfmathsetmacro{\rstep}{0.54};
\pgfmathsetmacro{\vstep}{0.36};
\pgfmathsetmacro{\dotsstep}{0.66};
\pgfmathsetmacro{\rdots}{0.30};
\pgfmathsetmacro{\rdotsbelow}{0.96};
\pgfmathsetmacro{\vdots}{\rdots*\vstep/\rstep};
\pgfmathsetmacro{\finalrstep}{0.84};
\pgfmathsetmacro{\legstep}{0.30};
\pgfmathsetmacro{\bm}{0.60};
\pgfmathsetmacro{\bedge}{0.42};
\node[out1](m0)at(0,0){};
\node[inner sep=0](m2)at($(m0)+(2*\rstep,0)$){};
\node[inner sep=0](m4)at($(m2)+(2*\rstep,0)$){};
\node[inner sep=0](m6)at($(m4)+(2*\rstep,0)$){};
\node[inner sep=0](m8)at($(m6)+(2*\dotsstep,0)$){};
\node[inner sep=0](m10)at($(m8)+(2*\rstep,0)$){};
\node[inner sep=0](b2)at($(m2)+(0,-\bm)$){};
\node[inner sep=0](b4)at($(m4)+(0,-\bm)$){};
\node[inner sep=0](b6)at($(m6)+(0,-\bm)$){};
\node[inner sep=0](b8)at($(m8)+(0,-\bm)$){};
\node[inner sep=0](b10)at($(m10)+(0,-\bm)$){};
\node[inner sep=0](b7l)at($(b6)+(\rdots,0)$){};
\node[inner sep=0](b7r)at($(b8)+(-\rdots,0)$){};
\node[in1](t1)at($(m0)+(\rstep,\vstep)$){};
\node[out1](t3)at($(m2)+(\rstep,\vstep)$){};
\node[out1](t5)at($(m4)+(\rstep,\vstep)$){};
\node[inner sep=0](t7l)at($(m6)+(\rdots,\vdots)$){};
\node[inner sep=0](t7r)at($(m8)+(-\rdots,\vdots)$){};
\node[inner sep=0](t7)at($(m6)+(\dotsstep,\vstep/2+\legstep)$){};
\node[out1](t9)at($(m8)+(\rstep,\vstep)$){};
\node[out1](t11)at($(m10)+(\finalrstep,\vstep)$){};
\node[inner sep=0](legl)at($(b10)+(\legstep,-\legstep)$){};
\node[inner sep=0](legr)at($(b10)+(\finalrstep,-\legstep)$){};
\node[inner sep=0](leglb)at($(legl)+(0,-\rdotsbelow)$){};
\node[inner sep=0](legrb)at($(legr)+(0,-\rdotsbelow)$){};
\node[inner sep=0](bound0)at($(m0)+(-\bedge,0)$){};
\node[inner sep=0](bound1)at($(t1)+(0,\bedge)$){};
\node[inner sep=0](bound2)at($(t3)+(0,\bedge)$){};
\node[inner sep=0](bound3)at($(t5)+(0,\bedge)$){};
\node[inner sep=0](bound4)at($(bound3)+(\rstep+\dotsstep,0)$){};
\node[inner sep=0]at($(bound4)+(0.03,0)$){$\cdots$};
\node[inner sep=0](bound4l)at($(bound4)+(-0.36,0)$){};
\node[inner sep=0](bound4r)at($(bound4)+(0.36,0)$){};
\node[inner sep=0](bounds)at($(t9)+(0,\bedge)$){};
\node[inner sep=0](bounds+1)at($(t11)+(0,\bedge)$){};
\node[inner sep=0](bound01)at($(bound1)+(-0.60,-0.24)$){};
\node[inner sep=0](bound-1)at($(bound0)+(-0.12,-\bm-\rdotsbelow)$){};
\node[inner sep=0,rotate=-90]at($(bound-1)+(0,-0.30)$){$\cdots$};
\node[inner sep=0](bounds+2)at($(legrb)+(0.72,\legstep)$){};
\node[inner sep=0,rotate=-90]at($(bounds+2)+(0,-0.30)$){$\cdots$};
\node[inner sep=0](bounds+12)at($(bounds+1)+(0.66,-0.44)$){};
\node[inner sep=0]at($(bound0)+(-0.18,0)$){$p$};
\node[inner sep=0]at($(bound1)+(0,0.24)$){$p{+}1$};
\node[inner sep=0]at($(bound2)+(0,0.24)$){$p{+}2$};
\node[inner sep=0]at($(bound3)+(0,0.24)$){$p{+}3$};
\node[inner sep=0]at($(bounds)+(0,0.24)$){$p{+}s$};
\node[inner sep=0]at($(bounds+1)+(0,0.24)$){$p{+}s{+}1$};
\node[out1](newlegr)at(b2){};
\node[in1](newlegl)at($(newlegr)+(-\finalrstep+\legstep,0)$){};
\node[inner sep=0](newleglb)at($(newlegl)+(0,-\rdotsbelow)$){};
\node[inner sep=0](newlegrb)at($(newlegr)+(0,-\rdotsbelow)$){};
\node[inner sep=0]at($(newleglb)+(-0.21,0.27)$){$f_1$};
\node[inner sep=0]at($(newlegrb)+(0.23,0.27)$){$f_2$};
\node[inner sep=0,rotate=-90]at($(newleglb)+(0,-0.30)$){$\cdots$};
\node[inner sep=0,rotate=-90]at($(newlegrb)+(0,-0.30)$){$\cdots$};
\path[thick](m0)edge(bound0.center) edge(t1) edge(newlegl) (newlegr)edge(t1) edge(newlegl) edge(newlegrb) (newlegl)edge(newleglb) (t1)edge(bound1.center) (t3)edge(bound2.center) (t5)edge(bound3.center) (t9)edge(bounds.center) (t11)edge(bounds+1.center);
\draw[thick]plot[smooth,tension=0.5]coordinates{(bound-1) (bound0) (bound01) (bound1) (bound2) (bound3) (bound4l)};
\draw[thick]plot[smooth,tension=0.5]coordinates{(bound4r) (bounds) (bounds+1) (bounds+12) (bounds+2)};
\end{tikzpicture}\quad
\\ \qquad$T'$\rule[-8pt]{0pt}{0pt}\qquad & \qquad$G''$\qquad \\ \hline
\end{tabular}
\caption{The complete binary trees $T$ and $T'$, and the plabic graphs $G(T)$ and $G''$, in Case 2 of the proof of \cref{thm:shift_map_holds}. If $k=0$, then the rightmost vertex shown in $T$ and $T'$ is the root vertex, and the edges $f_1$ and $f_2$ are incident to boundary vertices $1$ and $n$.}
\label{fig:finalH-tree-transformation}
\end{center}
\end{figure}
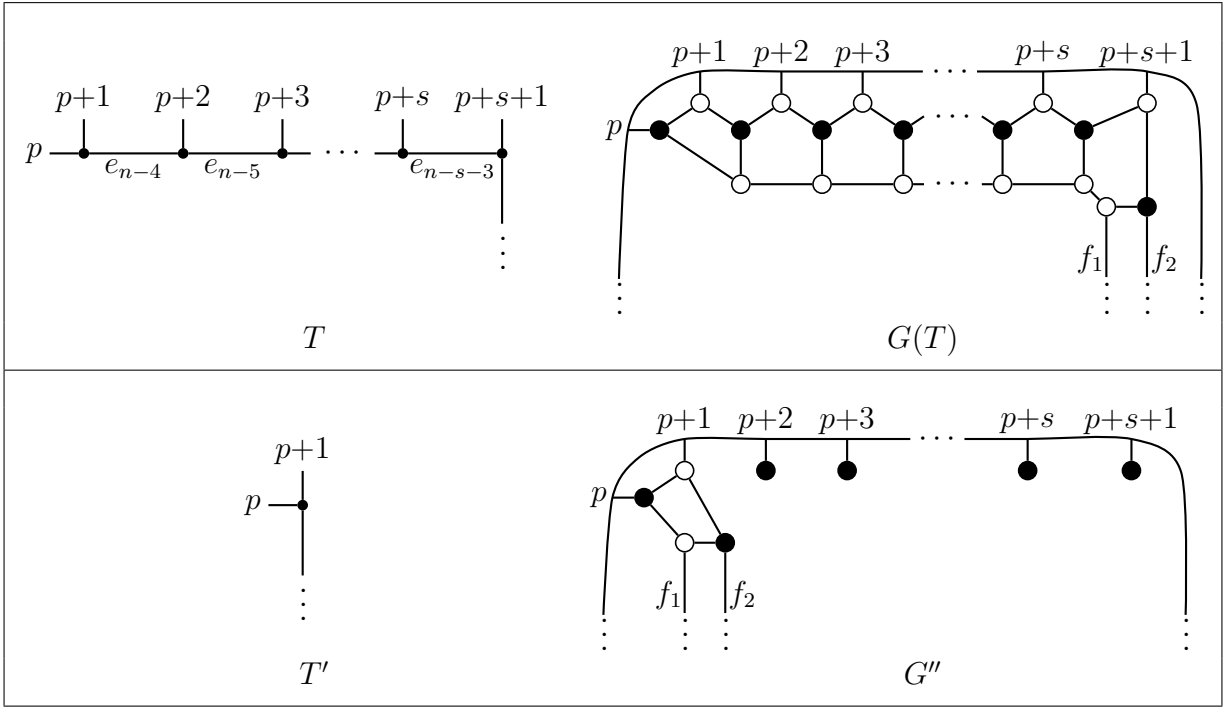
These edges appear in a horizontal path in $T$ as in 
\cref{fig:finalH-tree-transformation}, 
where the rightmost vertex is either the root vertex (in which case this picture is the entirety of $T$), or is joined to its parent vertex by a vertical edge. The children of the vertices on this path are all leaves, which are labeled by $p, p+1, \dots, p+s+1$ for some $p$. Proceeding in a similar manner to Case 1, we let $T'\in\mathcal{T}_{n-s,k,4}$ be obtained from $T$ by replacing the entirety of this horizontal path and its children by a vertex incident to two leaves labeled $p$ and $p+1$ (see \cref{fig:finalH-tree-transformation}). Note that $(W_U(T'), W_L(T'))$ is obtained from $(W_U(T), W_L(T))$ by deleting the last $s$ steps of each path (which are all horizontal). Let $D'\in\mathcal{D}_{n-s,k,4}$ be the associated $\oplus$-diagram. (See \cref{fig:finalH-example}.) By the induction hypothesis, we have $\pi_{D'} = c_{n-s}^2\pi_{G(T')}$. (As in Case $1$, we will regard permutations of $[n-s]$ as permutations of $[n]$ which fix $n-s+1, \dots, n$.)
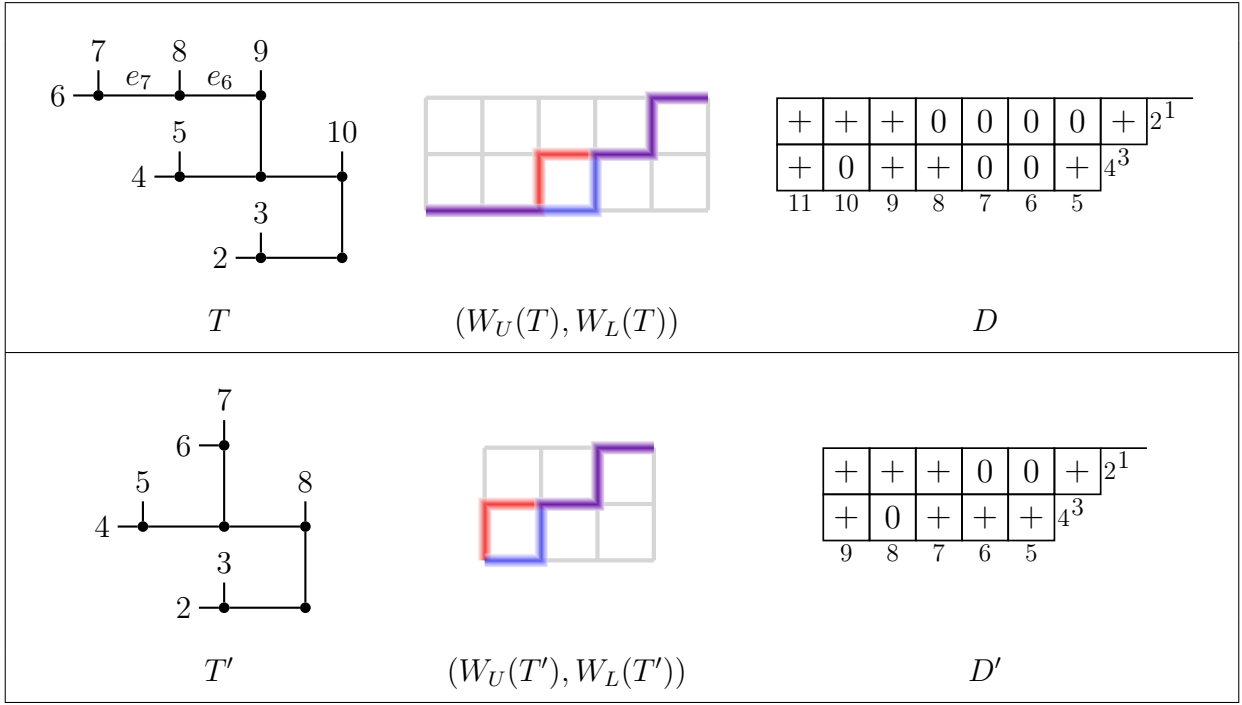
\begin{figure}[htb]
\begin{center}
\begin{tabular}{|ccc|}
\hline & & \\
\quad\begin{tikzpicture}[baseline=(current bounding box.center)]
\tikzstyle{out1}=[inner sep=0,minimum size=2.4mm,circle,draw=black,fill=black,semithick]
\tikzstyle{in1}=[inner sep=0,minimum size=2.4mm,circle,draw=black,fill=white,semithick]
\tikzstyle{vertex}=[inner sep=0,minimum size=1.2mm,circle,draw=black,fill=black,semithick]
\pgfmathsetmacro{\l}{0.36};
\pgfmathsetmacro{\d}{0.30};
\pgfmathsetmacro{\s}{1.08};
\node[vertex](r)at(0,0){};
\node[vertex](rh)at($(r)+(-\s,0)$){};
\node[vertex](rv)at($(r)+(0,\s)$){};
\node[vertex](rvh)at($(rv)+(-\s,0)$){};
\node[vertex](rvhh)at($(rvh)+(-\s,0)$){};
\node[vertex](rvhv)at($(rvh)+(0,\s)$){};
\node[vertex](rvhvh)at($(rvhv)+(-\s,0)$){};
\node[vertex](rvhvhh)at($(rvhvh)+(-\s,0)$){};
\node[inner sep=0](rhh)at($(rh)+(-\l,0)$){};
\node[inner sep=0](rhv)at($(rh)+(0,\l)$){};
\node[inner sep=0](rvv)at($(rv)+(0,\l)$){};
\node[inner sep=0](rvhhh)at($(rvhh)+(-\l,0)$){};
\node[inner sep=0](rvhhv)at($(rvhh)+(0,\l)$){};
\node[inner sep=0](rvhvv)at($(rvhv)+(0,\l)$){};
\node[inner sep=0](rvhvhv)at($(rvhvh)+(0,\l)$){};
\node[inner sep=0](rvhvhhv)at($(rvhvhh)+(0,\l)$){};
\node[inner sep=0](rvhvhhh)at($(rvhvhh)+(-\l,0)$){};
\node[inner sep=0]at($(rhh)+(-0.18,0)$){$2$};
\node[inner sep=0]at($(rhv)+(0,0.24)$){$3$};
\node[inner sep=0]at($(rvhhh)+(-0.18,0)$){$4$};
\node[inner sep=0]at($(rvhhv)+(0,0.24)$){$5$};
\node[inner sep=0]at($(rvhvhhh)+(-0.18,0)$){$6$};
\node[inner sep=0]at($(rvhvhhv)+(0,0.24)$){$7$};
\node[inner sep=0]at($(rvhvhv)+(0,0.24)$){$8$};
\node[inner sep=0]at($(rvhvv)+(0,0.24)$){$9$};
\node[inner sep=0]at($(rvv)+(0,0.24)$){$10$};
\path[thick](r)edge(rh) edge(rv) (rh)edge(rhh)edge(rhv) (rv)edge(rvh) edge(rvv) (rvh)edge(rvhh) edge(rvhv) (rvhh)edge(rvhhh) edge(rvhhv) (rvhv)edge node[above=-2pt]{$e_6$}(rvhvh) edge(rvhvv) (rvhvh)edge node[above=-2pt]{$e_7$}(rvhvhh) edge(rvhvhv) (rvhvhh)edge(rvhvhhh) edge(rvhvhhv);
\end{tikzpicture}\quad
&
\quad\;\tikzexternalenable\begin{tikzpicture}[baseline=(current bounding box.center)]
\tikzstyle{hu}=[top color=red!10,bottom color=red!10,middle color=red,opacity=0.70]
\tikzstyle{vu}=[left color=red!10,right color=red!10,middle color=red,opacity=0.70]
\tikzstyle{hl}=[top color=blue!10,bottom color=blue!10,middle color=blue,opacity=0.55]
\tikzstyle{vl}=[left color=blue!10,right color=blue!10,middle color=blue,opacity=0.55]
\pgfmathsetmacro{\u}{0.75};
\pgfmathsetmacro{\w}{0.12*\u};
\coordinate(h)at(-\u,0);
\coordinate(v)at(0,-\u);
\draw[step=\u,color=black!16,ultra thick](0,0)grid(5*\u,2*\u);
\node[inner sep=0](l1)at(5*\u,2*\u){};
\node[inner sep=0](l2)at($(l1)+(h)$){};
\node[inner sep=0](l3)at($(l2)+(v)$){};
\node[inner sep=0](l4)at($(l3)+(h)$){};
\node[inner sep=0](l5)at($(l4)+(v)$){};
\node[inner sep=0](l6)at($(l5)+(h)$){};
\node[inner sep=0](l7)at($(l6)+(h)$){};
\node[inner sep=0](l8)at($(l7)+(h)$){};
\node[inner sep=0](u1)at(5*\u,2*\u){};
\node[inner sep=0](u2)at($(u1)+(h)$){};
\node[inner sep=0](u3)at($(u2)+(v)$){};
\node[inner sep=0](u4)at($(u3)+(h)$){};
\node[inner sep=0](u5)at($(u4)+(h)$){};
\node[inner sep=0](u6)at($(u5)+(v)$){};
\node[inner sep=0](u7)at($(u6)+(h)$){};
\node[inner sep=0](u8)at($(u7)+(h)$){};
\begin{scope}
\clip($(u1.center)+(0,\w)$)--++($1*(h)+(-\w,0)$)--++(2*\w,-2*\w)--++($-1*(h)+(-\w,0)$);
\path[hu]($(u1.center)+(\w,\w)$)rectangle($(u2.center)+(-\w,-\w)$);
\end{scope}
\begin{scope}
\clip($(u2.center)+(-\w,\w)$)--++($1*(v)$)--++(2*\w,-2*\w)--++($-1*(v)$);
\path[vu]($(u2.center)+(\w,\w)$)rectangle($(u3.center)+(-\w,-\w)$);
\end{scope}
\begin{scope}
\clip($(u3.center)+(-\w,\w)$)--++($2*(h)$)--++(2*\w,-2*\w)--++($-2*(h)$);
\path[hu]($(u3.center)+(\w,\w)$)rectangle($(u5.center)+(-\w,-\w)$);
\end{scope}
\begin{scope}
\clip($(u5.center)+(-\w,\w)$)--++($1*(v)$)--++(2*\w,-2*\w)--++($-1*(v)$);
\path[vu]($(u5.center)+(\w,\w)$)rectangle($(u6.center)+(-\w,-\w)$);
\end{scope}
\begin{scope}
\clip($(u6.center)+(-\w,\w)$)--++($2*(h)+(\w,0)$)--++(0,-2*\w)--++($-2*(h)+(\w,0)$);
\path[hu]($(u6.center)+(\w,\w)$)rectangle($(u8.center)+(-\w,-\w)$);
\end{scope}
\begin{scope}
\clip($(l1.center)+(0,\w)$)--++($1*(h)+(-\w,0)$)--++(2*\w,-2*\w)--++($-1*(h)+(-\w,0)$);
\path[hl]($(l1.center)+(\w,\w)$)rectangle($(l2.center)+(-\w,-\w)$);
\end{scope}
\begin{scope}
\clip($(l2.center)+(-\w,\w)$)--++($1*(v)$)--++(2*\w,-2*\w)--++($-1*(v)$);
\path[vl]($(l2.center)+(\w,\w)$)rectangle($(l3.center)+(-\w,-\w)$);
\end{scope}
\begin{scope}
\clip($(l3.center)+(-\w,\w)$)--++($1*(h)$)--++(2*\w,-2*\w)--++($-1*(h)$);
\path[hl]($(l3.center)+(\w,\w)$)rectangle($(l4.center)+(-\w,-\w)$);
\end{scope}
\begin{scope}
\clip($(l4.center)+(-\w,\w)$)--++($1*(v)$)--++(2*\w,-2*\w)--++($-1*(v)$);
\path[vl]($(l4.center)+(\w,\w)$)rectangle($(l5.center)+(-\w,-\w)$);
\end{scope}
\begin{scope}
\clip($(l5.center)+(-\w,\w)$)--++($3*(h)+(\w,0)$)--++(0,-2*\w)--++($-3*(h)+(\w,0)$);
\path[hl]($(l5.center)+(\w,\w)$)rectangle($(l8.center)+(-\w,-\w)$);
\end{scope}
\end{tikzpicture}\tikzexternaldisable\;\quad
&
\quad\scalebox{0.8}{\begin{tikzpicture}[baseline=(current bounding box.center)]
\tikzstyle{out1}=[inner sep=0,minimum size=1.2mm,circle,draw=black,fill=black]
\tikzstyle{in1}=[inner sep=0,minimum size=1.2mm,circle,draw=black,fill=white]
\pgfmathsetmacro{\scalar}{4/3};
\pgfmathsetmacro{\unit}{\scalar*0.922/1.6};
\coordinate (vstep)at(0,-0.26*\unit);
\coordinate (hstep)at(0.20*\unit,0);
\draw[thick](8*\unit,0)--(9*\unit,0);
\foreach \x in {1,...,8}{
\draw[thick](\x*\unit-\unit,0)rectangle(\x*\unit,-\unit);}
\foreach \x in {1,...,7}{
\draw[thick](\x*\unit-\unit,-\unit)rectangle(\x*\unit,-2*\unit);}
\node[inner sep=0]at(0.5*\unit,-0.5*\unit){\scalebox{\scalar}{$+$}};
\node[inner sep=0]at(1.5*\unit,-0.5*\unit){\scalebox{\scalar}{$+$}};
\node[inner sep=0]at(2.5*\unit,-0.5*\unit){\scalebox{\scalar}{$+$}};
\node[inner sep=0]at(3.5*\unit,-0.5*\unit){\scalebox{\scalar}{$0$}};
\node[inner sep=0]at(4.5*\unit,-0.5*\unit){\scalebox{\scalar}{$0$}};
\node[inner sep=0]at(5.5*\unit,-0.5*\unit){\scalebox{\scalar}{$0$}};
\node[inner sep=0]at(6.5*\unit,-0.5*\unit){\scalebox{\scalar}{$0$}};
\node[inner sep=0]at(7.5*\unit,-0.5*\unit){\scalebox{\scalar}{$+$}};
\node[inner sep=0]at(0.5*\unit,-1.5*\unit){\scalebox{\scalar}{$+$}};
\node[inner sep=0]at(1.5*\unit,-1.5*\unit){\scalebox{\scalar}{$0$}};
\node[inner sep=0]at(2.5*\unit,-1.5*\unit){\scalebox{\scalar}{$+$}};
\node[inner sep=0]at(3.5*\unit,-1.5*\unit){\scalebox{\scalar}{$+$}};
\node[inner sep=0]at(4.5*\unit,-1.5*\unit){\scalebox{\scalar}{$0$}};
\node[inner sep=0]at(5.5*\unit,-1.5*\unit){\scalebox{\scalar}{$0$}};
\node[inner sep=0]at(6.5*\unit,-1.5*\unit){\scalebox{\scalar}{$+$}};
\node[inner sep=0]at($(8.5*\unit,0)+(vstep)$){$1$};
\node[inner sep=0]at($(8*\unit,-0.5*\unit)+(hstep)$){$2$};
\node[inner sep=0]at($(7.5*\unit,-1*\unit)+(vstep)$){$3$};
\node[inner sep=0]at($(7*\unit,-1.5*\unit)+(hstep)$){$4$};
\node[inner sep=0]at($(6.5*\unit,-2*\unit)+(vstep)$){$5$};
\node[inner sep=0]at($(5.5*\unit,-2*\unit)+(vstep)$){$6$};
\node[inner sep=0]at($(4.5*\unit,-2*\unit)+(vstep)$){$7$};
\node[inner sep=0]at($(3.5*\unit,-2*\unit)+(vstep)$){$8$};
\node[inner sep=0]at($(2.5*\unit,-2*\unit)+(vstep)$){$9$};
\node[inner sep=0]at($(1.5*\unit,-2*\unit)+(vstep)$){$10$};
\node[inner sep=0]at($(0.5*\unit,-2*\unit)+(vstep)$){$11$};
\end{tikzpicture}}\quad\qquad
\\ & & \\ \qquad$T$\rule[-8pt]{0pt}{0pt}\qquad & \quad$(W_U(T), W_L(T))$\quad & \quad$D$\quad\qquad \\ \hline & & \\
\quad\begin{tikzpicture}[baseline=(current bounding box.center)]
\tikzstyle{out1}=[inner sep=0,minimum size=2.4mm,circle,draw=black,fill=black,semithick]
\tikzstyle{in1}=[inner sep=0,minimum size=2.4mm,circle,draw=black,fill=white,semithick]
\tikzstyle{vertex}=[inner sep=0,minimum size=1.2mm,circle,draw=black,fill=black,semithick]
\pgfmathsetmacro{\l}{0.36};
\pgfmathsetmacro{\d}{0.30};
\pgfmathsetmacro{\s}{1.08};
\node[vertex](r)at(0,0){};
\node[vertex](rh)at($(r)+(-\s,0)$){};
\node[vertex](rv)at($(r)+(0,\s)$){};
\node[vertex](rvh)at($(rv)+(-\s,0)$){};
\node[vertex](rvhh)at($(rvh)+(-\s,0)$){};
\node[vertex](rvhv)at($(rvh)+(0,\s)$){};
\node[inner sep=0](rhh)at($(rh)+(-\l,0)$){};
\node[inner sep=0](rhv)at($(rh)+(0,\l)$){};
\node[inner sep=0](rvv)at($(rv)+(0,\l)$){};
\node[inner sep=0](rvhhh)at($(rvhh)+(-\l,0)$){};
\node[inner sep=0](rvhhv)at($(rvhh)+(0,\l)$){};
\node[inner sep=0](rvhvh)at($(rvhv)+(-\l,0)$){};
\node[inner sep=0](rvhvv)at($(rvhv)+(0,\l)$){};
\node[inner sep=0]at($(rhh)+(-0.18,0)$){$2$};
\node[inner sep=0]at($(rhv)+(0,0.24)$){$3$};
\node[inner sep=0]at($(rvhhh)+(-0.18,0)$){$4$};
\node[inner sep=0]at($(rvhhv)+(0,0.24)$){$5$};
\node[inner sep=0]at($(rvhvh)+(-0.18,0)$){$6$};
\node[inner sep=0]at($(rvhvv)+(0,0.24)$){$7$};
\node[inner sep=0]at($(rvv)+(0,0.24)$){$8$};
\path[thick](r)edge(rh) edge (rv) (rh)edge(rhh) edge(rhv) (rv)edge(rvh) edge(rvv) (rvh)edge(rvhh) edge(rvhv) (rvhh)edge(rvhhh) edge(rvhhv) (rvhv)edge(rvhvh) edge(rvhvv);
\end{tikzpicture}\quad
&
\quad\tikzexternalenable\begin{tikzpicture}[baseline=(current bounding box.center)]
\tikzstyle{hu}=[top color=red!10,bottom color=red!10,middle color=red,opacity=0.70]
\tikzstyle{vu}=[left color=red!10,right color=red!10,middle color=red,opacity=0.70]
\tikzstyle{hl}=[top color=blue!10,bottom color=blue!10,middle color=blue,opacity=0.55]
\tikzstyle{vl}=[left color=blue!10,right color=blue!10,middle color=blue,opacity=0.55]
\pgfmathsetmacro{\u}{0.75};
\pgfmathsetmacro{\w}{0.12*\u};
\coordinate(h)at(-\u,0);
\coordinate(v)at(0,-\u);
\draw[step=\u,color=black!16,ultra thick](0,0)grid(3*\u,2*\u);
\node[inner sep=0](l1)at(3*\u,2*\u){};
\node[inner sep=0](l2)at($(l1)+(h)$){};
\node[inner sep=0](l3)at($(l2)+(v)$){};
\node[inner sep=0](l4)at($(l3)+(h)$){};
\node[inner sep=0](l5)at($(l4)+(v)$){};
\node[inner sep=0](l6)at($(l5)+(h)$){};
\node[inner sep=0](u1)at(3*\u,2*\u){};
\node[inner sep=0](u2)at($(u1)+(h)$){};
\node[inner sep=0](u3)at($(u2)+(v)$){};
\node[inner sep=0](u4)at($(u3)+(h)$){};
\node[inner sep=0](u5)at($(u4)+(h)$){};
\node[inner sep=0](u6)at($(u5)+(v)$){};
\begin{scope}
\clip($(u1.center)+(0,\w)$)--++($1*(h)+(-\w,0)$)--++(2*\w,-2*\w)--++($-1*(h)+(-\w,0)$);
\path[hu]($(u1.center)+(\w,\w)$)rectangle($(u2.center)+(-\w,-\w)$);
\end{scope}
\begin{scope}
\clip($(u2.center)+(-\w,\w)$)--++($1*(v)$)--++(2*\w,-2*\w)--++($-1*(v)$);
\path[vu]($(u2.center)+(\w,\w)$)rectangle($(u3.center)+(-\w,-\w)$);
\end{scope}
\begin{scope}
\clip($(u3.center)+(-\w,\w)$)--++($2*(h)$)--++(2*\w,-2*\w)--++($-2*(h)$);
\path[hu]($(u3.center)+(\w,\w)$)rectangle($(u5.center)+(-\w,-\w)$);
\end{scope}
\begin{scope}
\clip($(u5.center)+(-\w,\w)$)--++($1*(v)+(0,-\w)$)--++(2*\w,0)--++($-1*(v)+(0,-\w)$);
\path[vu]($(u5.center)+(\w,\w)$)rectangle($(u6.center)+(-\w,-\w)$);
\end{scope}
\begin{scope}
\clip($(l1.center)+(0,\w)$)--++($1*(h)+(-\w,0)$)--++(2*\w,-2*\w)--++($-1*(h)+(-\w,0)$);
\path[hl]($(l1.center)+(\w,\w)$)rectangle($(l2.center)+(-\w,-\w)$);
\end{scope}
\begin{scope}
\clip($(l2.center)+(-\w,\w)$)--++($1*(v)$)--++(2*\w,-2*\w)--++($-1*(v)$);
\path[vl]($(l2.center)+(\w,\w)$)rectangle($(l3.center)+(-\w,-\w)$);
\end{scope}
\begin{scope}
\clip($(l3.center)+(-\w,\w)$)--++($1*(h)$)--++(2*\w,-2*\w)--++($-1*(h)$);
\path[hl]($(l3.center)+(\w,\w)$)rectangle($(l4.center)+(-\w,-\w)$);
\end{scope}
\begin{scope}
\clip($(l4.center)+(-\w,\w)$)--++($1*(v)$)--++(2*\w,-2*\w)--++($-1*(v)$);
\path[vl]($(l4.center)+(\w,\w)$)rectangle($(l5.center)+(-\w,-\w)$);
\end{scope}
\begin{scope}
\clip($(l5.center)+(-\w,\w)$)--++($1*(h)+(\w,0)$)--++(0,-2*\w)--++($-1*(h)+(\w,0)$);
\path[hl]($(l5.center)+(\w,\w)$)rectangle($(l6.center)+(-\w,-\w)$);
\end{scope}
\end{tikzpicture}\tikzexternaldisable\quad
& 
\quad\scalebox{0.8}{\begin{tikzpicture}[baseline=(current bounding box.center)]
\tikzstyle{out1}=[inner sep=0,minimum size=1.2mm,circle,draw=black,fill=black]
\tikzstyle{in1}=[inner sep=0,minimum size=1.2mm,circle,draw=black,fill=white]
\pgfmathsetmacro{\scalar}{4/3};
\pgfmathsetmacro{\unit}{\scalar*0.922/1.6};
\coordinate (vstep)at(0,-0.26*\unit);
\coordinate (hstep)at(0.20*\unit,0);
\draw[thick](6*\unit,0)--(7*\unit,0);
\foreach \x in {1,...,6}{
\draw[thick](\x*\unit-\unit,0)rectangle(\x*\unit,-\unit);}
\foreach \x in {1,...,5}{
\draw[thick](\x*\unit-\unit,-\unit)rectangle(\x*\unit,-2*\unit);}
\node[inner sep=0]at(0.5*\unit,-0.5*\unit){\scalebox{\scalar}{$+$}};
\node[inner sep=0]at(1.5*\unit,-0.5*\unit){\scalebox{\scalar}{$+$}};
\node[inner sep=0]at(2.5*\unit,-0.5*\unit){\scalebox{\scalar}{$+$}};
\node[inner sep=0]at(3.5*\unit,-0.5*\unit){\scalebox{\scalar}{$0$}};
\node[inner sep=0]at(4.5*\unit,-0.5*\unit){\scalebox{\scalar}{$0$}};
\node[inner sep=0]at(5.5*\unit,-0.5*\unit){\scalebox{\scalar}{$+$}};
\node[inner sep=0]at(0.5*\unit,-1.5*\unit){\scalebox{\scalar}{$+$}};
\node[inner sep=0]at(1.5*\unit,-1.5*\unit){\scalebox{\scalar}{$0$}};
\node[inner sep=0]at(2.5*\unit,-1.5*\unit){\scalebox{\scalar}{$+$}};
\node[inner sep=0]at(3.5*\unit,-1.5*\unit){\scalebox{\scalar}{$+$}};
\node[inner sep=0]at(4.5*\unit,-1.5*\unit){\scalebox{\scalar}{$+$}};
\node[inner sep=0]at($(6.5*\unit,0)+(vstep)$){$1$};
\node[inner sep=0]at($(6*\unit,-0.5*\unit)+(hstep)$){$2$};
\node[inner sep=0]at($(5.5*\unit,-1*\unit)+(vstep)$){$3$};
\node[inner sep=0]at($(5*\unit,-1.5*\unit)+(hstep)$){$4$};
\node[inner sep=0]at($(4.5*\unit,-2*\unit)+(vstep)$){$5$};
\node[inner sep=0]at($(3.5*\unit,-2*\unit)+(vstep)$){$6$};
\node[inner sep=0]at($(2.5*\unit,-2*\unit)+(vstep)$){$7$};
\node[inner sep=0]at($(1.5*\unit,-2*\unit)+(vstep)$){$8$};
\node[inner sep=0]at($(0.5*\unit,-2*\unit)+(vstep)$){$9$};
\end{tikzpicture}}\quad\qquad
\\ & & \\ \qquad$T'$\rule[-8pt]{0pt}{0pt}\qquad & \quad$(W_U(T'), W_L(T'))$\quad & \quad$D'$\quad\qquad \\ \hline
\end{tabular}
\caption{An example of Case 2 in the proof of \cref{thm:shift_map_holds}. Here $k=2$, $n=11$, $p=6$, $s=2$.}
\label{fig:finalH-example}
\end{center}
\end{figure}

In order to relate $\pi_{G(T)}$ and $\pi_{G(T')}$, we again introduce a plabic graph $G''$, obtained from $G(T')$ by inserting $s$ lollipops in between boundary vertices $p+1$ and $p+2$, and increasing the labels of the boundary vertices $p+2, p+3, \dots, n-s$ of $G(T')$ by $s$. We have
$$
\pi_{G''} = (p{+}2 \;\; p{+}3 \; \cdots \; n)^s\pi_{G(T')}(n \;\; n{-}1 \; \cdots \; p{+}2)^s,
$$
and we see from \cref{fig:finalH-tree-transformation} that $\pi_{G(T)} = \pi_{G''}(p \;\; p{+}1 \; \cdots \; p{+}s{+}1)^2$. This gives
\begin{align}\label{case_2_G_perm}
\begin{aligned}
\pi_{G(T)} &= (p{+}2 \;\; p{+}3 \; \cdots \; n)^s\pi_{G(T')}(n \;\; n{-}1 \; \cdots \; p{+}2)^s(p \;\; p{+}1 \; \cdots \; p{+}s{+}1)^2 \\
&= (p{+}2 \;\; p{+}3 \; \cdots \; n)^s\pi_{G(T')}(n \;\; n{-}1 \; \cdots \; p)^s.
\end{aligned}
\end{align}

On the other hand, by \cref{def:latticetocell}, $D$ is obtained from $D'$ by adding $s$ columns of all $0$'s; the bottom edges of these columns are labeled by $p, p+1, \dots, p+s-1$ when we label the southeast border of $D$ by $1, \dots, n$. Hence
$$
\pi_D = (p \;\; p{+}1 \; \cdots \; n)^s\pi_{D'}(n \;\; n{-}1 \; \cdots \; p)^s.
$$
Putting this together with the induction hypothesis $\pi_{D'} = c_{n-s}^2\pi_{G(T')}$ and \eqref{case_2_G_perm}, we obtain
\begin{align*}
\pi_D &= (p \;\; p{+}1 \; \cdots \; n)^sc_{n-s}^2\pi_{G(T')}(n \;\; n{-}1 \; \cdots \; p)^s \\
&= c_n^2(p{+}2 \;\; p{+}3 \; \cdots \; n)^s\pi_{G(T')}(n \;\; n{-}1 \; \cdots \; p)^s = c_n^2\pi_{G(T)}.\qedhere
\end{align*}
\end{pf}

\section{Number of cells in a decomposition of \texorpdfstring{$\mathcal{A}_{n,k,m}$}{A(n,k,m)} for arbitrary even \texorpdfstring{$m$}{m}}\label{sec:numerology}

\noindent Recall from \cref{sec:m2} that the $m=2$ amplituhedron $\mathcal{A}_{n,k,2}$ has a decomposition with $\binom{n-2}{k}$ top-dimensional cells, and from \cref{lem:bcfw-m=4} that the (conjectural) BCFW decomposition of the $m=4$ amplituhedron $\mathcal{A}_{n,k,4}$ contains $\frac{1}{n-3}\binom{n-3}{k+1}\binom{n-3}{k}$ top-dimensional cells. On the other hand, when $k=1$ the amplituhedron $\mathcal{A}_{n,1,m}$ is a cyclic polytope $C(n,m)$. Bayer \cite[Corollary 10]{bayer_93} (see also \cite[Corollary 1.2(ii)]{rambau_97}) showed that if $m$ is even, any triangulation of $C(n,m)$ has exactly $\binom{n-1-\frac{m}{2}}{\frac{m}{2}}$ top-dimensional simplices. In this section we conjecture a generalization of the above statements. We also give several families of combinatorial objects which are in bijection with the top-dimensional cells in our conjectural decomposition of $\mathcal{A}_{n,k,m}$.

For $a, b, c\in\mathbb{N}$, define
$$
M(a,b,c) := \prod_{i=1}^a\prod_{j=1}^b\prod_{k=1}^c\frac{i+j+k-1}{i+j+k-2}.
$$
Note that $M(a,b,c)$ is symmetric in $a,b,c$.

\begin{conj}\label{mainconj}
For even $m$, there is a cell decomposition of the amplituhedron $\mathcal{A}_{n,k,m}$, whose top-dimensional cells are the images of precisely $M(k, n-k-m, \frac{m}{2})$ cells of $\Gr_{k,n}^{\ge 0}$ of dimension $km$.
\end{conj}

We give a summary of all special cases in which \cref{mainconj} is known or conjecturally known:\\
\begin{center}
\setlength{\tabcolsep}{5pt}
\begin{tabular}{|c|c|c|}
\hline
{\bfseries special case} & $M(k, n-k-m, \frac{m}{2})$\rule{0pt}{13pt}\rule[-7pt]{0pt}{0pt} & {\bfseries explanation} \\ \hline
$m=0$ & $1$\rule{0pt}{13pt}\rule[-7pt]{0pt}{0pt} & $\mathcal{A}$ is a point \\ \hline
$m=2$ & $\displaystyle\binom{n-2}{k}$\rule{0pt}{21pt}\rule[-14pt]{0pt}{0pt} & \cite[Section 7]{ATT}\footnotemark \\ \hline
$m=4$ & $\displaystyle\frac{1}{n-3}\binom{n-3}{k+1}\binom{n-3}{k}$\rule{0pt}{21pt}\rule[-14pt]{0pt}{0pt} & conjectured in \cite{arkani-hamed_trnka} \\ \hline
$k=0$ & $1$\rule{0pt}{13pt}\rule[-7pt]{0pt}{0pt} & $\mathcal{A}$ is a point \\ \hline
$k=n-m$ & $1$\rule{0pt}{13pt}\rule[-7pt]{0pt}{0pt} & $\mathcal{A}\cong\Gr_{k,n}^{\ge 0}$  \\ \hline
$k=1$ & $\displaystyle\binom{n-1-\frac{m}{2}}{\frac{m}{2}}$\rule{0pt}{21pt}\rule[-14pt]{0pt}{0pt} & $\mathcal{A}\cong \text{cyclic polytope }C(n,m)$ \\ \hline
\end{tabular}\footnotetext{Arkani-Hamed, Thomas, and Trnka show that the $\binom{n-2}{k}$ top-dimensional cells in $\mathcal{A}_{n,k,2}$ are disjoint and cover a dense subset of the amplituhedron. It is not known whether this induces a cell decomposition.}
\end{center}

\begin{rmk}\label{oddm}
\cref{mainconj} only deals with the case of even $m$. For odd $m$, it is possible that a decomposition of $\mathcal{A}_{n+1,k,m+1}$ with $M(k,n-k-m,\frac{m+1}{2})$ top-dimensional cells could be used to give a decomposition of $\mathcal{A}_{n,k,m}$ with the same number of top-dimensional cells. This is the case when $m=1$, as two of us showed in \cite{karpwilliams}. It is also the case when $k=1$, since for odd $m$ there is a triangulation of the cyclic polytope $C(n,m)$ with $\binom{n-\frac{m+1}{2}}{\frac{m+1}{2}}$ top-dimensional simplices \cite[Corollary 1.2(ii)]{rambau_97}. When $m=3$, we came up with a natural construction of $M(k, n-k-3, 2)$ cells of $\Gr_{k,n}^{\ge 0}$ whose images we had hoped would give a decomposition of $\mathcal{A}_{n,k,3}$. However, these images inside $\mathcal{A}_{n,k,3}$ are not disjoint (see \cref{sec:m=3}). We mention that even in the case of cyclic polytopes, triangulations are not as well behaved in odd dimension: for even $m$ every triangulation of $C(n,m)$ has $\binom{n-1-\frac{m}{2}}{\frac{m}{2}}$ top-dimensional simplices, while for odd $m$ the number of top-dimensional simplices in a triangulation can lie anywhere between $\binom{n-1-\frac{m+1}{2}}{\frac{m-1}{2}}$ and $\binom{n-\frac{m+1}{2}}{\frac{m+1}{2}}$ \cite[Corollary 1.2(ii)]{rambau_97}.
\end{rmk}

The symmetry of $M(a,b,c)$ raises the following question.
\begin{question}
\cref{mainconj} suggests that there is a symmetry for amplituhedra $\mathcal{A}_{n,k,m}$ among the parameters $k$, $n-k-m$, and $\frac{m}{2}$. Is there an explanation for this symmetry?
\end{question}
In the case $m=4$, the symmetry between $k$ and $n-k-m$ comes from the well-known {\itshape parity} of the scattering amplitude (see \cite[Section 11]{ATT}). Below, we give a combinatorial explanation of this symmetry in terms of complete binary trees and decorated permutations. The possible symmetries between $\frac{m}{2}$ and the other two parameters are completely mysterious.
\begin{prop}\label{parity_symmetry}
Let $c_n := (n \;\; n{-}1 \; \cdots \; 2 \;\; 1)$ be the long cycle in the symmetric group on $[n]$, and $w_n$ the permutation given by $w_n(i) := n+1-i$ for $i\in [n]$. \\
(i) Given $T\in\mathcal{T}_{n,k,4}$, let $T'\in\mathcal{T}_{n,n-k-4,4}$ be obtained by reflecting $T$ through a line of slope $-1$ (i.e.\ by switching the roles of horizontal and vertical children). Then $\pi_{G(T')} = w_n\pi_{G(T)}w_n$.\\
(ii) The map $\pi\mapsto c_n^4w_n\pi w_n$ is an involution from $\Pi_{n,k,4}$ to $\Pi_{n,n-k-4,4}$.
\end{prop}

\begin{pf}
By \cref{def:tree-to-graph}, $G(T')$ is obtained from $G(T)$ by reflecting the plabic graph, interchanging black and white vertices, and reversing the order of the ground set. Reflecting and interchanging colors individually invert the decorated permutation (by \cref{def:rules}), and reversing the order of the ground set corresponds to conjugating by $w_n$. This proves part (i). Then using the fact that $w_nc_nw_n = c_n^{-1}$, we get $c_n^2\pi_{G(T')} = c_n^4w_n(c_n^2\pi_{G(T)})w_n$. Hence by \cref{trees_to_graphs}, the involution $T\mapsto T'$ from $\mathcal{T}_{n,k,4}$ to $\mathcal{T}_{n,n-k-4,4}$, which sends each tree to its reflection, induces the involution $\pi\mapsto c_n^4w_n\pi w_n$ from $\Pi_{n,k,4}$ to $\Pi_{n,n-k-4,4}$.
\end{pf}

Recall from \cref{sec:m2} that $\Pi_{n,k,2}$ is a set of decorated permutations which give a decomposition of $\mathcal{A}_{n,k,2}$. We can verify that the analogue of \cref{parity_symmetry}(ii) holds for $m=2$, i.e.\ $\pi\mapsto c_n^2w_n\pi w_n$ is an involution from $\Pi_{n,k,2}$ to $\Pi_{n,n-k-2,2}$.  This motivates the following question.
\begin{question}\label{parity_question}
Fix $n$ and $m$ with $m$ even, and assume for each $k$ we can find a collection $\Pi_{n,k,m}$ of decorated permutations corresponding to $km$-dimensional cells of $\Gr_{k,n}^{\ge 0}$, whose images induce a cell decomposition of $\mathcal{A}_{n,k,m}(Z)$. Can we choose $\Pi_{n,k,m}$ so that the map $\pi\mapsto c_n^mw_n\pi w_n$ is an involution from $\Pi_{n,k,m}$ to $\Pi_{n,n-k-m,m}$?
\end{question}
We do not expect \cref{parity_question} has a positive answer if $m$ is odd. Indeed, if we let $\Pi_{n,k,1}$ be the set of decorated permutations corresponding to the $m=1$ BCFW cells of $\Gr_{k,n}^{\ge 0}$ defined in \cite{karpwilliams}, then $\pi\mapsto c_nw_n\pi w_n$ does not in general take $\Pi_{n,k,1}$ to $\Pi_{n,n-k-1,1}$. (For example, $\Pi_{3,1,1} = \{(1 \;\; 2), (2 \;\; 3)\}$, and $c_3w_3(1 \;\; 2)w_3 = 321\notin\Pi_{3,1,1}$.)

\subsection{Combinatorial and geometric interpretations of $M(a,b,c)$}\label{sec:Nabc}

The number $M(a,b,c)$ has many interpretations which appear in the literature. We present some of them below. We refer to the article of Propp \cite{propp_99} for further background on this subject.
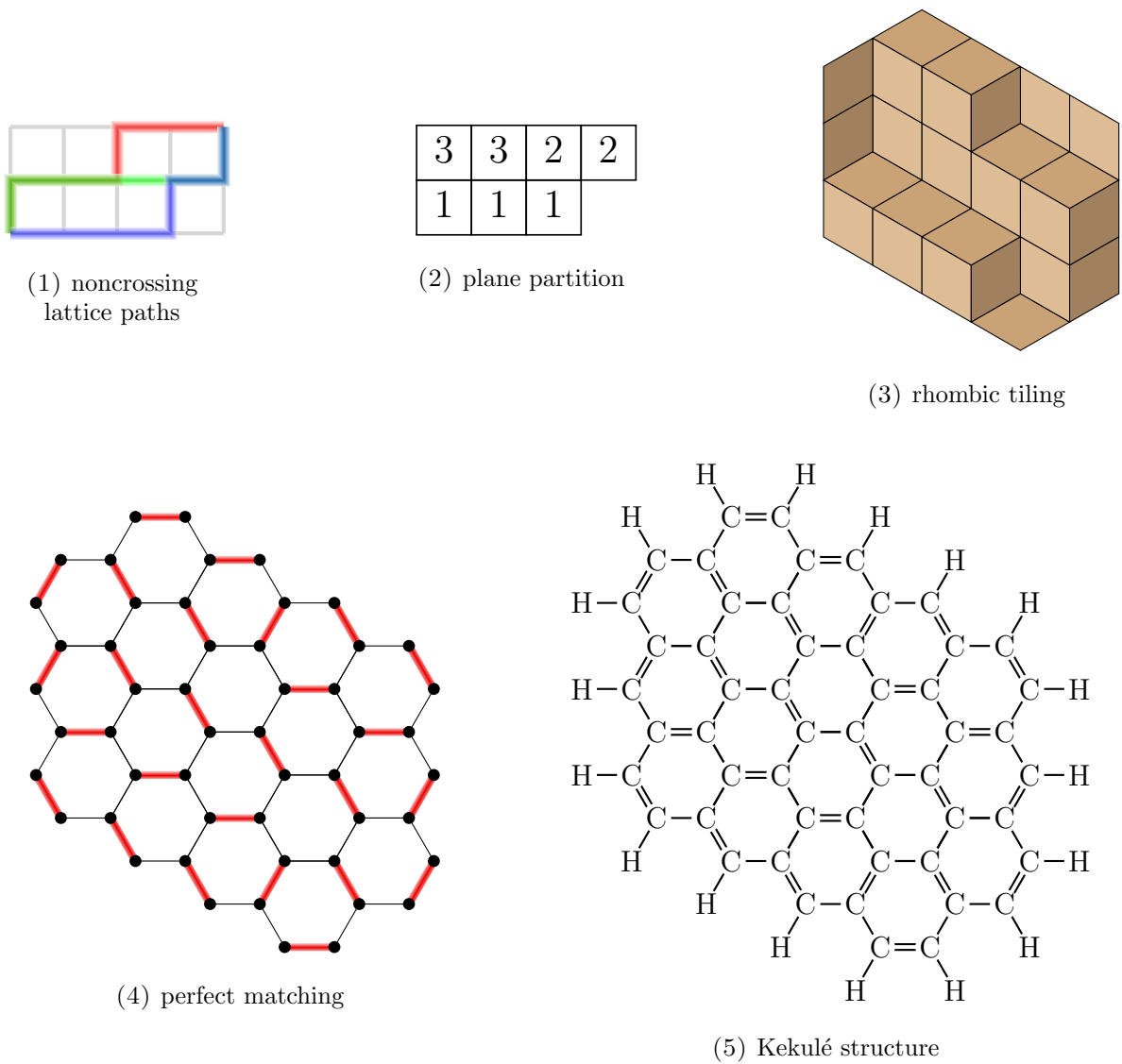
\begin{figure}[htb]
\captionsetup[subfigure]{margin=0pt,captionskip=12pt}
\renewcommand*\thesubfigure{\arabic{subfigure}}
\centering
\subfloat[][noncrossing lattice paths]{
\tikzexternalenable\begin{tikzpicture}[baseline=(current bounding box.center)]
\tikzstyle{hu}=[top color=red!10,bottom color=red!10,middle color=red,opacity=0.70]
\tikzstyle{vu}=[left color=red!10,right color=red!10,middle color=red,opacity=0.70]
\tikzstyle{hm}=[top color=green!10,bottom color=green!10,middle color=green,opacity=0.62]
\tikzstyle{vm}=[left color=green!10,right color=green!10,middle color=green,opacity=0.62]
\tikzstyle{hl}=[top color=blue!10,bottom color=blue!10,middle color=blue,opacity=0.55]
\tikzstyle{vl}=[left color=blue!10,right color=blue!10,middle color=blue,opacity=0.55]
\pgfmathsetmacro{\u}{0.75};
\pgfmathsetmacro{\w}{0.12*\u};
\coordinate(h)at(-\u,0);
\coordinate(v)at(0,-\u);
\draw[step=\u,color=black!16,ultra thick](0,0)grid(4*\u,2*\u);
\node[inner sep=0](l1)at(4*\u,2*\u){};
\node[inner sep=0](l2)at($(l1)+(v)$){};
\node[inner sep=0](l3)at($(l2)+(h)$){};
\node[inner sep=0](l4)at($(l3)+(v)$){};
\node[inner sep=0](l5)at($(l4)+(h)$){};
\node[inner sep=0](l6)at($(l5)+(h)$){};
\node[inner sep=0](l7)at($(l6)+(h)$){};
\node[inner sep=0](m1)at(4*\u,2*\u){};
\node[inner sep=0](m2)at($(m1)+(v)$){};
\node[inner sep=0](m3)at($(m2)+(h)$){};
\node[inner sep=0](m4)at($(m3)+(h)$){};
\node[inner sep=0](m5)at($(m4)+(h)$){};
\node[inner sep=0](m6)at($(m5)+(h)$){};
\node[inner sep=0](m7)at($(m6)+(v)$){};
\node[inner sep=0](u1)at(4*\u,2*\u){};
\node[inner sep=0](u2)at($(u1)+(h)$){};
\node[inner sep=0](u3)at($(u2)+(h)$){};
\node[inner sep=0](u4)at($(u3)+(v)$){};
\node[inner sep=0](u5)at($(u4)+(h)$){};
\node[inner sep=0](u6)at($(u5)+(h)$){};
\node[inner sep=0](u7)at($(u6)+(v)$){};
\begin{scope}
\clip($(u1.center)+(0,\w)$)--++($2*(h)+(-\w,0)$)--++(2*\w,-2*\w)--++($-2*(h)+(-\w,0)$);
\path[hu]($(u1.center)+(\w,\w)$)rectangle($(u3.center)+(-\w,-\w)$);
\end{scope}
\begin{scope}
\clip($(u3.center)+(-\w,\w)$)--++($1*(v)$)--++(2*\w,-2*\w)--++($-1*(v)$);
\path[vu]($(u3.center)+(\w,\w)$)rectangle($(u4.center)+(-\w,-\w)$);
\end{scope}
\begin{scope}
\clip($(u4.center)+(-\w,\w)$)--++($2*(h)$)--++(2*\w,-2*\w)--++($-2*(h)$);
\path[hu]($(u4.center)+(\w,\w)$)rectangle($(u6.center)+(-\w,-\w)$);
\end{scope}
\begin{scope}
\clip($(u6.center)+(-\w,\w)$)--++($1*(v)+(0,-\w)$)--++(2*\w,0)--++($-1*(v)+(0,-\w)$);
\path[vu]($(u6.center)+(\w,\w)$)rectangle($(u7.center)+(-\w,-\w)$);
\end{scope}
\begin{scope}
\clip($(m1.center)+(-\w,0)$)--++($1*(v)+(0,\w)$)--++(2*\w,-2*\w)--++($-1*(v)+(0,\w)$);
\path[vm]($(m1.center)+(\w,\w)$)rectangle($(m2.center)+(-\w,-\w)$);
\end{scope}
\begin{scope}
\clip($(m2.center)+(-\w,\w)$)--++($4*(h)$)--++(2*\w,-2*\w)--++($-4*(h)$);
\path[hm]($(m2.center)+(\w,\w)$)rectangle($(m6.center)+(-\w,-\w)$);
\end{scope}
\begin{scope}
\clip($(m6.center)+(-\w,\w)$)--++($1*(v)+(0,-\w)$)--++(2*\w,0)--++($-1*(v)+(0,-\w)$);
\path[vm]($(m6.center)+(\w,\w)$)rectangle($(m7.center)+(-\w,-\w)$);
\end{scope}
\begin{scope}
\clip($(l1.center)+(-\w,0)$)--++($1*(v)+(0,\w)$)--++(2*\w,-2*\w)--++($-1*(v)+(0,\w)$);
\path[vl]($(l1.center)+(\w,\w)$)rectangle($(l2.center)+(-\w,-\w)$);
\end{scope}
\begin{scope}
\clip($(l2.center)+(-\w,\w)$)--++($1*(h)$)--++(2*\w,-2*\w)--++($-1*(h)$);
\path[hl]($(l2.center)+(\w,\w)$)rectangle($(l3.center)+(-\w,-\w)$);
\end{scope}
\begin{scope}
\clip($(l3.center)+(-\w,\w)$)--++($1*(v)$)--++(2*\w,-2*\w)--++($-1*(v)$);
\path[vl]($(l3.center)+(\w,\w)$)rectangle($(l4.center)+(-\w,-\w)$);
\end{scope}
\begin{scope}
\clip($(l4.center)+(-\w,\w)$)--++($3*(h)+(\w,0)$)--++(0,-2*\w)--++($-3*(h)+(\w,0)$);
\path[hl]($(l4.center)+(\w,\w)$)rectangle($(l7.center)+(-\w,-\w)$);
\end{scope}
\end{tikzpicture}\tikzexternaldisable\label{fig:objects1}}\qquad\qquad\qquad
\subfloat[][plane partition]{
\begin{tikzpicture}[baseline=(current bounding box.center)]
\node[inner sep=0]at(0,0){\scalebox{1.33333333}{\begin{ytableau}
3 & 3 & 2 & 2 \\
1 & 1 & 1
\end{ytableau}}};
\end{tikzpicture}\label{fig:objects2}}\qquad\qquad\qquad
\subfloat[][rhombic tiling]{
\tikzexternalenable\begin{tikzpicture}[scale=0.8,rotate=-30,baseline=(current bounding box.center)]
\node[inner sep=0](a1)at($(0:-3)+(60:4)$){};
\node[inner sep=0](a2)at($(0:-2)+(60:4)$){};
\node[inner sep=0](a3)at($(0:0)+(60:3)$){};
\node[inner sep=0](a4)at($(0:1)+(60:3)$){};
\node[inner sep=0](a5)at($(0:-1)+(60:1)$){};
\node[inner sep=0](a6)at($(0:0)+(60:1)$){};
\node[inner sep=0](a7)at($(0:1)+(60:1)$){};
\node[inner sep=0](a8)at($(0:3)+(60:0)$){};
\node[inner sep=0](b1)at($(0:-1)+(60:5)$){};
\node[inner sep=0](b2)at($(0:0)+(60:5)$){};
\node[inner sep=0](b3)at($(0:-3)+(60:4)$){};
\node[inner sep=0](b4)at($(0:-2)+(60:4)$){};
\node[inner sep=0](b5)at($(0:-2)+(60:3)$){};
\node[inner sep=0](b6)at($(0:-1)+(60:3)$){};
\node[inner sep=0](b7)at($(0:0)+(60:3)$){};
\node[inner sep=0](b8)at($(0:1)+(60:3)$){};
\node[inner sep=0](b9)at($(0:2)+(60:2)$){};
\node[inner sep=0](b10)at($(0:-1)+(60:1)$){};
\node[inner sep=0](b11)at($(0:0)+(60:1)$){};
\node[inner sep=0](b12)at($(0:1)+(60:1)$){};
\node[inner sep=0](c1)at($(0:-1)+(60:5)$){};
\node[inner sep=0](c2)at($(0:-3)+(60:4)$){};
\node[inner sep=0](c3)at($(0:2)+(60:4)$){};
\node[inner sep=0](c4)at($(0:-2)+(60:3)$){};
\node[inner sep=0](c5)at($(0:3)+(60:3)$){};
\node[inner sep=0](c6)at($(0:2)+(60:2)$){};
\definecolor{bcolor}{rgb}{.867,0.737,0.588}\definecolor{acolor}{rgb}{.788,.639,.462}\definecolor{ccolor}{rgb}{.631,.502,.365}
\foreach \x in {1,...,8}{
\node[inner sep=0](au\x)at($(a\x)+(0:1)$){};
\node[inner sep=0](av\x)at($(au\x)+(60:1)$){};
\node[inner sep=0](aw\x)at($(av\x)+(0:-1)$){};
\draw[fill=acolor](a\x.center)--(au\x.center)--(av\x.center)--(aw\x.center)--(a\x.center);}
\foreach \x in {1,...,12}{
\node[inner sep=0](bu\x)at($(b\x)+(-60:1)$){};
\node[inner sep=0](bv\x)at($(bu\x)+(0:1)$){};
\node[inner sep=0](bw\x)at($(bv\x)+(-60:-1)$){};
\draw[fill=bcolor](b\x.center)--(bu\x.center)--(bv\x.center)--(bw\x.center)--(b\x.center);}
\foreach \x in {1,...,6}{
\node[inner sep=0](cu\x)at($(c\x)+(-120:1)$){};
\node[inner sep=0](cv\x)at($(cu\x)+(-60:1)$){};
\node[inner sep=0](cw\x)at($(cv\x)+(-120:-1)$){};
\draw[fill=ccolor](c\x.center)--(cu\x.center)--(cv\x.center)--(cw\x.center)--(c\x.center);}
\end{tikzpicture}\tikzexternaldisable\label{fig:objects3}}\\[8pt]
\subfloat[][perfect matching]{
\tikzexternalenable\begin{tikzpicture}[scale=0.70,rotate=-30,baseline=(current bounding box.center)]
\tikzstyle{vertex}=[inner sep=0,minimum size=4.8,circle,fill=black]
\definecolor{edgecolor}{rgb}{0.92,0,0}
\definecolor{egdecolor}{rgb}{0.92,0,0}
\tikzstyle{m-180}=[top color=egdecolor!10,bottom color=egdecolor!10,middle color=edgecolor,shading angle=50,opacity=1,fading speed=65]
\tikzstyle{m-60}=[top color=egdecolor!10,bottom color=egdecolor!10,middle color=edgecolor,shading angle=-180,opacity=1,fading speed=-20]
\tikzstyle{m60}=[top color=egdecolor!10,bottom color=egdecolor!10,middle color=edgecolor,shading angle=-50,opacity=1,fading speed=65]
\pgfmathsetmacro{\w}{0.08};
\foreach \i in {0,...,3}{
\draw[black](0:2*\i*cos{30}) -- ++(150:1) -- ++(-150:1) -- ++(-90:1) -- ++(-30:1) -- ++(30:1) -- ++(90:1);}
\foreach \i in {0,...,4}{
\draw[black]($(0:2*\i*cos{30})+(120:2*cos{30})$) -- ++(150:1) -- ++(-150:1) -- ++(-90:1) -- ++(-30:1) -- ++(30:1) -- ++(90:1);}
\foreach \i in {0,...,4}{
\draw[black]($(0:2*\i*cos{30})+(120:4*cos{30})$) -- ++(150:1) -- ++(-150:1) -- ++(-90:1) -- ++(-30:1) -- ++(30:1) -- ++(90:1);}
\foreach \i in {1,...,4}{
\draw[black]($(0:2*\i*cos{30})+(120:6*cos{30})$) -- ++(150:1) -- ++(-150:1) -- ++(-90:1) -- ++(-30:1) -- ++(30:1) -- ++(90:1);}
\path[m-60]($(0:2*3*cos{30}) + (120:0*2*cos{30}) + (-90:1)+(-60:\w)$)--++(-60:-2*\w)--++(-150:1)--++(-60:2*\w)--cycle;
\path[m-180]($(0:2*2*cos{30}) + (120:0*2*cos{30})+(-180:\w)$)--++(180:-2*\w)--++($(-90:1)$)--++(-180:2*\w)--cycle;
\path[m60]($(0:2*3*cos{30}) + (120:0*2*cos{30})+(60:\w)$)--++(60:-2*\w)--++($(150:1)$)--++(60:2*\w)--cycle;
\foreach \i in {0,...,2}{
\path[m60]($(0:\i*2*cos{30}) + (120:0*2*cos{30}) + (-150:1)+(-90:1)+(60:\w)$)--++(60:-2*\w)--++($(150:1)$)--++(60:2*\w)--cycle;}
\foreach \i in {0,...,2}{
\path[m-60]($(0:\i*2*cos{30}) + (120:1*2*cos{30}) + (-150:1) + (-90:1)+(-60:\w)$)--++(-60:-2*\w)--++(30:1)--++(-60:2*\w)--cycle;}
\foreach \i in {1,...,4}{
\path[m60]($(0:\i*2*cos{30}) + (120:2*2*cos{30}) + (-150:1)+(-90:1)+(60:\w)$)--++(60:-2*\w)--++($(150:1)$)--++(60:2*\w)--cycle;}
\path[m-180]($(180:2*cos{30}) + (120:1*2*cos{30})+(-180:\w)$)--++(180:-2*\w)--++($(-90:1)$)--++(-180:2*\w)--cycle;
\path[m-180]($(180:2*cos{30}) + (120:2*2*cos{30})+(-180:\w)$)--++(180:-2*\w)--++($(-90:1)$)--++(-180:2*\w)--cycle;
\path[m-180]($(0:2*4*cos{30}) + (120:0*2*cos{30}) + (150:1)+(-180:\w)$)--++(180:-2*\w)--++($(90:1)$)--++(-180:2*\w)--cycle;
\path[m-180]($(0:2*4*cos{30}) + (120:1*2*cos{30}) + (150:1)+(-180:\w)$)--++(180:-2*\w)--++($(90:1)$)--++(-180:2*\w)--cycle;
\path[m60]($(0:0*2*cos{30}) + (120:2*2*cos{30})+(60:\w)$)--++(60:-2*\w)--++($(150:1)$)--++(60:2*\w)--cycle;
\path[m60]($(0:1*2*cos{30}) + (120:2*2*cos{30})+(60:\w)$)--++(60:-2*\w)--++($(150:1)$)--++(60:2*\w)--cycle;
\path[m-60]($(0:2*2*cos{30}) + (120:2*2*cos{30})+(-60:\w)$)--++(-60:-2*\w)--++(30:1)--++(-60:2*\w)--cycle;
\path[m-60]($(0:3*2*cos{30}) + (120:2*2*cos{30})+(-60:\w)$)--++(-60:-2*\w)--++(30:1)--++(-60:2*\w)--cycle;
\path[m-180]($(0:2*2*cos{30}) + (120:2*2*cos{30}) + (150:1)+(-180:\w)$)--++(180:-2*\w)--++($(90:1)$)--++(-180:2*\w)--cycle;
\path[m-60]($(0:0*2*cos{30}) + (120:3*2*cos{30})+(-60:\w)$)--++(-60:-2*\w)--++(30:1)--++(-60:2*\w)--cycle;
\path[m-60]($(0:1*2*cos{30}) + (120:3*2*cos{30})+(-60:\w)$)--++(-60:-2*\w)--++(30:1)--++(-60:2*\w)--cycle;
\path[m60]($(0:3*2*cos{30}) + (120:3*2*cos{30})+(60:\w)$)--++(60:-2*\w)--++($(150:1)$)--++(60:2*\w)--cycle;
\path[m60]($(0:4*2*cos{30}) + (120:3*2*cos{30})+(60:\w)$)--++(60:-2*\w)--++($(150:1)$)--++(60:2*\w)--cycle;
\node[inner sep=0](s0)at($(-90:1)+(-150:1)$){};
\node[inner sep=0](s1)at($(s0)+(150:1)$){};
\node[inner sep=0](s2)at($(s1)+(90:1)$){};
\node[inner sep=0](s3)at($(s2)+(150:1)$){};
\node[inner sep=0](s4)at($(s3)+(90:1)$){};
\node[inner sep=0](s5)at($(s4)+(150:1)$){};
\node[inner sep=0](s6)at($(s5)+(90:1)$){};
\node[inner sep=0](s7)at($(s6)+(30:1)$){};
\node[inner sep=0](s8)at($(s7)+(90:1)$){};
\node[inner sep=0](s9)at($(s8)+(30:1)$){};
\foreach \i in {0,...,3}{
\node[vertex]at($(s0)+(-30:\i)+(30:\i)$){};}
\foreach \i in {0,...,4}{
\node[vertex]at($(s1)+(-30:\i)+(30:\i)$){};}
\foreach \i in {0,...,4}{
\node[vertex]at($(s2)+(-30:\i)+(30:\i)$){};}
\foreach \i in {0,...,5}{
\node[vertex]at($(s3)+(-30:\i)+(30:\i)$){};}
\foreach \i in {0,...,5}{
\node[vertex]at($(s4)+(-30:\i)+(30:\i)$){};}
\foreach \i in {0,...,5}{
\node[vertex]at($(s5)+(-30:\i)+(30:\i)$){};}
\foreach \i in {0,...,5}{
\node[vertex]at($(s6)+(-30:\i)+(30:\i)$){};}
\foreach \i in {0,...,4}{
\node[vertex]at($(s7)+(-30:\i)+(30:\i)$){};}
\foreach \i in {0,...,4}{
\node[vertex]at($(s8)+(-30:\i)+(30:\i)$){};}
\foreach \i in {0,...,3}{
\node[vertex]at($(s9)+(-30:\i)+(30:\i)$){};}
\end{tikzpicture}\tikzexternaldisable\label{fig:objects4}}\qquad\qquad
\subfloat[][Kekul\'{e} structure]{
\begin{tikzpicture}[scale=0.70,rotate=-30,baseline=(current bounding box.center)]
\pgfmathsetmacro{\scalar}{1};
\foreach \i in {0,...,3}{
\draw[thick] (0:2*\i*cos{30}) -- ++(150:1) -- ++(-150:1) -- ++(-90:1) -- ++(-30:1) -- ++(30:1) -- ++(90:1);}
\foreach \i in {0,...,4}{
\draw[thick] ($(0:2*\i*cos{30})+(120:2*cos{30})$) -- ++(150:1) -- ++(-150:1) -- ++(-90:1) -- ++(-30:1) -- ++(30:1) -- ++(90:1);}
\foreach \i in {0,...,4}{
\draw[thick] ($(0:2*\i*cos{30})+(120:4*cos{30})$) -- ++(150:1) -- ++(-150:1) -- ++(-90:1) -- ++(-30:1) -- ++(30:1) -- ++(90:1);}
\foreach \i in {1,...,4}{
\draw[thick] ($(0:2*\i*cos{30})+(120:6*cos{30})$) -- ++(150:1) -- ++(-150:1) -- ++(-90:1) -- ++(-30:1) -- ++(30:1) -- ++(90:1);}
\draw[style=double, double distance=1.5pt, thick] ($(0:2*3*cos{30}) + (120:0*2*cos{30}) + (-90:1)$) -- ++ (-150:1);
\draw[style=double, double distance=1.5pt, thick] ($(0:2*2*cos{30}) + (120:0*2*cos{30})$) -- ++ (-90:1);
\draw[style=double, double distance=1.5pt, thick] ($(0:2*3*cos{30}) + (120:0*2*cos{30})$) -- ++ (150:1);
\foreach \i in {0,...,2}{
\draw[style=double, double distance=1.5pt, thick] ($(0:\i*2*cos{30}) + (120:0*2*cos{30}) + (-150:1) + (-90:1)$) -- ++ (150:1); }
\foreach \i in {0,...,2}{
\draw[style=double, double distance=1.5pt, thick] ($(0:\i*2*cos{30}) + (120:1*2*cos{30}) + (-150:1) + (-90:1)$) -- ++ (30:1); }
\foreach \i in {1,...,4}{
\draw[style=double, double distance=1.5pt, thick] ($(0:\i*2*cos{30}) + (120:2*2*cos{30}) + (-150:1) + (-90:1)$) -- ++ (150:1); }
\draw[style=double, double distance=1.5pt, thick] ($(180:2*cos{30}) + (120:1*2*cos{30})$) -- ++ (-90:1);
\draw[style=double, double distance=1.5pt, thick] ($(180:2*cos{30}) + (120:2*2*cos{30})$) -- ++ (-90:1);
\draw[style=double, double distance=1.5pt, thick] ($(0:2*4*cos{30}) + (120:0*2*cos{30}) + (150:1)$) -- ++ (90:1);
\draw[style=double, double distance=1.5pt, thick] ($(0:2*4*cos{30}) + (120:1*2*cos{30}) + (150:1)$) -- ++ (90:1);
\draw[style=double, double distance=1.5pt, thick] ($(0:0*2*cos{30}) + (120:2*2*cos{30})$) -- ++ (150:1);
\draw[style=double, double distance=1.5pt, thick] ($(0:1*2*cos{30}) + (120:2*2*cos{30})$) -- ++ (150:1);
\draw[style=double, double distance=1.5pt, thick] ($(0:2*2*cos{30}) + (120:2*2*cos{30})$) -- ++ (30:1);
\draw[style=double, double distance=1.5pt, thick] ($(0:3*2*cos{30}) + (120:2*2*cos{30})$) -- ++ (30:1);
\draw[style=double, double distance=1.5pt, thick] ($(0:2*2*cos{30}) + (120:2*2*cos{30}) + (150:1)$) -- ++ (90:1);
\draw[style=double, double distance=1.5pt, thick] ($(0:0*2*cos{30}) + (120:3*2*cos{30})$) -- ++ (30:1);
\draw[style=double, double distance=1.5pt, thick] ($(0:1*2*cos{30}) + (120:3*2*cos{30})$) -- ++ (30:1);
\draw[style=double, double distance=1.5pt, thick] ($(0:3*2*cos{30}) + (120:3*2*cos{30})$) -- ++ (150:1);
\draw[style=double, double distance=1.5pt, thick] ($(0:4*2*cos{30}) + (120:3*2*cos{30})$) -- ++ (150:1);
\foreach \i in {0,...,2}{
\draw[thick] ($(-150:2) + (120:2*\i*cos{30})$) -- ++ (-150:1);
\node [circle, fill=white, inner sep=0] at ($(-150:2) + (120:2*\i*cos{30}) + (-150:1)$) {\scalebox{\scalar}{H}};
\draw[thick] ($(0:6*cos{30}) + (30:1) + (90:1) + (120:2*\i*cos{30})$) -- ++ (30:1);
\node [circle, fill=white, inner sep=0] at ($(0:6*cos{30}) + (30:2) + (90:1) + (120:2*\i*cos{30})$) {\scalebox{\scalar}{H}};
}
\draw[thick] ($(0:6*cos{30}) + (-90:1)$) -- ++ (-30:1);
\draw[thick] ($(0:6*cos{30}) + (30:1)$) -- ++ (-30:1);
\draw[thick] ($(120:6*cos{30})$) -- ++ (150:1);
\draw[thick] ($(120:6*cos{30}) + (-120:2*cos{30})$) -- ++ (150:1);
\node [circle, fill=white, inner sep=0] at ($(0:6*cos{30}) + (-90:1) + (-30:1)$) {\scalebox{\scalar}{H}};
\node [circle, fill=white, inner sep=0] at ($(0:6*cos{30}) + (30:1) + (-30:1)$) {\scalebox{\scalar}{H}};
\node [circle, fill=white, inner sep=0] at ($(120:6*cos{30}) + (150:1)$) {\scalebox{\scalar}{H}};
\node [circle, fill=white, inner sep=0] at ($(120:6*cos{30}) + (-120:2*cos{30}) + (150:1)$) {\scalebox{\scalar}{H}};

\foreach \i in {0,...,4}{
\node [circle, fill=white, inner sep=0] at ($(180:2*cos{30}) + (0:2*\i*cos{30})$) {\scalebox{\scalar}{C}};
\node [circle, fill=white, inner sep=0] at ($(180:2*cos{30}) + (-90:1) + (0:2*\i*cos{30})$) {\scalebox{\scalar}{C}};
\node [circle, fill=white, inner sep=0] at ($(180:2*cos{30}) + (90:3) + (150:1) + (0:2*\i*cos{30})$) {\scalebox{\scalar}{C}};
\node [circle, fill=white, inner sep=0] at ($(180:2*cos{30}) + (90:4) + (150:1) + (0:2*\i*cos{30})$) {\scalebox{\scalar}{C}};
}
\foreach \i in {0,...,3}{
\draw[thick] ($(180:2*cos{30}) + (90:5) + (0:2*\i*cos{30})$) -- ++ (90:1);  
\node [circle, fill=white, inner sep=0] at ($(180:2*cos{30}) + (90:5) + (0:2*\i*cos{30})$) {\scalebox{\scalar}{C}};
\node [circle, fill=white, inner sep=0] at ($(180:2*cos{30}) + (90:6) + (0:2*\i*cos{30})$) {\scalebox{\scalar}{H}};
\draw[thick] ($(180:2*cos{30}) + (-90:1) + (-30:1) + (0:2*\i*cos{30})$) -- ++ (-90:1);
\node [circle, fill=white, inner sep=0] at ($(180:2*cos{30}) + (-90:1) + (-30:1) + (0:2*\i*cos{30})$) {\scalebox{\scalar}{C}};
\node [circle, fill=white, inner sep=0] at ($(180:2*cos{30}) + (-90:2) + (-30:1) + (0:2*\i*cos{30})$) {\scalebox{\scalar}{H}};
}
\foreach \i in {0,...,5}{
\node [circle, fill=white, inner sep=0] at ($(180:2*cos{30}) + (150:1) + (0:2*\i*cos{30})$) {\scalebox{\scalar}{C}};
\node [circle, fill=white, inner sep=0] at ($(180:2*cos{30}) + (150:1) + (90:1)+ (0:2*\i*cos{30})$) {\scalebox{\scalar}{C}};
\node [circle, fill=white, inner sep=0] at ($(180:2*cos{30}) + (150:2) + (90:1)+ (0:2*\i*cos{30})$) {\scalebox{\scalar}{C}};
\node [circle, fill=white, inner sep=0] at ($(180:2*cos{30}) + (150:2) + (90:2)+ (0:2*\i*cos{30})$) {\scalebox{\scalar}{C}};
}
\end{tikzpicture}\label{fig:objects5}}
\caption{The objects appearing in \cref{prop:objects}. Here $(a,b,c) = (2,4,3)$.}
\label{fig:objects}
\end{figure}
\begin{prop}[{\cite[Sections IX and X]{macmahon}}]\label{prop:objects}
Fix $a,b,c\in\mathbb{N}$. Then $M(a,b,c)$ equals the number of the following objects (see \cref{fig:objects}):
\begin{enumerate}
\item collections of precisely $c$ noncrossing lattice paths inside an $a\times b$ rectangle;
\item \emph{plane partitions} which fit inside 
an $a \times b \times c$ box;
\item tilings of the hexagon $H(a,b,c) :=\begin{tikzpicture}[scale=0.4,baseline=(current bounding box.center)]
\pgfmathsetmacro{\a}{2};
\pgfmathsetmacro{\b}{4};
\pgfmathsetmacro{\c}{3};
\pgfmathsetmacro{\hspace}{0.9};
\pgfmathsetmacro{\vspace}{0.6};
\coordinate(v1)at(0,0);
\coordinate(v2)at($(v1)+(30:\a)$);
\coordinate(v3)at($(v2)+(-30:\b)$);
\coordinate(v4)at($(v3)+(-90:\c)$);
\coordinate(v5)at($(v4)+(-150:\a)$);
\coordinate(v6)at($(v5)+(150:\b)$);
\coordinate(tl)at($(v1)+(0,\a/2)+(-\hspace,0)$);
\coordinate(br)at($(v4)+(0,-\a/2)+(\hspace,-\vspace)$);
\useasboundingbox(tl)rectangle(br);
\path(v1)edge node[above left=-4pt]{$a$}(v2) (v2)edge node[above right=-4pt]{$b$}(v3) (v3)edge node[right=-2pt]{$c$}(v4) (v4)edge node[below right=-4pt]{$a$}(v5) (v5)edge node[below left=-4pt]{$b$}(v6) (v6)edge node[left=-2pt]{$c$}(v1);
\tkzMarkAngle[size=5mm,mark=||](v6,v1,v2)
\tkzMarkAngle[size=5mm,mark=||](v1,v2,v3)
\tkzMarkAngle[size=5mm,mark=||](v2,v3,v4)
\tkzMarkAngle[size=5mm,mark=||](v3,v4,v5)
\tkzMarkAngle[size=5mm,mark=||](v4,v5,v6)
\tkzMarkAngle[size=5mm,mark=||](v5,v6,v1)
\end{tikzpicture}$ by rhombi of the form $\,\begin{tikzpicture}[scale=0.92,baseline=(current bounding box.center)]
\coordinate(v1)at(0,0);
\coordinate(v2)at($(v1)+(30:1)$);
\coordinate(v3)at($(v2)+(-30:1)$);
\coordinate(v4)at($(v3)+(30:-1)$);
\path[use as bounding box](v1)edge node[above left=-4pt]{$1$}(v2) (v2)edge node[above right=-4pt]{$1$}(v3) (v3)edge node[below right=-4pt]{$1$}(v4) (v4)edge node[below left=-4pt]{$1$}(v1);
\tkzMarkAngle[size=2.08mm,mark=||](v1,v2,v3)
\tkzMarkAngle[size=2.208mm,mark=||](v3,v4,v1)
\tkzMarkAngle[size=2.8mm,mark=|](v2,v3,v4)
\tkzMarkAngle[size=2.8mm,mark=|](v4,v1,v2)
\end{tikzpicture}\,$;
\item perfect matchings of the {\itshape honeycomb lattice} $\mathcal{O}(a,b,c)$ defined below;
\item \emph{Kekul\'{e} structures} of a hexagon-shaped benzenoid with parameters $a,b,c$.
\end{enumerate}
\end{prop}
We briefly define the objects in \cref{prop:objects} and describe bijections between them, with reference to \cref{fig:objects}. Noncrossing lattice paths \subref{fig:objects1} were defined in \cref{def:noncrossing}. A {\itshape plane partition} \subref{fig:objects2} is a filling of the boxes of a Young diagram $Y$ with positive integers, such the numbers along each row (from left to right) and along each column (from top to bottom) are weakly decreasing. We get from a collection of noncrossing lattice paths \subref{fig:objects1} to a plane partition \subref{fig:objects2} by taking the Young diagram whose southeast border is the bottom lattice path, and writing in each box the number of lattice paths passing below it. We can depict a plane partition as a stacking of unit cubes in the nonnegative orthant of $\mathbb{R}^3$ \subref{fig:objects3}, by stacking $d$ unit cubes on top of each box of $Y$ filled with a $d$. We are considering stackings of cubes contained in an $a\times b\times c$ box, i.e.\ such that $Y$ is contained inside an $a\times b$ rectangle and the entries in its boxes are bounded above by $c$.\footnote{MacMahon was the first to enumerate any of the objects in \cref{prop:objects}, by showing that the number of plane partitions which fit inside an $a\times b\times c$ box equals $M(a,b,c)$ \cite[Sections IX and X]{macmahon}.} From a plane partition regarded as a stacking of cubes inside an $a\times b\times c$ box, we get a rhombic tiling of $H(a,b,c)$ \subref{fig:objects3} by orthogonally projecting the exterior surface of the stacking onto a symmetric affine plane.

We define the {\itshape honeycomb lattice} $\mathcal{O}(a,b,c)$ as the graph dual to the tiling of $H(a,b,c)$ by unit equilateral triangles; that is, the vertices of $\mathcal{O}(a,b,c)$ correspond to the triangles tiling $H(a,b,c)$, and two vertices of $\mathcal{O}(a,b,c)$ are adjacent precisely when the corresponding triangles share an edge. Alternatively, we obtain $\mathcal{O}(a,b,c)$ by gluing together regular hexagons into a hexagonal arrangement, as shown in \subref{fig:objects4}. A {\itshape perfect matching} of a graph is a subset of its edges which meets every vertex exactly once. From a rhombic tiling of $H(a,b,c)$ \subref{fig:objects3}, we obtain a perfect matching of $\mathcal{O}(a,b,c)$ \subref{fig:objects4} by including an edge in the matching if and only if the corresponding equilateral triangles in $H(a,b,c)$ are covered by the same rhombus.

Finally, from $\mathcal{O}(a,b,c)$, we obtain a {\itshape hexagon-shaped benzenoid} (with parameters $a,b,c$) by replacing each vertex by a carbon atom, and each edge by a bond between carbon atoms. Moreover, every carbon atom should be bonded to exactly three other atoms, so for each carbon atom bonded to only two others, we add a hydrogen atom bonded to it (by a single bond). A {\itshape Kekul\'{e} structure} of such a benzenoid specifies whether each bond between carbon atoms is a single bond or a double bond, subject to the condition that each carbon atom is tetravalent, i.e.\ it participates in exactly one double bond. (The tetravalency condition was first proposed by Kekul\'{e} \cite{kekule_65, kekule_66}.) We have a bijection from perfect matchings \subref{fig:objects4} to Kekul\'{e} structures \subref{fig:objects5}, which sends matched edges to double bonds and unmatched edges to single bonds.\footnote{In general, a {\itshape Kekul\'{e} structure} corresponds to a perfect matching of any graph formed by gluing together regular hexagons, not necessarily of $\mathcal{O}(a,b,c)$. Independently of work on plane partitions, the chemists Gordon and Davison \cite{gordon_davison_52} gave bijections between Kekul\'{e} structures \subref{fig:objects5}, perfect matchings \subref{fig:objects4}, and collections of noncrossing lattice paths \subref{fig:objects1}. (We thank Greg Kuperberg for bringing this to our attention.) They also state a formula suggested to them by Everett for the number of perfect matchings of $\mathcal{O}(a,b,b)$, and say that it is ``a special case of a more general equation established by Mr.\ M.\ Woodger'' in a forthcoming paper. Unfortunately, Woodger's work was never published. In later work on Kekul\'{e} structures, Cyvin \cite[(8)]{cyvin_1986} rediscovered MacMahon's formula for the number of perfect matchings of $\mathcal{O}(a,b,c)$, and it was reproven in \cite{bodroza_gutman_cyvin_tosic_1988}. We refer to \cite{kekule,klein_babic_trinajstic_02} for more on Kekul\'{e} structures.}
\begin{remark}
We note that there is a natural structure of \emph{distributive lattice} on the various combinatorial objects enumerated by $M(a,b,c)$ \cite[Theorem 2 and Example 2.1]{propp_02}. In terms of stackings of cubes inside an $a\times b\times c$ box, a cover relation in the lattice structure corresponds to adding a single unit cube.
It would be interesting to explore what this distributive lattice
structure tells us about the relative position of the corresponding cells in a decomposition of the amplituhedron.
Perhaps cover relations in the distributive lattice are related to whether the corresponding cells
are adjacent in the amplituhedron.
\end{remark}

\begin{remark}
In addition to the combinatorial interpretations given in \cref{prop:objects}, $M(a,b,c)$ also equals the dimension of the degree $c$
component of the homogeneous coordinate ring 
$\C[\Gr_{a,a+b}]$ \cite{hodge43}. So, the number of top-dimensional cells 
in the (conjectured) BCFW decomposition of $\mathcal{A}_{n,k,4}$ is equal to the dimension of the degree $2$ part of $\C[\Gr_{k,n-4}]$.  It would be interesting to 
give a geometric explanation of this statement.
\end{remark}

\section{Disjointness for BCFW cells when \texorpdfstring{$k=1$}{k=1}}\label{sec:warmup}

\noindent For completeness, and as a further warmup to proving that the images of the $k=2$ BCFW cells are disjoint in the amplituhedron (\cref{sec:k=2-disjointness}), we prove disjointness in the case $k=1$. This follows from the work of Rambau on triangulations of cyclic polytopes \cite{rambau_97}: any amplituhedron $\mathcal{A}_{n,1,4}(Z)$ is a $4$-dimensional cyclic polytope with $n$ vertices \cite{Sturmfels}, and the triangulation given by the BCFW recursion is among of the ones identified by Rambau. We give his argument below, rephrased in the language of sign variation and dominoes.

When $k=1$, the BCFW cells have the following explicit description, which we can verify from \cref{thm:BCFW} (see also \cite[(5.2)]{arkani-hamed_trnka}).

\begin{lem}\label{lem:bcfw-14}
For $k=1$, the positroids (see \cref{def:positroid}) of the $(k,n)$-BCFW cells are precisely $M = \{\{i\}, \{i+1\}, \{j\}, \{j+1\}, \{n\}\}$ for all $i,j\in [n-2]$ with $i+1 < j$.
\end{lem}
Then \cref{amp-conj} when $k=1$ says that every $4$-dimensional cyclic polytope with $n$ vertices is triangulated by the simplices whose vertex sets are of the form $\{i, i+1, j, j+1, n\}$. This is a special case of \cite[Theorem 4.2]{rambau_97}\footnote{In more detail, the BCFW triangulation is an iterated extension (in the sense of \cite[Definition 4.1]{rambau_97}) of the triangulation of a $1$-dimensional cyclic polytope into intervals between consecutive vertices.}.
\begin{prop}[\cite{rambau_97}]\label{prop:k=1-disjointness}
Let $Z\in\Mat_{5,n}^{>0}$. Then $\tilde{Z}$ maps the $(1,n)$-BCFW cells $\mathcal{C}_{n,1,4}$ injectively into the amplituhedron $\mathcal{A}_{n,1,4}(Z)$, and their images are pairwise disjoint.
\end{prop}

\begin{pf}[cf.\ {\cite[Remark 3.8]{rambau_97}}]
Let $V,V'\in\Gr_{1,n}^{\ge 0}$ be subspaces each contained in a $(1,n)$-BCFW cell, such that $\tilde{Z}(V) = \tilde{Z}(V')$. We must show that $V = V'$. Let $v\in\mathbb{R}^n$ be a basis vector of $V$, and $v'\in V'$ its matching vector as in \cref{criterion}. By \cref{lem:bcfw-14}, we can write $v-v'$ as a sum of $5$ or fewer dominoes (recall \cref{def:domino_short}). Hence $\var(v-v')\le 4$ by \cref{domino_sum}, whence $v = v'$ by \cref{criterion}.
\end{pf}

\begin{remark}
This argument generalizes to all $m$, to show that the images in the amplituhedron $\mathcal{A}_{n,1,m}(Z)$ of certain $m$-dimensional cells of $\Gr_{1,n}^{\ge 0}$ are mutually disjoint. When $m$ is even, these cells are indexed by the collection of positroids of the form $\{\{i_1\}, \{i_1+1\}, \dots, \{i_{m/2}\}, \{i_{m/2}+1\}, \{n\}\}$. When $m$ is odd, we can take the collection of positroids of the form $\{\{i_1\}, \{i_1+1\}, \dots, \{i_{(m+1)/2}\}, \{i_{(m+1)/2}+1\}\}$, or alternatively the collection of positroids of the form $\{\{1\}, \{i_1\}, \{i_1+1\}, \dots, \{i_{(m-1)/2}\}, \{i_{(m-1)/2}+1\}, \{n\}\}$. These in fact give triangulations of $\mathcal{A}_{n,1,m}(Z)$ by \cite[Theorem 4.2]{rambau_97}.
\end{remark}

\section{Domino bases for BCFW cells when \texorpdfstring{$k=2$}{k=2}}\label{sec:k=2-bases}

\noindent In this section we give a basis classification of the BCFW cells $\mathcal{C}_{n,2,4}$ (the case $k=2$) in terms of domino bases. We propose a generalization to all $k$ in the appendix (\cref{conj:domino}). We recall the definition of a domino (\cref{def:domino_short}), and introduce some new terminology.
\begin{defn}\label{def:domino}
We say that $d\in\mathbb{R}^n\setminus\{0\}$ is a {\itshape domino} if there exists $i\in [n]$ such that $d_j = 0$ for all $j\neq i, i+1$, and $d_i$ and $d_{i+1}$ are nonzero and have the same sign. (If $i = n$, then we require that $d_j = 0$ for all $j < n$, and that $d_n$ is nonzero.) In this case, we call $d$ an {\itshape $i$-domino}, we call $i$ the {\itshape index} of $d$, and we call the common sign of $d_i$ and $d_{i+1}$ the {\itshape sign} of $d$. We also regard the zero vector as a domino, with sign zero.
\end{defn}
For example, $(0,-1,-2,0,0)\in\mathbb{R}^5$ is a negative $2$-domino.
\begin{defn}\label{def:bar}
Given $v\in\mathbb{R}^n$ with $n\ge 1$, let $\overline{v}\in\mathbb{R}^n$ denote the vector obtained from $v$ by setting coordinate $n$ to $0$.
For $v\in\mathbb{R}^n$, we say that $v$ is \emph{orthodox} if $\overline{v}$ is a sum of two nonzero dominoes of the same sign with disjoint support, and \emph{deviant} if $\overline{v}$ is a difference of two such nonzero dominoes.
\end{defn}

\begin{figure}
\begin{align*}
& \textnormal{\bfseries Class 1.} &&
\begin{tikzpicture}[baseline=(current bounding box.center),x=16pt,y=16pt]
\useasboundingbox(-0.51,1)rectangle(0.51,-1);
\node at(0,0){\begin{ytableau}+ \\ +\end{ytableau}};
\end{tikzpicture}
\begin{tikzpicture}[baseline=(current bounding box.center),x=16pt,y=16pt]
\useasboundingbox(0,0)rectangle(0.3,-2.0635);
\draw[fill=black](0,0)rectangle(0.3,-2.0635);
\end{tikzpicture}
\begin{tikzpicture}[baseline=(current bounding box.center),x=16pt,y=16pt]
\useasboundingbox(-1.02,1)rectangle(1.02,-1);
\node at(0,0){\begin{ytableau}0 & 0 \\ + & + \end{ytableau}};
\end{tikzpicture}
\begin{tikzpicture}[baseline=(current bounding box.center),x=16pt,y=16pt]
\useasboundingbox(0,0)rectangle(0.3,-2.0635);
\draw[fill=black](0,0)rectangle(0.3,-2.0635);
\end{tikzpicture}
\begin{tikzpicture}[baseline=(current bounding box.center),x=16pt,y=16pt]
\useasboundingbox(-0.51,1)rectangle(0.51,-1);
\node at(0,0){\begin{ytableau}0 \\ +\end{ytableau}};
\end{tikzpicture}
\begin{tikzpicture}[baseline=(current bounding box.center),x=16pt,y=16pt]
\useasboundingbox(0,0)rectangle(0.3,-2.0635);
\draw[fill=black](0,0)rectangle(0.3,-1.0344);
\end{tikzpicture}
\begin{tikzpicture}[baseline=(current bounding box.center),x=16pt,y=16pt]
\useasboundingbox(-1.02,1)rectangle(1.02,-1);
\node at(0,0){\begin{ytableau}+ & + \\ \none\end{ytableau}};
\end{tikzpicture}
\begin{tikzpicture}[baseline=(current bounding box.center),x=16pt,y=16pt]
\useasboundingbox(0,0)rectangle(0.3,-2.0635);
\draw[fill=black](0,0)rectangle(0.3,-1.0344);
\end{tikzpicture}
\begin{tikzpicture}[baseline=(current bounding box.center),x=16pt,y=16pt]
\useasboundingbox(-0.51,1)rectangle(0.51,-1);
\node at(0,0){\begin{ytableau}+ \\ \none\end{ytableau}};
\end{tikzpicture} && \quad\begin{tabular}{|cc|cc|cc|cc|c}
+&+& +&+& 0&0& 0&0 &- \\
\hline
0&0& 0&0& +&+& +&+ &+ \\
\end{tabular} && \textnormal{orthodox} \\[4pt]
& \textnormal{\bfseries Class 2.} &&
\begin{tikzpicture}[baseline=(current bounding box.center),x=16pt,y=16pt]
\useasboundingbox(-0.51,1)rectangle(0.51,-1);
\node at(0,0){\begin{ytableau}+ \\ +\end{ytableau}};
\end{tikzpicture}
\begin{tikzpicture}[baseline=(current bounding box.center),x=16pt,y=16pt]
\useasboundingbox(0,0)rectangle(0.3,-2.0635);
\draw[fill=black](0,0)rectangle(0.3,-2.0635);
\end{tikzpicture}
\begin{tikzpicture}[baseline=(current bounding box.center),x=16pt,y=16pt]
\useasboundingbox(-1.02,1)rectangle(1.02,-1);
\node at(0,0){\begin{ytableau}0 & 0 \\ + & + \end{ytableau}};
\end{tikzpicture}
\begin{tikzpicture}[baseline=(current bounding box.center),x=16pt,y=16pt]
\useasboundingbox(0,0)rectangle(0.3,-2.0635);
\draw[fill=black](0,0)rectangle(0.3,-2.0635);
\end{tikzpicture}
\begin{tikzpicture}[baseline=(current bounding box.center),x=16pt,y=16pt]
\useasboundingbox(-1.02,1)rectangle(1.02,-1);
\node at(0,0){\begin{ytableau}+ & + \\ +\end{ytableau}};
\end{tikzpicture}
\begin{tikzpicture}[baseline=(current bounding box.center),x=16pt,y=16pt]
\useasboundingbox(0,0)rectangle(0.3,-2.0635);
\draw[fill=black](0,0)rectangle(0.3,-1.0344);
\end{tikzpicture}
\begin{tikzpicture}[baseline=(current bounding box.center),x=16pt,y=16pt]
\useasboundingbox(-0.51,1)rectangle(0.51,-1);
\node at(0,0){\begin{ytableau}+ \\ \none\end{ytableau}};
\end{tikzpicture} && \quad\begin{tabular}{|cc|ccc|cc|c}
+&+& +&+&0& 0&0 &- \\
\hline
0&0& 0&+&+& +&+ &+ \\
\end{tabular} && \textnormal{orthodox} \\[4pt]
& \textnormal{\bfseries Class 3.} &&
\begin{tikzpicture}[baseline=(current bounding box.center),x=16pt,y=16pt]
\useasboundingbox(-0.51,1)rectangle(0.51,-1);
\node at(0,0){\begin{ytableau}+ \\ +\end{ytableau}};
\end{tikzpicture}
\begin{tikzpicture}[baseline=(current bounding box.center),x=16pt,y=16pt]
\useasboundingbox(0,0)rectangle(0.3,-2.0635);
\draw[fill=black](0,0)rectangle(0.3,-2.0635);
\end{tikzpicture}
\begin{tikzpicture}[baseline=(current bounding box.center),x=16pt,y=16pt]
\useasboundingbox(-1.02,1)rectangle(1.02,-1);
\node at(0,0){\begin{ytableau}0 & + \\ + & + \end{ytableau}};
\end{tikzpicture}
\begin{tikzpicture}[baseline=(current bounding box.center),x=16pt,y=16pt]
\useasboundingbox(0,0)rectangle(0.3,-2.0635);
\draw[fill=black](0,0)rectangle(0.3,-2.0635);
\end{tikzpicture}
\begin{tikzpicture}[baseline=(current bounding box.center),x=16pt,y=16pt]
\useasboundingbox(-0.51,1)rectangle(0.51,-1);
\node at(0,0){\begin{ytableau}+ \\ +\end{ytableau}};
\end{tikzpicture}
\begin{tikzpicture}[baseline=(current bounding box.center),x=16pt,y=16pt]
\useasboundingbox(0,0)rectangle(0.3,-2.0635);
\draw[fill=black](0,0)rectangle(0.3,-1.0344);
\end{tikzpicture}
\begin{tikzpicture}[baseline=(current bounding box.center),x=16pt,y=16pt]
\useasboundingbox(-0.51,1)rectangle(0.51,-1);
\node at(0,0){\begin{ytableau}+ \\ \none\end{ytableau}};
\end{tikzpicture} && \quad\begin{tabular}{|cc|cc|cc|c}
+&+& +&+& 0&0 &- \\
\hline
0&0& +&+& +&+ &+ \\
\end{tabular} && \textnormal{orthodox} \\[4pt]
& \textnormal{\bfseries Class 4.} &&
\begin{tikzpicture}[baseline=(current bounding box.center),x=16pt,y=16pt]
\useasboundingbox(-0.51,1)rectangle(0.51,-1);
\node at(0,0){\begin{ytableau}+ \\ +\end{ytableau}};
\end{tikzpicture}
\begin{tikzpicture}[baseline=(current bounding box.center),x=16pt,y=16pt]
\useasboundingbox(0,0)rectangle(0.3,-2.0635);
\draw[fill=black](0,0)rectangle(0.3,-2.0635);
\end{tikzpicture}
\begin{tikzpicture}[baseline=(current bounding box.center),x=16pt,y=16pt]
\useasboundingbox(-1.02,1)rectangle(1.02,-1);
\node at(0,0){\begin{ytableau}+ & + \\ + & + \end{ytableau}};
\end{tikzpicture}
\begin{tikzpicture}[baseline=(current bounding box.center),x=16pt,y=16pt]
\useasboundingbox(0,0)rectangle(0.3,-2.0635);
\draw[fill=black](0,0)rectangle(0.3,-2.0635);
\end{tikzpicture}
\begin{tikzpicture}[baseline=(current bounding box.center),x=16pt,y=16pt]
\useasboundingbox(-0.51,1)rectangle(0.51,-1);
\node at(0,0){\begin{ytableau}0 \\ +\end{ytableau}};
\end{tikzpicture}
\begin{tikzpicture}[baseline=(current bounding box.center),x=16pt,y=16pt]
\useasboundingbox(0,0)rectangle(0.3,-2.0635);
\draw[fill=black](0,0)rectangle(0.3,-1.0344);
\end{tikzpicture}
\begin{tikzpicture}[baseline=(current bounding box.center),x=16pt,y=16pt]
\useasboundingbox(-0.51,1)rectangle(0.51,-1);
\node at(0,0){\begin{ytableau}+ \\ \none\end{ytableau}};
\end{tikzpicture} && \quad\begin{tabular}{|cc|cc|cc|c}
+&+& 0&0& -&- &- \\
\hline
+&+& +&+& +&+ &0 \\
\end{tabular} && \textnormal{deviant} \\[4pt]
& \textnormal{\bfseries Class 5.} &&
\begin{tikzpicture}[baseline=(current bounding box.center),x=16pt,y=16pt]
\useasboundingbox(-0.51,1)rectangle(0.51,-1);
\node at(0,0){\begin{ytableau}+ \\ +\end{ytableau}};
\end{tikzpicture}
\begin{tikzpicture}[baseline=(current bounding box.center),x=16pt,y=16pt]
\useasboundingbox(0,0)rectangle(0.3,-2.0635);
\draw[fill=black](0,0)rectangle(0.3,-2.0635);
\end{tikzpicture}
\begin{tikzpicture}[baseline=(current bounding box.center),x=16pt,y=16pt]
\useasboundingbox(-1.53,1)rectangle(1.53,-1);
\node at(0,0){\begin{ytableau}+ & + & 0 \\ 0 & + & + \end{ytableau}};
\end{tikzpicture}
\begin{tikzpicture}[baseline=(current bounding box.center),x=16pt,y=16pt]
\useasboundingbox(0,0)rectangle(0.3,-2.0635);
\draw[fill=black](0,0)rectangle(0.3,-2.0635);
\end{tikzpicture}
\begin{tikzpicture}[baseline=(current bounding box.center),x=16pt,y=16pt]
\useasboundingbox(-0.51,1)rectangle(0.51,-1);
\node at(0,0){\begin{ytableau}0 \\ +\end{ytableau}};
\end{tikzpicture}
\begin{tikzpicture}[baseline=(current bounding box.center),x=16pt,y=16pt]
\useasboundingbox(0,0)rectangle(0.3,-2.0635);
\draw[fill=black](0,0)rectangle(0.3,-1.0344);
\end{tikzpicture}
\begin{tikzpicture}[baseline=(current bounding box.center),x=16pt,y=16pt]
\useasboundingbox(-0.51,1)rectangle(0.51,-1);
\node at(0,0){\begin{ytableau}+ \\ \none\end{ytableau}};
\end{tikzpicture} && \quad\begin{tabular}{|cc|cc|ccc|c}
+&+& 0&0& 0&-&- &- \\
\hline
+&+& +&+& +&+&0 &0 \\
\end{tabular} && \textnormal{deviant} \\[4pt]
& \textnormal{\bfseries Class 6.} &&
\begin{tikzpicture}[baseline=(current bounding box.center),x=16pt,y=16pt]
\useasboundingbox(-0.51,1)rectangle(0.51,-1);
\node at(0,0){\begin{ytableau}+ \\ +\end{ytableau}};
\end{tikzpicture}
\begin{tikzpicture}[baseline=(current bounding box.center),x=16pt,y=16pt]
\useasboundingbox(0,0)rectangle(0.3,-2.0635);
\draw[fill=black](0,0)rectangle(0.3,-2.0635);
\end{tikzpicture}
\begin{tikzpicture}[baseline=(current bounding box.center),x=16pt,y=16pt]
\useasboundingbox(-1.02,1)rectangle(1.02,-1);
\node at(0,0){\begin{ytableau}+ & + \\ 0 & 0 \end{ytableau}};
\end{tikzpicture}
\begin{tikzpicture}[baseline=(current bounding box.center),x=16pt,y=16pt]
\useasboundingbox(0,0)rectangle(0.3,-2.0635);
\draw[fill=black](0,0)rectangle(0.3,-2.0635);
\end{tikzpicture}
\begin{tikzpicture}[baseline=(current bounding box.center),x=16pt,y=16pt]
\useasboundingbox(-1.02,1)rectangle(1.02,-1);
\node at(0,0){\begin{ytableau}0 & 0 \\ + & +\end{ytableau}};
\end{tikzpicture}
\begin{tikzpicture}[baseline=(current bounding box.center),x=16pt,y=16pt]
\useasboundingbox(0,0)rectangle(0.3,-2.0635);
\draw[fill=black](0,0)rectangle(0.3,-2.0635);
\end{tikzpicture}
\begin{tikzpicture}[baseline=(current bounding box.center),x=16pt,y=16pt]
\useasboundingbox(-0.51,1)rectangle(0.51,-1);
\node at(0,0){\begin{ytableau}0 \\ +\end{ytableau}};
\end{tikzpicture}
\begin{tikzpicture}[baseline=(current bounding box.center),x=16pt,y=16pt]
\useasboundingbox(0,0)rectangle(0.3,-2.0635);
\draw[fill=black](0,0)rectangle(0.3,-1.0344);
\end{tikzpicture}
\begin{tikzpicture}[baseline=(current bounding box.center),x=16pt,y=16pt]
\useasboundingbox(-0.51,1)rectangle(0.51,-1);
\node at(0,0){\begin{ytableau}+ \\ \none\end{ytableau}};
\end{tikzpicture} && \quad\begin{tabular}{|cc|cc|cc|cc|c}
+&+& 0&0& 0&0& -&- &- \\
\hline
+&+& +&+& +&+& 0&0 &0 \\
\end{tabular} && \textnormal{deviant} \\[4pt]
& \textnormal{\bfseries Class 7.} &&
\begin{tikzpicture}[baseline=(current bounding box.center),x=16pt,y=16pt]
\useasboundingbox(-0.51,1)rectangle(0.51,-1);
\node at(0,0){\begin{ytableau}+ \\ +\end{ytableau}};
\end{tikzpicture}
\begin{tikzpicture}[baseline=(current bounding box.center),x=16pt,y=16pt]
\useasboundingbox(0,0)rectangle(0.3,-2.0635);
\draw[fill=black](0,0)rectangle(0.3,-2.0635);
\end{tikzpicture}
\begin{tikzpicture}[baseline=(current bounding box.center),x=16pt,y=16pt]
\useasboundingbox(-1.02,1)rectangle(1.02,-1);
\node at(0,0){\begin{ytableau}+ & + \\ 0 & 0 \end{ytableau}};
\end{tikzpicture}
\begin{tikzpicture}[baseline=(current bounding box.center),x=16pt,y=16pt]
\useasboundingbox(0,0)rectangle(0.3,-2.0635);
\draw[fill=black](0,0)rectangle(0.3,-2.0635);
\end{tikzpicture}
\begin{tikzpicture}[baseline=(current bounding box.center),x=16pt,y=16pt]
\useasboundingbox(-1.02,1)rectangle(1.02,-1);
\node at(0,0){\begin{ytableau}0 & 0 \\ + & + \end{ytableau}};
\end{tikzpicture}
\begin{tikzpicture}[baseline=(current bounding box.center),x=16pt,y=16pt]
\useasboundingbox(0,0)rectangle(0.3,-2.0635);
\draw[fill=black](0,0)rectangle(0.3,-2.0635);
\end{tikzpicture}
\begin{tikzpicture}[baseline=(current bounding box.center),x=16pt,y=16pt]
\useasboundingbox(-0.51,1)rectangle(0.51,-1);
\node at(0,0){\begin{ytableau}+ \\ +\end{ytableau}};
\end{tikzpicture} && \quad\begin{tabular}{|ccc|cc|cc|c}
+&+&0& 0&0& -&- &- \\
\hline
+&+&+& +&+& 0&0 &0 \\
\end{tabular} && \textnormal{deviant} \\[4pt]
& \textnormal{\bfseries Class 8.} &&
\begin{tikzpicture}[baseline=(current bounding box.center),x=16pt,y=16pt]
\useasboundingbox(-0.51,1)rectangle(0.51,-1);
\node at(0,0){\begin{ytableau}+ \\ +\end{ytableau}};
\end{tikzpicture}
\begin{tikzpicture}[baseline=(current bounding box.center),x=16pt,y=16pt]
\useasboundingbox(0,0)rectangle(0.3,-2.0635);
\draw[fill=black](0,0)rectangle(0.3,-2.0635);
\end{tikzpicture}
\begin{tikzpicture}[baseline=(current bounding box.center),x=16pt,y=16pt]
\useasboundingbox(-1.53,1)rectangle(1.53,-1);
\node at(0,0){\begin{ytableau}+ & + & 0 \\ 0 & + & +\end{ytableau}};
\end{tikzpicture}
\begin{tikzpicture}[baseline=(current bounding box.center),x=16pt,y=16pt]
\useasboundingbox(0,0)rectangle(0.3,-2.0635);
\draw[fill=black](0,0)rectangle(0.3,-2.0635);
\end{tikzpicture}
\begin{tikzpicture}[baseline=(current bounding box.center),x=16pt,y=16pt]
\useasboundingbox(-0.51,1)rectangle(0.51,-1);
\node at(0,0){\begin{ytableau}+ \\ +\end{ytableau}};
\end{tikzpicture} && \quad\begin{tabular}{|ccc|ccc|c}
+&+&0& 0&-&- &- \\
\hline
+&+&+& +&+&0 &0 \\
\end{tabular} && \textnormal{deviant} \\[4pt]
& \textnormal{\bfseries Class 9.} &&
\begin{tikzpicture}[baseline=(current bounding box.center),x=16pt,y=16pt]
\useasboundingbox(-0.51,1)rectangle(0.51,-1);
\node at(0,0){\begin{ytableau}+ \\ +\end{ytableau}};
\end{tikzpicture}
\begin{tikzpicture}[baseline=(current bounding box.center),x=16pt,y=16pt]
\useasboundingbox(0,0)rectangle(0.3,-2.0635);
\draw[fill=black](0,0)rectangle(0.3,-2.0635);
\end{tikzpicture}
\begin{tikzpicture}[baseline=(current bounding box.center),x=16pt,y=16pt]
\useasboundingbox(-1.02,1)rectangle(1.02,-1);
\node at(0,0){\begin{ytableau}+ & + \\ + & + \end{ytableau}};
\end{tikzpicture}
\begin{tikzpicture}[baseline=(current bounding box.center),x=16pt,y=16pt]
\useasboundingbox(0,0)rectangle(0.3,-2.0635);
\draw[fill=black](0,0)rectangle(0.3,-2.0635);
\end{tikzpicture}
\begin{tikzpicture}[baseline=(current bounding box.center),x=16pt,y=16pt]
\useasboundingbox(-0.51,1)rectangle(0.51,-1);
\node at(0,0){\begin{ytableau}+ \\ +\end{ytableau}};
\end{tikzpicture} && \quad\begin{tabular}{|ccc|cc|c}
+&+&0& -&- &- \\
\hline
+&+&+& +&+ &0 \\
\end{tabular} && \textnormal{deviant}
\end{align*}
\caption{The $\oplus$-diagrams and standard basis vectors for the $9$ classes of BCFW cells for $k=2$. The black vertical bars in the $\oplus$-diagrams represent a (possibly empty) block of $0$'s of the appropriate height. Similarly, the vertical bars in the matrices represent a (possibly empty) block of $0$'s.}
\label{fig:cell-classification}
\end{figure}

\begin{thm}\label{thm:4-dominoes}
Each BCFW cell $S$ of $\Gr_{2,n}^{\geq 0}$ falls into one of $9$ classes, shown in \cref{fig:cell-classification}. Any element $V\in S$ can be written as the row span of a $2 \times n$ matrix with rows $d$ and $e$, whose sign patterns are specified precisely in \cref{fig:cell-classification}. In the figure, a vertical line represents a (possibly empty) block of $0$'s. We call $d$ and $e$ the \emph{standard basis vectors} of $V$. Note that $d$ is either orthodox or deviant; we call $V$ either \emph{orthodox} or \emph{deviant}, accordingly.

Moreover, we can write the row vectors $d$ and $e$ in terms of 
linearly independent positive dominoes $d^{(1)}, d^{(2)}, d^{(3)}, d^{(4)}\in\mathbb{R}^n$, such that the following holds (where $i_j$ is the index of $d^{(j)}$):
\begin{itemize}
\item if $V$ is orthodox, then $i_1 + 1 < i_2 \le i_3 < i_4 - 1$, $\overline{d} = d^{(1)} + d^{(2)}$ with $d_n < 0$, and $\overline{e} = d^{(3)} + d^{(4)}$ with $e_n > 0$;
\item if $V$ is deviant, then $i_1 + 1 < i_2 + 1 < i_3 \le i_4$, $\overline{d} = d^{(1)} - d^{(4)}$ with $d_n < 0$, and $e = d^{(1)} + d^{(2)} + d^{(3)}$ (with $e_n = 0$).
\end{itemize}
We call $d^{(1)}, d^{(2)}, d^{(3)}, d^{(4)}$ the \emph{fundamental dominoes} of $V$.
\end{thm}
For example, the matrix 
$\begin{bmatrix}0 & 0 & 3 & 7 & 0 & 1 & 2 & 0 & 0 & 0 & -5 \\ 0 & 0 & 0 & 0 & 0 & 1 & 4 & 1 & 2 & 0 & 2\end{bmatrix}$ represents an orthodox element of $\Gr_{2,11}^{\ge 0}$ in a Class $3$ BCFW cell of \cref{fig:cell-classification}.

\begin{pf}
By \cref{thm:BCFW}, the BCFW cells $\mathcal{C}_{n,2,4}$ of $\Gr_{2,n}^{\geq 0}$ correspond to the $\oplus$-diagrams $\mathcal{D}_{n,2,4}$, and a straightforward case analysis using \cref{def:latticetocell} implies that these are precisely the $\oplus$-diagrams in \cref{fig:cell-classification}. It remains to show that for each such $\oplus$-diagram, an arbitrary element of the corresponding positroid cell can be represented by the $2\times n$ matrix to its right in \cref{fig:cell-classification}, and that we can find dominoes $d^{(1)}, d^{(2)}, d^{(3)}, d^{(4)}$ satisfying the desired properties. We can ignore the black vertical bars in the $\oplus$-diagrams, which correspond exactly to the vertical bars in the matrices. Since the proofs for each of the $9$ classes are similar, we work out the details in two representative cases (Classes $3$ and $8$), and leave the others as an exercise. The idea is to use \Le -moves (\cref{le_moves}) to turn each $\oplus$-diagram into a \Le -diagram, from which we obtain a hook diagram (\cref{def:Le-plabic}) and matrix parameterization of the corresponding cell (\cref{network_param}). From here we can verify the required properties.

{\bfseries Class 3.} In this case the $\oplus$-diagram is a \Le -diagram, with the following hook diagram:
$$
\;\;\begin{tikzpicture}[baseline=(current bounding box.center)]
\pgfmathsetmacro{\scalar}{1.6};
\pgfmathsetmacro{\unit}{\scalar*0.922/1.6};
\coordinate (vstep)at(0,-0.24*\unit);
\coordinate (hstep)at(0.17*\unit,0);
\coordinate (vepsilon)at(0,-0.02*\unit);
\coordinate (hepsilon)at(0.02*\unit,0);
\useasboundingbox(0,0)rectangle(5.25*\unit,-2.5*\unit);
\foreach \x in {1,...,5}{
\draw[thick](\x*\unit-\unit,0)rectangle(\x*\unit,-\unit);}
\foreach \x in {1,...,4}{
\draw[thick](\x*\unit-\unit,-\unit)rectangle(\x*\unit,-2*\unit);}
\node[inner sep=0]at(0.5*\unit,-0.5*\unit){\scalebox{\scalar}{$+$}};
\node[inner sep=0]at(1.5*\unit,-0.5*\unit){\scalebox{\scalar}{$0$}};
\node[inner sep=0]at(2.5*\unit,-0.5*\unit){\scalebox{\scalar}{$+$}};
\node[inner sep=0]at(3.5*\unit,-0.5*\unit){\scalebox{\scalar}{$+$}};
\node[inner sep=0]at(4.5*\unit,-0.5*\unit){\scalebox{\scalar}{$+$}};
\node[inner sep=0]at(0.5*\unit,-1.5*\unit){\scalebox{\scalar}{$+$}};
\node[inner sep=0]at(1.5*\unit,-1.5*\unit){\scalebox{\scalar}{$+$}};
\node[inner sep=0]at(2.5*\unit,-1.5*\unit){\scalebox{\scalar}{$+$}};
\node[inner sep=0]at(3.5*\unit,-1.5*\unit){\scalebox{\scalar}{$+$}};
\node[inner sep=0]at($(5*\unit,-0.5*\unit)+(hstep)$){$1$};
\node[inner sep=0]at($(4.5*\unit,-1*\unit)+(vstep)$){$2$};
\node[inner sep=0]at($(4*\unit,-1.5*\unit)+(hstep)$){$3$};
\node[inner sep=0]at($(3.5*\unit,-2*\unit)+(vstep)$){$4$};
\node[inner sep=0]at($(2.5*\unit,-2*\unit)+(vstep)$){$5$};
\node[inner sep=0]at($(1.5*\unit,-2*\unit)+(vstep)$){$6$};
\node[inner sep=0]at($(0.5*\unit,-2*\unit)+(vstep)$){$7$};
\end{tikzpicture}\qquad\qquad\qquad
\begin{tikzpicture}[baseline=(current bounding box.center)]
\tikzstyle{out1}=[inner sep=0,minimum size=1.2mm,circle,draw=black,fill=black]
\tikzstyle{in1}=[inner sep=0,minimum size=1.2mm,circle,draw=black,fill=white]
\pgfmathsetmacro{\unit}{0.922};
\useasboundingbox(0,0)rectangle(5.25*\unit,-2.5*\unit);
\coordinate (vstep)at(0,-0.24*\unit);
\coordinate (hstep)at(0.17*\unit,0);
\coordinate (vepsilon)at(0,-0.02*\unit);
\coordinate (hepsilon)at(0.02*\unit,0);
\draw[thick](0,0)--(5*\unit,0)--(5*\unit,-1*\unit)--(4*\unit,-1*\unit)--(4*\unit,-2*\unit)--(0,-2*\unit)--(0,0);
\node[inner sep=0]at($(5*\unit,-0.5*\unit)+(hstep)$){$1$};
\node[inner sep=0]at($(4.5*\unit,-1*\unit)+(vstep)$){$2$};
\node[inner sep=0]at($(4*\unit,-1.5*\unit)+(hstep)$){$3$};
\node[inner sep=0]at($(3.5*\unit,-2*\unit)+(vstep)$){$4$};
\node[inner sep=0]at($(2.5*\unit,-2*\unit)+(vstep)$){$5$};
\node[inner sep=0]at($(1.5*\unit,-2*\unit)+(vstep)$){$6$};
\node[inner sep=0]at($(0.5*\unit,-2*\unit)+(vstep)$){$7$};
\node[inner sep=0](b1)at($(5*\unit,-0.5*\unit)+(hepsilon)$){};
\node[inner sep=0](b2)at($(4.5*\unit,-1*\unit)+(vepsilon)$){};
\node[inner sep=0](b3)at($(4*\unit,-1.5*\unit)+(hepsilon)$){};
\node[inner sep=0](b4)at($(3.5*\unit,-2*\unit)+(vepsilon)$){};
\node[inner sep=0](b5)at($(2.5*\unit,-2*\unit)+(vepsilon)$){};
\node[inner sep=0](b6)at($(1.5*\unit,-2*\unit)+(vepsilon)$){};
\node[inner sep=0](b7)at($(0.5*\unit,-2*\unit)+(vepsilon)$){};
\node[out1](i15)at($(5*\unit,-1*\unit)+(-0.5*\unit,0.5*\unit)$){};
\node[out1](i14)at($(4*\unit,-1*\unit)+(-0.5*\unit,0.5*\unit)$){};
\node[out1](i13)at($(3*\unit,-1*\unit)+(-0.5*\unit,0.5*\unit)$){};
\node[out1](i11)at($(1*\unit,-1*\unit)+(-0.5*\unit,0.5*\unit)$){};
\node[out1](i24)at($(4*\unit,-2*\unit)+(-0.5*\unit,0.5*\unit)$){};
\node[out1](i23)at($(3*\unit,-2*\unit)+(-0.5*\unit,0.5*\unit)$){};
\node[out1](i22)at($(2*\unit,-2*\unit)+(-0.5*\unit,0.5*\unit)$){};
\node[out1](i21)at($(1*\unit,-2*\unit)+(-0.5*\unit,0.5*\unit)$){};
\path[-latex',thick](b1)edge node[above=-3pt]{$a_1$}(i15) (i15)edge node[above=-3pt]{$a_2$}(i14) (i14)edge node[above=-3pt]{$a_3$}(i13) (i13)edge node[above=-3pt]{$a_4$}(i11) (b3)edge node[above=-3pt]{$b_1$}(i24) (i24)edge node[above=-3pt]{$b_2$}(i23) (i23)edge node[above=-3pt]{$b_3$}(i22) (i22)edge node[above=-3pt]{$b_4$}(i21) (i15)edge(b2) (i14)edge(i24) (i24)edge(b4) (i13)edge(i23) (i23)edge(b5) (i22)edge(b6) (i11)edge(i21) (i21)edge(b7);
\end{tikzpicture}\;\;.
$$
By \cref{network_param}, an arbitrary element of the corresponding cell of $\Gr_{2,7}^{\ge 0}$ is represented by
$$
\begin{bmatrix}
1 & a_1 & 0 & -a_1 a_2 & -a_1 a_2 (a_3+b_2) & -a_1 a_2 b_3 (a_3+b_2) & -a_1 a_2
(a_3a_4 + b_3b_4(a_3+b_2)) \\ 
0 & 0 & 1 & b_1 & b_1 b_2 & b_1 b_2 b_3 & b_1 b_2 b_3 b_4
\end{bmatrix},
$$
where the $a_i$'s and $b_i$'s are positive real numbers. Let $v,w\in\mathbb{R}^7$ denote the first and second rows of this matrix. If we replace $v$ by $v' := v + \frac{a_1 a_2 (a_3+b_2)}{b_1 b_2}w$, then we obtain a matrix whose sign pattern matches the matrix in \cref{fig:cell-classification}. Therefore any element $V$ in this cell has basis vectors $d := v'$ and $e := w$, and we can write $\overline{d} = d^{(1)}+d^{(2)}$ and $\overline{e} = d^{(3)}+d^{(4)}$ for some positive dominoes $d^{(1)}, d^{(2)}, d^{(3)}, d^{(4)}$. We can check that their indices satisfy the stated inequalities, and they are linearly independent (the fact that $\Delta_{\{3,4\}}\neq 0$ implies that $d^{(2)}$ and $d^{(3)}$ are linearly independent).

{\bfseries Class 8.}  
After applying a \Le -move to the $\oplus$-diagram, we obtain the following \Le -diagram and hook diagram:
$$
\;\;
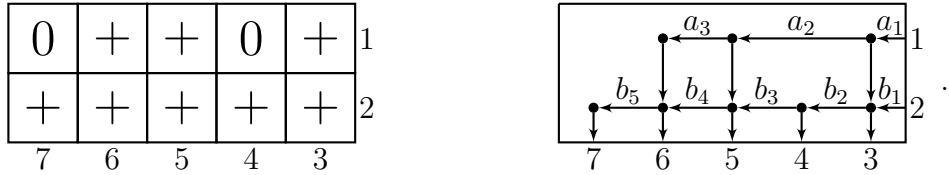
\begin{tikzpicture}[baseline=(current bounding box.center)]
\pgfmathsetmacro{\scalar}{1.6};
\pgfmathsetmacro{\unit}{\scalar*0.922/1.6};
\coordinate (vstep)at(0,-0.24*\unit);
\coordinate (hstep)at(0.17*\unit,0);
\coordinate (vepsilon)at(0,-0.02*\unit);
\coordinate (hepsilon)at(0.02*\unit,0);
\useasboundingbox(0,0)rectangle(5.25*\unit,-2.5*\unit);
\foreach \x in {1,...,5}{
\draw[thick](\x*\unit-\unit,0)rectangle(\x*\unit,-\unit);}
\foreach \x in {1,...,5}{
\draw[thick](\x*\unit-\unit,-\unit)rectangle(\x*\unit,-2*\unit);}
\node[inner sep=0]at(0.5*\unit,-0.5*\unit){\scalebox{\scalar}{$0$}};
\node[inner sep=0]at(1.5*\unit,-0.5*\unit){\scalebox{\scalar}{$+$}};
\node[inner sep=0]at(2.5*\unit,-0.5*\unit){\scalebox{\scalar}{$+$}};
\node[inner sep=0]at(3.5*\unit,-0.5*\unit){\scalebox{\scalar}{$0$}};
\node[inner sep=0]at(4.5*\unit,-0.5*\unit){\scalebox{\scalar}{$+$}};
\node[inner sep=0]at(0.5*\unit,-1.5*\unit){\scalebox{\scalar}{$+$}};
\node[inner sep=0]at(1.5*\unit,-1.5*\unit){\scalebox{\scalar}{$+$}};
\node[inner sep=0]at(2.5*\unit,-1.5*\unit){\scalebox{\scalar}{$+$}};
\node[inner sep=0]at(3.5*\unit,-1.5*\unit){\scalebox{\scalar}{$+$}};
\node[inner sep=0]at(4.5*\unit,-1.5*\unit){\scalebox{\scalar}{$+$}};
\node[inner sep=0]at($(5*\unit,-0.5*\unit)+(hstep)$){$1$};
\node[inner sep=0]at($(5*\unit,-1.5*\unit)+(hstep)$){$2$};
\node[inner sep=0]at($(4.5*\unit,-2*\unit)+(vstep)$){$3$};
\node[inner sep=0]at($(3.5*\unit,-2*\unit)+(vstep)$){$4$};
\node[inner sep=0]at($(2.5*\unit,-2*\unit)+(vstep)$){$5$};
\node[inner sep=0]at($(1.5*\unit,-2*\unit)+(vstep)$){$6$};
\node[inner sep=0]at($(0.5*\unit,-2*\unit)+(vstep)$){$7$};
\end{tikzpicture}\qquad\qquad\qquad
\begin{tikzpicture}[baseline=(current bounding box.center)]
\tikzstyle{out1}=[inner sep=0,minimum size=1.2mm,circle,draw=black,fill=black]
\tikzstyle{in1}=[inner sep=0,minimum size=1.2mm,circle,draw=black,fill=white]
\pgfmathsetmacro{\unit}{0.922};
\useasboundingbox(0,0)rectangle(5.25*\unit,-2.5*\unit);
\coordinate (vstep)at(0,-0.24*\unit);
\coordinate (hstep)at(0.17*\unit,0);
\coordinate (vepsilon)at(0,-0.02*\unit);
\coordinate (hepsilon)at(0.02*\unit,0);
\draw[thick](0,0)--(5*\unit,0)--(5*\unit,-2*\unit)--(0,-2*\unit)--(0,0);
\node[inner sep=0]at($(5*\unit,-0.5*\unit)+(hstep)$){$1$};
\node[inner sep=0]at($(5*\unit,-1.5*\unit)+(hstep)$){$2$};
\node[inner sep=0]at($(4.5*\unit,-2*\unit)+(vstep)$){$3$};
\node[inner sep=0]at($(3.5*\unit,-2*\unit)+(vstep)$){$4$};
\node[inner sep=0]at($(2.5*\unit,-2*\unit)+(vstep)$){$5$};
\node[inner sep=0]at($(1.5*\unit,-2*\unit)+(vstep)$){$6$};
\node[inner sep=0]at($(0.5*\unit,-2*\unit)+(vstep)$){$7$};
\node[inner sep=0](b1)at($(5*\unit,-0.5*\unit)+(hepsilon)$){};
\node[inner sep=0](b2)at($(5*\unit,-1.5*\unit)+(hepsilon)$){};
\node[inner sep=0](b3)at($(4.5*\unit,-2*\unit)+(vepsilon)$){};
\node[inner sep=0](b4)at($(3.5*\unit,-2*\unit)+(vepsilon)$){};
\node[inner sep=0](b5)at($(2.5*\unit,-2*\unit)+(vepsilon)$){};
\node[inner sep=0](b6)at($(1.5*\unit,-2*\unit)+(vepsilon)$){};
\node[inner sep=0](b7)at($(0.5*\unit,-2*\unit)+(vepsilon)$){};
\node[out1](i15)at($(5*\unit,-1*\unit)+(-0.5*\unit,0.5*\unit)$){};
\node[out1](i13)at($(3*\unit,-1*\unit)+(-0.5*\unit,0.5*\unit)$){};
\node[out1](i12)at($(2*\unit,-1*\unit)+(-0.5*\unit,0.5*\unit)$){};
\node[out1](i25)at($(5*\unit,-2*\unit)+(-0.5*\unit,0.5*\unit)$){};
\node[out1](i24)at($(4*\unit,-2*\unit)+(-0.5*\unit,0.5*\unit)$){};
\node[out1](i23)at($(3*\unit,-2*\unit)+(-0.5*\unit,0.5*\unit)$){};
\node[out1](i22)at($(2*\unit,-2*\unit)+(-0.5*\unit,0.5*\unit)$){};
\node[out1](i21)at($(1*\unit,-2*\unit)+(-0.5*\unit,0.5*\unit)$){};
\path[-latex',thick](b1)edge node[above=-3pt]{$a_1$}(i15) (i15)edge node[above=-3pt]{$a_2$}(i13) (i13)edge node[above=-3pt]{$a_3$}(i12)  (b2)edge node[above=-3pt]{$b_1$}(i25) (i25)edge node[above=-3pt]{$b_2$}(i24) (i24)edge node[above=-3pt]{$b_3$}(i23) (i23)edge node[above=-3pt]{$b_4$}(i22)
(i22)edge node[above=-3pt]{$b_5$}(i21)
 (i15)edge(i25) (i13)edge(i23) (i12)edge(i22) (i25)edge(b3) (i24)edge(b4) (i23)edge(b5) (i22)edge(b6) (i21)edge(b7);
\end{tikzpicture}\;\;.
$$
Therefore an arbitrary element of the corresponding cell of $\Gr_{2,7}^{\ge 0}$ is represented by
$$
\scalebox{0.98}{$\begin{bmatrix}
1 & 0 & -a_1 & -a_1 b_2 & -a_1 (a_2+b_2 b_3) & -a_1 (a_2(a_3+b_4)+b_2 b_3 b_4) & 
-a_1 b_5(a_2 (a_3+b_4)+b_2 b_3 b_4) \\
0 & 1 & b_1 & b_1 b_2 & b_1 b_2 b_3 & b_1 b_2 b_3 b_4 & b_1 b_2 b_3 b_4 b_5\end{bmatrix}$},
$$
where the $a_i$'s and $b_i$'s are positive real numbers. Let $v,w\in\mathbb{R}^7$ denote the first and second rows of this matrix. If we replace $v$ by $v' := v + \frac{a_1}{b_1}w$ and then $w$ by $w' := \frac{a_1 a_2 (a_3+b_4)}{b_1 b_2 b_3 b_4}w + v'$, we obtain a matrix whose sign pattern matches the matrix in \cref{fig:cell-classification}. Therefore any element $V$ in this cell has basis vectors $d := v'$ and $e := w'$, and we can write $\overline{d} = d^{(1)}-d^{(4)}$ and $\overline{e} = d^{(1)}+d^{(2)}+d^{(3)}$ for some positive dominoes $d^{(1)}, d^{(2)}, d^{(3)}, d^{(4)}$. We can again check that their indices satisfy the stated inequalities and that they are linearly independent.
\end{pf}

\begin{defn}\label{def:dom}
Let $V\in\Gr_{2,n}^{\ge 0}$ come from a BCFW cell, with fundamental dominoes $d^{(1)}, d^{(2)}, d^{(3)}, d^{(4)}$ from \cref{thm:4-dominoes}. Note that for any $v\in V$, we can write $\overline{v} = \sum_{j=1}^4 \alpha_jd^{(j)}$ for some unique $\alpha_1, \alpha_2, \alpha_3, \alpha_4\in\mathbb{R}$. We let $\Dom_V(v)$ denote the sequence of dominoes $(\alpha_1d^{(1)}, \alpha_2d^{(2)}, \alpha_3d^{(3)}, \alpha_4d^{(4)})$, and $\dom_V(v)\in\{0,+,-\}^4$ the sign vector of $(\alpha_1, \alpha_2, \alpha_3, \alpha_4)$.
\end{defn}

We state two corollaries of \cref{thm:4-dominoes} that we will use in \cref{sec:k=2-disjointness}.
\begin{cor}\label{cor:4-dominoes-signs}
Suppose that $V\in\Gr_{2,n}^{\ge 0}$ comes from a BCFW cell, and $v\in V$ such that $\dom_V(v)$ has no zero components. Then we have the following:
\begin{itemize}
\item if $V$ is orthodox, $\dom_V(v)$ equals $\pm(+,+,+,+)$ or $\pm(+,+,-,-)$;
\item if $V$ is deviant, $\dom_V(v)$ equals $\pm(+,+,+,-)$ or $\pm(-,+,+,+)$ or $\pm(+,+,+,+)$.
\end{itemize}
\end{cor}

\begin{cor}\label{cor:no-support-4}
Suppose that $V\in\Gr_{2,n}^{\ge 0}$ comes from a BCFW cell with standard basis vectors $d,e\in\mathbb{R}^n$ from \cref{thm:4-dominoes}, and $v\in V$ such that $\overline{v}$ has support size at most $4$ (in particular, this includes orthodox and deviant $v$). Then we have the following:
\begin{itemize}
\item if $V$ is orthodox, then $v$ is a scalar multiple of $d$ or $e$;
\item if $V$ is deviant, then $v$ is a scalar multiple of $d$ or $d-e$.
\end{itemize}
\end{cor}

\section{Disjointness for BCFW cells when \texorpdfstring{$k=2$}{k=2}}
\label{sec:k=2-disjointness}

\noindent This section is devoted to the proof of the following result.
\begin{thm}\label{thm:disjointness=k=2}
For $m=4$, $\tilde{Z}$ maps the BCFW cells $\mathcal{C}_{n,2,4}$ of $\Gr_{2,n}^{\ge 0}$ injectively into the amplituhedron $\mathcal{A}_{n,2,4}(Z)$, and their images are pairwise disjoint.
\end{thm}

\subsection{Lemmas on dominoes}
We begin by proving some useful results on dominoes. We already have \cref{domino_sum}, but we will need more powerful tools.

\begin{defn}\label{def:domino-sequence}
Let $D\subseteq\mathbb{R}^n$ be a finite multiset of dominoes, and $v\in\mathbb{R}^n$ the sum of the dominoes in $D$. Given $I = \{i_1 < \dots < i_k\}\subseteq [n]$, an {\itshape $I$-alternating domino sequence for $v$} (with respect to $D$) is a sequence $(d^{(1)}, \dots, d^{(k)})$ of distinct nonzero dominoes in $D$ such that
\begin{itemize}
\item $d^{(j)}_{i_j}$ has the same sign as $v_{i_j}$ for all $j\in [k]$;
\item the {\itshape sign sequence} $(\sign(d^{(1)}), \dots, \sign(d^{(l)}))\in\{+,-\}^l$ alternates in sign.
\end{itemize}
We call $k$ the {\itshape length} of $(d^{(1)}, \dots, d^{(k)})$.
\end{defn}

\begin{lem}\label{lem:obo}
Suppose that $D\subseteq\mathbb{R}^n$ is a finite multiset of dominoes, $v\in\mathbb{R}^n$ is the sum of the dominoes in $D$, and $I\subseteq [n]$. Then $v$ alternates in sign restricted to $I$ if and only if $v$ has an $I$-alternating domino sequence. The indices of dominoes in any $I$-alternating domino sequence for $v$ are weakly increasing, with each index appearing at most twice.
\end{lem}
We think of an $I$-alternating domino sequence $(d^{(1)}, \dots, d^{(k)})$ for $v$ as a witness for the corresponding alternating subsequence of $v$.
\begin{eg}
Let $D$ be the set of dominoes
$d^{(1)} := (1,1,0,0,0)$, 
$d^{(2)} := (0,-3,-1,0,0)$,
$d^{(3)} := (0,1,2,0,0)$,
$d^{(4)} := (0,0,0,-2,-4)$,
and let $v = (1,-1,1,-2,-4)$ be their sum. Note that $v$ alternates in sign restricted to $I := \{1,2,3,5\}$. The unique corresponding $I$-alternating domino sequence is $(d^{(1)}, d^{(2)}, d^{(3)}, d^{(4)})$, with corresponding indices $1,2,2,4$.
\end{eg}

\begin{pf}
Write $I = \{i_1 < \dots < i_k\}$.

($\Rightarrow$): Suppose that $v$ alternates in sign on $I$. Then for each $j\in [k]$, we can find $d^{(j)}\in D$ such that $d^{(j)}_{i_j}$ is nonzero and has the same sign as $v_{i_j}$. Letting $a_j$ be the index of $d^{(j)}$ for $j\in [k]$, we have $a_j\in\{i_j,i_j-1\}$, so $a_j < a_{j'}$ for $j'\ge j+2$. Moreover, $d^{(j)}$ and $d^{(j+1)}$ have opposite signs for $j\in [k-1]$, so $d^{(1)}, \dots, d^{(k)}$ are all distinct. Hence $(d^{(1)}, \dots, d^{(k)})$ is an $I$-alternating domino sequence for $v$.

($\Leftarrow$): Suppose that $(d^{(1)}, \dots, d^{(k)})$ is an $I$-alternating domino sequence for $v$, with corresponding indices $a_1, \dots, a_k$. Then $\sign(v|_I)$ equals the sign sequence of $(d^{(1)}, \dots, d^{(k)})$, which alternates in sign. Moreover, since $a_j\in\{i_j,i_j-1\}$ for all $j\in [k]$, we have $a_1 \le \dots \le a_k$ and $a_j < a_j'$ for $j' \ge j+2$. Since $d^{(j)}$ and $d^{(j+1)}$ have opposite signs for $j\in [k-1]$, this also shows that $d^{(1)}, \dots, d^{(k)}$ are all distinct.
\end{pf}

\begin{lem}\label{lem:delete-from-sum}
Suppose that $(d^{(1)}, \dots, d^{(k)})$ is an $I$-alternating domino sequence for $v\in\mathbb{R}^n$, where $I = \{i_1 < \dots < i_k\}$. Choose $j'\in [k]$ and any subset  $J\subseteq [k]\setminus \{j'\}$. Then the vector $v-\sum_{j\in J}d^{(j)}$ has the same sign in coordinate $i_{j'}$ as the domino $d^{(j')}$.
\end{lem}

\begin{pf}
We know from \cref{def:domino-sequence} that $v_{i_{j'}}$ has the same sign as $d^{(j')}$. It now suffices to observe that for all $j\in J$, $d^{(j)}_{i_{j'}}$ is either zero or has the opposite sign as $d^{(j')}$.
\end{pf}

We will also use the following facts about dominoes and sign variation, which are straightforward to verify.
\begin{lem}\label{lem:add-single}
Suppose that $v\in\mathbb{R}^n$, and $d\in\mathbb{R}^n$ is an $i$-domino. \\
(i) We have $\var(v+d)\le\var(v) + 2$. \\
(ii) If $v_j = 0$ for all $j\le i$ or $v_j = 0$ for all $j > i$, then $\var(v+d)\le\var(v) + 1$.
\end{lem}

\begin{defn}\label{def:shuffle}
A {\itshape shuffle} of two sequences 
 $(s_1, \dots, s_k)$ and $(t_1, \dots, t_l)$ 
is a sequence of length $k+l$ formed by permuting $s_1, \dots, s_k, t_1, \dots, t_l$, such that $s_1, \dots, s_k$ appear in the same relative order, and $t_1, \dots, t_l$ appear in the same relative order.
\end{defn}
For example, the shuffles of $(a,b)$ and $(c,d)$ are precisely $(a,b,c,d)$, $(a,c,b,d)$, $(a,c,d,b)$, $(c,a,b,d)$, $(c,a,d,b)$, and $(c,d,a,b)$.
\begin{lem}\label{lem:no-domino-left-behind}
Suppose that $V,V'\in\Gr_{2,n}^{\ge 0}$ come from BCFW cells, and $v\in V$, $v'\in V$. Note that $\overline{v+v'}$ is the sum of the dominoes in the multiset $D$ of nonzero dominoes appearing in $\Dom_V(v)$ and $\Dom_{V'}(v')$. Then any alternating domino sequence for $\overline{v+v'}$ with respect to $D$ in which all the dominoes in $D$ appear is obtained by shuffling $\Dom_V(v)$ and $\Dom_V(v')$ and deleting all zero dominoes.
\end{lem}

\begin{pf}
Suppose that $(f^{(1)}, \dots, f^{(l)})$ is an $I$-alternating domino sequence for $\overline{v+v'}$ with respect to $D$ which uses all the dominoes in $D$, where $I = \{i_1 < \dots < i_l\}$. Then the indices of $(f^{(1)}, \dots, f^{(l)})$ are weakly increasing, so we must show that if two dominoes in $\Dom_V(v)$ or two dominoes in $\Dom_{V'}(v')$ have the same index and both appear in $(f^{(1)}, \dots, f^{(l)})$, then they appear in the same relative order. It suffices to prove this for $\Dom_V(v)$. We proceed by contradiction and suppose that there exist dominoes $f^{(j)}, f^{(j+1)}\in D$ with the same index (they are necessarily adjacent in the domino sequence, by \cref{lem:obo}), but appear in $\Dom_V(v)$ in the opposite order, i.e.\ $f^{(j+1)}$ before $f^{(j)}$. We will deduce that $\var(v)\ge 2$, which contradicts \cref{thm:G-K}(i).

Perhaps after multiplying all of $v, v', f^{(1)}, \dots, f^{(l)}$ by $-1$, we may assume that $f^{(j)}$ is negative and $f^{(j+1)}$ is positive. Note that $\overline{v}$ equals $\overline{v+v'}$ minus the sum of the dominoes appearing in $\Dom_{V'}(v')$. Since $(f^{(1)}, \dots, f^{(l)})$ are precisely all dominoes in $D$, by \cref{lem:delete-from-sum}, $v_{i_j} < 0$ and $v_{i_{j+1}} > 0$. We now let $d^{(1)}, d^{(2)}, d^{(3)}, d^{(4)}$ be the fundamental dominoes of $V$, and refer to \cref{thm:4-dominoes}. 
First we suppose that $V$ is orthodox. 
Since $f^{(j)}$ and $f^{(j+1)}$ have the same index, $V$ 
must be of Class 3. 
Then $f^{(j+1)}$ is a positive scalar multiple of $d^{(2)}$, and $f^{(j)}$ is a negative scalar multiple of $d^{(3)}$. Hence $\overline{v} = \alpha(d^{(1)} + d^{(2)}) - \beta(d^{(3)} + d^{(4)})$ for some $\alpha,\beta > 0$. 
We see that if $i$ is  the index of $d^{(1)}$, then $i < i_j$ and $v_i > 0$, which gives $\var(v)\ge 2$, as desired. Instead suppose that $V$ is deviant.
Then it must be of Class 4 or 9. Then $f^{(j+1)}$ is a positive scalar multiple of $d^{(3)}$, and $f^{(j)}$ is a negative scalar multiple of $d^{(4)}$. Hence $\overline{v} = \gamma(d^{(1)} - d^{(4)}) + \delta(d^{(1)} + d^{(2)} + d^{(3)})$ for some $\gamma,\delta > 0$. Again $v$ is positive at the index of $d^{(1)}$, giving $\var(v)\ge 2$.
\end{pf}

\begin{cor}\label{cor:no-domino-left-behind}
Suppose that $V,V'\in\Gr_{2,n}^{\ge 0}$ come from BCFW cells, and $v\in V$, $v'\in V$ are distinct such that $Z(v) = Z(v')$ and the sum of the support sizes of $\dom_V(v)$ and $\dom_{V'}(-v')$ is at most $6$, i.e.\ at most $6$ dominoes
appear with nonzero coefficient in 
$\Dom_V(v)$ and $\Dom_{V'}(-v')$. Then the sum of the support sizes is exactly $6$, and some shuffle of $\dom_V(v)$ and $\dom_{V'}(-v')$ alternates in sign (ignoring the zero components).
\end{cor}
This will be useful for us, since if $d,e$ are the standard basis vectors of $V$ from \cref{thm:4-dominoes}, then $\dom_V(d)$ has support size $2$, and also $\dom_V(e)$ has support size $2$ if $V$ is orthodox.
\begin{pf}
By \cref{thm:G-K}(ii), we have $\var(v-v')\ge 6$, and so $\var(\overline{v-v'})\ge 5$. By \cref{lem:obo}, $\overline{v-v'}$ has an alternating domino sequence of length $6$, which necessarily uses all the (nonzero) dominoes in $\Dom_V(v)$ and $\Dom_{V'}(-v')$, since by assumption at most $6$ dominoes appear. By \cref{lem:no-domino-left-behind}, the sequence is obtained by shuffling $\Dom_V(v)$ and $\Dom_{V'}(-v')$ and deleting all zero dominoes. The corresponding sign sequence (which alternates in sign) is therefore obtained by shuffling $\dom_V(v)$ and $\dom_{V'}(-v')$ and deleting all zero components.
\end{pf}

\subsection{Orthodox vs.\ orthodox cells}

\begin{lem}
\label{lem:orthodox-vs-orthodox}
Suppose that $V, V'\in\Gr_{2,n}^{\ge 0}$ are distinct and orthodox. Then $\tilde{Z}(V)\neq\tilde{Z}(V')$.
\end{lem}
\begin{pf}
Suppose otherwise that $\tilde{Z}(V) = \tilde{Z}(V')$. Let $d,e\in\mathbb{R}^n$ be the standard basis vectors of $V$ from \cref{thm:4-dominoes}, which are orthodox, and take matching vectors $v,w\in V'$ for $d,e$, i.e.\ $Z(d) = Z(v)$ and $Z(e) = Z(w)$. We claim that $d = v$. Otherwise, by \cref{cor:no-domino-left-behind}, some shuffle of $\dom_V(d)$ and $\dom_{V'}(-v)$ alternates in sign after deleting its two zero components, and in particular contains exactly three $+$'s and three $-$'s. But $\dom_V(d)$ contains two $+$'s, and by \cref{cor:4-dominoes-signs} $\dom_{V'}(v)$ contains an even number of $+$'s. This contradiction shows $d = v$. Similarly, $e = w$. Since $d$ and $e$ span $V$, we get $V\subseteq V'$, a contradiction.
\end{pf}

\subsection{Deviant vs.\ deviant cells}

\begin{lem}
\label{lem:deviant-vs-deviant}
Suppose that $V,V'\in\Gr_{2,n}^{\ge 0}$ are distinct and deviant. Then $\tilde{Z}(V)\neq\tilde{Z}(V')$.
\end{lem}
\begin{pf}
Suppose otherwise that $\tilde{Z}(V) = \tilde{Z}(V')$. Denote the standard basis vectors and fundamental dominoes of $V$ and $V'$ by $d$, $e$, $d^{(1)}$, $d^{(2)}$, $d^{(3)}$, $d^{(4)}$, and $d'$, $e'$, ${d^{(1)}}'$, ${d^{(2)}}'$, ${d^{(3)}}'$, ${d^{(4)}}'$, respectively. Take matching vectors $v,w\in V'$ for $d,e$. We claim that $d = v$. Otherwise, by \cref{cor:no-domino-left-behind}, some shuffle of $\dom_V(d)$ and $\dom_{V'}(-v)$ alternates in sign after deleting its two zero components, and in particular contains exactly three $+$'s and three $-$'s. But $\dom_V(d)$ contains one $+$, and by \cref{cor:4-dominoes-signs} $\dom_{V'}(v)$ contains either at least three $+$'s or at most one $+$. This contradiction shows $d = v$.

Since $e$ has support size $3$, we cannot apply \cref{cor:no-domino-left-behind} to $e$ and $w$, but we can deduce that $e \neq w$, since otherwise $V = V'$. Hence $\var(e-w)\ge 6$ by \cref{criterion}. Since $d = v$, we may rescale our fundamental dominoes and standard basis vectors so that $d = d'$, $d^{(1)} = {d^{(1)}}'$, and $d^{(4)} = {d^{(4)}}'$. We can write $d = d' = d^{(1)} - d^{(4)} -f$ for some positive $n$-domino $f$, and $w = \alpha d' + \beta e'$ for some $\alpha,\beta\in\mathbb{R}$. Then
$$
e-w = (1-\alpha-\beta)d^{(1)} + d^{(2)} + d^{(3)} + \alpha d^{(4)} - \beta{d^{(2)}}' - \beta{d^{(3)}}' + \alpha f.
$$
Since $\var(e-w)\ge 6$, by \cref{lem:obo} we can arrange the $7$ dominoes summed above into an alternating domino sequence $(f^{(1)}, \dots, f^{(7)})$. The indices of the dominoes weakly increase, so $f^{(1)} = (1-\alpha-\beta)d^{(1)}$ and $f^{(7)} = \alpha f$. Since $f^{(1)}$ and $f^{(7)}$ have the same sign, we get $\sign(1 - \alpha - \beta) = \sign(\alpha)$. Hence the sign sequence of $(f^{(1)}, \dots, f^{(7)})$ is a permutation of
$$
\sign(\alpha), +, +, \sign(\alpha), -\sign(\beta), -\sign(\beta), \sign(\alpha).
$$
This sequence must contain either $3$ or $4$ $+$'s, and therefore $\alpha$ and $\beta$ are negative. But this contradicts $\sign(1-\alpha-\beta) = \sign(\alpha)$.
\end{pf}

\subsection{Orthodox vs.\ deviant cells}

This is the hardest case. Until now, our strategy has been to pick a standard basis vector, look at its matching vector, and obtain a contradiction. It turns out this is not sufficient here, and we have to work harder to prove disjointness.

We start by examining the supports $\supp(v) := \{i\in [n] : v_i\neq 0\}$ of vectors $v\in\mathbb{R}^n$.

\begin{lem}\label{lem:equal-support-range}
Suppose that $V,V'\in\Gr_{2,n}^{\ge 0}$ come from BCFW cells, and $\tilde{Z}(V) = \tilde{Z}(V')$. Then $\min_{v\in V\setminus\{0\}}\supp(v) = \min_{v\in V'\setminus\{0\}}\supp(v)$, $\max_{v\in V\setminus\{0\}}\supp(\overline{v}) = \max_{v\in V'\setminus\{0\}}\supp(\overline{v})$.
\end{lem}

\begin{pf}
Suppose otherwise, so that without loss of generality, we may assume that either $\min_{v\in V\setminus\{0\}}\supp(v) < \min_{v\in V'\setminus\{0\}}\supp(v)$ or $\max_{v\in V\setminus\{0\}}\supp(\overline{v}) > \max_{v\in V'\setminus\{0\}}\supp(\overline{v})$. Let $d$, $e$, $d^{(1)}$, $d^{(2)}$, $d^{(3)}$, $d^{(4)}$ be the standard basis vectors and fundamental dominoes of $V$. First we consider the case $\min_{v\in V\setminus\{0\}}\supp(v) < \min_{v\in V'\setminus\{0\}}\supp(v)$. Write $\overline{d} = d^{(1)} + f$, where $f := d^{(2)}$ if $V$ is orthodox and $f := -d^{(4)}$ if $V$ is deviant, and let $v'\in V'$ be the matching vector for $d$. Note that the index of $d^{(1)}$ is strictly less than $\min(\supp(v'))$. Hence by applying \cref{lem:add-single}(i) to $f$, \cref{lem:add-single}(ii) to $d^{(1)}$, and \cref{thm:G-K}(i) to $v'$, we obtain
$$
\var(\overline{d-v'}) = \var(\overline{-v'}+f+d^{(1)}) \le \var(\overline{v'}) + 2 + 1 \le 1 + 2 + 1 = 4.
$$
But we also have $d\neq v'$, so $\var(\overline{d-v'})\ge \var(d-v') - 1 \ge 5$ by \cref{criterion}, a contradiction. We can treat the case $\max_{v\in V\setminus\{0\}}\supp(\overline{v}) > \max_{v\in V'\setminus\{0\}}\supp(\overline{v})$ by a similar argument, where if $V$ is orthodox we replace $d$ with $e$.
\end{pf}

\begin{lem}
\label{lem:deviant-vs-orthodox-share-vector}
Suppose that $V\in\Gr_{k,n}^{\ge 0}$ is orthodox and $V'\in\Gr_{k,n}^{\ge 0}$ is deviant with $\tilde{Z}(V) = \tilde{Z}(V')$. Then $V\cap V'$ contains no orthodox vectors.
\end{lem}

\begin{pf}
Suppose otherwise that there exists an orthodox $v\in V\cap V'$. Denote the standard basis vectors and fundamental dominoes of $V$ and $V'$ by $d$, $e$, $d^{(1)}$, $d^{(2)}$, $d^{(3)}$, $d^{(4)}$, and $d'$, $e'$, ${d^{(1)}}'$, ${d^{(2)}}'$, ${d^{(3)}}'$, ${d^{(4)}}'$, respectively. By \cref{cor:no-support-4}, $v$ is a scalar multiple of both $e'-d'$ and either $d$ or $e$. After rescaling the vectors appropriately, we may assume that $\overline{v}$ equals both $e'-d'$ and either $d$ or $e$.

Let $w\in V$ be the matching vector for $d'$, and write $-w = \alpha d + \beta e$ for some $\alpha,\beta\in\mathbb{R}$. By \cref{cor:no-support-4} we have $w \neq d'$, so $\var(\overline{d'-w}) \geq 5$ by \cref{criterion}, i.e.\ $\overline{d'-w}$ alternates in sign on some $I = \{i_1 < \dots < i_6\}\subseteq [n-1]$. By \cref{lem:obo} and \cref{lem:no-domino-left-behind}, $\overline{d'-w}$ has an $I$-alternating domino sequence obtained by shuffling $({d^{(1)}}', -{d^{(4)}}')$ and $(\alpha d^{(1)}, \alpha d^{(2)}, \beta d^{(3)}, \beta d^{(4)})$, which we see must equal
$$
(\alpha d^{(1)}, {d^{(1)}}', \alpha d^{(2)}, \beta d^{(3)}, -{d^{(4)}}', \beta d^{(4)}),
$$
with $\beta > 0$. In particular, the index of $d^{(2)}$ is strictly less than that of ${d^{(4)}}'$, which implies $e' - d'\neq d$. Hence $e'-d' = e$. Now let $x := (1-\beta)d' + \beta e'\in V'$, so that $\var(x)\le 1$ by \cref{thm:G-K}(i). We can write $\overline{x} = \overline{d' - w} - \alpha d^{(1)} - \alpha d^{(2)}$, so by \cref{lem:delete-from-sum},
$$
\sign(x|_{\{i_2, i_4, i_5, i_6\}}) = \sign({d^{(1)}}', \beta d^{(3)}, -{d^{(4)}}', \beta d^{(4)}) = (+,+,-,+).
$$
This implies $\var(x)\ge 2$.
\end{pf}

\begin{lem}\label{lem:deviant-vs-orthodox}
Let $V\in\Gr_{2,n}^{\ge 0}$ be orthodox and $V'\in\Gr_{2,n}^{\ge 0}$ be deviant. Then $\tilde{Z}(V)\neq\tilde{Z}(V')$.
\end{lem}

\begin{pf}
Suppose otherwise that $\tilde{Z}(V) = \tilde{Z}(V')$. Denote the standard basis vectors and fundamental dominoes of $V$ and $V'$ by $d$, $e$, $d^{(1)}$, $d^{(2)}$, $d^{(3)}$, $d^{(4)}$, and $d'$, $e'$, ${d^{(1)}}'$, ${d^{(2)}}'$, ${d^{(3)}}'$, ${d^{(4)}}'$, respectively. Let $v, w\in V'$ be the matching vectors for $d,e$, whence $v\neq d$ and $w\neq e$ by \cref{lem:deviant-vs-orthodox-share-vector}. Hence $\var(\overline{d-v}), \var(\overline{e-w})\ge 5$ by \cref{criterion}, so we can take $I,J\in\binom{[n]}{6}$ such that $\overline{d-v}$ alternates on $I$ and $\overline{e-w}$ alternates on $J$.

Let us write $-v = \alpha d' + \beta e'$ and $-w = \gamma d' + \delta e'$ for some $\alpha,\beta,\gamma,\delta\in\mathbb{R}$. By \cref{lem:obo} and \cref{lem:no-domino-left-behind}, $\overline{d-v}$ has an $I$-alternating domino sequence obtained by shuffling $(d^{(1)}, d^{(2)})$ and $((\alpha + \beta){d^{(1)}}', \beta{d^{(2)}}', \beta{d^{(3)}}', -\alpha{d^{(4)}}')$. By \cref{lem:equal-support-range}, $d^{(1)}$ and ${d^{(1)}}'$ have the same support, say $\{i, i+1\}$. Hence this shuffle equals
$$
((\alpha + \beta){d^{(1)}}', d^{(1)}, \beta{d^{(2)}}', d^{(2)}, \beta{d^{(3)}}', -\alpha{d^{(4)}}')
$$
with $\alpha < 0$ and $\beta < 0$. Since $(d-v)_i$ has the same sign as $(\alpha + \beta){d^{(1)}}'$, we have $(d-v)_i < 0$. Similarly, $\overline{e-w}$ has a $J$-alternating domino sequence obtained by shuffling $(d^{(3)}, d^{(4)})$ and $((\gamma + \delta){d^{(1)}}', \delta{d^{(2)}}', \delta{d^{(3)}}', -\gamma{d^{(4)}}')$. By \cref{lem:equal-support-range}, $d^{(4)}$ and ${d^{(4)}}'$ have the same support, so this shuffle equals
$$
((\gamma + \delta){d^{(1)}}', \delta{d^{(2)}}', d^{(3)}, \delta{d^{(3)}}', d^{(4)}, -\gamma{d^{(4)}}'),
$$
with $\gamma > 0$, $\delta < 0$, and $\gamma + \delta > 0$. Since $(\gamma + \delta){d^{(1)}}'$ is the only domino above whose support contains $i$, we have $(e-w)_i > 0$.

Now let
\begin{align}\label{x-orthodox-deviant-proof}
x := \delta(d-v) - \beta(e-w) = \delta d - \beta e + (\alpha\delta - \beta\gamma)d',
\end{align}
so that $x\in\ker(Z)$. Note that $x\neq 0$, since otherwise $d'\in V$, contradicting \cref{cor:no-support-4}. Hence $\var(\overline{x})\ge 5$ by \cref{thm:G-K}(ii), i.e. $x$ alternates in sign on $K$ for some $K\in\binom{[n]}{6}$. By \cref{lem:obo} and \cref{lem:no-domino-left-behind}, $\overline{x}$ has a $K$-alternating domino sequence obtained by shuffling $(\delta d^{(1)}, \delta d^{(2)}, -\beta d^{(3)}, -\beta d^{(4)})$ and $((\alpha\delta - \beta\gamma){d^{(1)}}', -(\alpha\delta - \beta\gamma){d^{(4)}}')$, which necessarily equals
$$
(\delta d^{(1)}, (\alpha\delta - \beta\gamma){d^{(1)}}', \delta d^{(2)}, -\beta d^{(3)}, -(\alpha\delta - \beta\gamma){d^{(4)}}', -\beta d^{(4)}).
$$
Hence $x_i$ has the same sign as $\delta d^{(1)}$, so $x_i < 0$. But $\beta < 0$ , $\delta < 0$, $(d-v)_i < 0$, and $(e-w)_i > 0$, so we see from \eqref{x-orthodox-deviant-proof} that $x_i > 0$.
\end{pf}

We can now deduce \cref{thm:disjointness=k=2} from \cref{thm:4-dominoes}, \cref{lem:orthodox-vs-orthodox}, \cref{lem:deviant-vs-deviant}, and \cref{lem:deviant-vs-orthodox}.

\section{A non-triangulation for \texorpdfstring{$m=3$}{m=3}}\label{sec:m=3}

\noindent In this section we show that the $m=3$ amplituhedron is {\itshape not} triangulated by a seemingly natural collection of cells coming from the BCFW cells for $m=4$. As background, we recall from \cref{sec:m2} (in particular, see \cref{m=2_to_m=1}) that we have two sets of \Le -diagrams $\mathcal{D}_{n,k,2}$ and $\mathcal{D}_{n,k,1}$, which give triangulations of $\mathcal{A}_{n,k,2}(Z)$ and $\mathcal{A}_{n,k,1}(Z)$, respectively. Moreover, we have a bijection $\mathcal{D}_{n+1,k,2} \to \mathcal{D}_{n,k,1}$, which takes a \Le -diagram $D$ and deletes its leftmost column.

Analogously, let $\mathcal{D}_{n,k,3}$ be the set of \Le -diagrams formed from $\oplus$-diagrams in $\mathcal{D}_{n+1,k,4}$ by deleting the leftmost column, so that we have a bijection $\mathcal{D}_{n+1,k,4}\to\mathcal{D}_{n,k,3}$ given by deleting the leftmost column. (The fact that this is a well-defined bijection follows from the proof of \cref{lem:reduced}.) To our surprise, we found that the images under $\tilde{Z}$ of the cells of $\Gr_{k,n}^{\ge 0}$ corresponding to $\mathcal{D}_{n,k,3}$ do not triangulate the $m=3$ amplituhedron $\mathcal{A}_{n,k,3}(Z)$, since their images are not mutually disjoint.

For example, consider the $\oplus$-diagrams
$$
D_1 := \begin{ytableau}
+ & + & + & 0 & 0 & 0 & 0 & 0 & 0 & 0 & + \\
+ & 0 & 0 & + & + & 0 & 0 & +
\end{ytableau}\;,\quad\;\;
D_2 := \begin{ytableau}
+ & + & + & 0 & 0 & 0 & 0 & 0 & 0 & 0 & + \\
+ & 0 & 0 & 0 & 0 & + & + & 0 & 0 & +
\end{ytableau}\;
$$
in $\mathcal{D}_{13,2,4}$, which are both Class 6 $\oplus$-diagrams from \cref{fig:cell-classification}. Let $D_1'$ and $D_2'$ be the \Le -diagrams of type $(2,12)$ formed from $D_1$ and $D_2$ by deleting the leftmost column. Given $Z\in\Mat_{5,12}^{>0}$, by the claim in \cite[Lemma 4.1]{karp}, the sign patterns of the nonzero vectors in $\ker(Z)$ are precisely those with at least $5$ sign changes. In particular $\ker(Z)$ contains a vector $v\in\mathbb{R}^{12}$ with $\sign(v) = (+, +, -, -, +, +, -, -, +, +, -, -)$.

Now let $V_1, V_2\in\Gr_{2,12}^{\ge 0}$ be represented by the matrices
$$
\scalebox{0.86}{$\begin{bmatrix}
v_1 & v_2 & 0 & 0 & 0 & 0 & 0 & 0 & 0 & 0 & v_{11} & v_{12} \\
v_1 & v_2 & 0 & 0 & v_5 & v_6 & 0 & 0 & v_9 & v_{10} & v_{11} & v_{12}
\end{bmatrix}$}, \quad
\scalebox{0.86}{$\begin{bmatrix}
v_1 & v_2 & 0 & 0 & 0 & 0 & 0 & 0 & 0 & 0 & v_{11} & v_{12} \\
0 & 0 & -v_3 & -v_4 & 0 & 0 & -v_7 & -v_8 & 0 & 0 & 0 & 0
\end{bmatrix}$},
$$
respectively. We can check using the network parameterizations coming from the \Le -diagrams (\cref{network_param}) that $V_1\in S_{D_1'}$, $V_2\in S_{D_2'}$. The difference of the two matrices above is
$$
\begin{bmatrix}
0 & 0 & 0 & 0 & 0 & 0 & 0 & 0 & 0 & 0 & 0 & 0 \\
v_1 & v_2 & v_3 & v_4 & v_5 & v_6 & v_7 & v_8 & v_9 & v_{10} & v_{11} & v_{12}
\end{bmatrix},
$$
whose rows are both in $\ker(Z)$. Hence $\tilde{Z}(V_1) = \tilde{Z}(V_2)$, showing that the images of $S_{D_1'}$ and $S_{D_2'}$ in $\mathcal{A}_{12,2,3}(Z)$ intersect.

\begin{prob}\label{prob:m3}
Can we find $3k$-dimensional cells of $\Gr_{k,n}^{\ge 0}$, naturally in bijection with the BCFW cells $\mathcal{C}_{n+1,k,4}$, whose images under $\tilde{Z}$ `triangulate' the $m=3$ amplituhedron $\mathcal{A}_{n,k,3}(Z)$?
\end{prob}

\section*{Appendix. Dyck paths and BCFW domino bases (with Hugh Thomas\except{toc}{\protect\footnotemark})}
\footnotetext{H.T.\ is supported by an NSERC Discovery Grant and the Canada Research Chairs program. He gratefully acknowledges the hospitality of the Fields Institute.  His work on this project was in part supported by the Munich Institute for Astro- and Particle Physics (MIAPP) of the DFG cluster of excellence ``Origin and Structure of the Universe''.}
\setcounter{section}{1}
\setcounter{equation}{0}
\renewcommand{\thesection}{\Alph{section}}
\tikzexternalenable

\begin{defn}\label{defn:dyck}
A {\itshape Dyck path} $P$ is a path in the plane from $(0,0)$ to $(2n,0)$ for some $n\ge 0$, formed by $n$ up steps $(1,1)$ and $n$ down steps $(1,-1)$, which never passes below the $x$-axis. A local maximum of $P$ is called a {\itshape peak}. Let $\mathcal{P}_{n,k,4}$ denote the set of Dyck paths with $2(n-3)$ steps and precisely $n-3-k$ peaks. For example, the Dyck path $P$ shown in \cref{fig:dyck_bijection} is in $\mathcal{P}_{12,3,4}$.
\end{defn}

The cardinality of $\mathcal{P}_{n,k,4}$ 
equals the Narayana number $N_{n-3,k+1}$ \cite[A46]{stanley_catalan}, the number of $(k,n)$-BCFW cells. We will give a bijection $\mathcal{L}_{n,k,4}\longleftrightarrow\mathcal{P}_{n,k,4}$, which thereby allows us to label the $(k,n)$-BCFW cells by the Dyck paths $\mathcal{P}_{n,k,4}$. We then provide a way to conjecturally obtain $k$ basis vectors for any element of a BCFW cell from its Dyck path.

\subsection{From pairs of lattice paths to Dyck paths}

\begin{figure}[ht]
$$
\begin{tikzpicture}[baseline=(current bounding box.center)]
\tikzstyle{hu}=[top color=red!10,bottom color=red!10,middle color=red,opacity=0.70]
\tikzstyle{vu}=[left color=red!10,right color=red!10,middle color=red,opacity=0.70]
\tikzstyle{hl}=[top color=blue!10,bottom color=blue!10,middle color=blue,opacity=0.55]
\tikzstyle{vl}=[left color=blue!10,right color=blue!10,middle color=blue,opacity=0.55]
\pgfmathsetmacro{\u}{0.62};
\pgfmathsetmacro{\w}{0.12*\u};
\coordinate(h)at(-\u,0);
\coordinate(v)at(0,-\u);
\draw[step=\u,color=black!16,ultra thick](0,0)grid(5*\u,3*\u);
\node[inner sep=0](l1)at(5*\u,3*\u){};
\node[inner sep=0](l2)at($(l1)+(v)$){};
\node[inner sep=0](l3)at($(l2)+(v)$){};
\node[inner sep=0](l4)at($(l3)+(h)$){};
\node[inner sep=0](l5)at($(l4)+(v)$){};
\node[inner sep=0](l6)at($(l5)+(h)$){};
\node[inner sep=0](l7)at($(l6)+(h)$){};
\node[inner sep=0](l8)at($(l7)+(h)$){};
\node[inner sep=0](l9)at($(l8)+(h)$){};
\node[inner sep=0](u1)at(5*\u,3*\u){};
\node[inner sep=0](u2)at($(u1)+(h)$){};
\node[inner sep=0](u3)at($(u2)+(h)$){};
\node[inner sep=0](u4)at($(u3)+(v)$){};
\node[inner sep=0](u5)at($(u4)+(h)$){};
\node[inner sep=0](u6)at($(u5)+(h)$){};
\node[inner sep=0](u7)at($(u6)+(v)$){};
\node[inner sep=0](u8)at($(u7)+(v)$){};
\node[inner sep=0](u9)at($(u8)+(h)$){};
\begin{scope}
\clip($(u1.center)+(0,\w)$)--++($2*(h)+(-\w,0)$)--++(2*\w,-2*\w)--++($-2*(h)+(-\w,0)$);
\path[hu]($(u1.center)+(\w,\w)$)rectangle($(u3.center)+(-\w,-\w)$);
\end{scope}
\begin{scope}
\clip($(u3.center)+(-\w,\w)$)--++($1*(v)$)--++(2*\w,-2*\w)--++($-1*(v)$);
\path[vu]($(u3.center)+(\w,\w)$)rectangle($(u4.center)+(-\w,-\w)$);
\end{scope}
\begin{scope}
\clip($(u4.center)+(-\w,\w)$)--++($2*(h)$)--++(2*\w,-2*\w)--++($-2*(h)$);
\path[hu]($(u4.center)+(\w,\w)$)rectangle($(u6.center)+(-\w,-\w)$);
\end{scope}
\begin{scope}
\clip($(u6.center)+(-\w,\w)$)--++($2*(v)$)--++(2*\w,-2*\w)--++($-2*(v)$);
\path[vu]($(u6.center)+(\w,\w)$)rectangle($(u8.center)+(-\w,-\w)$);
\end{scope}
\begin{scope}
\clip($(u8.center)+(-\w,\w)$)--++($1*(h)+(\w,0)$)--++(0,-2*\w)--++($-1*(h)+(\w,0)$);
\path[hu]($(u8.center)+(\w,\w)$)rectangle($(u9.center)+(-\w,-\w)$);
\end{scope}
\begin{scope}
\clip($(l1.center)+(-\w,0)$)--++($2*(v)+(0,\w)$)--++(2*\w,-2*\w)--++($-2*(v)+(0,\w)$);
\path[vl]($(l1.center)+(\w,\w)$)rectangle($(l3.center)+(-\w,-\w)$);
\end{scope}
\begin{scope}
\clip($(l3.center)+(-\w,\w)$)--++($1*(h)$)--++(2*\w,-2*\w)--++($-1*(h)$);
\path[hl]($(l3.center)+(\w,\w)$)rectangle($(l4.center)+(-\w,-\w)$);
\end{scope}
\begin{scope}
\clip($(l4.center)+(-\w,\w)$)--++($1*(v)$)--++(2*\w,-2*\w)--++($-1*(v)$);
\path[vl]($(l4.center)+(\w,\w)$)rectangle($(l5.center)+(-\w,-\w)$);
\end{scope}
\begin{scope}
\clip($(l5.center)+(-\w,\w)$)--++($4*(h)+(\w,0)$)--++(0,-2*\w)--++($-4*(h)+(\w,0)$);
\path[hl]($(l5.center)+(\w,\w)$)rectangle($(l9.center)+(-\w,-\w)$);
\end{scope}
\end{tikzpicture}\quad\longleftrightarrow\quad
\begin{tikzpicture}[baseline=(current bounding box.center)]
\pgfmathsetmacro{\u}{0.60};
\coordinate(u)at(\u,\u);
\coordinate(d)at(\u,-\u);
\tikzstyle{out1}=[inner sep=0,minimum size=1.2mm,circle,draw=black,fill=black]
\coordinate(p1)at(0,0){};
\coordinate(p2)at($(p1)+(u)$){};
\coordinate(p3)at($(p2)+(u)$){};
\coordinate(p4)at($(p3)+(u)$){};
\coordinate(p5)at($(p4)+(d)$){};
\coordinate(p6)at($(p5)+(d)$){};
\coordinate(p7)at($(p6)+(u)$){};
\coordinate(p8)at($(p7)+(u)$){};
\coordinate(p9)at($(p8)+(d)$){};
\coordinate(p10)at($(p9)+(u)$){};
\coordinate(p11)at($(p10)+(d)$){};
\coordinate(p12)at($(p11)+(d)$){};
\coordinate(p13)at($(p12)+(u)$){};
\coordinate(p14)at($(p13)+(d)$){};
\coordinate(p15)at($(p14)+(d)$){};
\coordinate(p16)at($(p15)+(u)$){};
\coordinate(p17)at($(p16)+(d)$){};
\coordinate(p18)at($(p17)+(u)$){};
\coordinate(p19)at($(p18)+(d)$){};
\draw[thick](p1.center)--(p2.center)--(p3.center)--(p4.center)--(p5.center)--(p6.center)--(p7.center)--(p8.center)--(p9.center)--(p10.center)--(p11.center)--(p12.center)--(p13.center)--(p14.center)--(p15.center)--(p16.center)--(p17.center)--(p18.center)--(p19.center);
\foreach \x in {1,...,19}{
\node[out1]at(p\x){};}
\end{tikzpicture}
$$
\caption{The pair of noncrossing lattice paths in $\mathcal{L}_{12,3,4}$ from \cref{fig:PathsToLeDiagram}, and its corresponding Dyck path in $\mathcal{P}_{12,3,4}$.}
\label{fig:dyck_bijection}
\end{figure}
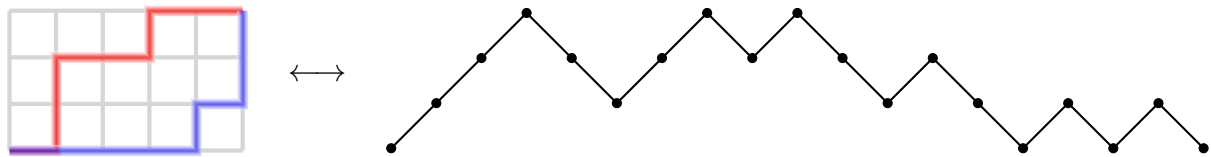

\begin{defn}\label{lattice_to_dyck}
To any pair of noncrossing lattice paths $(W_U, W_L) \in \mathcal{L}_{n,k,4}$, we associate a Dyck path $P(W_U, W_L)\in \mathcal{P}_{n,k,4}$ by the following recursive definition. (We use $+$ and $-$ to denote up and down steps of a Dyck path, and $\cdot$ to denote concatenation of paths.)
\begin{itemize}
\item If $W_U$ and $W_L$ are the trivial paths of length zero, then $P(W_U, W_L) := +-$.
\item If $W_U$ and $W_L$ both begin with a vertical step, then we can write $W_U = V\cdot W_U'$, $W_L = V\cdot W_L'$. We set
$$
P(W_U, W_L) := +\cdot P(W_U', W_L')\cdot -.
$$
\item Otherwise, let $(W_U'', W_L'')$ be the final portion of $(W_U, W_L)$ starting at the first overlapping vertical step of $W_U$ and $W_L$. (If $W_U$ and $W_L$ have no overlapping vertical steps, then we let $W_U''$ and $W_L''$ be the trivial paths.) Then we can write $W_U = H\cdot W_U'\cdot W_U''$, $W_L = W_L'\cdot H\cdot W_L''$. We set
$$
P(W_U, W_L) := P(W_U', W_L') \cdot P(W_U'', W_L'').
$$
\end{itemize}
For example, see \cref{fig:dyck_bijection}. This gives us a map $\Omega_{\mathcal{L}\mathcal{P}}:\mathcal{L}_{n,k,4}\to\mathcal{P}_{n,k,4}$, which sends $(W_U, W_L)$ to $P(W_U, W_L)$.
\end{defn}

In order to show that $\Omega_{\mathcal{L}\mathcal{P}}$ is a bijection, we will define its inverse, which has an elegant description. As far as we know, these inverse bijections have not appeared in the literature.
\begin{defn}\label{dyck_to_lattice}
Let $P$ be a Dyck path. For any point $p$ on $P$, we let $\shadow_p(P)$ be the waterline when we turn $P$ upside-down and fill it with as much water as possible without submerging $p$; this  is a piecewise linear curve between the endpoints of $P$. Let $\touch_p(P)$ be the number of down steps of $P$ whose right endpoint lies on $\shadow_p(P)$ to the right of $p$, but otherwise does not intersect $\shadow_p(P)$. Then we associate a pair of noncrossing lattice paths $(W_U(P), W_L(P))$ to $P$ as follows.
\begin{itemize}
\item $W_L(P)$ is obtained by reading the up steps of $P$ from left to right, recording $V$ for every up step followed by an up step and $H$ for every up step followed by a down step (except for the final up step).
\item Let $p_1, \dots, p_k$ be the left endpoints of the up steps of $P$ which are followed by an up step, ordered from left to right. We let $W_U(P)$ be such that the distance between the vertical edges of $W_U(P)$ and $W_L(P)$ in row $i$ equals $\touch_{p_i}(P)-1$, for $1 \le i \le k$.
\end{itemize}
Note that this gives a map $\Omega_{\mathcal{P}\mathcal{L}}:\mathcal{P}_{n,k,4}\to\mathcal{L}_{n,k,4}$, which sends $P$ to $(W_U(P), W_L(P))$.
\end{defn}

\begin{eg}
Let $P\in\mathcal{P}_{12,3,4}$ be the Dyck path in \cref{fig:dyck_bijection}. Reading the up steps of $P$ from left to right (and ignoring the final up step), we get $W_L(P) = VVHVHHHH$. To find $W_U(P)$, we let $p_1, p_2, p_3$ denote the left endpoints of the three up steps of $P$ which are followed by an up step. (See \cref{fig:shadows}.) Then $\shadow_{p_1}(P)$ is the portion of the $x$-axis between the endpoints of $P$, and $\touch_{p_1}(P) = 3$. We have $\shadow_{p_2}(P) = \shadow_{p_3}(P)$, and we can calculate that $\touch_{p_2}(P) = 5$, $\touch_{p_3}(P) = 4$. (Note that $\touch_{p_2}(P) = 1+\touch_{p_3}(P)$, since $p_3$ contributes to $\touch_{p_2}(P)$ but not to $\touch_{p_3}(P)$. Also, we point out that the first down step of $P$ which ends on the $x$-axis does not contribute to $\touch_{p_2}(P)$ or $\touch_{p_3}(P)$, since the entire step lies on the shadow.) Therefore $W_U(P)$ is such that the distances between $W_U(P)$ and $W_L(P)$ in rows $1$, $2$, and $3$ are, respectively, $3-1$, $5-1$, and $4-1$. This agrees with \cref{fig:dyck_bijection}.
\end{eg}
\begin{figure}[ht]
\begin{center}
\begin{tabular}{cl}
\begin{tikzpicture}[baseline=(current bounding box.center)]
\pgfmathsetmacro{\u}{0.60};
\pgfmathsetmacro{\r}{0.18};
\pgfmathsetmacro{\vstep}{0.42};
\coordinate(u)at(\u,\u);
\coordinate(d)at(\u,-\u);
\tikzstyle{out1}=[inner sep=0,minimum size=1.2mm,circle,draw=black,fill=black]
\coordinate(p1)at(0,0){};
\coordinate(p2)at($(p1)+(u)$){};
\coordinate(p3)at($(p2)+(u)$){};
\coordinate(p4)at($(p3)+(u)$){};
\coordinate(p5)at($(p4)+(d)$){};
\coordinate(p6)at($(p5)+(d)$){};
\coordinate(p7)at($(p6)+(u)$){};
\coordinate(p8)at($(p7)+(u)$){};
\coordinate(p9)at($(p8)+(d)$){};
\coordinate(p10)at($(p9)+(u)$){};
\coordinate(p11)at($(p10)+(d)$){};
\coordinate(p12)at($(p11)+(d)$){};
\coordinate(p13)at($(p12)+(u)$){};
\coordinate(p14)at($(p13)+(d)$){};
\coordinate(p15)at($(p14)+(d)$){};
\coordinate(p16)at($(p15)+(u)$){};
\coordinate(p17)at($(p16)+(d)$){};
\coordinate(p18)at($(p17)+(u)$){};
\coordinate(p19)at($(p18)+(d)$){};
\fill[blue!10](p1)--(p15)--(p13)--(p12)--(p10)--(p9)--(p8)--(p6)--(p4)--cycle (p15)--(p17)--(p16)--cycle (p17)--(p19)--(p18)--cycle;
\draw[thick,dashed](p1)--(p19);
\draw[thick](p1.center)--(p2.center)--(p3.center)--(p4.center)--(p5.center)--(p6.center)--(p7.center)--(p8.center)--(p9.center)--(p10.center)--(p11.center)--(p12.center)--(p13.center)--(p14.center)--(p15.center)--(p16.center)--(p17.center)--(p18.center)--(p19.center);
\foreach \x in {1,...,19}{
\node[out1]at(p\x){};}
\node[inner sep=0]at($(p1)+(0,-\vstep)$){$p_1$};
\node[inner sep=0]at(9*\u,-1*\u){$\shadow_{p_1}(P)$};
\draw[thick](p1)circle(\r);
\draw[thick,red](p15)circle(\r) (p17)circle(\r) (p19)circle(\r);
\end{tikzpicture} & $\begin{array}{l}\touch_{p_1}(P) = 3 \\ ~\end{array}$ \\
~ & ~ \\
~ & ~ \\
\begin{tikzpicture}[baseline=(current bounding box.center)]
\pgfmathsetmacro{\u}{0.60};
\pgfmathsetmacro{\r}{0.18};
\pgfmathsetmacro{\vstep}{0.42};
\coordinate(u)at(\u,\u);
\coordinate(d)at(\u,-\u);
\tikzstyle{out1}=[inner sep=0,minimum size=1.2mm,circle,draw=black,fill=black]
\coordinate(p1)at(0,0){};
\coordinate(p2)at($(p1)+(u)$){};
\coordinate(p3)at($(p2)+(u)$){};
\coordinate(p4)at($(p3)+(u)$){};
\coordinate(p5)at($(p4)+(d)$){};
\coordinate(p6)at($(p5)+(d)$){};
\coordinate(p7)at($(p6)+(u)$){};
\coordinate(p8)at($(p7)+(u)$){};
\coordinate(p9)at($(p8)+(d)$){};
\coordinate(p10)at($(p9)+(u)$){};
\coordinate(p11)at($(p10)+(d)$){};
\coordinate(p12)at($(p11)+(d)$){};
\coordinate(p13)at($(p12)+(u)$){};
\coordinate(p14)at($(p13)+(d)$){};
\coordinate(p15)at($(p14)+(d)$){};
\coordinate(p16)at($(p15)+(u)$){};
\coordinate(p17)at($(p16)+(d)$){};
\coordinate(p18)at($(p17)+(u)$){};
\coordinate(p19)at($(p18)+(d)$){};
\fill[blue!10](p2)--(p6)--(p4)--cycle (p6)--(p12)--(p10)--(p9)--(p8)--cycle (p12)--(p14)--(p13)--cycle (p15)--(p17)--(p16)--cycle (p17)--(p19)--(p18)--cycle;
\draw[thick,dashed](p2)--(p14) (p15)--(p19);
\draw[thick](p1.center)--(p2.center)--(p3.center)--(p4.center)--(p5.center)--(p6.center)--(p7.center)--(p8.center)--(p9.center)--(p10.center)--(p11.center)--(p12.center)--(p13.center)--(p14.center)--(p15.center)--(p16.center)--(p17.center)--(p18.center)--(p19.center);
\foreach \x in {1,...,19}{
\node[out1]at(p\x){};}
\node[inner sep=0]at($(p2)+(0,-\vstep)$){$p_2$};
\node[inner sep=0]at($(p6)+(0,-\vstep)$){$p_3$};
\node[inner sep=0]at(9*\u,-1*\u){$\shadow_{p_2}(P) = \shadow_{p_3}(P)$};
\draw[thick](p2)circle(\r);
\draw[thick,red](p6)circle(\r) (p12)circle(\r) (p14)circle(\r) (p17)circle(\r) (p19)circle(\r);
\end{tikzpicture} & $\begin{array}{l}\touch_{p_2}(P) = 5 \\[2pt] \touch_{p_3}(P) = 4 \\ ~\end{array}$
\end{tabular}
\end{center}
\caption{Calculating $\shadow_{p_i}(P)$ and $\touch_{p_i}(P)$ for the Dyck path $P$ from \cref{fig:dyck_bijection}.}
\label{fig:shadows}
\end{figure}

\begin{prop}\label{prop:dyck_bijection}
The map $\Omega_{\mathcal{L}\mathcal{P}}:\mathcal{L}_{n,k,4}\to\mathcal{P}_{n,k,4}$ is a bijection with inverse $\Omega_{\mathcal{P}\mathcal{L}}$.
\end{prop}

\begin{pf}
We can show that $\Omega_{\mathcal{P}\mathcal{L}}\circ\Omega_{\mathcal{L}\mathcal{P}}$ is the identity on $\mathcal{L}_{n,k,4}$ by induction, using the recursive nature of \cref{lattice_to_dyck}. The result then follows from the fact that $|\mathcal{L}_{n,k,4}| = N_{n-3,k+1} = |\mathcal{P}_{n,k,4}|$. (Alternatively, we can verify that $\Omega_{\mathcal{L}\mathcal{P}}\circ\Omega_{\mathcal{P}\mathcal{L}}$ is the identity on $\mathcal{P}_{n,k,4}$.)
\end{pf}

By \cref{thm:BCFW}, the $(k,n)$-BCFW cells are labeled by the $\oplus$-diagrams $\mathcal{D}_{n,k,4}$. Therefore the bijection $\Omega_{\mathcal{L}\mathcal{P}}\circ\Omega_{\mathcal{L}\mathcal{D}}^{-1} : \mathcal{D}_{n,k,4}\to\mathcal{P}_{n,k,4}$ allows us to label the $(k,n)$-BCFW cells by $\mathcal{P}_{n,k,4}$.

\subsection{From Dyck paths to domino bases}

Recall the definition of an $i$-domino from \cref{def:domino}.

\begin{defn}\label{dyck_to_basis}
Given $P\in\mathcal{P}_{n,k,4}$, label the up steps of $P$ by $1, \dots, n-3$ from left to right, and label each down step so that it has the same label as the next up step. (We label the final down step by $n-2$.) We match each up step of $P$ to the next down step at the same height. Let $\up_1 < \dots < \up_k$ be the labels of the up steps of $P$ which are followed by an up step, and for $i\in [k]$ let $\down_i$ be the label of the matching down step of the up step $\up_i$.

Now for $i= 1, \dots, k$, we define the following dominoes in $\mathbb{R}^n$. Let $d^{(i)}$ be a positive $\up_i$-domino, and $e^{(i)}$ a $\down_i$-domino which has sign $(-1)^{|\{j\in [k] : \up_i < \up_j < \down_i\}|}$. If the up step $\up_i$ of $P$ begins on the $x$-axis, then we let $f^{(i)}$ be the $n$-domino $(0,\dots, 0, (-1)^{k-i})$. Otherwise, take $i' < i$ so that $\up_{i'}$ labels the last up step of $P$ before $\up_i$ which finishes at the same height that $\up_i$ begins, and let $f^{(i)} := (-1)^{i-i'-1}d^{(i')}$. We set
$$
v^{(i)} := d^{(i)} + e^{(i)} + f^{(i)}\in\mathbb{R}^n.
$$
We call any $k$-tuple of linearly independent vectors $(v^{(1)}, \dots, v^{(k)})$ which we can obtain in this way a {\itshape $P$-domino basis}.
\end{defn}

\begin{conj}\label{conj:domino}
Let $S\subseteq\Gr_{k,n}^{\ge 0}$ be a $(k,n)$-BCFW cell labeled by the Dyck path $P\in\mathcal{P}_{n,k,4}$. Then any $V\in S$ has a $P$-domino basis.
\end{conj}

\begin{eg}\label{eg_domino_basis}
Let $P\in\mathcal{P}_{12,3,4}$ be the Dyck path from \cref{fig:dyck_bijection}. Then the edges of $P$ are labeled as follows:
$$
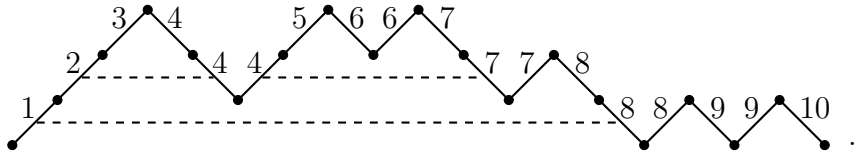
\begin{tikzpicture}[baseline=(p1)]
\pgfmathsetmacro{\u}{0.60};
\pgfmathsetmacro{\r}{0.18};
\pgfmathsetmacro{\vstep}{0.20};
\pgfmathsetmacro{\hstepup}{-0.08};
\pgfmathsetmacro{\hstepdown}{0.08};
\coordinate(u)at(\u,\u);
\coordinate(d)at(\u,-\u);
\tikzstyle{out1}=[inner sep=0,minimum size=1.2mm,circle,draw=black,fill=black]
\coordinate(p1)at(0,0){};
\coordinate(p2)at($(p1)+(u)$){};
\coordinate(p3)at($(p2)+(u)$){};
\coordinate(p4)at($(p3)+(u)$){};
\coordinate(p5)at($(p4)+(d)$){};
\coordinate(p6)at($(p5)+(d)$){};
\coordinate(p7)at($(p6)+(u)$){};
\coordinate(p8)at($(p7)+(u)$){};
\coordinate(p9)at($(p8)+(d)$){};
\coordinate(p10)at($(p9)+(u)$){};
\coordinate(p11)at($(p10)+(d)$){};
\coordinate(p12)at($(p11)+(d)$){};
\coordinate(p13)at($(p12)+(u)$){};
\coordinate(p14)at($(p13)+(d)$){};
\coordinate(p15)at($(p14)+(d)$){};
\coordinate(p16)at($(p15)+(u)$){};
\coordinate(p17)at($(p16)+(d)$){};
\coordinate(p18)at($(p17)+(u)$){};
\coordinate(p19)at($(p18)+(d)$){};
\foreach \x in {1,...,18}{
\pgfmathsetmacro{\y}{\x +1}
\node[inner sep=0](m\x)at($1/2*(p\x)+1/2*(p\y)$){};}
\draw[thick,dashed](m1)--(m14) (m2)--(m5) (m6)--(m11);
\draw[thick](p1.center)--(p2.center)--(p3.center)--(p4.center)--(p5.center)--(p6.center)--(p7.center)--(p8.center)--(p9.center)--(p10.center)--(p11.center)--(p12.center)--(p13.center)--(p14.center)--(p15.center)--(p16.center)--(p17.center)--(p18.center)--(p19.center);
\foreach \x in {1,...,19}{
\node[out1]at(p\x){};}
\node[inner sep=0]at($(m1)+(\hstepup-0.01,\vstep)$){$1$};
\node[inner sep=0]at($(m2)+(\hstepup-0.01,\vstep)$){$2$};
\node[inner sep=0]at($(m3)+(\hstepup,\vstep)$){$3$};
\node[inner sep=0]at($(m4)+(\hstepdown-0.02,\vstep)$){$4$};
\node[inner sep=0]at($(m5)+(\hstepdown-0.02,\vstep)$){$4$};
\node[inner sep=0]at($(m6)+(\hstepup,\vstep)$){$4$};
\node[inner sep=0]at($(m7)+(\hstepup,\vstep)$){$5$};
\node[inner sep=0]at($(m8)+(\hstepdown,\vstep)$){$6$};
\node[inner sep=0]at($(m9)+(\hstepup,\vstep)$){$6$};
\node[inner sep=0]at($(m10)+(\hstepdown,\vstep)$){$7$};
\node[inner sep=0]at($(m11)+(\hstepdown,\vstep)$){$7$};
\node[inner sep=0]at($(m12)+(\hstepup+0.02,\vstep)$){$7$};
\node[inner sep=0]at($(m13)+(\hstepdown,\vstep)$){$8$};
\node[inner sep=0]at($(m14)+(\hstepdown,\vstep)$){$8$};
\node[inner sep=0]at($(m15)+(\hstepup,\vstep)$){$8$};
\node[inner sep=0]at($(m16)+(\hstepdown,\vstep)$){$9$};
\node[inner sep=0]at($(m17)+(\hstepup+0.01,\vstep)$){$9$};
\node[inner sep=0]at($(m18)+(\hstepdown+0.09,\vstep)$){$10$};
\end{tikzpicture}\;\;.
$$
We see that
$$
\up_1 = 1,\; \up_2 = 2,\; \up_3 = 4,\quad \down_1 = 8,\; \down_2 = 4,\; \down_3 = 7.
$$
Then according to \cref{dyck_to_basis}, the $P$-domino bases $(v^{(1)}, v^{(2)}, v^{(3)})$ are precisely those given by the rows of the matrix\vspace*{-4pt}
\begin{align}\label{eqn:domino_matrix}
\begin{bmatrix}v^{(1)} \\ v^{(2)} \\ v^{(3)}\end{bmatrix} =\hspace*{-7pt}
\kbordermatrix{
& 1 & 2 & 3 & 4 & 5 & 6 & 7 & 8 & 9 & 10 & 11 & 12 \cr
& \alpha & \beta & 0 & 0 & 0 & 0 & 0 & \gamma & \delta & 0 & 0 & 1 \cr
& \alpha & \beta+\varepsilon & \zeta & \eta & \theta & 0 & 0 & 0 & 0 & 0 & 0 & 0 \cr
& -\alpha & -\beta & 0 & \iota & \kappa & 0 & \lambda & \mu & 0 & 0 & 0 & 0},
\end{align}
where $\alpha, \beta, \gamma, \delta, \varepsilon, \zeta, \eta, \theta, \iota, \kappa, \lambda, \mu > 0$. Let us explain in detail how the third row, $v^{(3)}$, is obtained. First, $d^{(3)}$ is a positive $4$-domino, which we take to be $(\dots, 0, \iota, \kappa, 0, \dots)$. Now, $e^{(3)}$ is a 7-domino, and is positive since none of the elements of $\{1,2,4\}$ are strictly between $4$ and $7$. We take $e^{(3)}$ to be $(\dots, 0, \lambda, \mu, 0, \dots)$. Finally, the last up step of $P$ before $\up_3$ which finishes at the same height that $\up_3$ begins is $\up_1$. Therefore $f^{(3)} = (-1)^{3-1-1}d^{(1)}$, and we have already taken $d^{(1)} = (\alpha, \beta, 0, \dots)$.

We can check that every element of the BCFW cell in $\Gr_{3,12}^{\ge 0}$ labeled by $P$ can be represented by a unique matrix of the form \eqref{eqn:domino_matrix}, by taking the $\oplus$-diagram in \cref{fig:PathsToLeDiagram}, performing \Le -moves to obtain a \Le -diagram, constructing the network parameterization matrix from \cref{network_param}, and carrying out appropriate row operations. This verifies \cref{conj:domino} in this particular case. It seems likely that the same method can be used to prove \cref{conj:domino} in general. To do so, one would need to figure out a systematic way to  perform the required \Le -moves and row operations. We also observe that constraining all parameters to be positive is not sufficient for the element of $\Gr_{3,12}$ represented by \eqref{eqn:domino_matrix} to lie in the corresponding BCFW cell; we also have the nontrivial inequality $\eta\kappa > \theta\iota$.
\end{eg}

We note in the cases $k=1$ and $k=2$, \cref{conj:domino} follows from the same arguments used to establish \cref{lem:bcfw-14} and \cref{thm:4-dominoes}. When $k=2$, the $P$-domino basis vectors $v^{(1)}$, $v^{(2)}$ of an element $V$ of a BCFW cell are the rows of the corresponding matrix in \cref{fig:cell-classification}, up to rescaling the rows by positive constants.
\begin{rmk}
If we take $P\in\mathcal{P}_{n,k,4}$ and delete the $2(n-3-k)$ edges incident to a peak, we obtain a Dyck path $P'$ with $2k$ steps. When $k=2$, $P'$ equals either $\begin{tikzpicture}[scale=0.35]
\pgfmathsetmacro{\u}{0.60};
\coordinate(u)at(\u,\u);
\coordinate(d)at(\u,-\u);
\tikzstyle{out1}=[inner sep=0,minimum size=1.2mm,circle,draw=black,fill=black]
\node[out1](p1)at(0,0){};
\node[out1](p2)at($(p1)+(u)$){};
\node[out1](p3)at($(p2)+(d)$){};
\node[out1](p4)at($(p3)+(u)$){};
\node[out1](p5)at($(p4)+(d)$){};
\draw[thick](p1.center)--(p2.center)--(p3.center)--(p4.center)--(p5.center);
\end{tikzpicture}$ or $\begin{tikzpicture}[scale=0.35]
\pgfmathsetmacro{\u}{0.60};
\coordinate(u)at(\u,\u);
\coordinate(d)at(\u,-\u);
\tikzstyle{out1}=[inner sep=0,minimum size=1.2mm,circle,draw=black,fill=black]
\node[out1](p1)at(0,0){};
\node[out1](p2)at($(p1)+(u)$){};
\node[out1](p3)at($(p2)+(u)$){};
\node[out1](p4)at($(p3)+(d)$){};
\node[out1](p5)at($(p4)+(d)$){};
\draw[thick](p1.center)--(p2.center)--(p3.center)--(p4.center)--(p5.center);
\end{tikzpicture}$, depending on whether the BCFW cell of $P$ is orthodox or deviant, respectively. In general, it may make sense to divide the $(k,n)$-BCFW cells into $C_k = \frac{1}{k+1}\binom{2k}{k}$ classes, based on $P'$.
\end{rmk}

\bibliographystyle{alpha}
\bibliography{bibliography}

\newcommand{\etalchar}[1]{$^{#1}$}
\begin{thebibliography}{AHBC{\etalchar{+}}16}

\bibitem[AHBC{\etalchar{+}}16]{abcgpt}
Nima Arkani-Hamed, Jacob Bourjaily, Freddy Cachazo, Alexander Goncharov,
  Alexander Postnikov, and Jaroslav Trnka.
\newblock {\em Grassmannian Geometry of Scattering Amplitudes}.
\newblock Cambridge University Press, 2016.
\newblock Preliminary version titled ``Scattering {A}mplitudes and the
  {P}ositive {G}rassmannian'' on the arXiv at
  \href{https://arxiv.org/abs/1212.5605}{\texttt{https://arxiv.org/abs/1212.5605}}.

\bibitem[AHT14]{arkani-hamed_trnka}
Nima Arkani-Hamed and Jaroslav Trnka.
\newblock The amplituhedron.
\newblock {\em J. High Energy Phys.}, (10):33, 2014.

\bibitem[AHTT]{ATT}
Nima Arkani-Hamed, Hugh Thomas, and Jaroslav Trnka.
\newblock Unwinding the amplituhedron in binary.
\newblock Preprint,
  \href{https://arxiv.org/abs/1704.05069}{\texttt{https://arxiv.org/abs/1704.05069}}.

\bibitem[Bay93]{bayer_93}
Margaret~M. Bayer.
\newblock Equidecomposable and weakly neighborly polytopes.
\newblock {\em Israel J. Math.}, 81(3):301--320, 1993.

\bibitem[BCF05]{BCF}
Ruth Britto, Freddy Cachazo, and Bo~Feng.
\newblock New recursion relations for tree amplitudes of gluons.
\newblock {\em Nuclear Phys. B}, 715(1-2):499--522, 2005.

\bibitem[BCFW05]{BCFW}
Ruth Britto, Freddy Cachazo, Bo~Feng, and Edward Witten.
\newblock Direct proof of the tree-level scattering amplitude recursion
  relation in {Y}ang-{M}ills theory.
\newblock {\em Phys. Rev. Lett.}, 94(18):181602, 4, 2005.

\bibitem[BGCT88]{bodroza_gutman_cyvin_tosic_1988}
O.~Bodro\v{z}a, I.~Gutman, S.~J. Cyvin, and R.~To\v{s}i\'c.
\newblock Number of {K}ekul\'e structures of hexagon-shaped benzenoids.
\newblock {\em J. Math. Chem.}, 2(3):287--298, 1988.

\bibitem[BH15]{bai_he_15}
Yuntao Bai and Song He.
\newblock The amplituhedron from momentum twistor diagrams.
\newblock {\em J. High Energy Phys.}, (2):065, front matter+38, 2015.

\bibitem[CG88]{kekule}
S.~J. Cyvin and I.~Gutman.
\newblock {\em Kekul\'e structures in benzenoid hydrocarbons}, volume~46 of
  {\em Lecture notes in chemistry}.
\newblock Springer, 1988.

\bibitem[Cyv86]{cyvin_1986}
Sven~J. Cyvin.
\newblock The number of \emph{Kekul\'{e}} structures of hexagon-shaped
  benzenoids and members of other related classes.
\newblock {\em Monatshefte f{\"u}r Chemie/Chemical Monthly}, 117(1):33--45,
  1986.

\bibitem[FZ99]{FZ}
Sergey Fomin and Andrei Zelevinsky.
\newblock Double {B}ruhat cells and total positivity.
\newblock {\em J. Amer. Math. Soc.}, 12(2):335--380, 1999.

\bibitem[GD52]{gordon_davison_52}
M.~Gordon and W.~H.~T. Davison.
\newblock Theory of resonance topology of fully aromatic hydrocarbons. {I}.
\newblock {\em J. Chem. Phys.}, 20(3):428--435, 1952.

\bibitem[GK50]{gantmakher_krein_50}
F.~R. Gantmaher and M.~G. Kre{\u\i}n.
\newblock {\em Oscillyacionye matricy i yadra i malye kolebaniya mehani\v
  ceskih sistem}.
\newblock Gosudarstv. Isdat. Tehn.-Teor. Lit., Moscow-Leningrad, 1950.
\newblock 2d ed. Translated into English by A. Eremenko
  \cite{gantmakher_krein_translation}.

\bibitem[GK02]{gantmakher_krein_translation}
F.~P. Gantmacher and M.~G. Krein.
\newblock {\em Oscillation matrices and kernels and small vibrations of
  mechanical systems}.
\newblock AMS Chelsea Publishing, Providence, RI, revised edition, 2002.
\newblock Translation based on the 1950 Russian original. Edited and with a
  preface by Alex Eremenko.

\bibitem[Hod43]{hodge43}
W.~V.~D. Hodge.
\newblock Some enumerative results in the theory of forms.
\newblock {\em Proc. Cambridge Philos. Soc.}, 39:22--30, 1943.

\bibitem[Kar17]{karp}
Steven~N. Karp.
\newblock Sign variation, the {G}rassmannian, and total positivity.
\newblock {\em J. Combin. Theory Ser. A}, 145:308--339, 2017.

\bibitem[KBT02]{klein_babic_trinajstic_02}
D.~J. Klein, D.~Babi\'{c}, and N.~Trinajsti\'{c}.
\newblock Enumeration in chemistry.
\newblock In {\em Chemical Modelling: Applications and Theory}, volume~2, pages
  56--95. The Royal Society of Chemistry, 2002.

\bibitem[Kek65]{kekule_65}
Aug. Kekul\'{e}.
\newblock Sur la constitution des substances aromatiques.
\newblock {\em Bulletin de la Soci\'{e}t\'{e} chimique de Paris, nouvelle
  s\'{e}rie}, 3:98--110, 1865.

\bibitem[Kek66]{kekule_66}
Aug. Kekul\'{e}.
\newblock Untersuchungen \"{u}ber aromatische {V}erbindungen.
\newblock {\em Annalen der Chemie und Pharmacie}, 137(2):129--196, 1866.

\bibitem[KW14]{KodamaWilliams}
Yuji Kodama and Lauren Williams.
\newblock K{P} solitons and total positivity for the {G}rassmannian.
\newblock {\em Invent. Math.}, 198(3):637--699, 2014.

\bibitem[KW17]{karpwilliams}
Steven~N. Karp and Lauren~K. Williams.
\newblock The {$m=1$} amplituhedron and cyclic hyperplane arrangements.
\newblock {\em Int. Math. Res. Not. IMRN \textup{(to appear)}}, 2017.

\bibitem[Lus94]{lusztig}
G.~Lusztig.
\newblock Total positivity in reductive groups.
\newblock In {\em Lie theory and geometry}, volume 123 of {\em Progr. Math.},
  pages 531--568. Birkh\"auser Boston, Boston, MA, 1994.

\bibitem[LW08]{LamWilliams}
Thomas Lam and Lauren Williams.
\newblock Total positivity for cominuscule {G}rassmannians.
\newblock {\em New York J. Math.}, 14:53--99, 2008.

\bibitem[Mac16]{macmahon}
Major Percy~A. MacMahon.
\newblock {\em Combinatory Analysis. {V}ol. {II}}.
\newblock Cambridge University Press, Cambridge, 1916.

\bibitem[Mor]{morales}
Alejandro~H. Morales.
\newblock Positive {G}rassmannian, lectures by {A}.\ {P}ostnikov.
\newblock Lecture notes from 18.318: Topics in Combinatorics at MIT,
  \url{http://www.math.ucla.edu/~ahmorales/18.318lecs/lectures.pdf}.

\bibitem[MR]{MarshRietsch}
Robert Marsh and Konstanze Rietsch.
\newblock The {$B$}-model connection and mirror symmetry for {G}rassmannians.
\newblock Preprint,
  \href{https://arxiv.org/abs/1307.1085}{\texttt{https://arxiv.org/abs/1307.1085}}.

\bibitem[Pos]{postnikov}
Alexander Postnikov.
\newblock Total positivity, {G}rassmannians, and networks.
\newblock Preprint,
  \href{http://math.mit.edu/~apost/papers/tpgrass.pdf}{\texttt{http://math.mit.edu/\textasciitilde
  apost/papers/tpgrass.pdf}}.

\bibitem[Pro]{propp_02}
James Propp.
\newblock Lattice structure for orientations of graphs.
\newblock Preprint, \url{https://arxiv.org/abs/math/0209005}.

\bibitem[Pro99]{propp_99}
James Propp.
\newblock Enumeration of matchings: problems and progress.
\newblock In {\em New perspectives in algebraic combinatorics ({B}erkeley,
  {CA}, 1996--97)}, volume~38 of {\em Math. Sci. Res. Inst. Publ.}, pages
  255--291. Cambridge Univ. Press, Cambridge, 1999.

\bibitem[PSW09]{PSW}
Alexander Postnikov, David Speyer, and Lauren Williams.
\newblock Matching polytopes, toric geometry, and the totally non-negative
  {G}rassmannian.
\newblock {\em J. Algebraic Combin.}, 30(2):173--191, 2009.

\bibitem[Ram97]{rambau_97}
J{\"o}rg Rambau.
\newblock Triangulations of cyclic polytopes and higher {B}ruhat orders.
\newblock {\em Mathematika}, 44(1):162--194, 1997.

\bibitem[Sta15]{stanley_catalan}
Richard~P. Stanley.
\newblock {\em Catalan numbers}.
\newblock Cambridge University Press, New York, 2015.

\bibitem[Stu88]{Sturmfels}
Bernd Sturmfels.
\newblock Totally positive matrices and cyclic polytopes.
\newblock In {\em Proceedings of the {V}ictoria {C}onference on {C}ombinatorial
  {M}atrix {A}nalysis ({V}ictoria, {BC}, 1987)}, volume 107, pages 275--281,
  1988.

\end{thebibliography}

\end{document}